\documentclass[3p]{elsarticle}

\usepackage{amssymb}
\usepackage{amsmath}
\usepackage{amsfonts}
\usepackage{epsfig}
\usepackage[mathscr]{eucal}
\usepackage[]{algorithm2e}
\usepackage{algorithmic}
\usepackage{setspace,color}
\usepackage{tikz}

\usetikzlibrary{snakes}
\usetikzlibrary{arrows,shapes}

\numberwithin{equation}{section}

\newtheorem{remark}{Remark}[section]

\def\R{{\mathbb R}}
\def\C{{\mathbb C}}
\def\DD{{\mathbb D}}
\def\N{{\mathbb N}}

\def\Z{{\mathbb Z}}

\def\KK{{\mathbb K}}
\def\MM{{\mathbb M}}

\def\II{{\mathbb I}}
\def\JJ{{\mathbb J}}

\def\BB{{\mathbb B}}
\def\OO{{\mathbb O}}
\def\su{\mathrm{supp}}

\def\diam{\mathrm{diam}}

\def\Per{\mathrm{Per}}

\def\div{\mathrm{div}}
\def\Om{\Omega}

\def\f{\frac}
\def\p{\partial}

\def\na{\nabla}
\def\la{\langle}
\def\ra{\rangle}

\def\bb{{\boldsymbol b}}

\def\rd{{\mathrm d}}
\def\bsd{{\boldsymbol d}}

\def\bsf{{\boldsymbol f}}

\def\bsg{{\boldsymbol g}}
\def\bsh{{\boldsymbol h}}
\def\bi{{\mathbf i}}

\def\bsi{{\boldsymbol i}}
\def\bsj{{\boldsymbol j}}
\def\bk{{\boldsymbol k}}

\def\bn{\boldsymbol{n}}

\def\bq{\boldsymbol{q}}

\def\bu{\boldsymbol{u}}
\def\bsv{\boldsymbol{v}}

\def\w{{\boldsymbol w}}

\def\x{\boldsymbol{x}}

\def\bsA{{\boldsymbol A}}

\def\bsI{{\boldsymbol I}}

\def\bsW{{\boldsymbol W}}

\def\bsT{\boldsymbol{T}}

\def\mB{{\mathcal B}}

\def\mI{{\mathcal I}}

\def\mN{{\mathcal N}}
\def\mO{{\mathcal O}}

\def\mT{{\mathcal T}}

\def\msF{{\mathscr F}}

\def\msS{{\mathscr S}}
\def\msT{{\mathscr T}}

\def\msX{{\mathscr X}}
\def\msY{{\mathscr Y}}

\def\im{\mathfrak{Im}}
\def\aal{\boldsymbol{\alpha}}
\def\bbe{\boldsymbol{\beta}}
\def\gga{\boldsymbol{\gamma}}

\def\lam{\boldsymbol{\lambda}}

\def\nnu{\boldsymbol{\nu}}
\def\rrh{\boldsymbol{\rho}}

\def\bi{\begin{itemize}} \def\ei{\end{itemize}}
\def\be{\begin{eqnarray*}}
\def\ee{\end{eqnarray*}}
\def\diam{\mathop{\rm diam}\nolimits}

\def\0{{\mathbf 0}}
\def\one{{\mathbf 1}}

\newcommand{\beq}{\begin{equation}}
\newcommand{\eeq}{\end{equation}}

\def\xxi{{\boldsymbol{\xi}}}
\def\eet{{\boldsymbol\eta}}

\def\La{{\boldsymbol\Lambda}}

\def\Ph{\boldsymbol{\Phi}}
\def\Sig{\boldsymbol{\Sigma}}

\def\Ps{{\boldsymbol \Psi}}
\def\wt{\widetilde}
\def\wh{\widehat}

\newcommand{\eps}{\varepsilon}
\def\la{\langle}
\def\ra{\rangle}

\def\XXint#1#2#3{{\setbox0=\hbox{$#1{#2#3}{\int}$ }
\vcenter{\hbox{$#2#3$ }}\kern-.55\wd0}}

\newtheorem{proposition}{Proposition}[section]

\newtheorem{thm}{Theorem}[section]
\newtheorem{corollary}{Corollary}[section]
\newdefinition{definition}{Definition}[section]
\newdefinition{rmk}{Remark}[section]
\newdefinition{notation}{Notation}[section]
\newproof{pf}{Proof}
\newproof{pop}{Proof of Proposition \ref{Prop1}}
\newproof{pop1}{Proof of Theorem \ref{Th1}}
\newproof{pop2}{Proof of Proposition \ref{Prop2}}

\title{An Edge Driven Wavelet Frame Model for Image Restoration}

\author[INS]{Jae Kyu Choi}
\ead{jaycjk@sjtu.edu.cn}
\author[BICMR]{Bin Dong\fnref{BD}\corref{cor}}
\ead{dongbin@math.pku.edu.cn}
\author[INS,MOELSC]{Xiaoqun Zhang\fnref{XQZ}}
\ead{xqzhang@sjtu.edu.cn}

\cortext[cor]{Corresponding author}

\address[INS]{Institute of Natural Sciences, Shanghai Jiao Tong University, 200240 Shanghai China}
\address[BICMR]{Beijing International Center for Mathematical Research, Peking University, 100871 Beijing China}
\address[MOELSC]{School of Mathematical Sciences, and MOE-LSC, Shanghai Jiao Tong University, 200240 Shanghai China}

\fntext[BD]{B. Dong is supported in part by the Thousand Talents Plan of China.}
\fntext[XQZ]{J. Choi and X.  Zhang are partially supported by the Young Top-notch Talent program of China, 973 program (No. 2015CB856004) and Sino-German center grant (GZ1025).}

\begin{document}

\begin{abstract}
Wavelet frame systems are known to be effective in capturing singularities from noisy and degraded images. In this paper, we introduce a new edge driven wavelet frame  model for image restoration by approximating images as piecewise smooth functions. With an implicit representation of image singularities sets, the proposed model inflicts different strength of regularization on smooth and singular image regions and edges. The proposed edge driven model is robust to both image approximation and singularity estimation. The implicit formulation also enables an asymptotic analysis of the proposed models and a rigorous connection between the discrete model and a general continuous variational model. Finally, numerical results on image inpainting and deblurring show that the proposed model is compared favorably against several popular image restoration models.
\end{abstract}

\begin{keyword}
Image restoration \sep (tight) wavelet frames \sep framelets \sep edge estimation \sep variational method \sep pointwise convergence \sep $\Gamma$-convergence
\end{keyword}

\maketitle


\pagestyle{myheadings}
\thispagestyle{plain}
\markboth{Jae Kyu Choi, Bin Dong, and Xiaoqun Zhang}{An Edge Driven Wavelet Frame Model for Image Restoration}

\section{Introduction}\label{Introduction}

Image restoration, including image denoising, deblurring, inpainting, computed tomography, etc., is one of the most important areas in imaging science. It aims at recovering an image of high-quality from a given measurement which is degraded during the process of imaging, acquisition, and communication. An image restoration problem is typically modeled as the following linear inverse problem:
\begin{align}\label{Linear_IP}
\bsf=\bsA\bu+\eet,
\end{align}
where $\bsf$ is the degraded measurement or the observed image, $\eet$ is a certain additive noise, and $\bsA$ is some linear operator which takes different forms for different image restoration problems. Note that this paper involves both functions (operators) and their discrete counterparts. We shall use regular characters to denote functions or operators and use bold-faced characters to denote their discrete analogs. For example, we use $A$ to denote a linear operator between two function spaces and $u$ as an element in a function space, while we use $\bsA$ and $\bu$ to denote their corresponding discretized versions (the type of discretization will be made clear later).

The operator $\bsA$ is in general ill-conditioned (e.g. for deblurring) or non-invertible (e.g. for inpainting). Naive inversions of \eqref{Linear_IP} in the presence of noise $\eet$ will inevitably lead to significant noise amplification. Hence, in order to obtain a high quality recovery from the ill-posed linear inverse problem \eqref{Linear_IP}, a proper regularization on the images to be recovered is needed. Successful regularization based methods include the Rudin-Osher-Fatemi model \cite{L.I.Rudin1992} and its nonlocal variants \cite{G.Gilboa2008,X.Zhang2010}, the inf-convolution model \cite{A.Chambolle1997}, the total generalized variation (TGV) model \cite{K.Bredies2014,K.Bredies2010}, the combined first and second order total variation model \cite{M.Bergounioux2010,J.Liang2015,K.Papafitsoros2014}, and the applied harmonic analysis approach such as curvelets \cite{E.Candes2006}, Gabor frames \cite{Daubechies1992,Groechenig2001,H.Ji2016a,Mallat2008}, shearlets \cite{G.Kutyniok2011}, complex tight framelets \cite{B.Han2014}, wavelet frames \cite{C.Bao2016,J.F.Cai2009,J.F.Cai2010,J.F.Cai2009/10,R.H.Chan2003,I.Daubechies2007,B.Dong2013b,M.Elad2005,M.J.Fadili2009,M.A.T.Figueiredo2003,J.L.Starck2005,Y.Zhang2013}, etc.  The common concept of these methods  is  to find  sparse approximation of images using a properly designed linear transformation together with a sparsity promoting regularization term (such as the widely used $\ell_1$ norm). A typical $\ell_1$ norm based regularization model takes the following form
\begin{align}\label{Model:L1}
\min_{\bu}~\lambda\big\|\Ph\bu\big\|_1+\f{1}{2}\big\|\bsA\bu-\bsf\big\|_2^2
\end{align}
where $\Phi$ is some sparsifying linear transform (such as wavelet transform or $\nabla$). This general formulation is widely applied in image restoration for regularizing designed smooth image components while preserving image singularities.

Meanwhile, the idea of explicitly taking image singularities into consideration was first explored in the pioneer work \cite{D.Mumford1989}, where the following model, known as the Mumford-Shah model, was introduced:
\begin{align}\label{MumfordShah}
\min_{u,\Sigma}~\lambda\int_{\Om\setminus\Sigma}\big|\na u\big|^2 \rd\x+\nu\big|\Sigma\big|+\f{1}{2}\big\|u-f\big\|_{L_2(\Om)}^2.
\end{align}
Here, $\big|\Sigma\big|$ denotes the length of one-dimensional curve $\Sigma$ representing edges. Due to the smoothness promoting property of $\ell_2$ norm, the above Mumford-Shah functional encourages $u$ to be smooth except along $\Sigma$ (see \cite{G.Aubert2006,Chambolle1999,D.Mumford1989} for detailed surveys on the Mumford-Shah model and \cite{L.Bar2011,M.Jiang2014} for the applications to image restoration).
In a discrete setting, if we know the exact locations of image singularities, then  we can recover the image $\bu$ with sharp edges by solving the following minimization problem:
\begin{align}\label{Tikhonov_Oracle}
\min_{\bu}~\lambda\left\|\big(\Ph\bu\big)_{\Sig^c}\right\|_2^2+\frac{1}{2}\big\|\bsA\bu-\bsf\big\|_2^2
\end{align}
where $\Sig$ is the index set of pixels corresponding to image singularities. The problem \eqref{Tikhonov_Oracle} is easy to solve once we know $\Sig$. However, the restoration result of \eqref{Tikhonov_Oracle} can be highly sensitive to the estimation of $\Sig$, and the main challenge lies in how to identify $\Sig$ as accurately as possible from degraded observed images.

Sparse regularization with wavelet frame transforms \eqref{Model:L1} is successfully applied in various imaging problems, due to its effectiveness of capturing multiscale singularities using compactly supported wavelet frame functions of varied vanishing moments. In connection with Mumford-Shah model, the authors in \cite{J.F.Cai2016} exploited the favorable properties of wavelet frames, and proposed the following piecewise smooth wavelet frame image restoration model:
\begin{align}\label{CDS}
\min_{\bu,\Sig}~\left\|\big(\lam\cdot\bsW\bu\big)_{\Sig^c}\right\|_2^2+\left\|\big(\gga\cdot\bsW\bu\big)_{\Sig}\right\|_1+\f{1}{2}\big\|\bsA\bu-\bsf\big\|_2^2,
\end{align}
where $\bsW$ is a wavelet frame transform and $\Sig$ is the image singularities set to be estimated. As image singularities can be well approximated by wavelet frame coefficients of large magnitude, \eqref{CDS} uses the $\ell_2$ norm to promote the smoothness of image away from $\Sig$, and uses the $\ell_1$ norm to recover sharp features lying in $\Sig$ \cite{J.F.Cai2016}. The authors proved that under the assumption of a fixed index set $\Sig$, the discrete model \eqref{CDS} converges to a new variational model as the resolution goes to infinity. A special case of the variational model is related to (and yet significantly different from) the Mumford-Shah functional \eqref{MumfordShah}. As a byproduct of the analysis in \cite{J.F.Cai2016}, it demonstrated that the model \eqref{CDS} is more computationally tractable than the Mumford-Shah model \eqref{MumfordShah}. Interested readers should consult \cite{J.F.Cai2016} for more details.

Another model that exploits the similar idea is the following constrained minimization model proposed in \cite{H.Ji2016}:
\begin{align}\label{JLS}
\begin{split}
&\min_{\bu,\Sig}~\lambda\left\|\big(\bsW\bu\big)_{\Sig^c}\right\|_2^2+\frac 1 2\big\|\bsA\bu-\bsf\big\|_2^2\\
&\text{subject to}~~~|\Sig^c|\geq t~~\text{and}~~\Sig\in\mO,
\end{split}
\end{align}
where $\mO$ is the feasible set for $\Sig$, and the constraint on $|\Sig^c|$ is imposed to promote the regularity of the singularity set, by implication, the sparsity of the wavelet frame coefficients $\bsW\bu$. Unlike \eqref{CDS} which directly updates $\Sig$ by comparing the $\ell_1$ norm and $\ell_2$ norm of $\bsW\bu$ at each step, additional geometric constraints on $\Sig$ in \eqref{JLS} are utilized to regularize image singularities.

Even though both \eqref{CDS} and \eqref{JLS} showed significant improvements over the typical wavelet frame sparsity based image restoration model \eqref{Model:L1}, the above two models have their own drawbacks. For \eqref{CDS}, since $\Sig$ is estimated solely depending on the wavelet frame coefficients, the estimated $\Sig$ may capture the unwanted isolated singularities when the measurement $\bsf$ is severely noisy. In addition, since $\bsW\bu$ is split into the $\ell_1$  and the $\ell_2$ part, the reconstructed image may suffer from the staircase effect on the interface of $\Sig$ and $\Sig^c$. For \eqref{JLS}, as the coefficients $\bsW\bu$ on $\Sig$ are not directly penalized, \eqref{JLS} may introduce overly sharpened singularities compared to \eqref{CDS}, especially in the case of deblurring with a severely degraded $\bsf$. In addition, it is difficult to rigorously analyze the model and its solutions with the presence of the singularities set $\Sig$.

In this paper, we propose a new edge driven wavelet frame based image restoration model. We use the term ``edge driven'' as the proposed model continues to exploit the idea of alternate recovery of the image and the estimation of its singularities set in a different form. Here, we provide a first glance of the model as follows:
\begin{align}\label{Proposed}
\min_{\bu,0\leq\bsv\leq 1}~\left\|(\one-\bsv)\cdot\big(\lam\cdot\bsW\bu\big)\right\|_1+\left\|\bsv\cdot\big(\gga\cdot\bsW'\bu\big)\right\|_1+\big\|\rrh\cdot\bsW''\bsv\big\|_1+\f{1}{2}\big\|\bsA\bu-\bsf\big\|_2^2,
\end{align}
where $\bu$ is the image to be reconstructed, $\bsv$ denotes a relaxed set indicator of the singularities set, and $\bsW$, $\bsW'$, and $\bsW''$ are three wavelet frame transforms applied to different components of the images. For the clarity of presentation, the detailed definition and the analysis of the model in a multi-level decomposition form are postponed until Section \ref{WaveletModelAlg}-\ref{VariationAsymAnal}.

Our model is closely related to the piecewise smooth wavelet frame models \eqref{CDS} and \eqref{JLS}. In fact,  \eqref{Proposed} can be viewed as a relaxation of
\begin{align}\label{Proposed_Original}
\min_{\bu,\Sig}~\left\|\big(\lam\cdot\bsW\bu\big)_{\Sig^c}\right\|_1+\left\|\big(\gga\cdot\bsW'\bu\big)_{\Sig}\right\|_1+\big\|\rrh\cdot\bsW''\one_{\Sig}\big\|_1+\f{1}{2}\big\|\bsA\bu-\bsf\big\|_2^2
\end{align}
where $\Sig$ is the estimated singularities of $\bu$  and $\one_{\Sig}$ is its set indicator. The first term is used to restore smooth regions of an image, while the second term preserves singularities, and the third term provides the regularization on singularities to enhance sharp image features. In other words, our model inflicts a different strength of regularization in smooth image regions and near image singularities such as edges, and actively restores/enhances sharp image features at the same time. As the first two terms are exchangeable, an appropriate choice of the wavelet frame transforms as well as the associated parameters is needed to obtain desired effects. The details of the properties of the three transforms will be detailed in Section \ref{ModelAlg}.

Compared to the two existing models \eqref{CDS} and \eqref{JLS}, it should be noted that instead of using $\ell_2$ norm, $\ell_1$ norm is used to promote regularity in the smooth region, as the image singularities can be better protected if the singularity set $\Sig$ is not accurate. This leads to a more robust image approximation that is less sensitive to the estimation of the singularities of the unknown true image from the degraded measurement. In addition, an implicit and relaxed representation of the singularity set allows continuous overlap between the smooth and the sharp image regions in the transform domain. We expect that such overlap helps to suppress the staircase effects near the interface. Finally, representing the singularity set implicitly enables us to provide an asymptotic analysis of the model with respect to both $u$ and $v$, in contrast to that of \eqref{CDS} where the singularity set is assumed to be fixed.


To facilitate a better understanding of the proposed model \eqref{Proposed} and its relation to some existing variational models, we will present an asymptotic analysis of the proposed model. We discover that the continuum limit of the proposed model (after a reformulation) takes the following form
\begin{align}\label{Variational}
\begin{split}
\min_{u,0\leq v\leq1}~\lambda\int_{\Om}(1-v)\left(\sum_{\aal\in\II}|\p^{\aal}u|^2\right)^{\f{1}{2}}\rd\x&+\gamma\int_{\Om}v\left(\sum_{\aal\in\II'}|\p^{\aal}u|^2\right)^{\f{1}{2}}\rd\x\\
&+\rho\int_{\Om}\left(\sum_{\aal\in\II''}|\p^{\aal}v|^2\right)^{\f{1}{2}}\rd\x+\f{1}{2}\big\|Au-f\big\|_{L_2(\Om)}^2,
\end{split}
\end{align}
which is an edge driven variational model that includes several existing variational and partial differential equation (PDE) models as special cases (see Subsection \ref{CorrespondingVariational} for more details).

The rest of this paper is organized as follows. In Section \ref{WaveletFrameReview}, we introduce some basics of wavelet frame that will be used in later sections. We propose the discrete edge driven wavelet frame based model and its associated algorithm in Section \ref{WaveletModelAlg}. Numerical simulations of our proposed model and comparisons with some of the existing models are conducted at the end of this section. In Section \ref{VariationAsymAnal}, we present the continuum limit of the the proposed discrete model and provide a rigorous asymptotic analysis. All technical proofs will be postponed to the appendix.

\section{Preliminaries on Wavelet Frame}\label{WaveletFrameReview}

In this section, we present some basics of wavelet frame theory and some preliminary results.

\subsection{Tight Wavelet Frames}\label{WaveletFrame} In this subsection, we briefly introduce the concept of tight frames and wavelet tight frames. For the details, one may consult \cite{Daubechies1992,I.Daubechies2003,A.Ron1997} for theories of frames and wavelet frames, \cite{Shen2010} for a short survey on the theory and applications of frames, and \cite{B.Dong2013,B.Dong2015} for more detailed surveys.

A countable set $\msX\subseteq L_2(\R^d)$ with $d\in\N$ is called a tight frame of $L_2(\R^d)$ if
\begin{align}\label{TightFrame}
\|f\|_{L_2(\R^d)}^2=\sum_{\varphi\in\msX}|\la f,\varphi\ra|^2~~~~~\text{for all}~~~f\in L_2(\R^d),
\end{align}
where $\la\cdot,\cdot\ra$ is the inner product on $L_2(\R^d)$, and $\la f,\varphi\ra$ is called the canonical coefficient of $f$.

For given $\Ps=\big\{\psi_l:l=1,\cdots,r\big\}\subseteq L_2(\R^d)$ and $N\in\N$, the corresponding quasi-affine system $\msX^N(\Ps)$ generated by $\Ps$ is defined by the collection of the dilations and the shifts of the members in $\Ps$:
\begin{align}\label{QASystem}
\msX^N(\Ps)=\big\{\psi_{l,n,\bk}:1\leq l\leq r,~n\in\Z,~\bk\in\Z^d\big\}
\end{align}
where $\psi_{l,n,\bk}$ is defined as
\begin{align}\label{QAFramelet}
\psi_{l,n,\bk}(\x)=\left\{\begin{array}{cl}
2^{\f{nd}{2}}\psi_l(2^n\x-\bk)~~&~~n\geq N;\vspace{0.5em}\\
2^{\left(n-\f{N}{2}\right)d}\psi_l(2^n\x-2^{n-N}\bk)~~&~~n<N.
\end{array}\right.
\end{align}
When $\msX^N(\Ps)$ forms a tight frame of $L_2(\R^d)$, each $\psi_1,\cdots,\psi_r$ is called a (tight) framelet and the entire system $\msX^N(\Ps)$ is called a (tight) wavelet frame. In particular when $N=0$, we simply write $\msX(\Ps)=\msX^0(\Ps)$. Note that in the literature, the affine system is widely used, which corresponds to the decimate wavelet (frame) transform. The quasi-affine system, which corresponds to the undecimated wavelet (frame) transformation, was first introduced and analyzed in \cite{A.Ron1997}. Throughout this paper, we only discuss the quasi-affine system \eqref{QAFramelet} because it generally performs better in image restoration and the connection to PDE is more natural than the widely used affine system \cite{J.F.Cai2012,J.F.Cai2016,B.Dong2016}. The interested reader can find further details on the affine wavelet frame systems and its connections to the quasi-affine frames in \cite{A.Chai2007,B.Dong2013,A.Ron1997}.

The constructions of framelets $\Ps$, which are desirably (anti-)symmetric and compactly supported functions, are usually based on a multiresolution analysis (MRA) generated by some refinable function $\phi$ with a refinement mask $\bq_0$ such that
\begin{align}\label{MRA-RF}
\phi(\x)=2^d\sum_{\bk\in\Z^d}\bq_0[\bk]\phi(2\x-\bk).
\end{align}
The idea of an MRA based construction of $\Ps=\big\{\psi_1,\cdots,\psi_r\big\}\subseteq L_2(\R^d)$ is to find finitely supported masks $\bq_l$ such that
\begin{align}\label{MRA-Fra}
\psi_l(\x)=2^d\sum_{\bk\in\Z^d}\bq_l[\bk]\phi(2\x-\bk)~~~~~l=1,\cdots,r.
\end{align}
The sequences $\bq_1,\cdots,\bq_r$ are called wavelet frame mask or the high pass filters of the system, and the refinement mask $\bq_0$ is also called the low pass filter.

The unitary extension principle (UEP) of \cite{A.Ron1997} provides a general theory of the construction of MRA based tight wavelet frames. Briefly speaking, as long as $\big\{\bq_0,\bq_1,\cdots,\bq_r\big\}$ are compactly supported and their Fourier series
\begin{align*}
\wh{\bq}_l(\xxi)=\sum_{\bk\in\Z^d}\bq_l[\bk]e^{-i\xxi\cdot\bk},~~~~~~l=0,\cdots,r,~~~\xxi\in\R^d
\end{align*}
satisfy
\begin{align}\label{UEP}
\sum_{l=0}^r\left|\wh{\bq}_l(\xxi)\right|^2=1~~~~\text{and}~~~~\sum_{l=0}^r\wh{\bq}_l(\xxi)\overline{\wh{\bq}_l(\xxi+\nnu)}=0
\end{align}
for all $\nnu\in\big\{0,\pi\big\}^d\setminus\big\{\0\big\}$ and $\xxi\in[-\pi,\pi]^d$, the quasi-affine system $\msX(\Ps)$ with $\Ps=\big\{\psi_1,\cdots,\psi_r\big\}$ defined by \eqref{MRA-Fra} forms a tight frame of $L_2(\R^d)$, and the filters $\big\{\bq_0,\bq_1,\cdots,\bq_r\big\}$ form a discrete tight frame on $\ell_2(\Z^d)$ \cite{B.Dong2013}.

One of the most widely used examples is the piecewise linear B-spline \cite{I.Daubechies2003} for $L_2(\R)$, which has one refinable function and two framelets with the associated filters
\begin{align*}
\bq_0=\f{1}{4}\big[\begin{array}{ccc}
1&2&1
\end{array}\big],~~~~\bq_1=\f{\sqrt{2}}{4}\big[\begin{array}{ccc}
1&0&-1
\end{array}\big],~~\text{and}~~\bq_2=\f{1}{4}\big[\begin{array}{ccc}
-1&2&-1
\end{array}\big].
\end{align*}
Indeed, it can be shown that the above $\big\{\bq_0,\bq_1,\bq_2\big\}$ satisfies \eqref{UEP}, so that $\msX(\Ps)$ with $\Ps=\big\{\psi_1,\psi_2\big\}$ defined by \eqref{MRA-Fra} forms a tight frame on $L_2(\R)$.

For the practical concern, we need to construct tight frames for $L_2(\R^d)$ with $d\geq2$, because the discrete image is two or three dimensional array. One possible way is by taking tensor products of univariate tight frames \cite{J.F.Cai2012,J.F.Cai2016,Daubechies1992,B.Dong2013}. Throughout this paper, we will only consider two-dimensional case. Given a set of univariate masks $\big\{\bq_0,\bq_1,\cdots,\bq_r\big\}$, we define two-dimensional masks $\bq_{\aal}[\bk]$ with $\aal=(\alpha_1,\alpha_2)$ and $\bk=(k_1,k_2)$ as
\begin{align*}
\bq_{\aal}[\bk]=\bq_{\alpha_1}[k_1]\bq_{\alpha_2}[k_2],~~~~~~0\leq\alpha_1,\alpha_2\leq r,~~\bk=(k_1,k_2)\in\Z^2
\end{align*}
so that the corresponding $2$D refinable function and framelets are defined as
\begin{align*}
\psi_{\aal}(\x)=\psi_{\alpha_1}(x_1)\psi_{\alpha_2}(x_2),~~~~~0\leq\alpha_1,\alpha_2\leq r,~~\x=(x_1,x_2)\in\R^2
\end{align*}
with $\psi_0=\phi$ for convenience. If the univariate masks $\big\{\bq_l:l=1,\cdots,r\big\}$ are constructed from UEP, then it can be verified that $\big\{\bq_{\aal}:\aal\in\{0,\cdots,r\big\}^2\setminus\{\0\}\big\}$ satisfies \eqref{UEP} and thus $\msX(\Ps)$ with
\begin{align*}
\Ps=\big\{\psi_{\aal}:\aal\in\{0,\cdots,r\}^2\setminus\{\0\}\big\}
\end{align*}
forms a tight frame for $L_2(\R^2)$.

In the discrete setting, $\bu\in\mI_2$, where $\mI_2\simeq\R^{N_1\times N_2}$ denotes the space of two-dimensional discrete images. Throughout this paper, we assume for simplicity that all images are square images; $N_1=N_2=N$, and we only consider the MRA based tensor product wavelet frame system. We denote the two-dimensional fast (discrete) framelet transform, or the analysis operator (see, e.g., \cite{B.Dong2013}) with $L$ levels of decomposition as
\begin{align}
\bsW\bu=\left\{\bsW_{l,\aal}\bu:(l,\aal)\in\big(\big\{0,\cdots,L-1\big\}\times\BB\big)\cup\big\{(L-1,\0)\big\}\right\}~~~~~~~\bu\in\mI_2
\end{align}
where $\BB=\big\{0,\cdots,r\big\}^2\setminus\big\{\0\big\}$ is the framelet band. Then $\bsW$ is a linear operator with the frame coefficients $\bsW_{l,\aal}\bu\in\mI_2$ of $\bu$ at level $l$ and band $\aal$ being defined as
\begin{align*}
\bsW_{l,\aal}\bu=\bq_{l,\aal}[-\cdot]\circledast\bu.
\end{align*}
Here, $\circledast$ denotes the discrete convolution with a certain boundary condition (e.g., periodic boundary condition), and $\bq_{l,\aal}$ is defined as
\begin{align}
\bq_{l,\aal}=\wt{\bq}_{l,\aal}\circledast\wt{\bq}_{l-1,\0}\circledast\cdots\circledast\wt{\bq}_{0,\0}~~\text{with}~~\wt{\bq}_{l,\bsi}[\bk]=\left\{\begin{array}{rl}
\bq_{\aal}[2^{-l}\bk],&\bk\in 2^l\Z^2;\vspace{0.5em}\\
0,&\bk\notin 2^l\Z^2.
\end{array}\right.
\end{align}
Notice that $\bq_{0,\aal}=\bq_{\aal}$ and $\bsW_{0,\aal}\bu=\bsW_{\aal}\bu=\bq_{\aal}[-\cdot]\circledast\bu$.

The synthesis framelet transform is denoted as $\bsW^T$, the adjoint of $\bsW$. Since we consider a tight wavelet frame, we have the following perfect reconstruction formula
\begin{align*}
\bu=\bsW^T\bsW\bu
\end{align*}
for all $\bu\in\mI_2$.

\subsection{Vanishing Moments and Related Theory}\label{VanishingMoments}

The vanishing moments of framelets are closely related to the orders of differential operators and their corresponding finite difference operators. It is a crucial observation first made in \cite{J.F.Cai2012} and was further explored in \cite{,J.F.Cai2016,B.Dong2013a,B.Dong2016}, and will be vital to our analysis as well.

Throughout this paper, for a given multi-index $\aal=(\alpha_1,\alpha_2)\in\N_0^2$, we denote $|\aal|=\alpha_1+\alpha_2$. For two multi-indices $\aal$ and $\bbe$, we say $\bbe\leq\aal$ if $\beta_j\leq\alpha_j$ for all $j=1,2$. For $\x\in\R^2$ and a multi-index $\aal\in\N_0^2$, we denote $\x^{\aal}=x_1^{\alpha_1}x_2^{\alpha_2}$. We also define the mixed partial differential operator $\p^{\aal}$ as
\begin{align*}
\p^{\aal}=\p_2^{\alpha_2}\p_1^{\alpha_1}~~\text{where}~~\p_j=\f{\p}{\p x_j}.
\end{align*}
In particular, we use $\p_{\x}^{\aal}$ and $\p_{\xxi}^{\aal}$ to highlight the variable whenever it is needed to avoid confusion. For one-dimensional case, we will use the standard notation $f'$, $f''$, $f^{(\alpha)}$ etc.

Recall that the vanishing moments of a univariate function is the order of zeros of its Fourier transform at the origin. More precisely, $\psi\in L_2(\R)$ has \emph{vanishing moments of order} $\alpha\in\N_0$ if
\begin{align*}
\int_{-\infty}^{\infty}x^{\beta}\psi(x)\rd x=i^{\beta}\wh{\psi}^{(\beta)}(0)=0
\end{align*}
for all $\beta<\alpha$ but $\int_{-\infty}^{\infty}x^{\alpha}\psi(x)\rd x=i^{\alpha}\wh{\psi}^{(\alpha)}(0)\neq0$. Here, $\wh{\psi}$ is the Fourier transform of $\psi$ defined as
\begin{align*}
\wh{\psi}(\xi)=\msF(\psi)(\xi)=\int_{-\infty}^{\infty}\psi(x)e^{-i\xi x}\rd x,~~~~~~~\xi\in\R.
\end{align*}
We say that $\psi$ has the vanishing moment of order $0$ if $\int_{-\infty}^{\infty}\psi\rd x\neq0$. Likewise, we can define the vanishing moments of two-dimensional framelet function $\psi\in L_2(\R^2)$. We say that $\psi\in L_2(\R^2)$ has \emph{vanishing moments of order} $\aal=(\alpha_1,\alpha_2)\in\N_0^2$ provided that $\int_{\R^2}\x^{\aal}\psi(\x)\rd\x=i^{|\aal|}\p_{\xxi}^{\aal}\wh{\psi}(\0)\neq0$, and
\begin{align*}
\int_{\R^2}\x^{\bbe}\psi(\x)\rd\x=i^{|\bbe|}\p_{\xxi}^{\bbe}\wh{\psi}(\0)=0
\end{align*}
for all $\bbe\in\N_0^2$ with $|\bbe|<|\aal|$ and for all $\bbe\in\N_0^2$ with $|\bbe|=|\aal|$ but $\bbe\neq\aal$. Here, $\wh{\psi}$ is the Fourier transform of $\psi$ defined as
\begin{align*}
\wh{\psi}(\xxi)=\msF(\psi)(\xxi)=\int_{\R^2}\psi(\x)e^{-i\xxi\cdot\x}\rd\x,~~~~~~\xxi\in\R^2.
\end{align*}
We say that $\psi$ has a vanishing moment of order $\0$ if $\wh{\psi}(\0)=\int_{\R^2}\psi(\x)\rd\x\neq0$. Note that if $\psi$ is a tensor product framelet function, then its vanishing moments are determined by the vanishing moments of constituent univariate framelet functions.

In the literature of wavelet frame, we interpret the digital image $\bu$ as discrete sampling of an underlying function $u$ via the inner product with the corresponding refinable function $\phi$:
\begin{align*}
\bu[\bk]=2^n\big\la u,\phi_{n,\bk}\big\ra
\end{align*}
where for two-dimensional cases, $\phi_{n,\bk}$ (as well as $\psi_{n,\bk}$, etc.) takes the form
\begin{align}\label{Frameletn}
\phi_{n,\bk}=2^n\phi(2^n\cdot-\bk).
\end{align}
When discrete wavelet transform is applied on $\bu$, the underlying quasi-affine system we use is $\msX^n(\Ps)$, and we have
\begin{align}\label{Frameletn-1}
\phi_{n-1,\bk}=2^{n-2}\phi(2^{n-1}\cdot-2^{-1}\bk).
\end{align}
By \eqref{MRA-Fra}, we have $\psi_{\aal,n-1,\bk}=\sum_{\bsj\in\Z^2}\bq_{\aal}[\bsj-\bk]\phi_{n,\bsj}$, and the coefficients in the $\aal$th band satisfy
\begin{align*}
\big(\bsW_{\aal}\bu\big)[\bk]=\big(\bq_{\aal}[-\cdot]\ast\bu\big)[\bk]&=\sum_{\bsj\in\Z^2}\bq_{\aal}[\bsj-\bk]\bu[\bsj]=2^n\sum_{\bsj\in\Z^2}\bq_{\aal}[\bsj-\bk]\big\la u,\phi_{n,\bsj}\big\ra\\
&=2^n\left\la u,\sum_{\bsj\in\Z^2}\bq_{\aal}[\bsj-\bk]\phi_{n,\bsj}\right\ra=2^n\big\la u,\psi_{\aal,n-1,\bk}\big\ra,
\end{align*}
where $\ast$ denotes the discrete convolution. The key observation made by \cite{J.F.Cai2012} is that for the piecewise B-spline framelets $\psi_{\aal}$, there exists a function $\varphi_{\aal}$ associated to $\psi_{\aal}$ such that $\int_{\R^2}\varphi_{\aal}\rd\x\neq0$, $\psi_{\aal}=\p^{\aal}\varphi_{\aal}$ and $\su(\psi_{\aal})=\su(\varphi_{\aal})$, and the explicit formulae of $\varphi_{\aal}$ are given in \cite{Z.Shen2013}. With the aid of the theory of distribution \cite{Hormander1983,Rudin1991}, Proposition \ref{Prop1} generalizes the same result to any tensor product framelet. The proof can be found in \ref{ProofProp1}.

\begin{proposition}\label{Prop1} Assume that a framelet function $\psi_{\aal}\in L_2(\R^2)$ has vanishing moments of order $\aal$, and it is generated by the tensor product of univariate framelet functions. If its support is a two-dimensional box $[a_1,b_1]\times[a_2,b_2]$, then there exists the unique $\varphi_{\aal}\in L_2(\R^2)$ such that $\varphi_{\aal}$ is differentiable up to order $\aal$ a.e.,
\begin{align*}
c_{\aal}=\int_{\R^2}\varphi_{\aal}(\x)\rd\x\neq0~~~~\text{and}~~~~\psi_{\aal}=\p^{\aal}\varphi_{\aal}.
\end{align*}
Moreover, $\su(\varphi_{\aal})=\su(\psi_{\aal})$.
\end{proposition}

Let $\Om=(0,1)^2\subseteq\R^2$. Recall that $W_1^s(\Om)$ with $s\in\N$ is the Sobolev space defined as
\begin{align}\label{W_1^sSobolev}
W_1^s(\Om)=\big\{u\in L_1(\Om):\p^{\aal}u\in L_1(\Om)~~\text{for}~~|\aal|\leq s\big\},
\end{align}
where $\p^{\aal}u$ denotes the $\aal$th weak derivative of $u$. Then $W_1^s(\Om)$ equipped with the norm defined as
\begin{align*}
\|u\|_{W_1^s(\Om)}=\sum_{|\aal|\leq s}\|\p^{\aal}u\|_{L_1(\Om)}
\end{align*}
is a Banach space. Note that $W_1^s(\Om)\subseteq L_2(\Om)$ by Sobolev imbedding theorem \cite{Adams1975,H.Attouch2014}, and the proof of Proposition \ref{Prop1} implies that $\varphi_{\aal}$ is at least bounded and continuous. Hence, both $\big\la u,\psi_{\aal,n-1,\bk}\big\ra$ and $\big\la\p^{\aal}u,\varphi_{\aal,n-1,\bk}\big\ra$ are always well-defined for $u\in W_1^s(\Om)$ whenever $\su(\psi_{\aal,n-1,\bk})\subseteq\overline{\Om}$, and we arrive at the following proposition which provides a connection between $\big\la u,\psi_{\aal,n-1,\bk}\big\ra$ and $\big\la\p^{\aal}u,\varphi_{\aal,n-1,\bk}\big\ra$.

\begin{proposition}\label{Prop4} Let a tensor product framelet function $\psi_{\aal}\in L_2(\R^2)$ have vanishing moments of order $\aal$ with $|\aal|\leq s$, and let $\su(\psi_{\aal})=[a_1,a_2]\times[b_1,b_2]$. For $n\in\N$ and $\bk\in\Z^2$ with $\su(\psi_{\aal,n-1,\bk})\subseteq\overline\Om$, we have
\begin{align}\label{WCMeaning}
\big\la u,\psi_{\aal,n-1,\bk}\big\ra=(-1)^{|\aal|}2^{|\aal|(1-n)}\big\la\p^{\aal}u,\varphi_{\aal,n-1,\bk}\big\ra
\end{align}
for every $u\in W_1^s(\Om)$.
\end{proposition}

\begin{pf} By Proposition \ref{Prop1}, there exists the unique $\varphi_{\aal}$ corresponding to $\psi_{\aal}$ such that $\int_{\R^2}\varphi_{\aal}\rd\x\neq0$, $\su(\varphi_{\aal})=\su(\psi_{\aal})$, and $\psi_{\aal}=\p^{\aal}\varphi_{\aal}$ a.e. Then by the chain rule
\begin{align*}
\p^{\aal}\varphi_{\aal,n-1,\bk}=2^{(|\aal|+1)(n-1)-1}\psi_{\aal}(2^{n-1}\cdot-2^{-1}\bk)=2^{|\aal|(n-1)}\psi_{\aal,n-1,\bk}
\end{align*}
where $\varphi_{\aal,n-1,\bk}$ and $\psi_{\aal,n-1,\bk}$ are defined as in \eqref{Frameletn-1}. This means that
\begin{align*}
\big\la u,\psi_{\aal,n-1,\bk}\big\ra=2^{|\aal|(1-n)}\big\la u,\p^{\aal}\varphi_{\aal,n-1,\bk}\big\ra.
\end{align*}
The proof is completed by the integration by parts formula \cite[Proposition 4.2]{J.F.Cai2016}:
\begin{align*}
\big\la u,\p^{\aal}\varphi_{\aal,n-1,\bk}\big\ra=\sum_{\substack{\bbe_l\in\DD_{\aal}\\
1\leq l\leq|\aal|}}(-1)^{l-1}\int_{\p\Om}\msT(\p^{\bbe_l}u)(\p^{\aal-\bbe_{l+1}}\varphi_{\aal,n-1,\bk})\bn_{\bbe_{l+1}-\bbe_l}ds+(-1)^{|\aal|}\big\la\p^{\aal}u,\varphi_{\aal,n-1,\bk}\big\ra.
\end{align*}
Here, $\msT(\p^{\bbe_l}u)=\p^{\bbe_l}u\big|_{\p\Om}$, $\DD_{\aal}$ is the index set defined as
\begin{align*}
\DD_{\aal}=\big\{\bbe_l<\aal:|\bbe_l|=l-1,~\bbe_l<\bbe_{l+1},~\text{for}~l=1,2,\cdots,|\aal|\big\},
\end{align*}
and $\bn_{\bbe}=n_1$ if $\bbe=(1,0)$ and $\bn_{\bbe}=n_2$ if $\bbe=(0,1)$ with $\bn=(n_1,n_2)$ being the outward unit normal of $\p\Om$. Note that every integration on $\p\Om$ vanishes, because from the proof of Proposition \ref{Prop1}, it can be easily verified that $\su(\p^{\bbe}\varphi_{\aal})=\su(\varphi_{\aal})$ for $\0\leq\bbe\leq\aal$, so that $\su(\p^{\bbe}\varphi_{\aal,n-1,\bk})\subseteq\overline\Om$ for $\0\leq\bbe\leq\aal$. Hence,
\begin{align*}
\big\la u,\p^{\aal}\varphi_{\aal,n-1,\bk}\big\ra=(-1)^{|\aal|}\la\p^{\aal}u,\varphi_{\aal,n-1,\bk}\big\ra,
\end{align*}
which completes the proof.\qquad$\square$
\end{pf}

\begin{remark} Recently in \cite{B.Dong2016}, the authors obtained a similar result as Proposition \ref{Prop1} for generic tensor product framelets. However, it was not clear from their analysis that $\su(\varphi_{\aal})=\su(\psi_{\aal})$ as well as the regularity of $\varphi_{\aal}$. Therefore, the conclusion of Proposition \ref{Prop1} is stronger.
\end{remark}

\section{Edge Driven Wavelet Frame Based Image Restoration}\label{WaveletModelAlg}

In this section, we present our edge driven wavelet frame based image restoration model with full details. We also present an alternating optimization algorithm which iteratively updates the image to be recovered and the set of singularities. The proposed model and algorithm are all in discrete settings, where all variables are discrete arrays.

\subsection{Image Restoration Model}\label{ModelAlg}

 We denote by $\OO=\big\{0,1,\cdots,N-1\big\}^2$ the set of indices of the $N\times N$ Cartesian grid which discretizes the domain $\Om=(0,1)^2$. Recall that the space of all two-dimensional array on the grid $\OO$ is denoted as $\mI_2$. Let $\bsA$ be some linear operator mapping $\mI_2$ into itself, so that both the (unknown) true image $\bu$ and the degraded measurement (or the observed image) $\bsf$ are the elements of $\mI_2$.

We propose our wavelet frame based image restoration model as
\begin{align}\label{OurModel}
\min_{\bu,~0\leq\bsv\leq1}~\left\|(\one-\bsv)\cdot\big(\lam\cdot\bsW\bu\big)\right\|_1+\left\|\bsv\cdot\big(\gga\cdot\bsW'\bu\big)\right\|_1+\big\|\rrh\cdot\bsW''\bsv\big\|_1+\f{1}{2}\big\|\bsA\bu-\bsf\big\|_2^2,
\end{align}
where
\begin{align*}
\left\|(\one-\bsv)\cdot\big(\lam\cdot\bsW\bu\big)\right\|_1&=\sum_{\bk\in\OO}\sum_{l=0}^{L-1}\big(\one-\bsv_l[\bk]\big)\left(\sum_{\aal\in\BB}\lambda_{l,\aal}[\bk]\bigg|\big(\bsW_{l,\aal}\bu\big)[\bk]\bigg|^2\right)^{\f{1}{2}},\\
\left\|\bsv\cdot\big(\gga\cdot\bsW'\bu\big)\right\|_1&=\sum_{\bk\in\OO}\sum_{l=0}^{L-1}\bsv_l[\bk]\left(\sum_{\aal\in\BB'}\gamma_{l,\aal}[\bk]\bigg|\big(\bsW_{l,\aal}'\bu\big)[\bk]\bigg|^2\right)^{\f{1}{2}},\\
\big\|\rrh\cdot\bsW''\bsv\big\|_1&=\sum_{l=0}^{L-1}\underbrace{\sum_{\bk\in\OO}\sum_{m=0}^{L''-1}\left(\sum_{\aal\in\BB''}\rho_{l,m,\aal}[\bk]\bigg|\big(\bsW_{m,\aal}''\bsv_l\big)[\bk]\bigg|^2\right)^{\f{1}{2}}}_{:=\big\|\rrh_l\cdot\bsW''\bsv_l\big\|_1},
\end{align*}
and $\BB$, $\BB'$, and $\BB''$ denote the framelet bands of $\bsW$, $\bsW'$, and $\bsW''$ respectively:
\begin{align*}
\BB&=\big\{0,1,\cdots,r\big\}^2\setminus\big\{\0\big\},\\
\BB'&=\big\{0,1,\cdots,r'\big\}^2\setminus\big\{\0\big\},\\
\BB''&=\big\{0,1,\cdots,r''\big\}^2\setminus\big\{\0\big\}.
\end{align*}
To better understand the proposed model \eqref{OurModel}, we observe that it can be regarded as a relaxation of the following model:
\begin{align}\label{OurModelOriginal}
\min_{\bu,\Sig}~\left\|\big(\lam\cdot\bsW\bu\big)_{\Sig^c}\right\|_1+\left\|\big(\gga\cdot\bsW'\bu\big)_{\Sig}\right\|_1+\big\|\rrh\cdot\bsW''\one_{\Sig}\big\|_1+\f{1}{2}\big\|\bsA\bu-\bsf\big\|_2^2
\end{align}
where $\Sig=\big(\Sig_0,\cdots,\Sig_{L-1}\big)$ with $\Sig_l$ being the estimated singularity region for $l=0,\cdots,L-1$, which will be denoted as the $(l+1)$st level singularity in what follows, and $\one_{\Sig}=\big(\one_{\Sig_0},\cdots,\one_{\Sig_{L-1}}\big)$ with $\one_{\Sig_l}$ being the labelling binary image of $\Sig_l$: $\one_{\Sig_l}[\bk]=1$ if $\bk\in\Sig_l$, and $0$ otherwise. The first two terms in \eqref{OurModelOriginal} are defined as
\begin{align*}
\left\|\big(\lam\cdot\bsW\bu\big)_{\Sig^c}\right\|_1&=\sum_{l=0}^{L-1}\sum_{\bk\in\OO\setminus\Sig_l}\left(\sum_{\aal\in\BB}\lambda_{l,\aal}[\bk]\bigg|\big(\bsW_{l,\aal}\bu\big)[\bk]\bigg|^2\right)^{\f{1}{2}},\\
\left\|\big(\gga\cdot\bsW'\bu\big)_{\Sig}\right\|_1&=\sum_{l=0}^{L-1}\sum_{\bk\in\Sig_l}\left(\sum_{\aal\in\BB'}\gamma_{l,\aal}[\bk]\bigg|\big(\bsW_{l,\aal}'\bu\big)[\bk]\bigg|^2\right)^{\f{1}{2}}
\end{align*}
respectively.

Comparing \eqref{OurModel} to \eqref{OurModelOriginal}, we can see that the first term restores the smooth regions of image, while the second term preserves the singularities, and the third term provides the regularization on the singularities to enhance sharp image features. In other words, our model takes  different regularization in smooth image regions and near image singularities such as edges, and actively restores sharp image features at the same time. However, since the first two terms are exchangeable, an appropriate choice of the wavelet frame transforms as well as the associated parameters is necessary to enforce desired effects with the two terms. From Propositions \ref{Prop1} and \ref{Prop4}, we can see that image singularities (i.e. jumps and jumps after lower order differentiations) can be well captured by framelets of lower order vanishing moments. In our model \eqref{OurModel}, $\bsW$ consists of filters whose vanishing moments of the highest order is higher than those of $\bsW'$ (i.e. $r>r'$). Besides, since the magnitudes of the wavelet frame coefficients have to be as small as possible in smooth image regions, we choose the parameters so that $\lam$ is overall larger than $\gga$.

Compared to the existing models \eqref{CDS} and \eqref{JLS} which also treat images as piecewise smooth functions, our model \eqref{OurModel} uses $\ell_1$ norm to promote smoothness rather than $\ell_2$ norm. By doing so, we can better protect the singularities that are not captured by $\bsv$ than using the $\ell_2$ norm which can smear these singularities out. This leads to a  more robust  estimation of the singularities of the unknown true image from the degraded measurement than \eqref{CDS} and \eqref{JLS}. In addition, unlike \eqref{CDS} and \eqref{JLS} which explicitly takes the singularity set into account, our model adopts an implicit representation of the singularity set by relaxing the binary image $\one_{\Sig}$ into $\bsv$ taking values in $[0,1]$. This relaxation allows an overlap between the smooth and the sharp image regions in the transform domain, which will be helpful to suppress the staircase effects near the interface. Furthermore, as will be rigorously analyzed in Section \ref{VariationAsymAnal}, this implicit representation of the singularity set enables us to provide an asymptotic analysis of the model with respect to both $u$ and $v$, in contrast to that of \eqref{CDS} where the singularity set is assumed to be fixed.

We would like to mention that our model mainly focuses on the restoration of images which can be well approximated by piecewise smooth functions. Therefore, our model may not be suitable for images having textures. Indeed, textures can be sparsely approximated by systems with oscillating patterns such as local cosine systems \cite{J.F.Cai2016,Meyer2001}, rather than piecewise smooth functions. However, we can easily modify the proposed model by adopting the idea of a two system model (e.g. \cite{J.F.Cai2010,J.F.Cai2012,J.F.Cai2009/10,B.Dong2012,B.Dong2013,M.Elad2005,J.L.Starck2005}) to better handle images with textures. Nonetheless, we will not discuss details on such variant of our model, as it is beyond the scope of this paper. We will focus on recovering images that are piecewise smoothness.

\subsection{Algorithm for Image Restoration Model}\label{Algorithm}

\begin{algorithm}[ht!]
\textbf{Step 0.} $\bu^0$, $\bsv^0$\;
\For{$k=0,1,2,\cdots$}
{\textbf{Step 1.} Given $\bsv^k$, solve
\begin{align}\label{usubprob}
\bu^{k+1}=\arg\min_{\bu}~\left\|\big(\one-\bsv^k\big)\cdot\big(\lam\cdot\bsW\bu\big)\right\|_1+\left\|\bsv^k\cdot\big(\gga\cdot\bsW'\bu\big)\right\|_1+\f{1}{2}\big\|\bsA\bu-\bsf\big\|_2^2.
\end{align}
\textbf{Step 2.} Given $\bu^{k+1}$, solve
\begin{align}\label{vsubprob}
\bsv^{k+1}=\arg\min_{0\leq\bsv\leq1}~\left\|(\one-\bsv)\cdot\big(\lam\cdot\bsW\bu^{k+1}\big)\right\|_1+\left\|\bsv\cdot\big(\gga\cdot\bsW'\bu^{k+1}\big)\right\|_1+\big\|\rrh\cdot\bsW''\bsv\big\|_1.
\end{align}}
\caption{Alternating Minimization Algorithm for \eqref{OurModel}}\label{Alg1}
\end{algorithm}

The proposed alternating minimization algorithm for \eqref{OurModel} is given by Algorithm \ref{Alg1}. To solve the $\bu$ subproblem \eqref{usubprob}, we use the split Bregman algorithm \cite{J.F.Cai2009/10,J.Eckstein1992,T.Goldstein2009}, which is a widely used method for solving various convex sparse optimization problems in variational image restoration. For completeness, we present the full details of the split Bregman algorithm solving the subproblem \eqref{usubprob} as follows: let $\bsd_1^0=\bb_1^0=\bsd_2^0=\bb_2^0=\0$. For $j=0,1,2,\cdots$
\begin{align}\label{usubprobSB}
\begin{split}
\bu^{j+1}&=\arg\min_{\bu}~\f{1}{2}\big\|\bsA\bu-\bsf\big\|_2^2+\f{\mu_1}{2}\big\|\bsW\bu-\bsd_1^j+\bb_1^j\big\|_2^2+\f{\mu_2}{2}\big\|\bsW'\bu-\bsd_2^j+\bb_2^j\big\|_2^2\\
\bsd_1^{j+1}&=\arg\min_{\bsd_1}~\left\|\big(\one-\bsv\big)\cdot\big(\lam\cdot\bsd_1\big)\right\|_1+\f{\mu_1}{2}\big\|\bsd_1-\bsW\bu^{j+1}-\bb_1^j\big\|_2^2\\
\bsd_2^{j+1}&=\arg\min_{\bsd_2}~\left\|\bsv\cdot\big(\gga\cdot\bsd_2\big)\right\|_1+\f{\mu_2}{2}\big\|\bsd_2-\bsW'\bu^{j+1}-\bb_2^j\big\|_2^2\\
\bb_1^{j+1}&=\bb_1^j+\bsW\bu^{j+1}-\bsd_1^{j+1}\\
\bb_2^{j+1}&=\bb_2^j+\bsW'\bu^{j+1}-\bsd_2^{j+1},
\end{split}
\end{align}
where we omit the outer iteration superscript $k$ for notational simplicity. Note that each of the subproblem of \eqref{usubprobSB} has a closed-form solution and it can be rewritten as
\begin{align}\label{usubprobexplicit}
\begin{split}
\bu^{j+1}&=\left[\bsA^T\bsA+\big(\mu_1+\mu_2\big)\bsI\right]^{-1}\left[\bsA^T\bsf+\mu_1\bsW^T\big(\bsd_1^j-\bb_1^j\big)+\mu_2\big(\bsW'\big)^T\big(\bsd_2^j-\bb_2^j\big)\right]\\
\bsd_1^{j+1}&=\mT_{(\one-\bsv)\cdot\lam/\mu_1}\big(\bsW\bu^{j+1}+\bb_1^j\big)\\
\bsd_2^{j+1}&=\mT_{\bsv\cdot\gga/\mu_2}\big(\bsW'\bu^{j+1}+\bb_2^j\big)\\
\bb_1^{j+1}&=\bb_1^j+\bsW\bu^{j+1}-\bsd_1^{j+1}\\
\bb_2^{j+1}&=\bb_2^j+\bsW'\bu^{j+1}-\bsd_2^{j+1}.
\end{split}
\end{align}
Here, the isotropic shrinkage $\mT_{\bsv\cdot\lam}$ is defined as
\begin{align*}
\bigg(\mT_{\bsv\cdot\lam}\big(\w\big)\bigg)_{l,\aal}[\bk]=\left\{\begin{array}{ll}
\w_{l,\aal}[\bk],&\aal=\0,\vspace{0.5em}\\
\f{\w_{l,\aal}[\bk]}{R_l[\bk]}\max\big\{R_l[\bk]-\bsv_l[\bk]\lam_{l,\aal}[\bk],0\big\},&\aal\in\BB,
\end{array}\right.
\end{align*}
with $R_l[\bk]=\left(\sum_{\aal\in\BB}\big|\w_{l,\aal}[\bk]\big|^2\right)^{1/2}$ for $l=0,\cdots,L-1$ and $\bk\in\OO$. If $\bsv\equiv\one$, then we write it as $\mT_{\lam}$.

The subproblem \eqref{vsubprob} for variable $\bsv$ can be reformulated as
\begin{align*}
\min_{0\leq\bsv\leq1}~\sum_{l=0}^{L-1}\bigg(\big\la\one-\bsv_l,\bsg_{1,l}\big\ra+\big\la\bsv_l,\bsg_{2,l}\big\ra+\big\|\rrh_l\cdot\bsW''\bsv_l\big\|_1\bigg)
\end{align*}
where $\bsg_{1,l}$ and $\bsg_{2,l}$ for $l=0,\cdots,L-1$ are respectively defined as
\begin{align*}
\bsg_{1,l}[\bk]&=\left(\sum_{\aal\in\BB}\lambda_{l,\aal}[\bk]\bigg|\big(\bsW_{l,\aal}\bu\big)[\bk]\bigg|^2\right)^{\f{1}{2}}\\
\bsg_{2,l}[\bk]&=\left(\sum_{\aal\in\BB'}\gamma_{l,\aal}[\bk]\bigg|\big(\bsW_{l,\aal}'\bu\big)[\bk]\bigg|^2\right)^{\f{1}{2}}.
\end{align*}
This subproblem can also be solved using the split Bregman algorithm. Since each $\bsv_0,\cdots,\bsv_{L-1}$ can be computed separately in the same way, we omit the subscript $l$ and the outer iteration superscript $k$. The algorithm solving the subproblem \eqref{vsubprob} is as follows: let $\bsd^0=\bb^0=\0$. For $j=0,1,2,\cdots$
\begin{align}\label{vsubprobSB}
\begin{split}
\bsv^{j+\f{1}{2}}&=\arg\min_{\bsv}~\big\la\one-\bsv,\bsg_1\big\ra+\big\la\bsv,\bsg_2\big\ra+\f{\mu}{2}\big\|\bsW''\bsv-\bsd^j+\bb^j\big\|_2^2\\
\bsv^{j+1}&=\min\left\{\max\big(\bsv^{j+\f{1}{2}},0\big),1\right\}\\
\bsd^{j+1}&=\arg\min_{\bsd}~\big\|\rrh\cdot\bsd\big\|_1+\f{\mu}{2}\big\|\bsd-\bsW''\bsv^{j+1}-\bb^j\big\|_2^2\\
\bb^{j+1}&=\bb^j+\bsW''\bsv^{j+1}-\bsd^{j+1}.
\end{split}
\end{align}
Note that each step of \eqref{vsubprobSB} has a closed-form solution. Thus, \eqref{vsubprobSB} can be rewritten as
\begin{align}\label{vsubprobexplicit}
\begin{split}
\bsv^{j+\f{1}{2}}&=\big(\bsW''\big)^T\big(\bsd^j-\bb^j)+\big(\bsg_1-\bsg_2\big)/\mu\\
\bsv^{j+1}&=\min\left\{\max\big(\bsv^{j+\f{1}{2}},0\big),1\right\}\\
\bsd^{j+1}&=\mT_{\rrh/\mu}\big(\bsW''\bsv^{j+1}+\bb^j\big)\\
\bb^{j+1}&=\bb^j+\bsW''\bsv^{j+1}-\bsd^{j+1}.
\end{split}
\end{align}
From the reconstructed $\bsv$, we obtain the estimated $(l+1)$st level singularity by $\Sig_l=\big\{\bk\in\OO:\bsv_l>t_l\big\}$ for $t_l\in[0,1]$ and $l=0,\cdots,L-1$. In our numerical simulations, we set $t_l=t=0.5$ for $0\leq l\leq L-1$.


\subsection{Numerical Results}\label{SimulationResults}

In this subsection, we conduct some numerical simulations on image inpainting and image deblurring using Algorithm \ref{Alg1}. In all of the numerical simulations, we will use the piecewise cubic B-spline wavelet frame for $\bsW$, and the piecewise linear B-spline for $\bsW'$ and $\bsW''$. The levels of decomposition, i.e. $L$ and $L''$ are chosen differently depending on the image restoration problems. We compare the results obtained from our proposed model \eqref{OurModel} with the piecewise smooth (PS) model \eqref{CDS} in \cite{J.F.Cai2016}, and the geometric structure (GS) model \eqref{JLS} in \cite{H.Ji2016}. We also compare with the total generalized variation (TGV) model \cite{K.Bredies2010}:
\begin{align}\label{TGV}
\min_{u,w}~\alpha\left(\big\|\na u-w\big\|_{L_1(\Om)}+\beta\big\|\na_s w\big\|_{L_1(\Om)}\right)+\f{1}{2}\big\|Au-f\big\|_{L_2(\Om)}^2
\end{align}
which is solved by the modified primal-dual hybrid gradient method \cite{A.Chambolle2011,E.Esser2010}. Here, $\na_s=\f{1}{2}\big(\na+\na^T\big)$, and we use forward difference with periodic boundary condition to discretize \eqref{TGV}.

In all image restoration problems, the true image $\bu$ takes the integer values in $[0,255]$. For the image inpainting, $\bsA=\one_{\La}$ with a known $\La\subsetneq\OO$ and the measurement $\bsf$ is designed as
\begin{align*}
\bsf[\bk]=\left\{\begin{array}{cl}
\bu[\bk]+\eet[\bk],~&\bk\in\La,\vspace{0.5em}\\
\text{arbitrary},~&\bk\notin\La.
\end{array}\right.
\end{align*}
In particular, we focus on the task of removing texts and scratches. For the image deblurring, $\bsA$ is taken to be the convolution operator with the kernel generated in MATLAB by ``fspecial(`gaussian',2,15)''. In any case, the additive noise $\eet$ with standard deviation $4$ is also added. For the quantitative comparison on each model, we calculate the peak signal to noise ratio (PSNR) value defined by
\begin{align*}
\mathrm{PSNR}:=-20\log_{10}\f{\|\bu-\wt{\bu}\|_2}{255N}
\end{align*}
where $\wt{\bu}$ is the recovered image.

\subsubsection{Image Inpainting}\label{SimulationInpaint}

\begin{figure}[ht]
\begin{center}
\begin{minipage}{4.2cm}
\includegraphics[width=4.2cm]{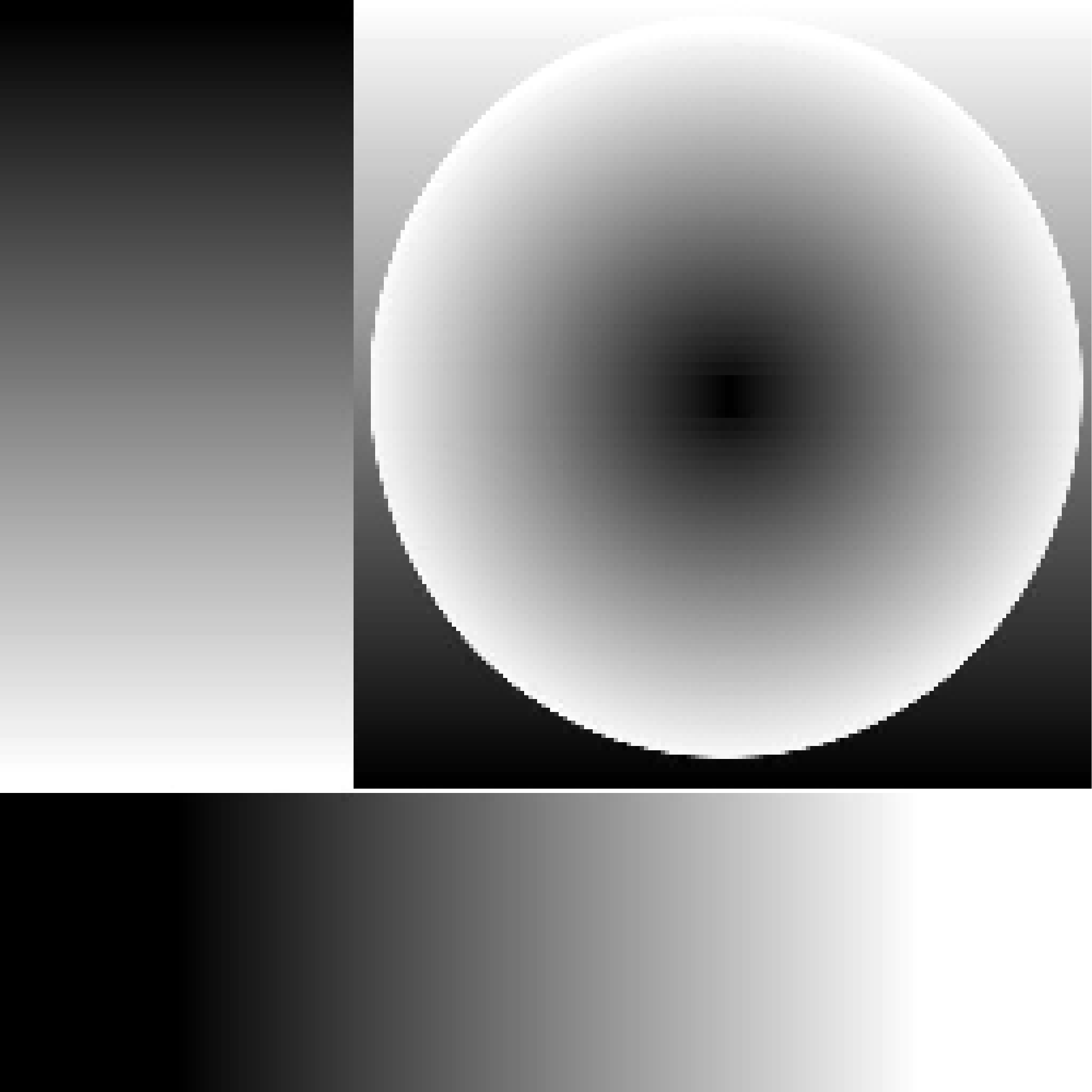}
\end{minipage}
\begin{minipage}{4.2cm}
\includegraphics[width=4.2cm]{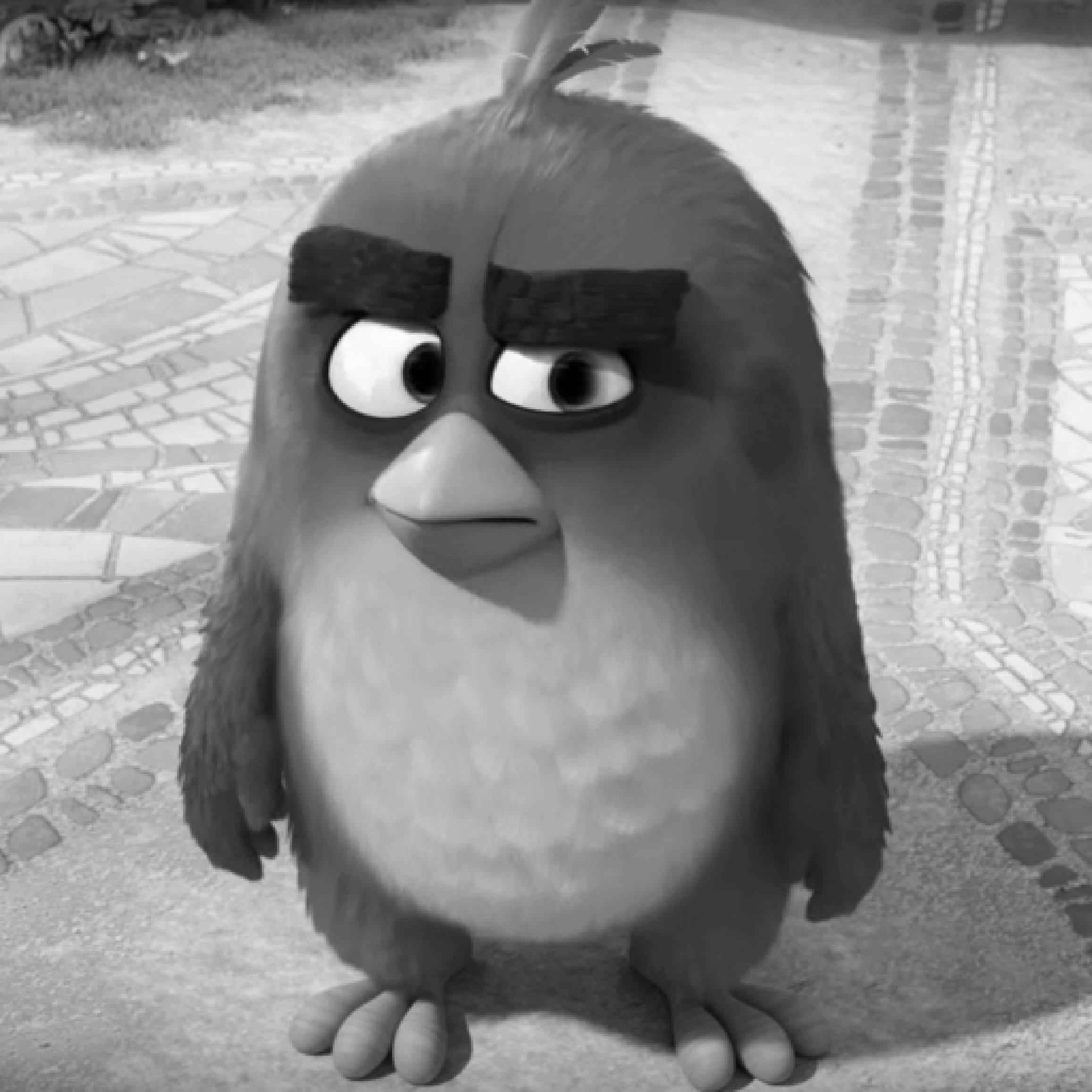}
\end{minipage}
\begin{minipage}{4.2cm}
\includegraphics[width=4.2cm]{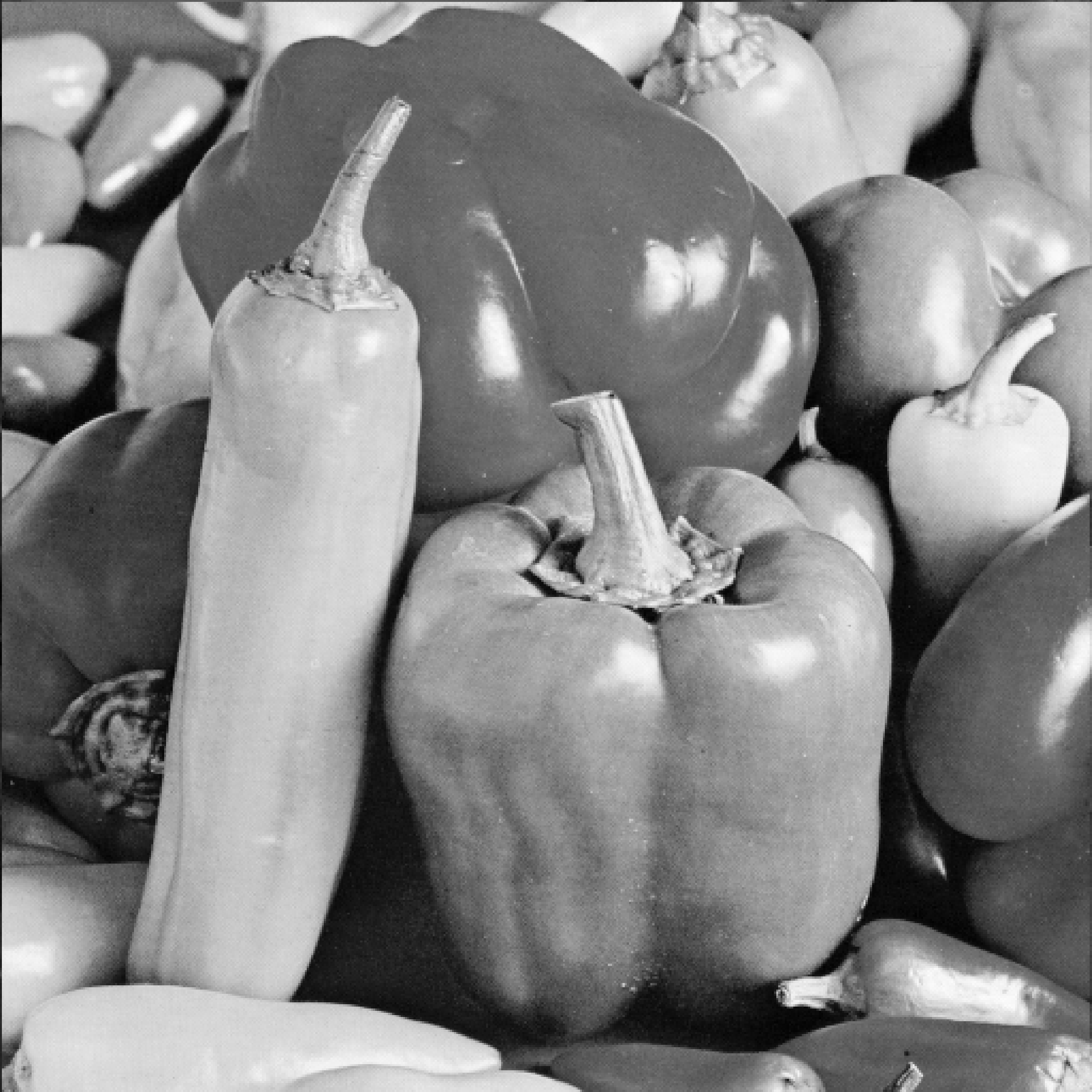}
\end{minipage}\vspace{0.25em}\\
\begin{minipage}{4.2cm}
\includegraphics[width=4.2cm]{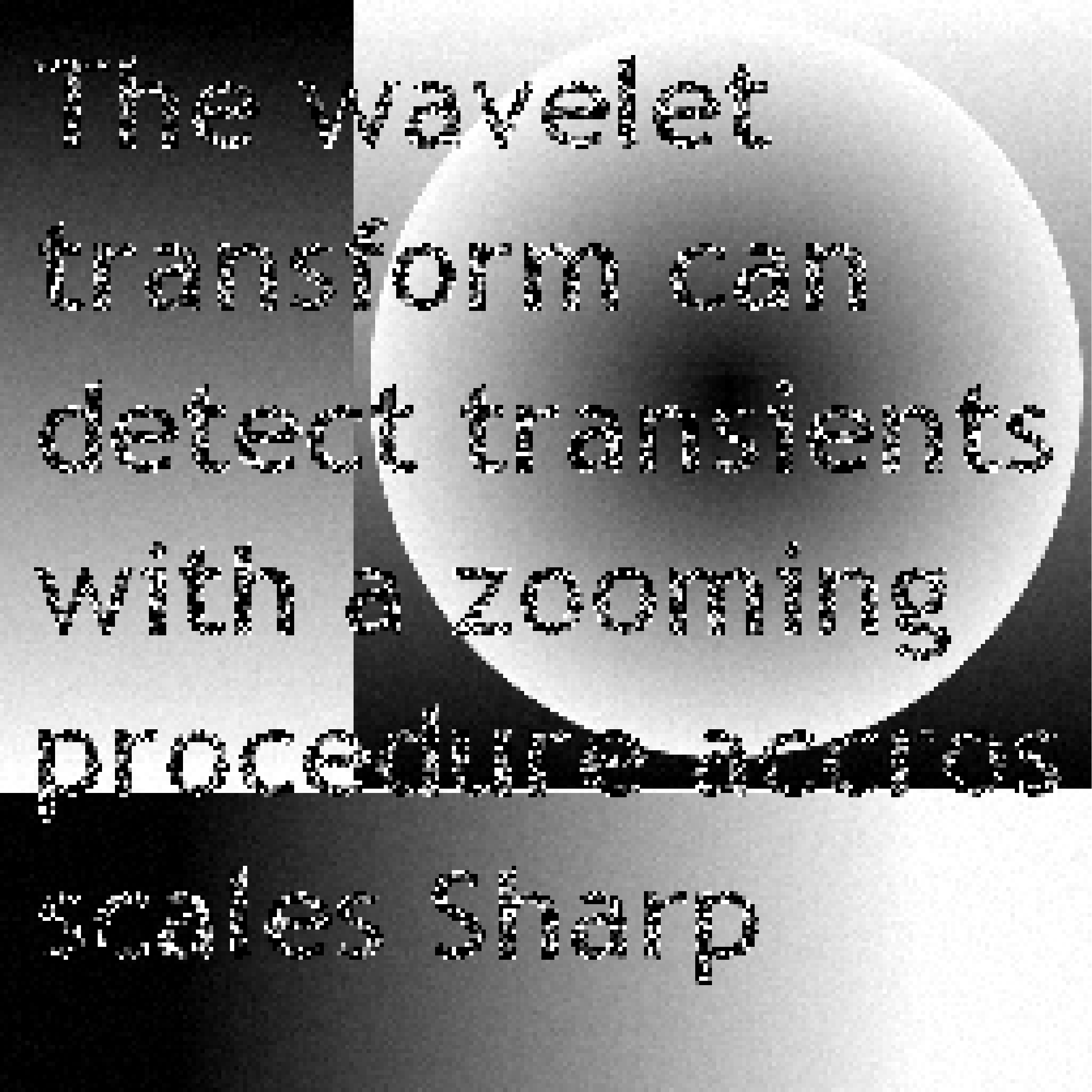}
\end{minipage}
\begin{minipage}{4.2cm}
\includegraphics[width=4.2cm]{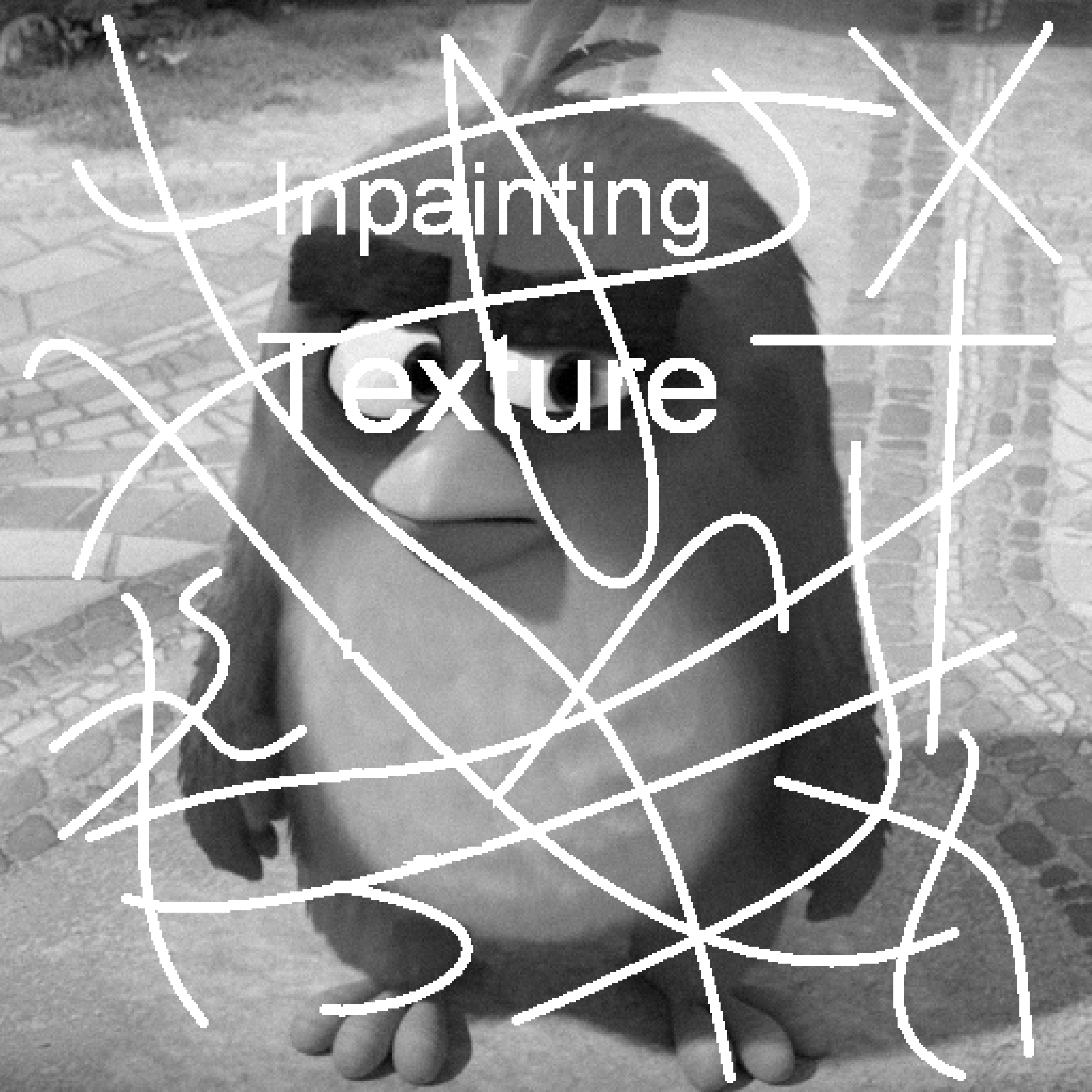}
\end{minipage}
\begin{minipage}{4.2cm}
\includegraphics[width=4.2cm]{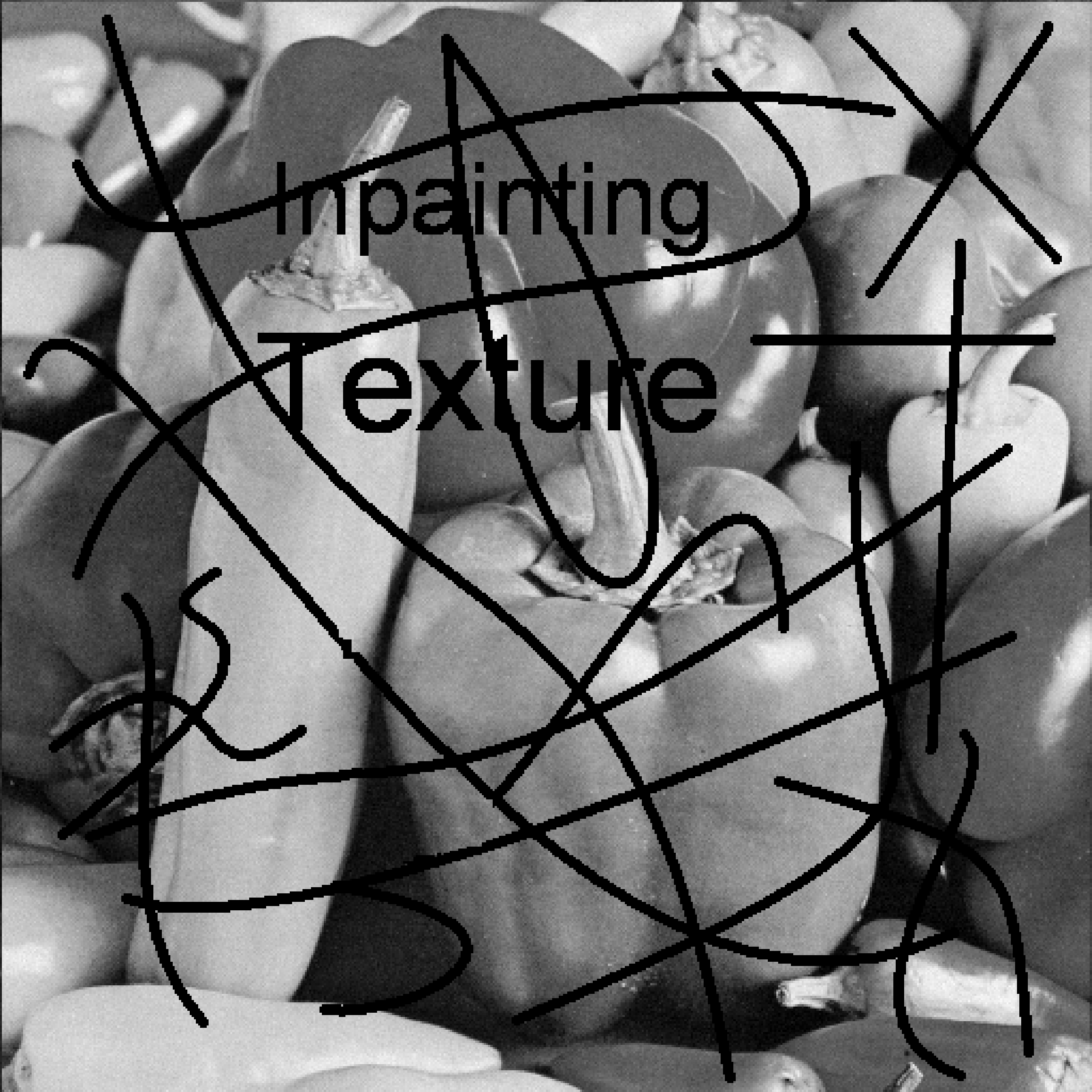}
\end{minipage}\vspace{0.25em}\\
\begin{minipage}{4.2cm}\begin{center}{\small{Slope}}\end{center}\end{minipage}\begin{minipage}{4.2cm}\begin{center}{\small{Angry Birds}}\end{center}\end{minipage}\begin{minipage}{4.2cm}\begin{center}{\small{Peppers}}\end{center}\end{minipage}
\caption{Visualization of original images and the observed images. Throughout this paper, all figures are shown in the window level $[0,255]$ for the fair comparison.}\label{fig:InpaintOriginalMeasure}
\end{center}
\end{figure}

For image inpainting, we test three images as shown in Figure \ref{fig:InpaintOriginalMeasure}, which will be denoted as ``Slope'', ``Angry Birds'', and ``Peppers'' respectively. We initialize our algorithm by choosing $\bu^0=\0$ and $\bsv^0=\0$. The level of decomposition for $\bsW$ and $\bsW'$ is chosen to be $1$. For $\bsW''$, the level of decomposition is chosen to be $4$. For the PS model \eqref{CDS} and the GS model \eqref{JLS}, we use the piecewise linear B-spline wavelet frame with $1$ level of decomposition for ``Slope'', and the piecewise cubic B-spline wavelet frame with $1$ level of decomposition for the others. The parameters $\lam$, $\gga$, $\rrh$ in our model \eqref{OurModel} are chosen as $\lambda_{l,\aal}=\lambda$, $\gamma_{l,\aal}=\gamma$, and $\rho_{l,m,\aal}=\rho$. In addition, the parameters in the PS model \eqref{CDS}, the GS model \eqref{JLS}, and the TGV model \eqref{TGV} as well as our model \eqref{OurModel} are manually chosen to achieve optimal results. (Empirically, we observe that choosing parameters in our model \eqref{OurModel} so that $\big\|\lam\big\|_1>\big\|\gga\big\|_1$ is a good choice.)

\begin{table}[ht]
\begin{center}
\begin{tabular}{|c||c|c|c|c|c|}\hline
Image&Observed&TGV Model \cite{K.Bredies2010}&PS Model \cite{J.F.Cai2016}&GS Model \cite{H.Ji2016}&Our Model \eqref{OurModel}\\ \hline
Slope&$13.8916$&$32.8187$&$30.6408$&$31.4253$&$\textbf{33.7157}$\\ \hline
Angry Birds&$14.1856$&$35.3697$&$35.3974$&$35.3537$&$\textbf{36.0355}$\\ \hline
Peppers&$14.8327$&$34.1675$&$34.0387$&$34.0219$&$\textbf{34.4252}$\\ \hline
\end{tabular}
\caption{Comparison of the PSNR values of four models for inpainting.}\label{tab:PSNRCompareInpaint}
\end{center}
\end{table}

\begin{figure}[ht]
\begin{center}
\begin{minipage}{3.15cm}
\includegraphics[width=3.15cm]{SlopeCorrupted.pdf}
\end{minipage}
\begin{minipage}{3.15cm}
\includegraphics[width=3.15cm]{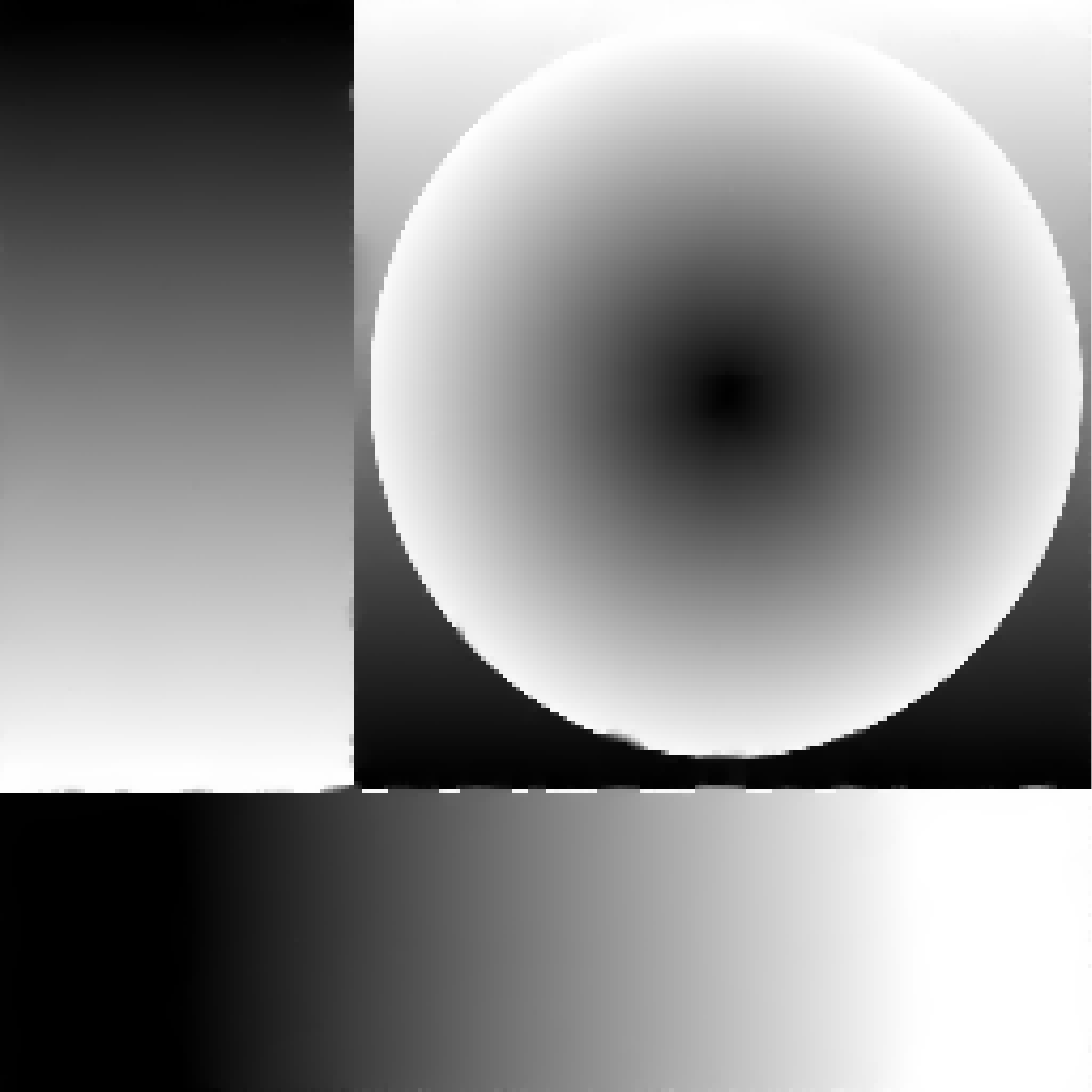}
\end{minipage}
\begin{minipage}{3.15cm}
\includegraphics[width=3.15cm]{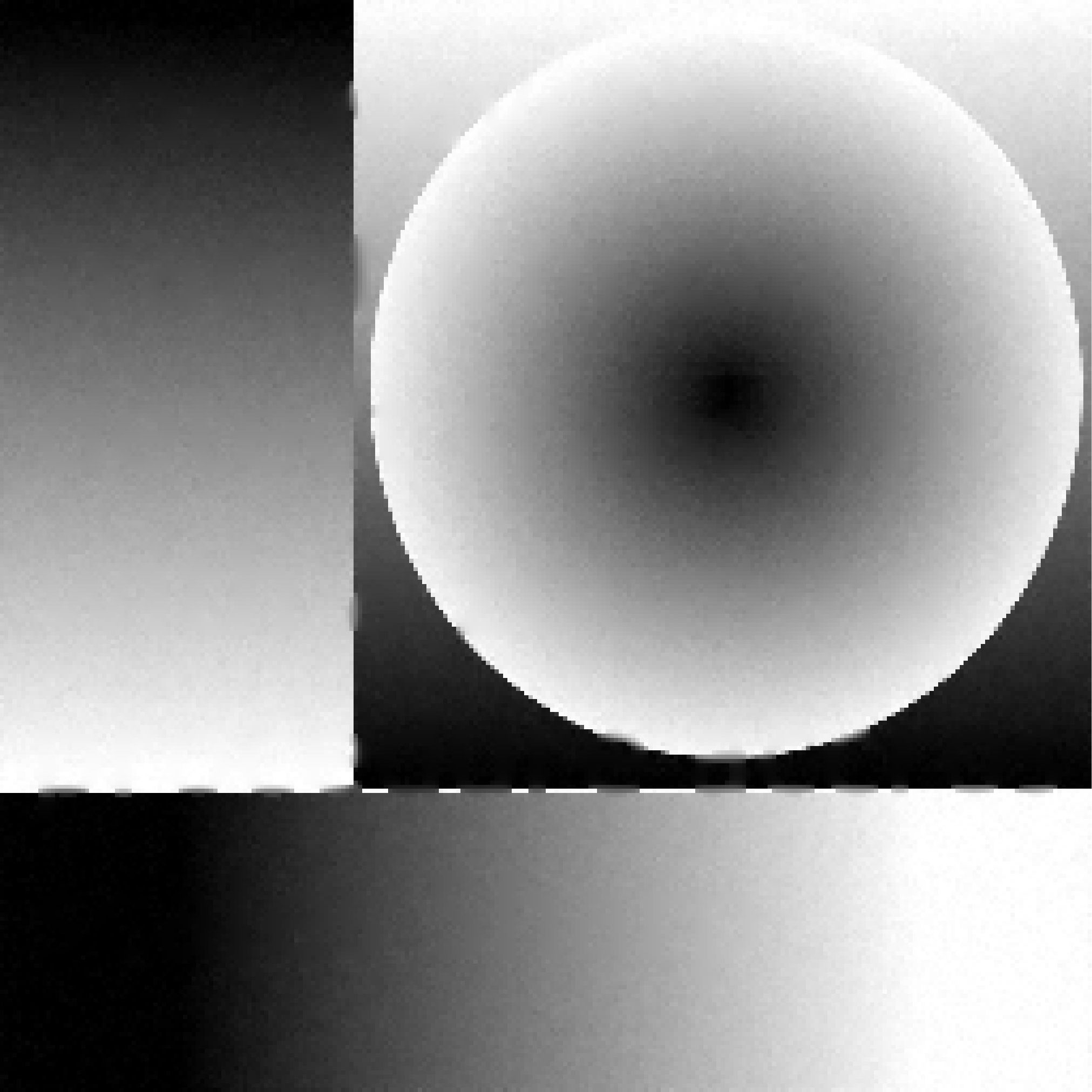}
\end{minipage}
\begin{minipage}{3.15cm}
\includegraphics[width=3.15cm]{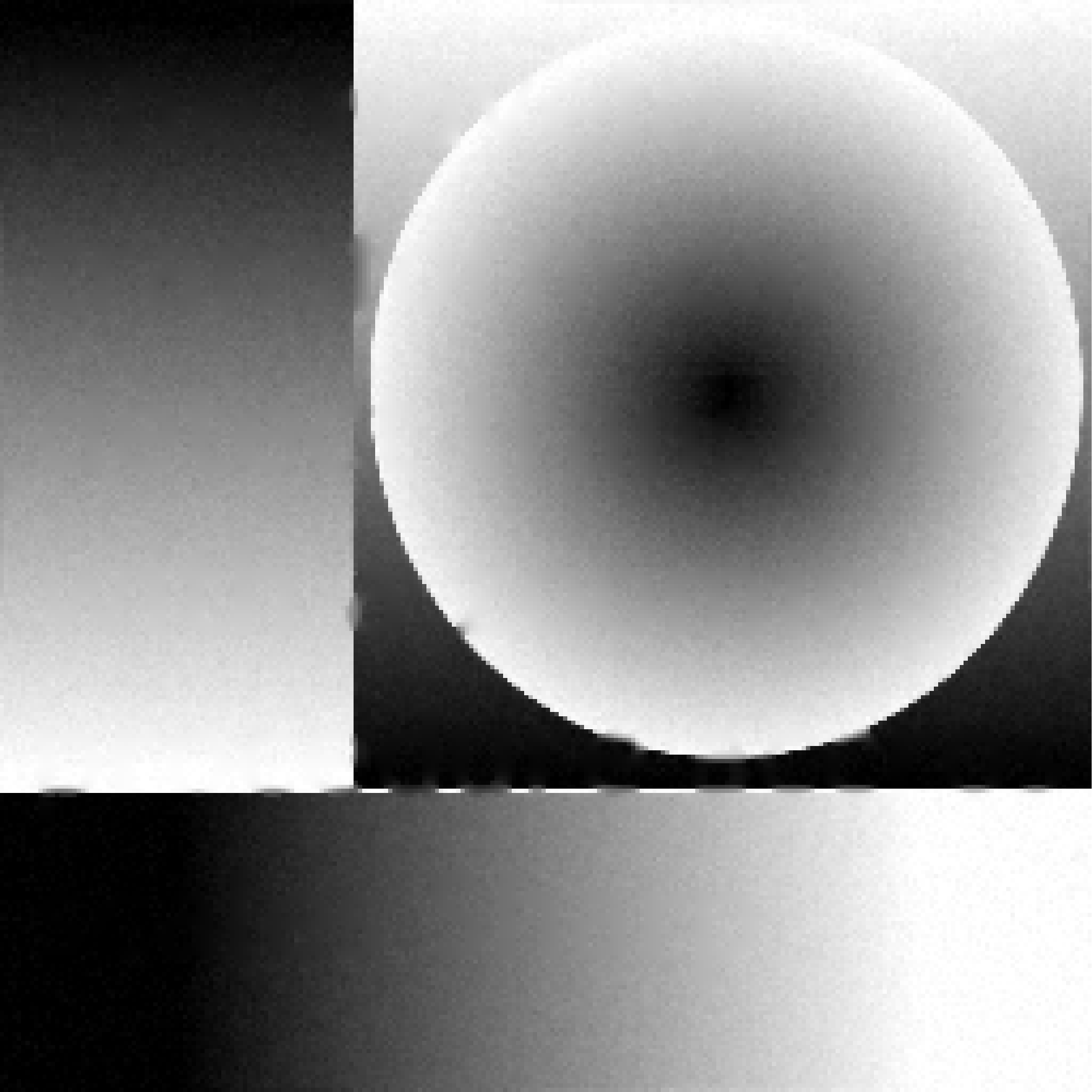}
\end{minipage}
\begin{minipage}{3.15cm}
\includegraphics[width=3.15cm]{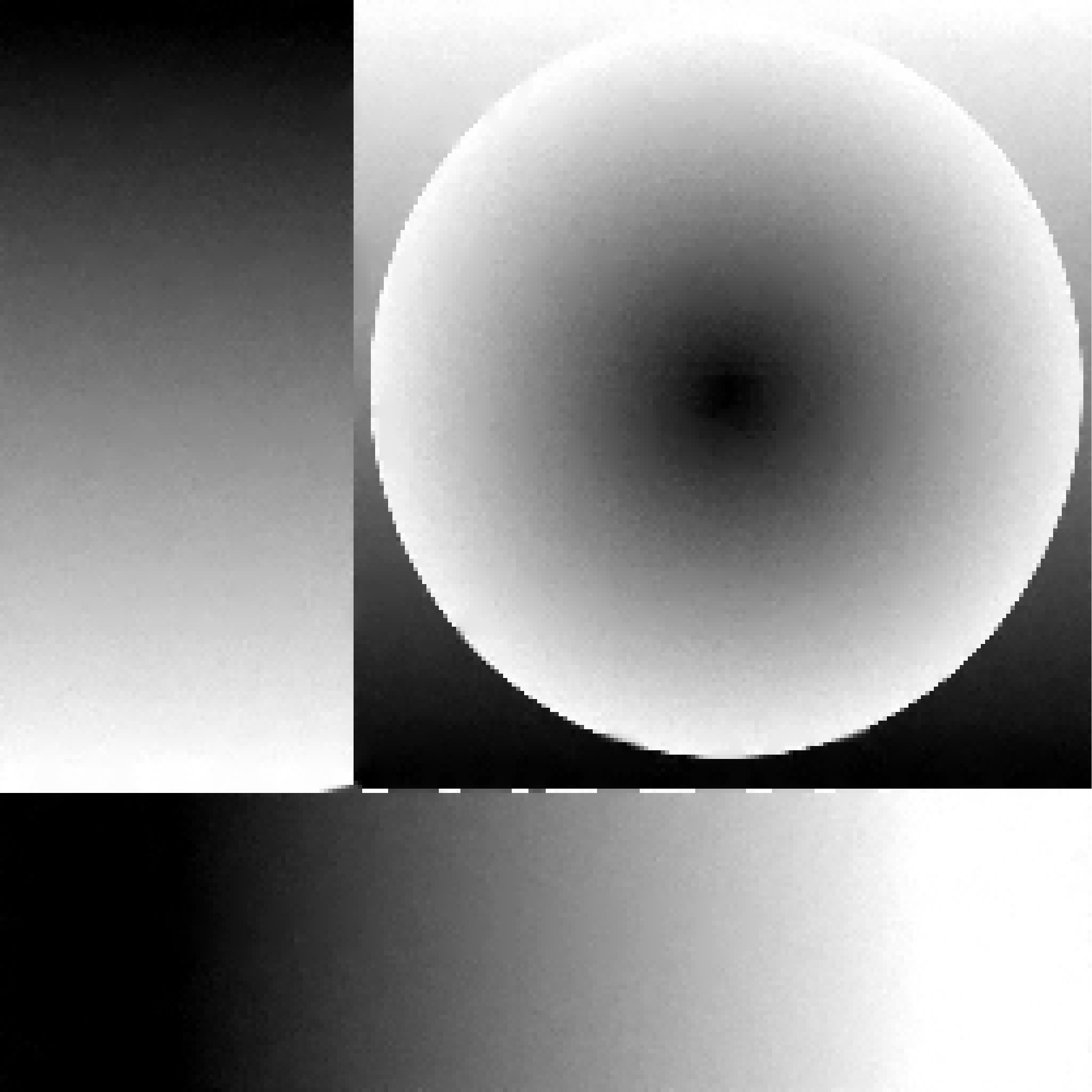}
\end{minipage}\vspace{0.25em}\\
\begin{minipage}{3.15cm}
\includegraphics[width=3.15cm]{AngryBirdsCorrupted.pdf}
\end{minipage}
\begin{minipage}{3.15cm}
\includegraphics[width=3.15cm]{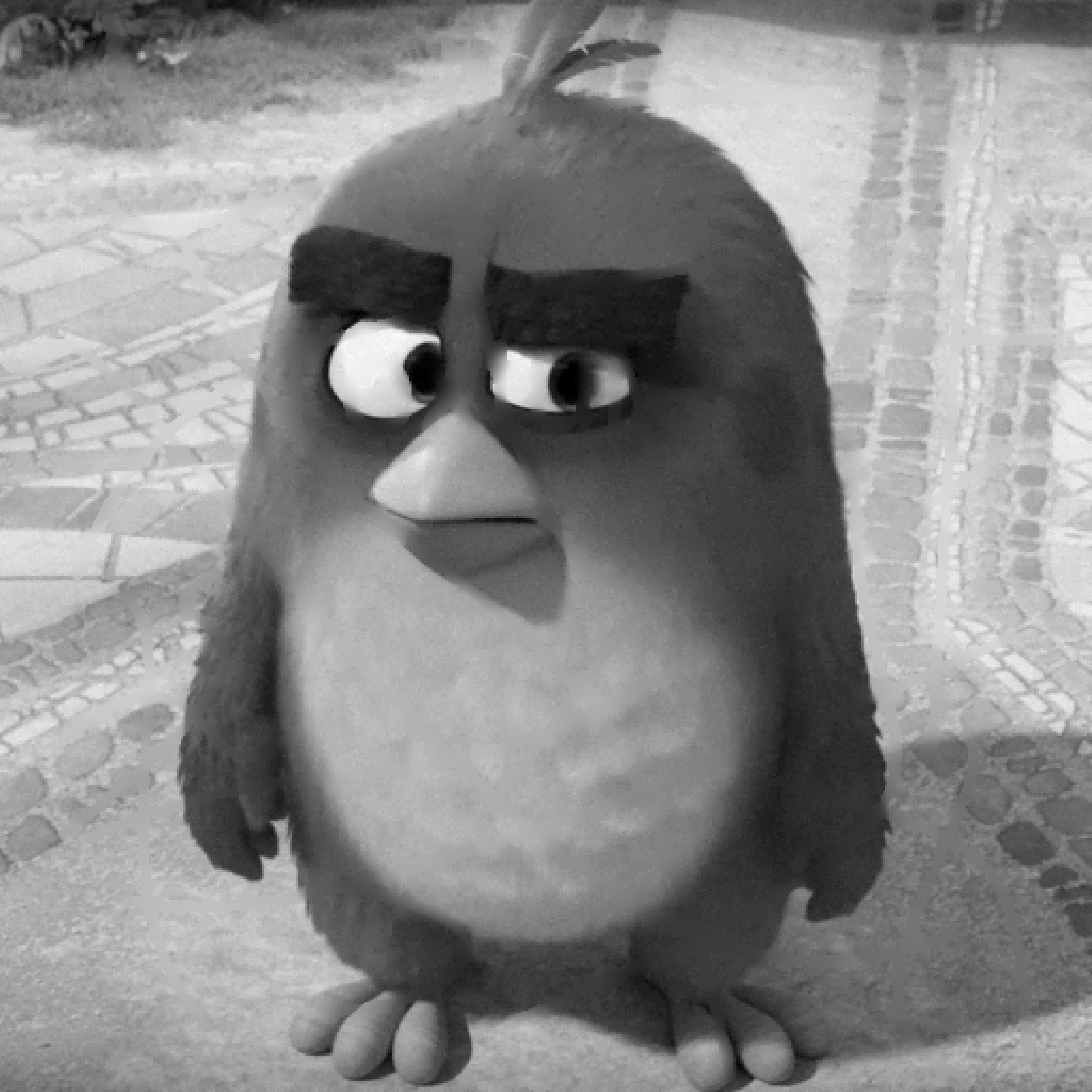}
\end{minipage}
\begin{minipage}{3.15cm}
\includegraphics[width=3.15cm]{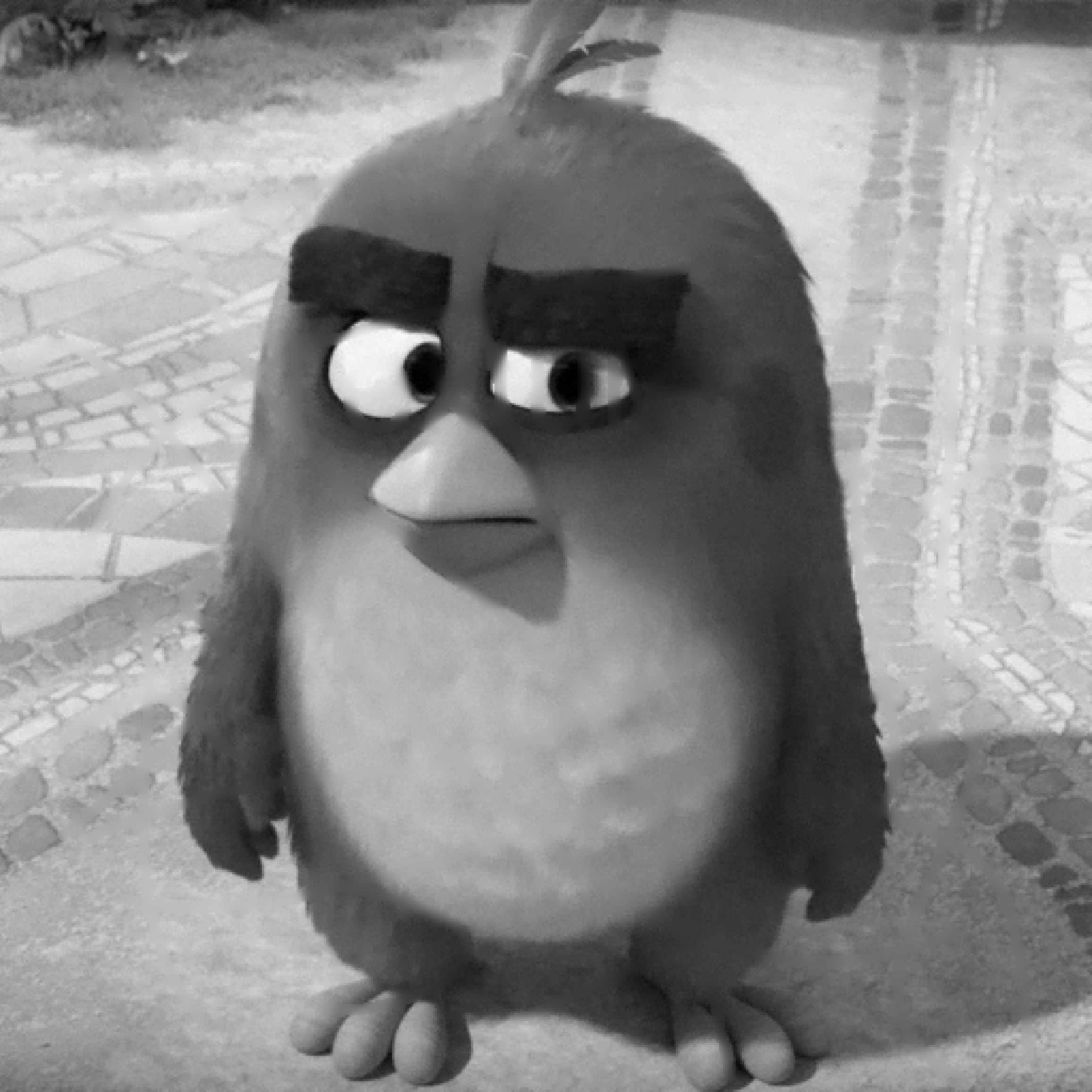}
\end{minipage}
\begin{minipage}{3.15cm}
\includegraphics[width=3.15cm]{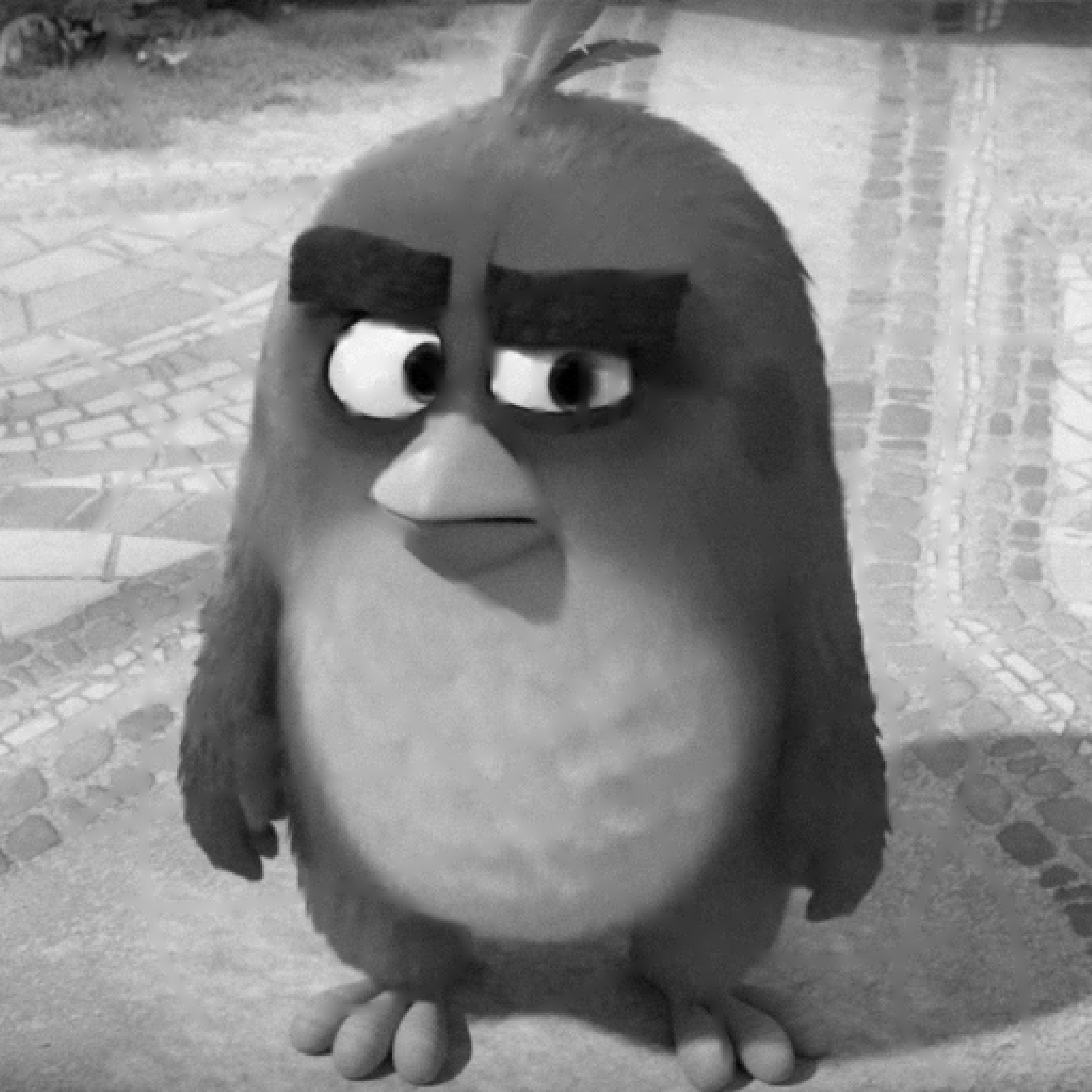}
\end{minipage}
\begin{minipage}{3.15cm}
\includegraphics[width=3.15cm]{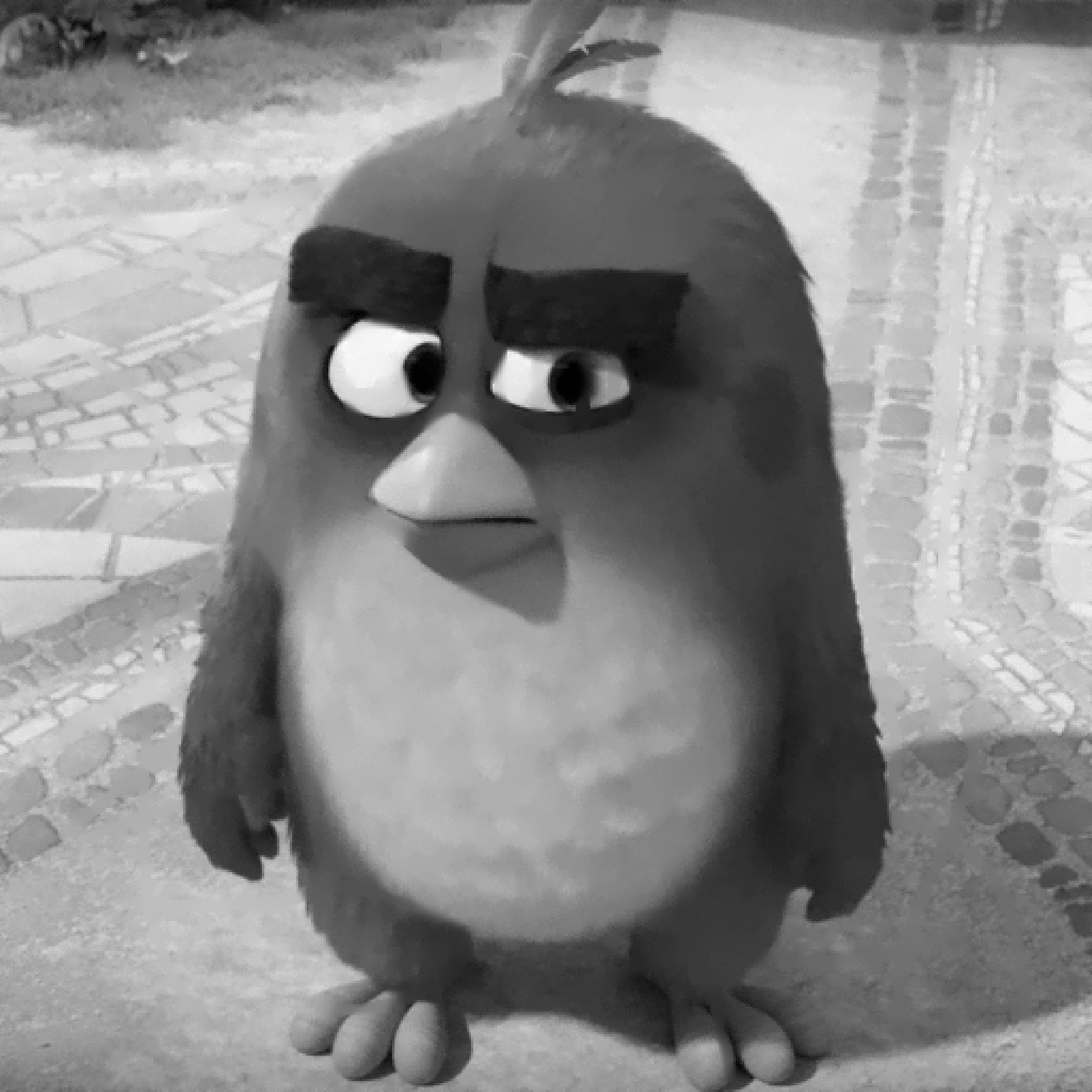}
\end{minipage}\vspace{0.25em}\\
\begin{minipage}{3.15cm}
\includegraphics[width=3.15cm]{PeppersCorrupted.pdf}
\end{minipage}
\begin{minipage}{3.15cm}
\includegraphics[width=3.15cm]{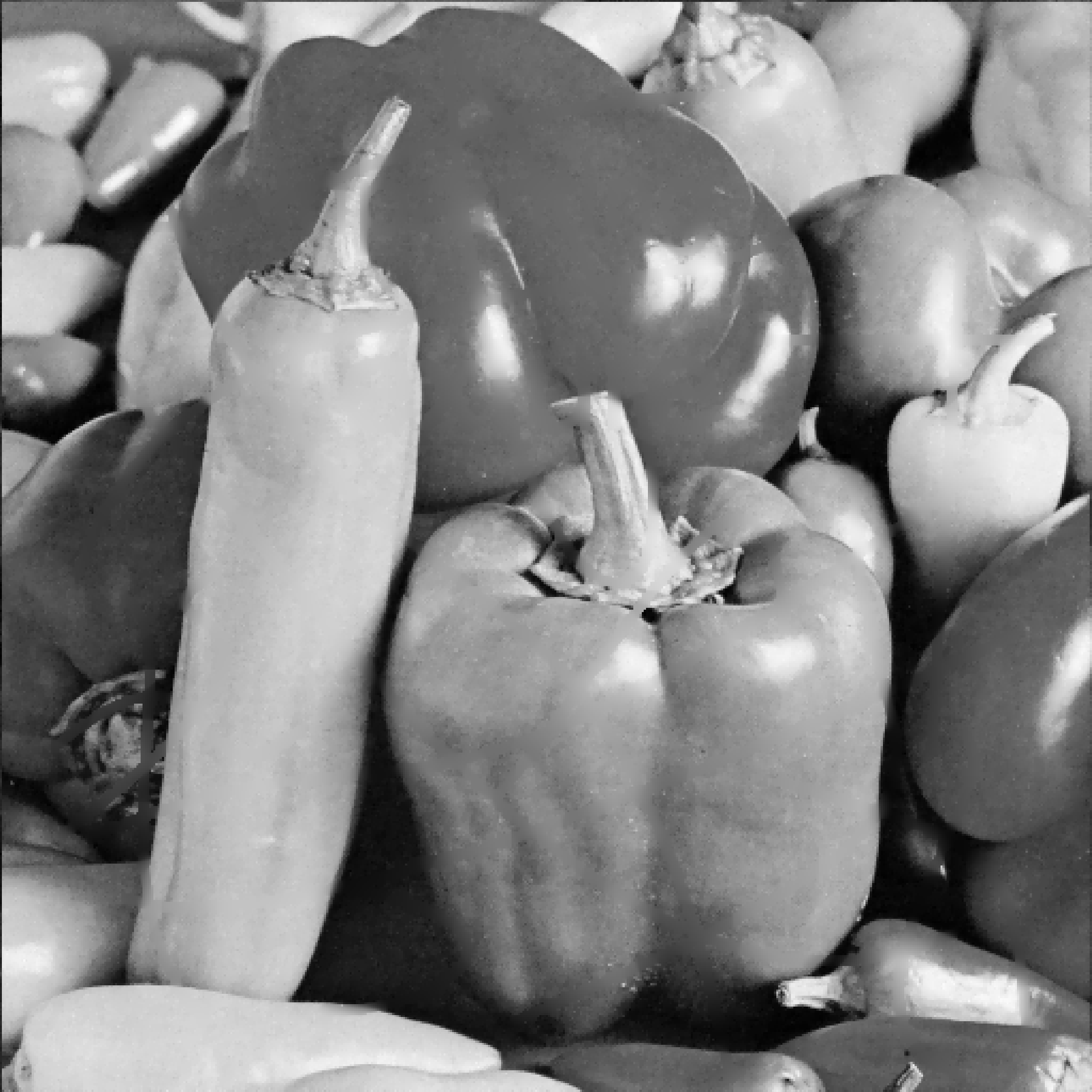}
\end{minipage}
\begin{minipage}{3.15cm}
\includegraphics[width=3.15cm]{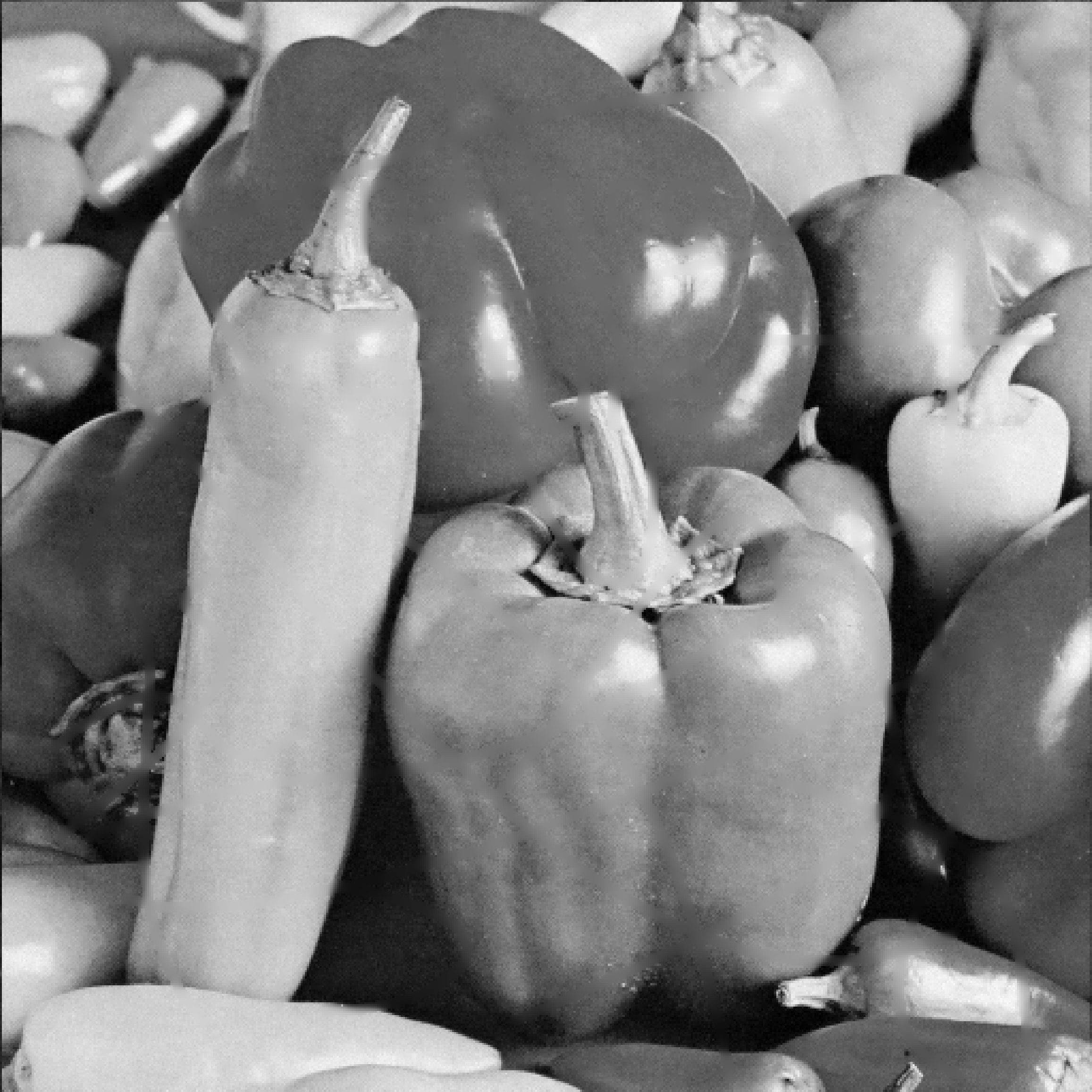}
\end{minipage}
\begin{minipage}{3.15cm}
\includegraphics[width=3.15cm]{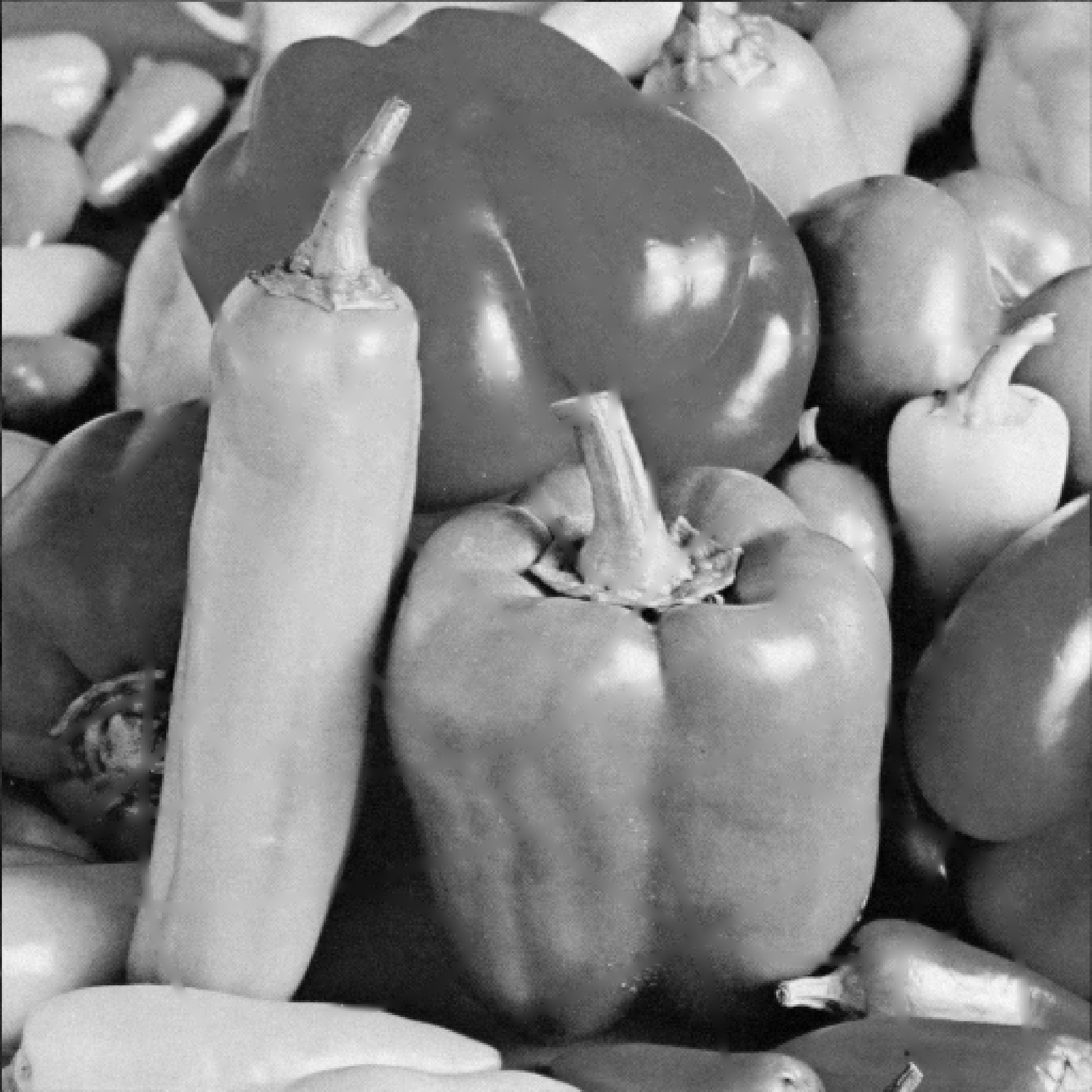}
\end{minipage}
\begin{minipage}{3.15cm}
\includegraphics[width=3.15cm]{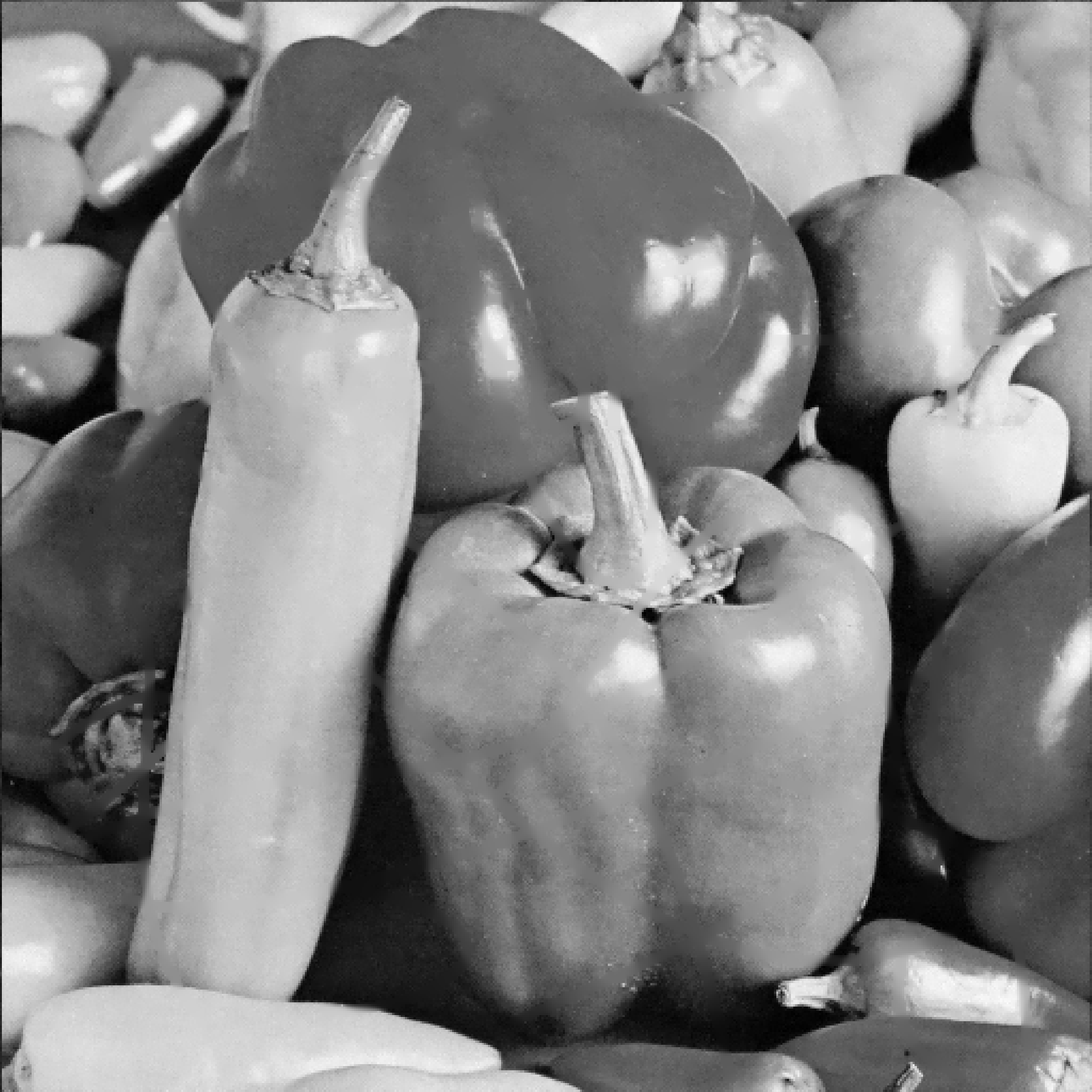}
\end{minipage}\vspace{0.25em}\\
\begin{minipage}{3.15cm}\begin{center}{\small{Observed}}\end{center}\end{minipage}\begin{minipage}{3.15cm}\begin{center}{\small{TGV Model \cite{K.Bredies2010}}}\end{center}\end{minipage}\begin{minipage}{3.15cm}\begin{center}{\small{PS Model \cite{J.F.Cai2016}}}\end{center}\end{minipage}\begin{minipage}{3.15cm}\begin{center}{\small{GS Model \cite{H.Ji2016}}}\end{center}\end{minipage}\begin{minipage}{3.15cm}\begin{center}{\small{Our Model \eqref{OurModel}}}\end{center}\end{minipage}
\caption{Visual comparisons of inpainted images for removing texts and scratches by four methods.}\label{fig:InpaintingResults}
\end{center}
\end{figure}

\begin{figure}[ht]
\begin{center}
\begin{minipage}{3.15cm}
\includegraphics[width=3.15cm]{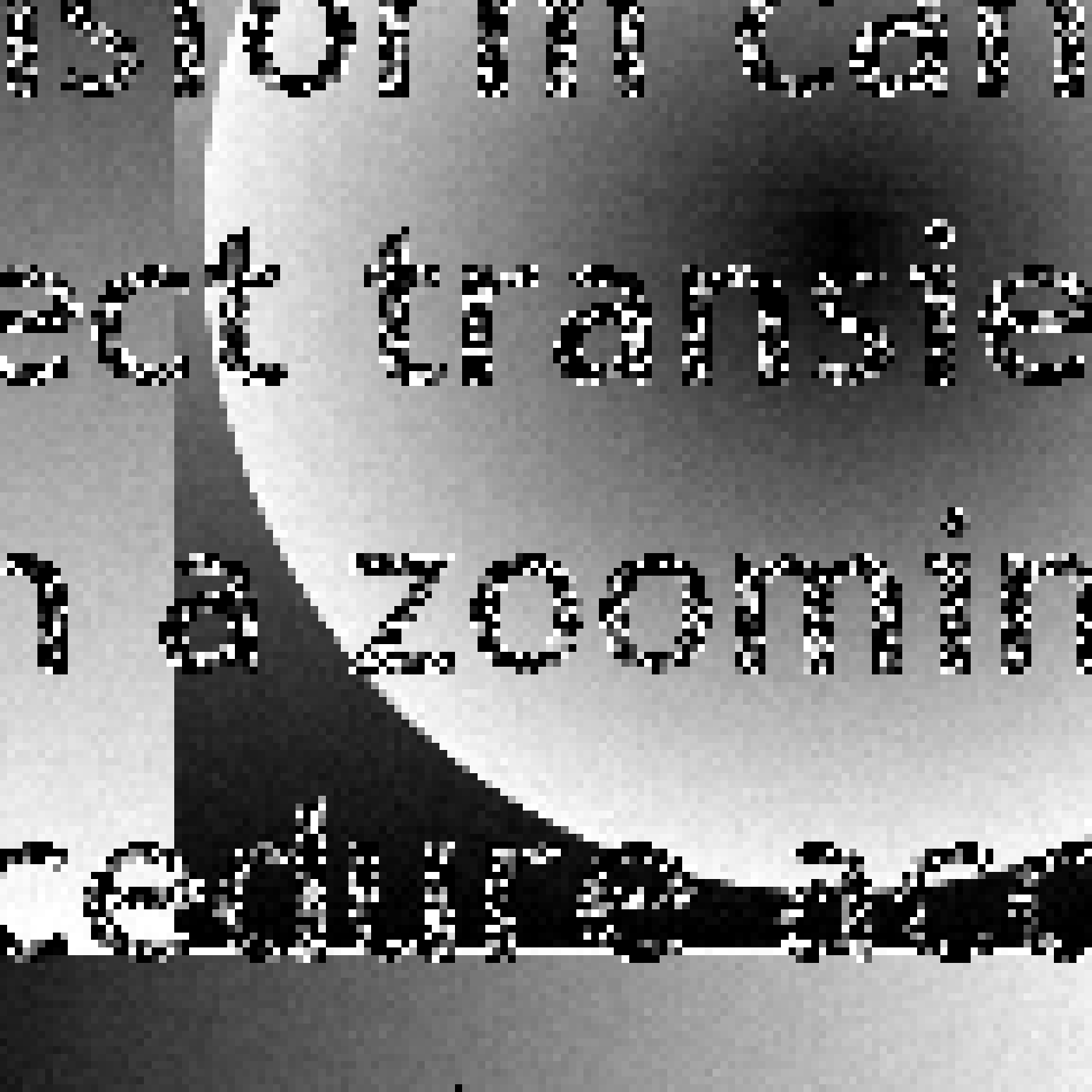}
\end{minipage}
\begin{minipage}{3.15cm}
\includegraphics[width=3.15cm]{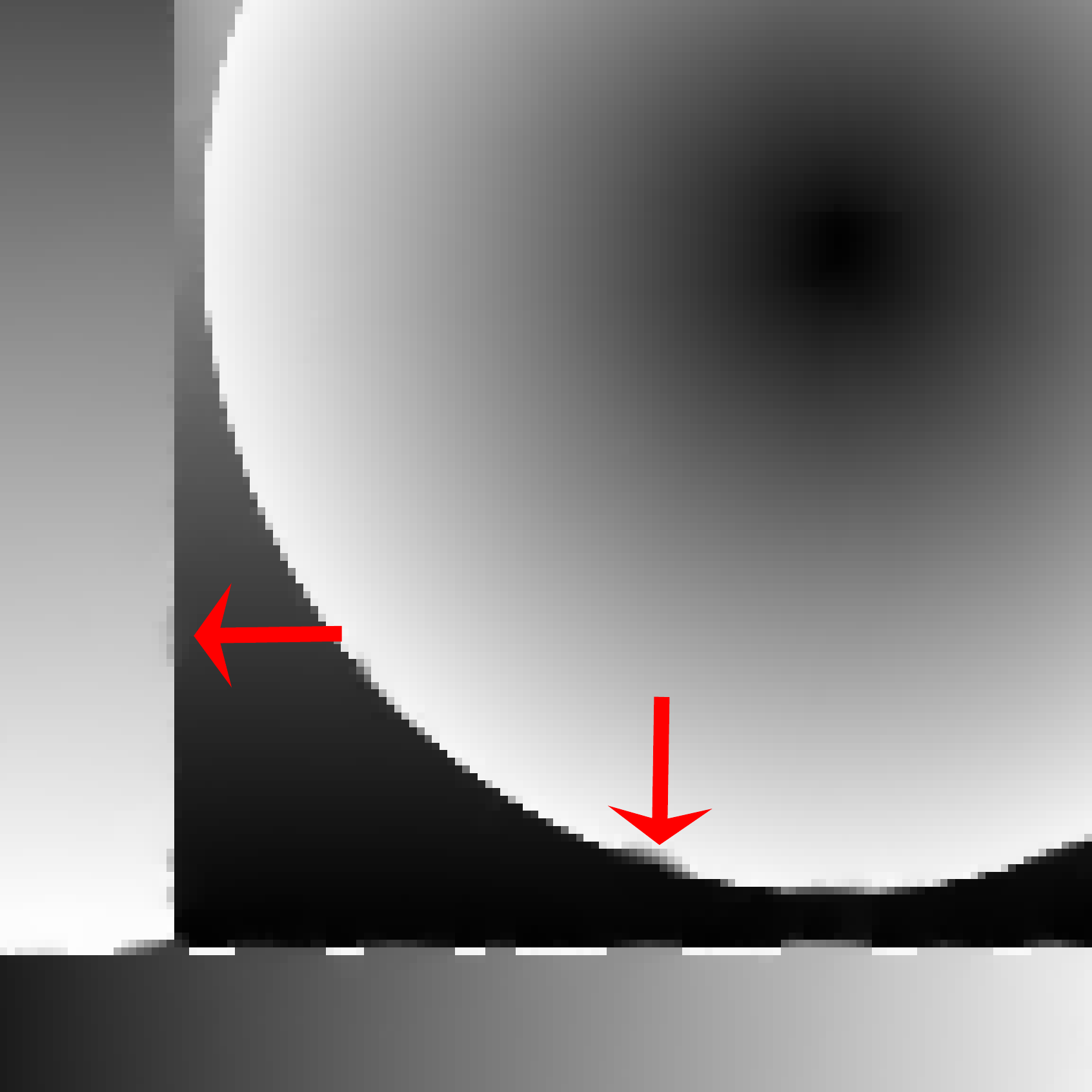}
\end{minipage}
\begin{minipage}{3.15cm}
\includegraphics[width=3.15cm]{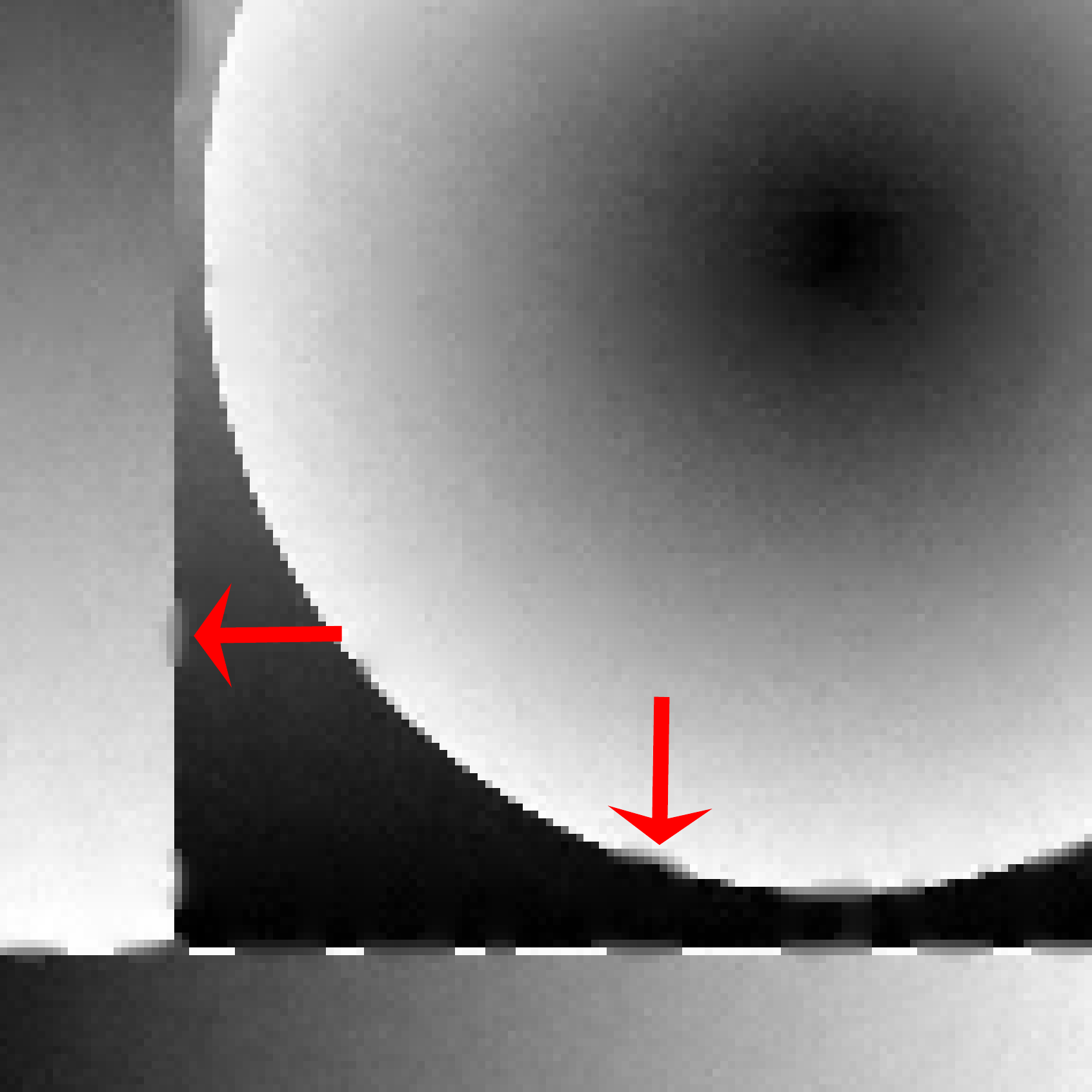}
\end{minipage}
\begin{minipage}{3.15cm}
\includegraphics[width=3.15cm]{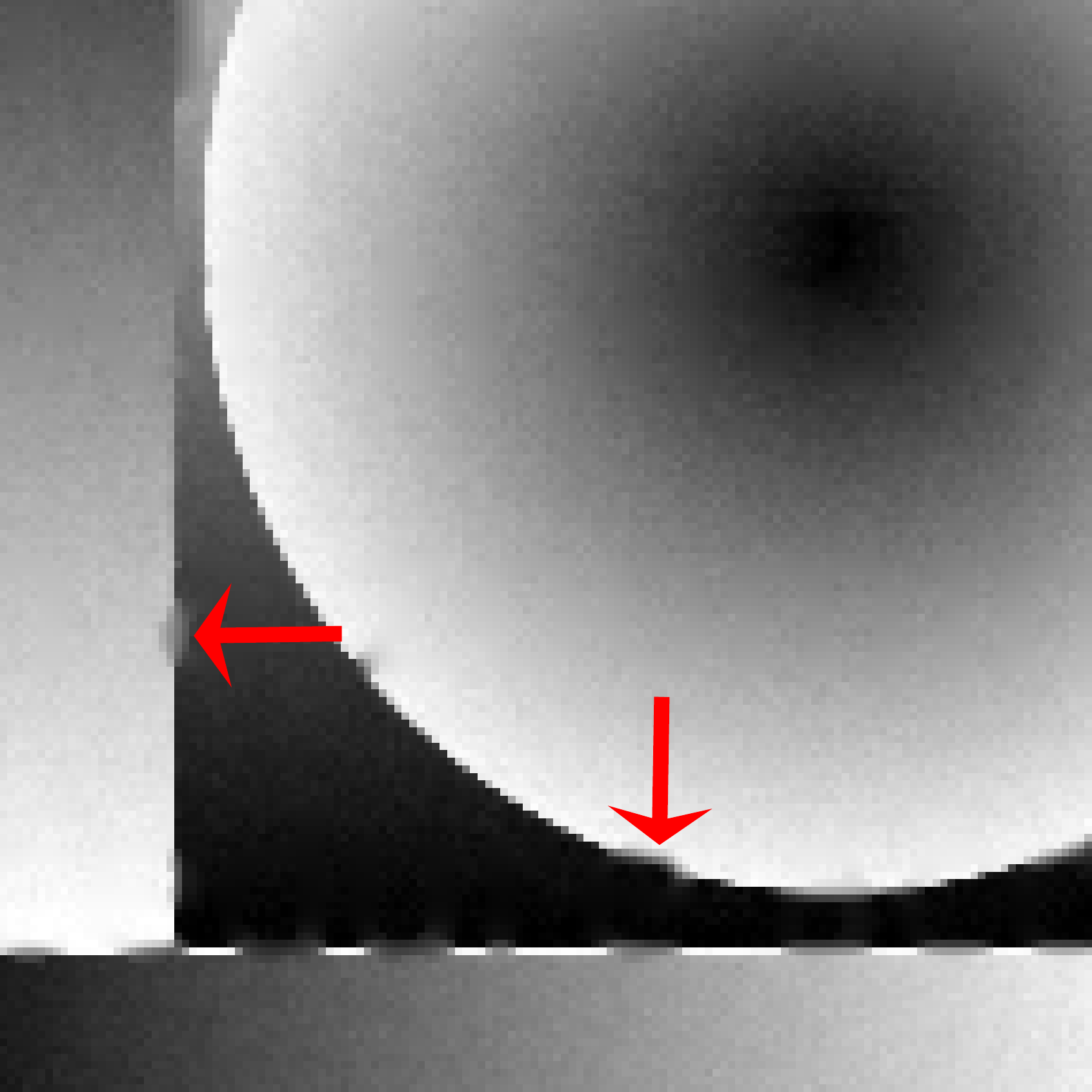}
\end{minipage}
\begin{minipage}{3.15cm}
\includegraphics[width=3.15cm]{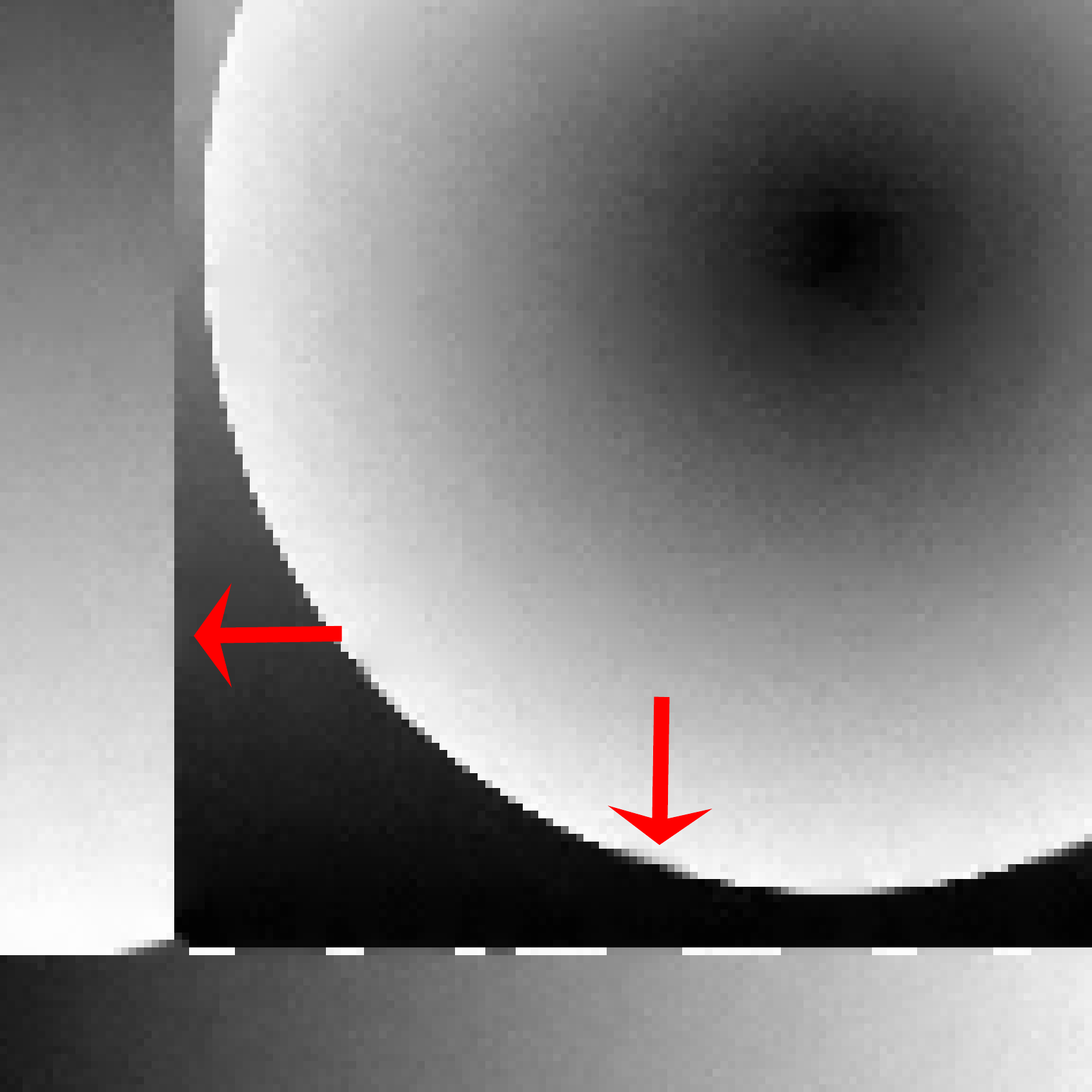}
\end{minipage}\vspace{0.25em}\\
\begin{minipage}{3.15cm}
\includegraphics[width=3.15cm]{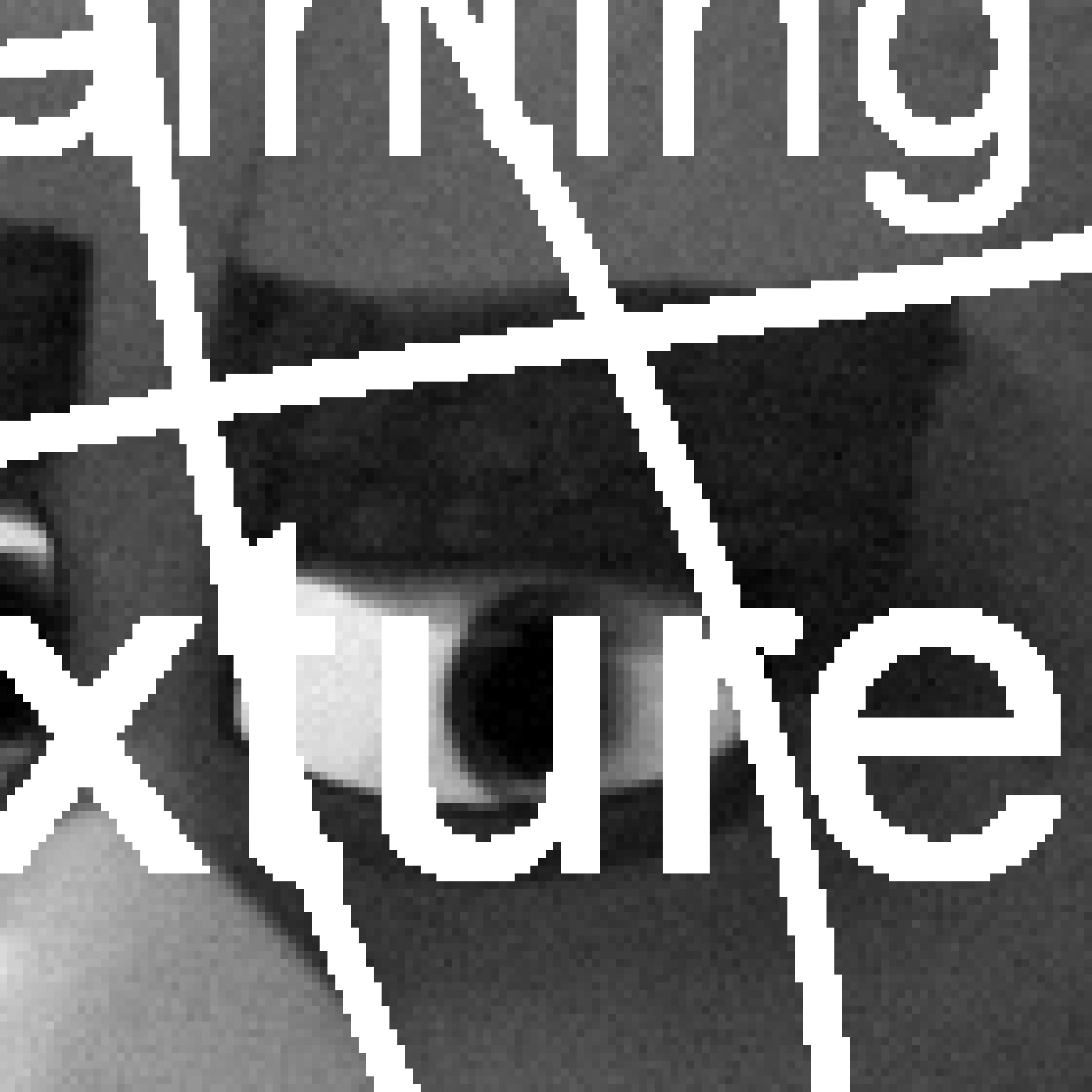}
\end{minipage}
\begin{minipage}{3.15cm}
\includegraphics[width=3.15cm]{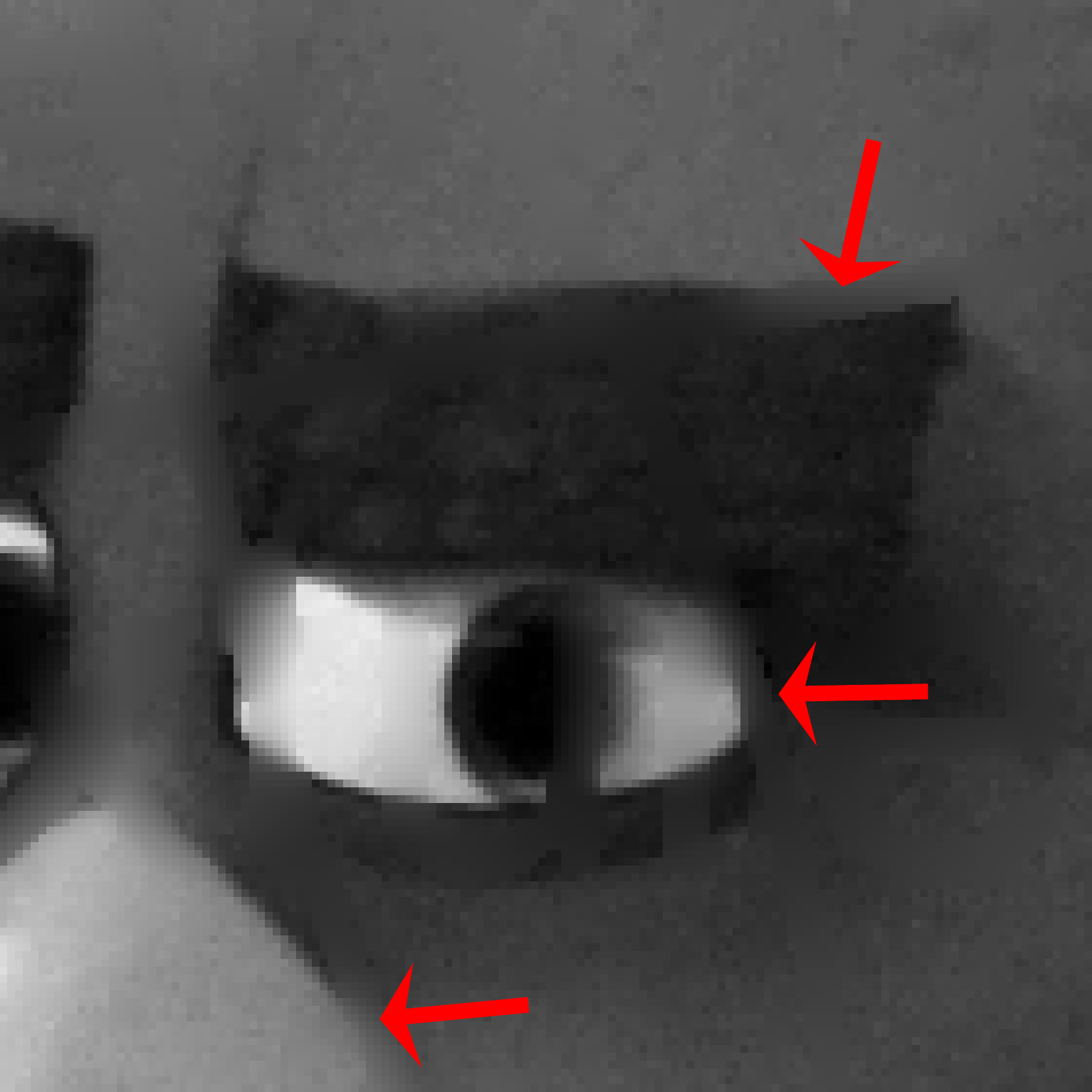}
\end{minipage}
\begin{minipage}{3.15cm}
\includegraphics[width=3.15cm]{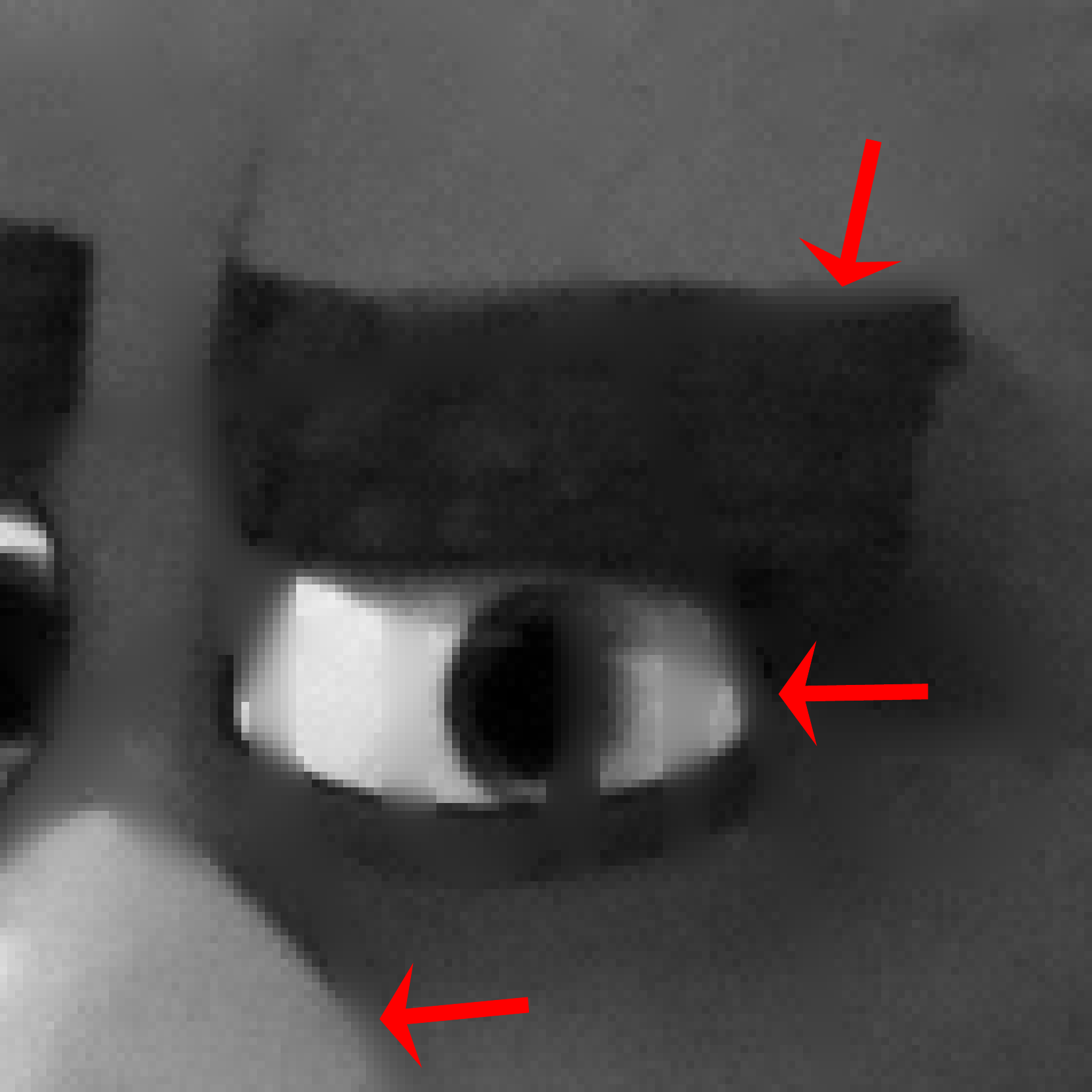}
\end{minipage}
\begin{minipage}{3.15cm}
\includegraphics[width=3.15cm]{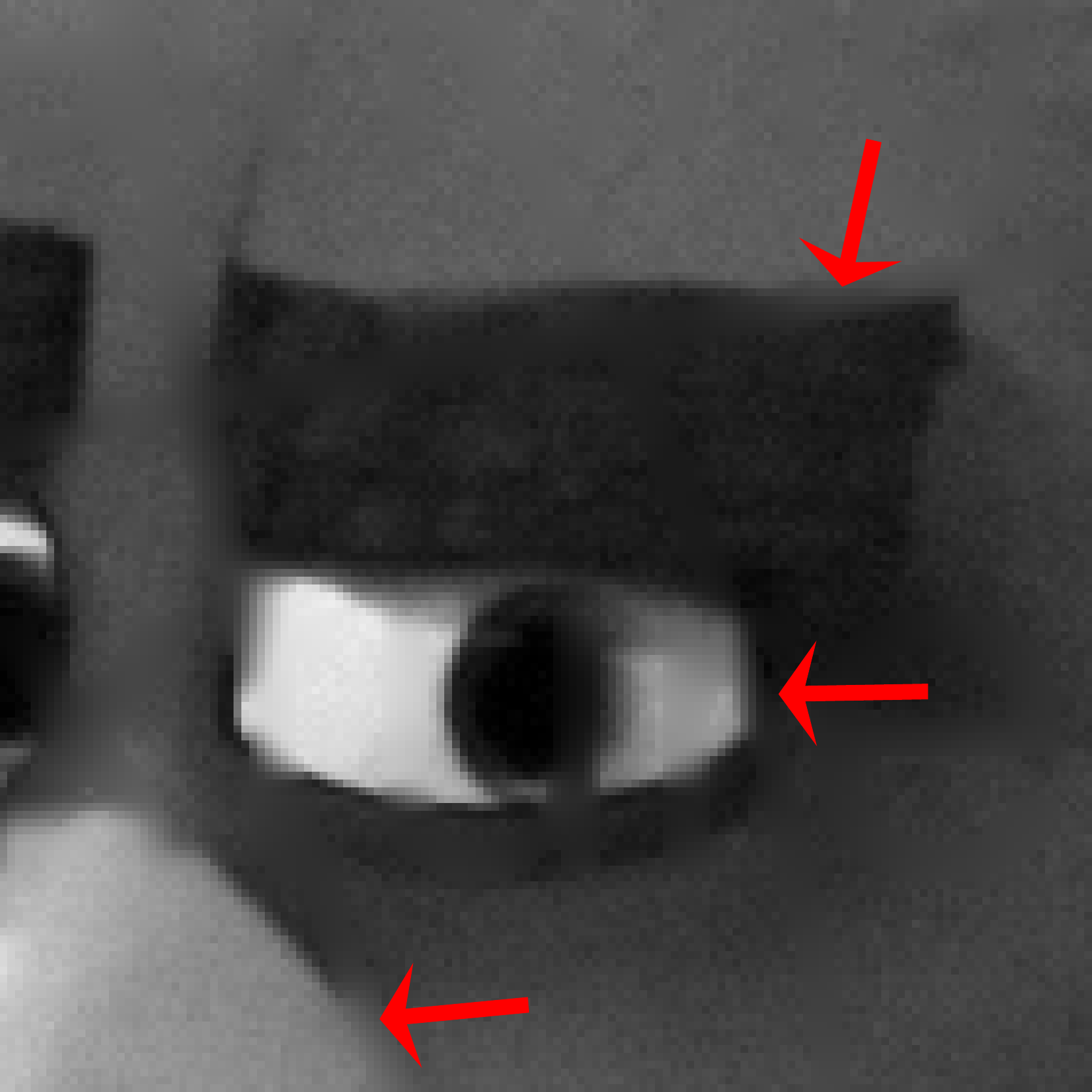}
\end{minipage}
\begin{minipage}{3.15cm}
\includegraphics[width=3.15cm]{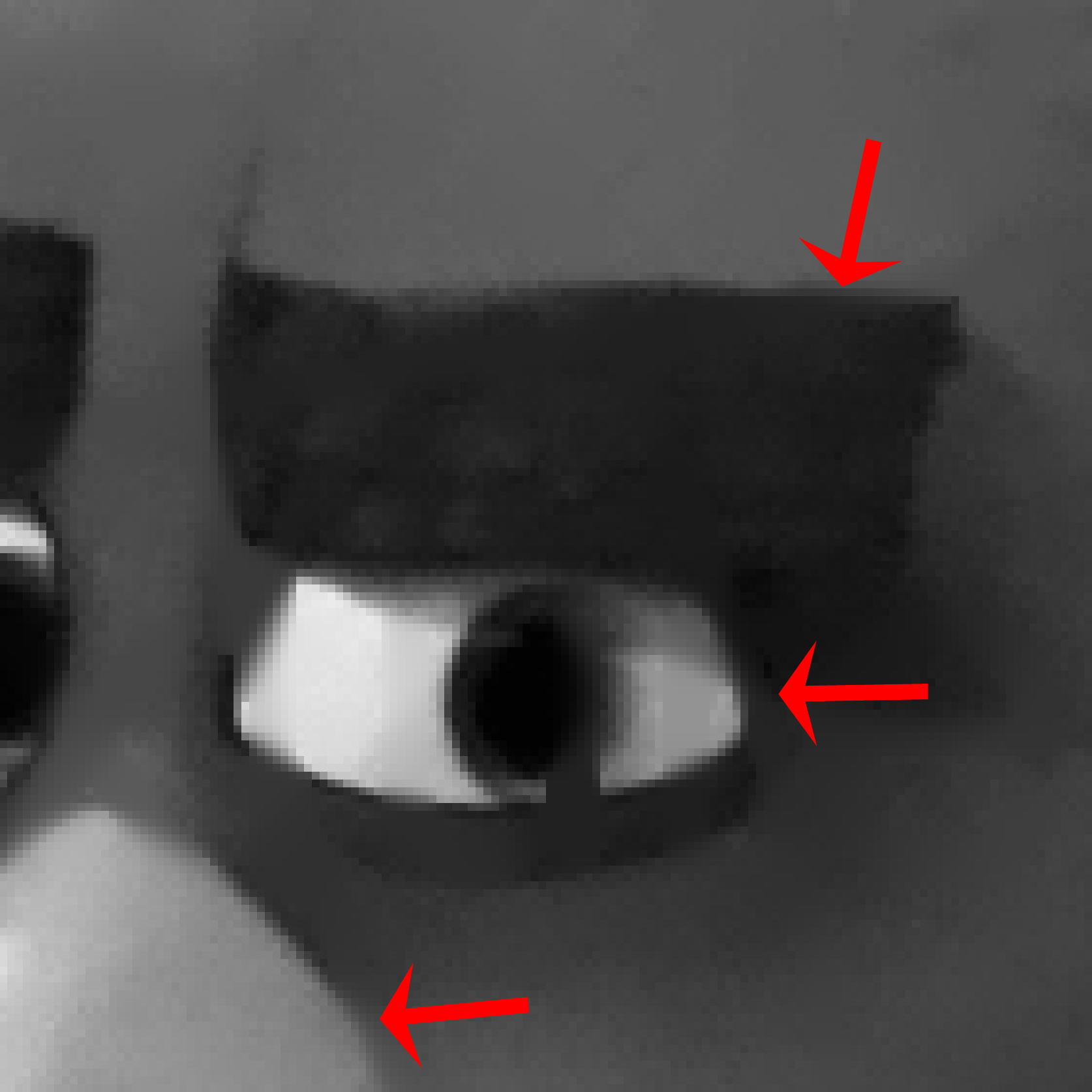}
\end{minipage}\vspace{0.25em}\\
\begin{minipage}{3.15cm}
\includegraphics[width=3.15cm]{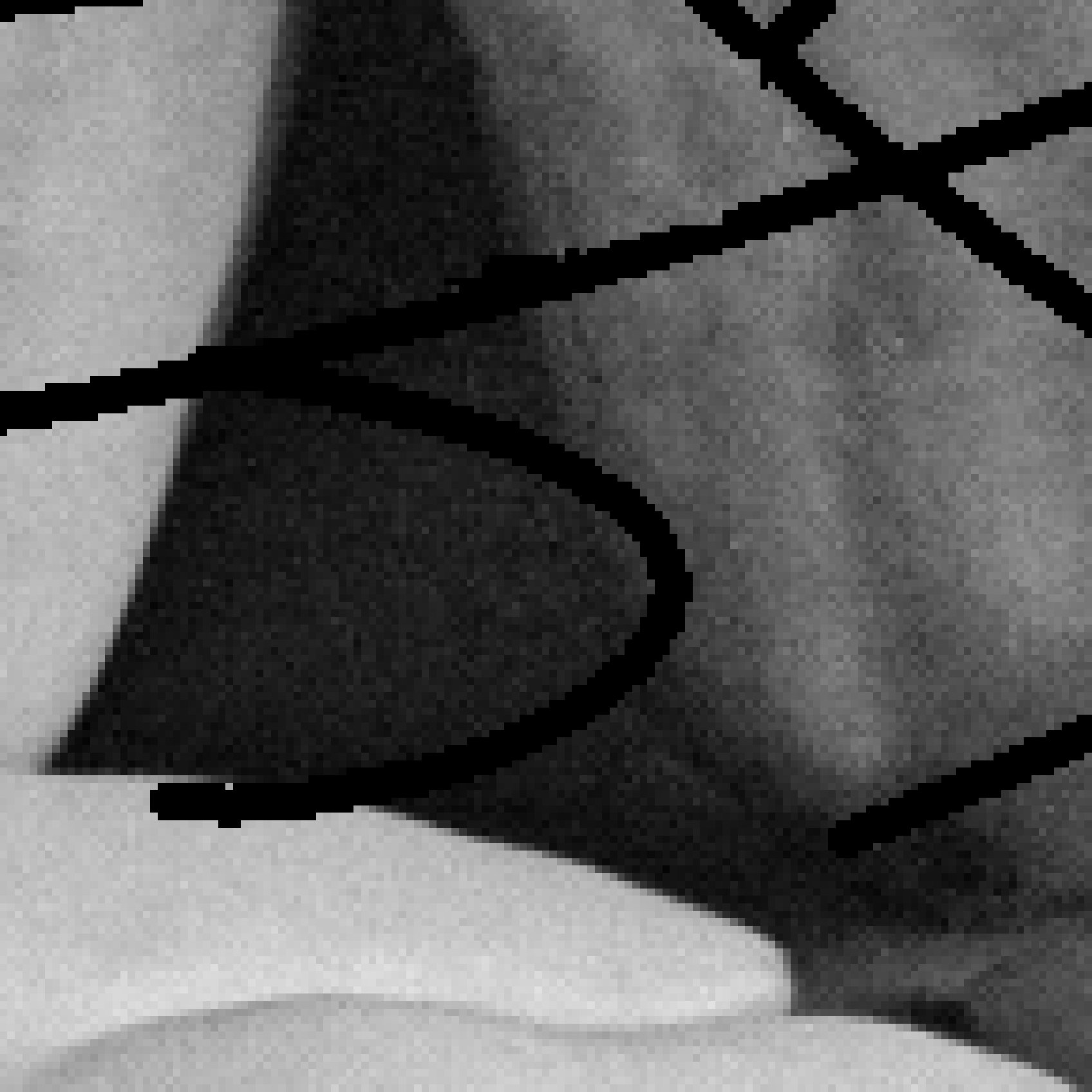}
\end{minipage}
\begin{minipage}{3.15cm}
\includegraphics[width=3.15cm]{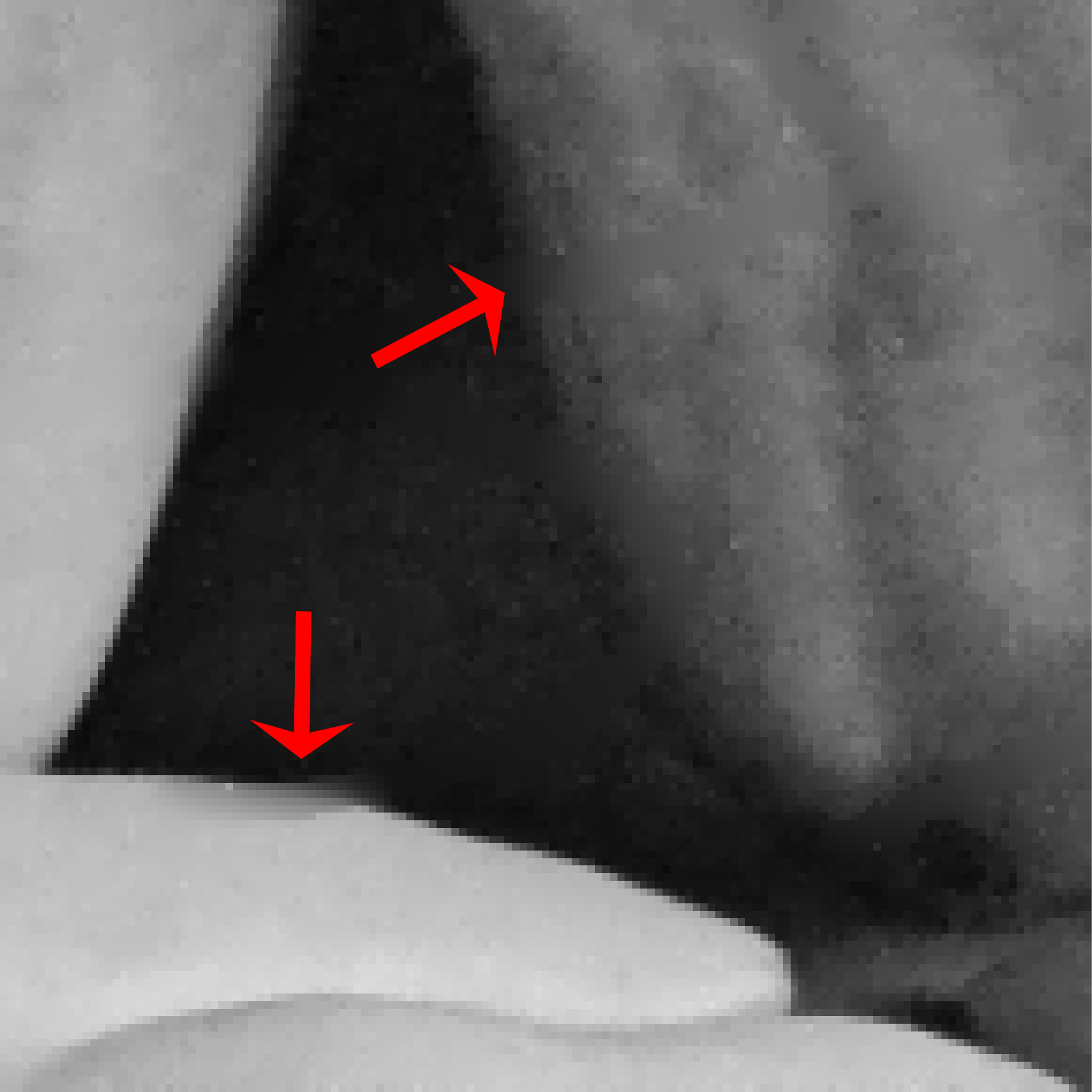}
\end{minipage}
\begin{minipage}{3.15cm}
\includegraphics[width=3.15cm]{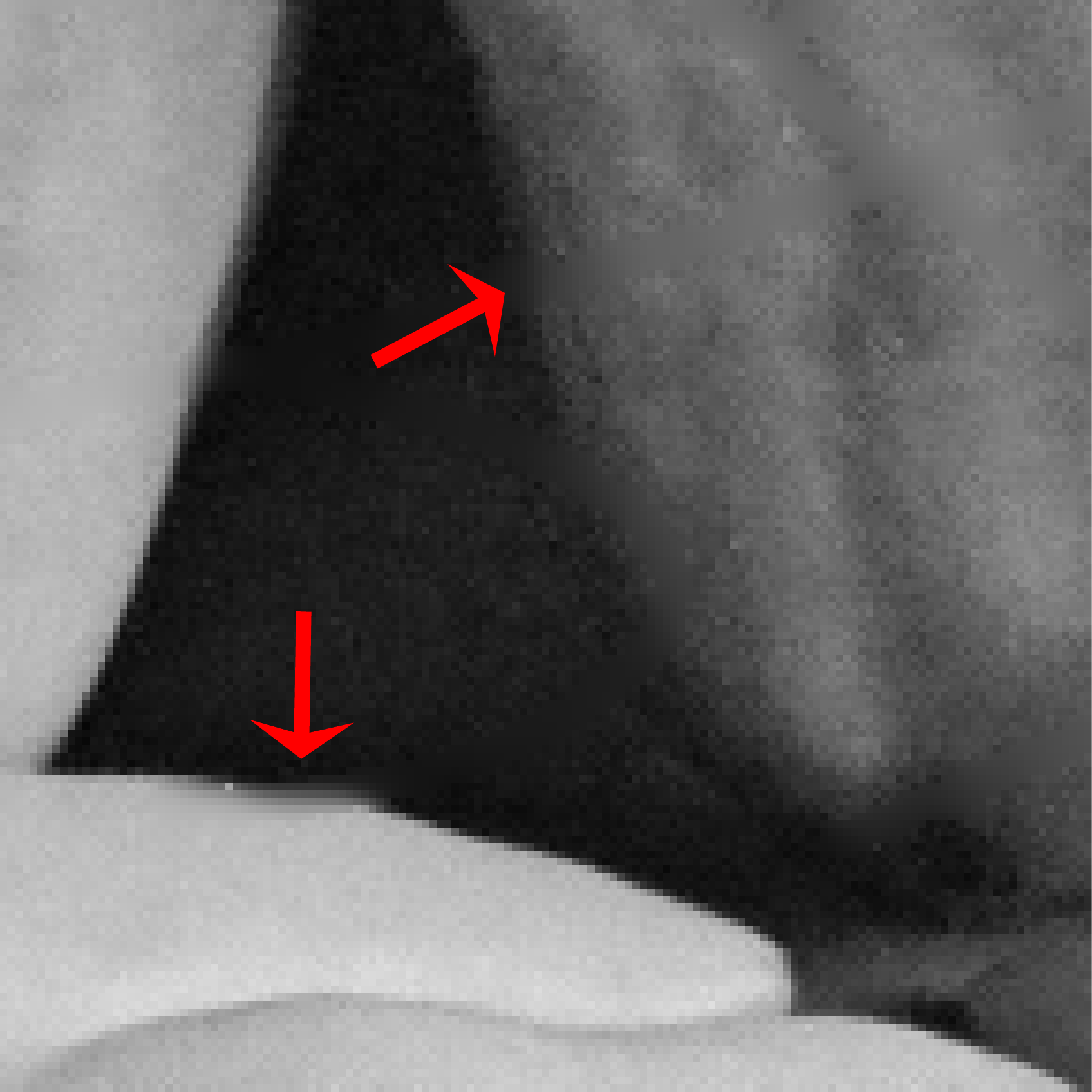}
\end{minipage}
\begin{minipage}{3.15cm}
\includegraphics[width=3.15cm]{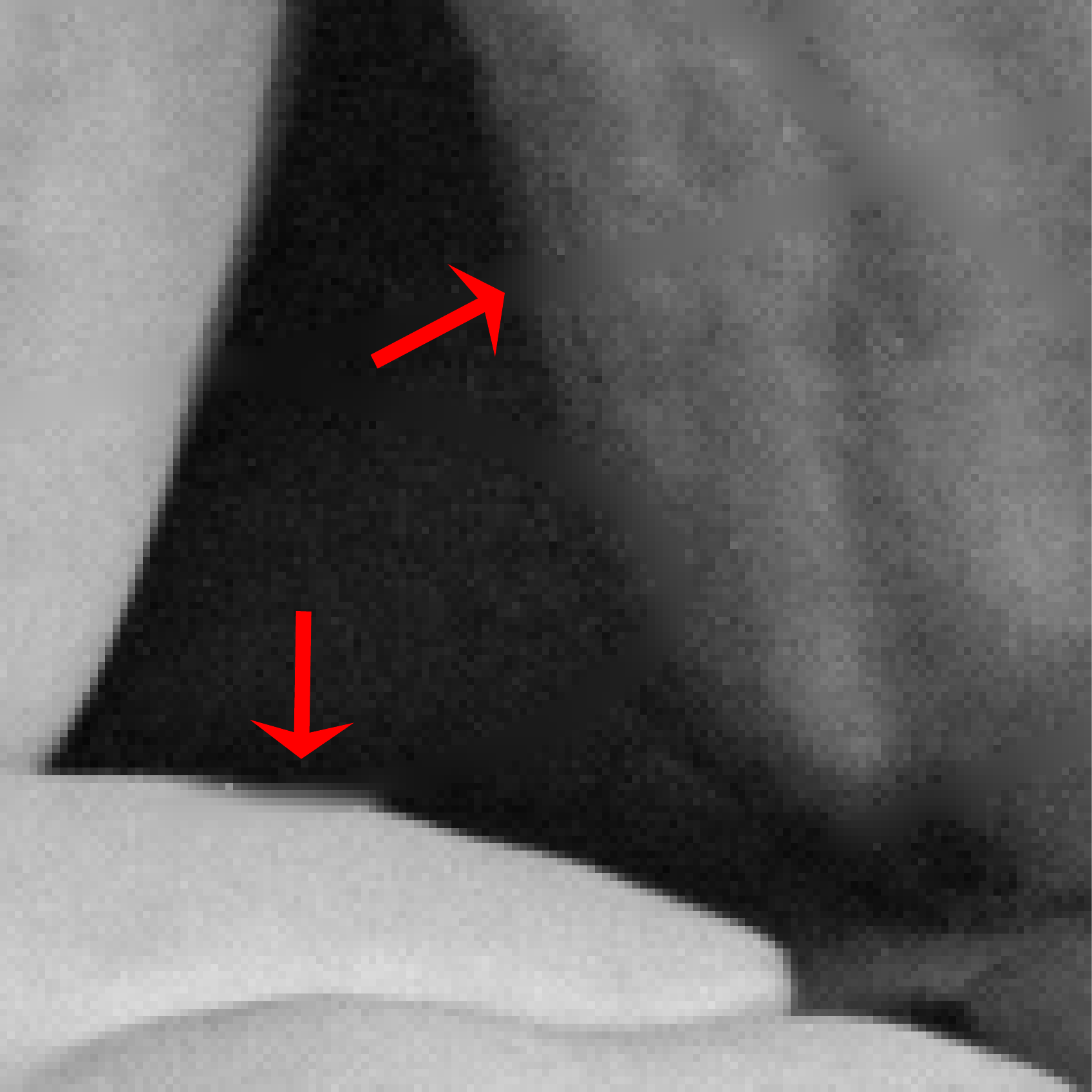}
\end{minipage}
\begin{minipage}{3.15cm}
\includegraphics[width=3.15cm]{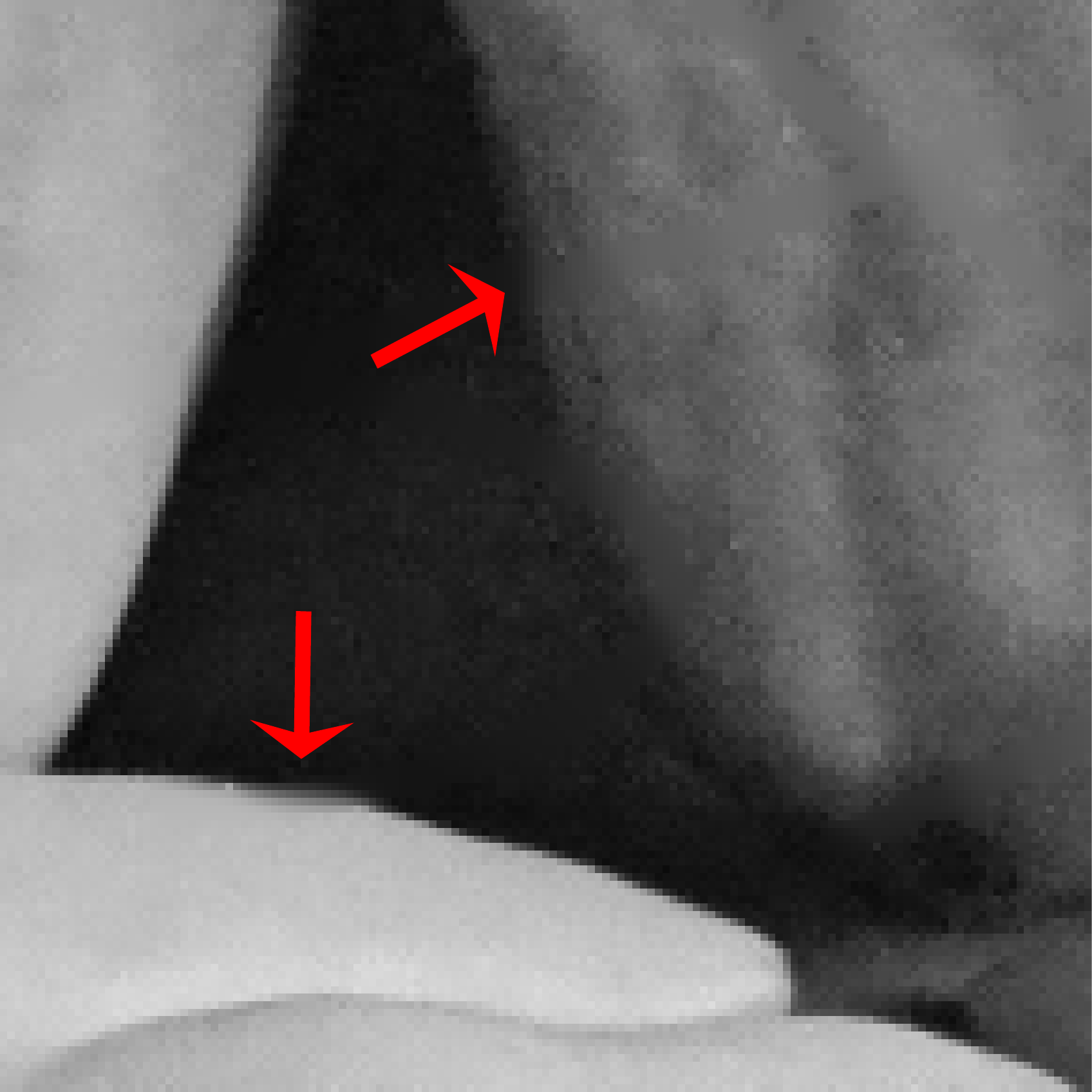}
\end{minipage}\vspace{0.25em}\\
\begin{minipage}{3.15cm}\begin{center}{\small{Observed}}\end{center}\end{minipage}\begin{minipage}{3.15cm}\begin{center}{\small{TGV Model \cite{K.Bredies2010}}}\end{center}\end{minipage}\begin{minipage}{3.15cm}\begin{center}{\small{PS Model \cite{J.F.Cai2016}}}\end{center}\end{minipage}\begin{minipage}{3.15cm}\begin{center}{\small{GS Model \cite{H.Ji2016}}}\end{center}\end{minipage}\begin{minipage}{3.15cm}\begin{center}{\small{Our Model \eqref{OurModel}}}\end{center}\end{minipage}
\caption{Zoom-in views of Figure \ref{fig:InpaintingResults}. The red arrows indicate the region worth noticing.}\label{fig:InpaintingResultsZoom}
\end{center}
\end{figure}

Table \ref{tab:PSNRCompareInpaint} summarizes the results of the aforementioned four models for image inpainting, and Figure \ref{fig:InpaintingResults} and Figure \ref{fig:InpaintingResultsZoom} present visual comparisons of the results. It can be seen from Table \ref{tab:PSNRCompareInpaint} that our model \eqref{OurModel} consistently outperforms other image restoration models. Compared to the $\ell_2$ norm based PS model \eqref{CDS} and GS model \eqref{JLS}, we can see that our model does not smear out the singularities that are not captured by $\bsv$, leading to the visual improvements that are consistent with the improvements in PSNR values.

\begin{figure}[ht]
\begin{center}
\begin{minipage}{4.2cm}
\includegraphics[width=4.2cm]{SlopeCorrupted.pdf}
\end{minipage}
\begin{minipage}{4.2cm}
\includegraphics[width=4.2cm]{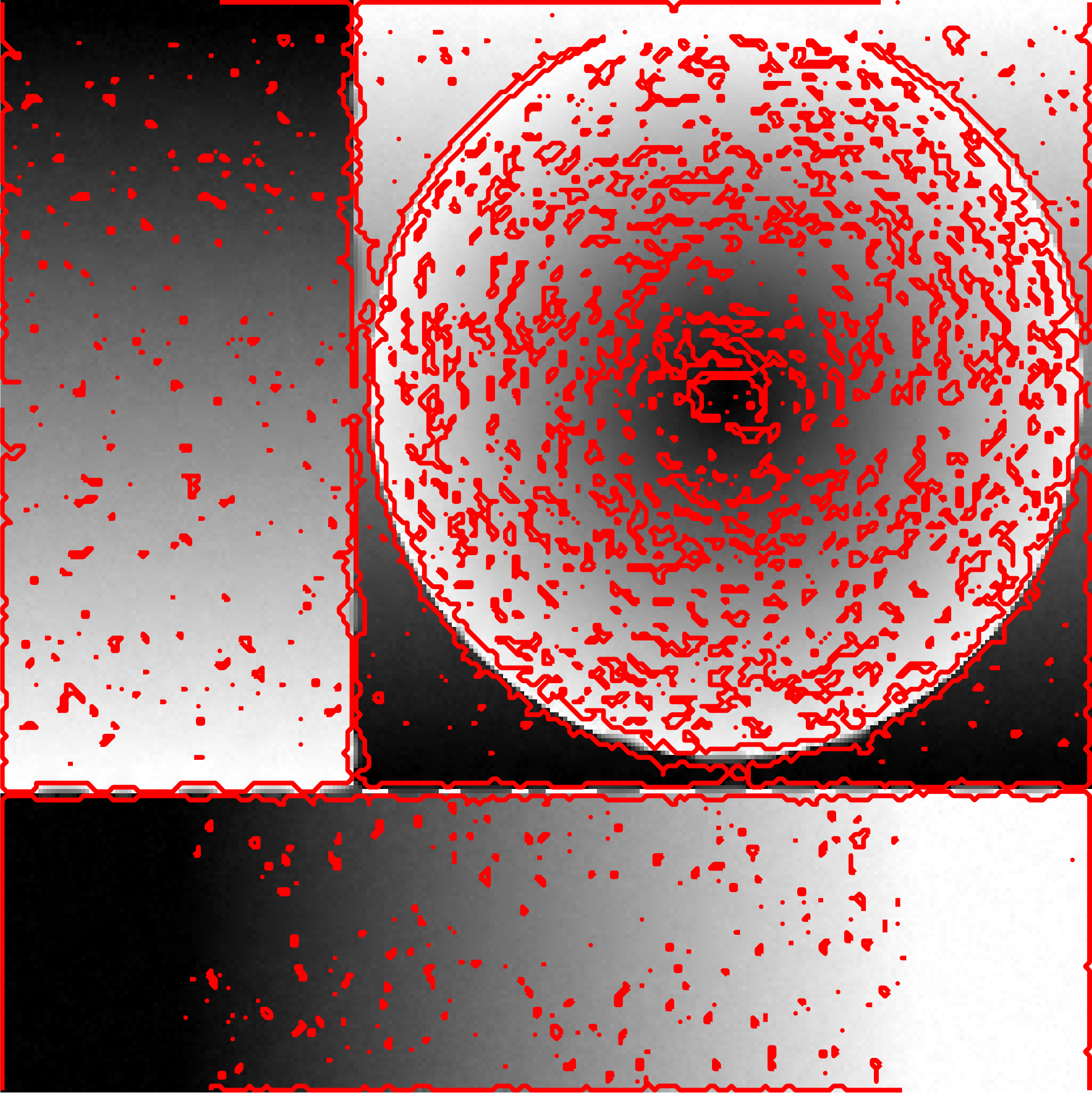}
\end{minipage}
\begin{minipage}{4.2cm}
\includegraphics[width=4.2cm]{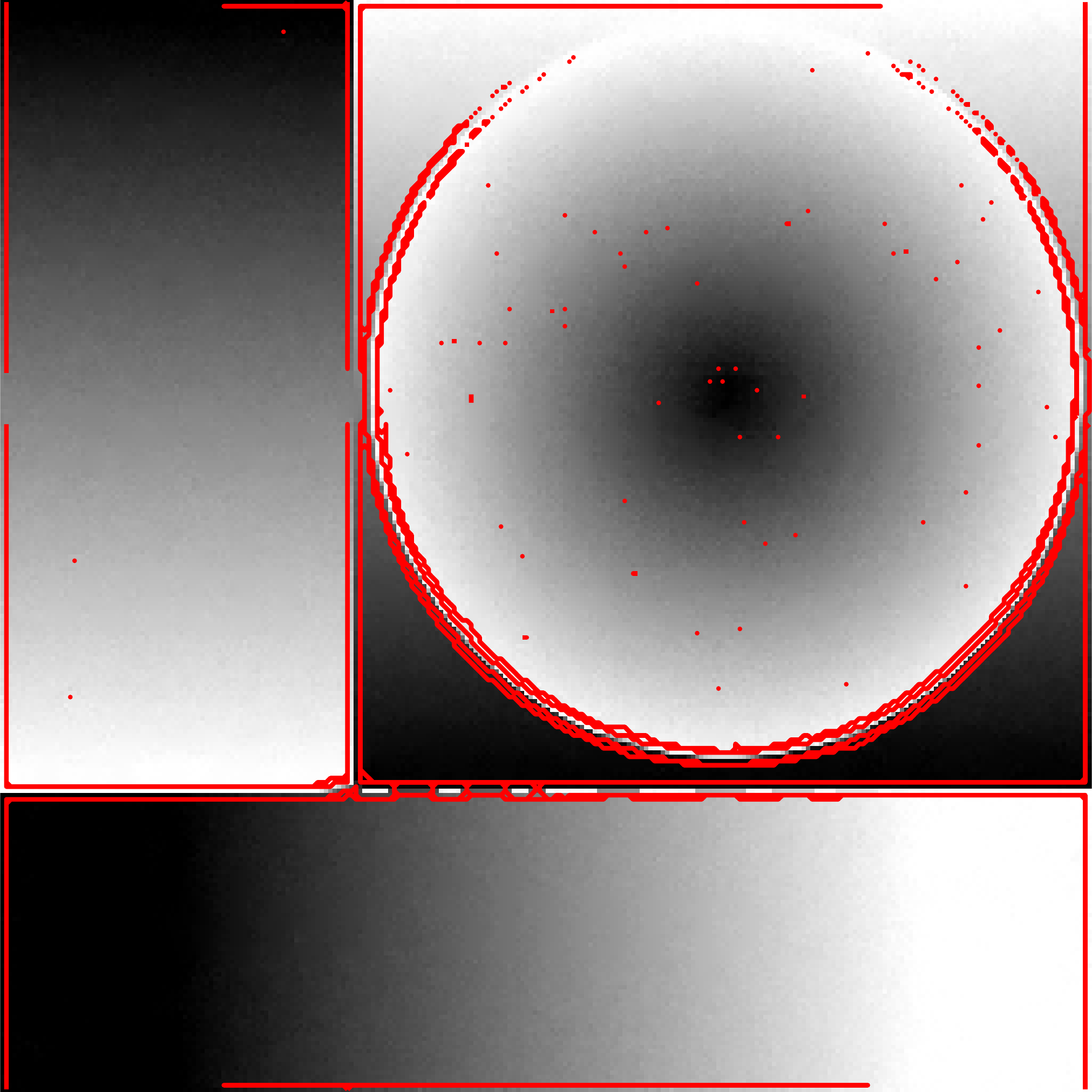}
\end{minipage}\vspace{0.25em}\\
\begin{minipage}{4.2cm}
\includegraphics[width=4.2cm]{AngryBirdsCorrupted.pdf}
\end{minipage}
\begin{minipage}{4.2cm}
\includegraphics[width=4.2cm]{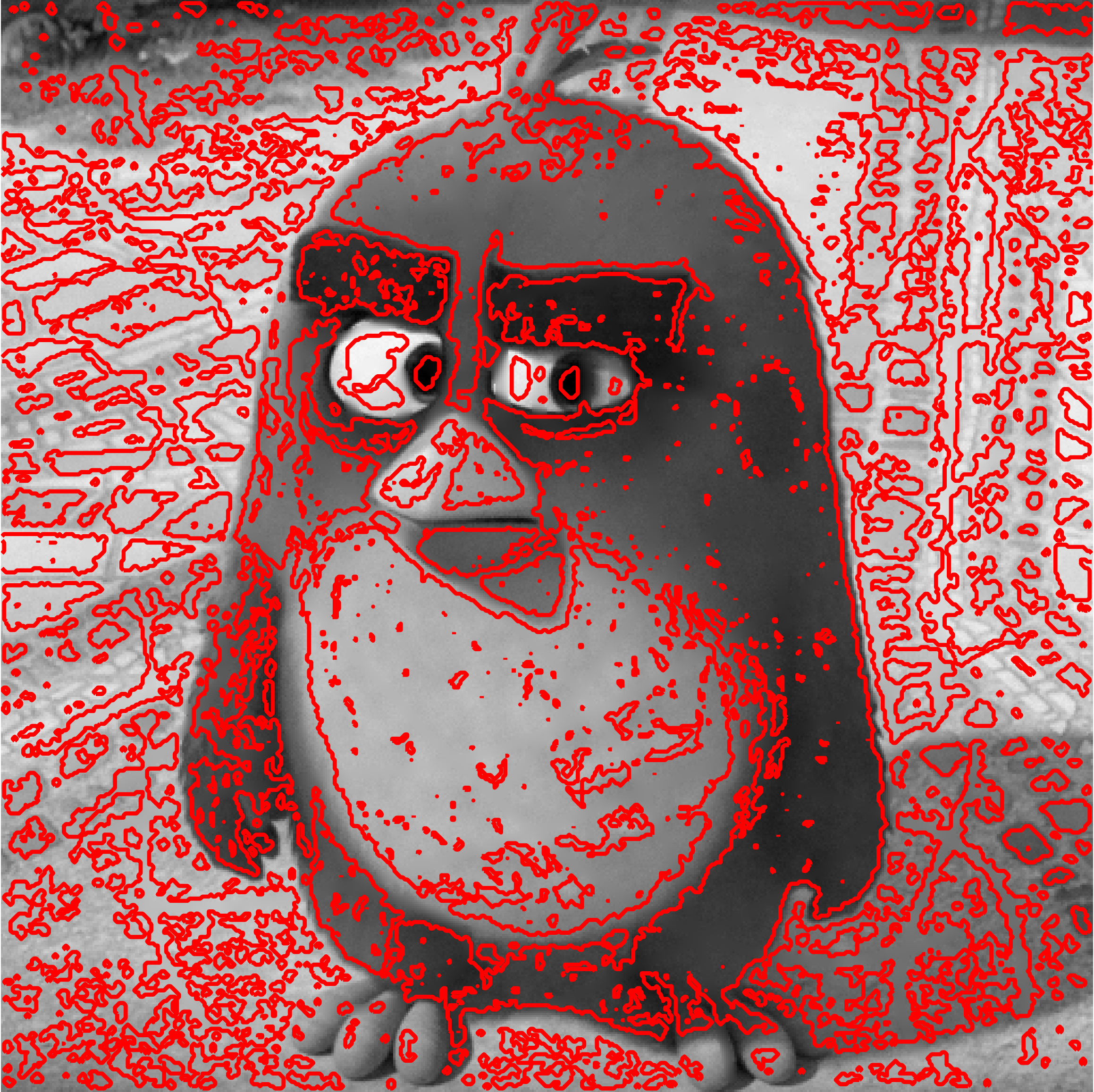}
\end{minipage}
\begin{minipage}{4.2cm}
\includegraphics[width=4.2cm]{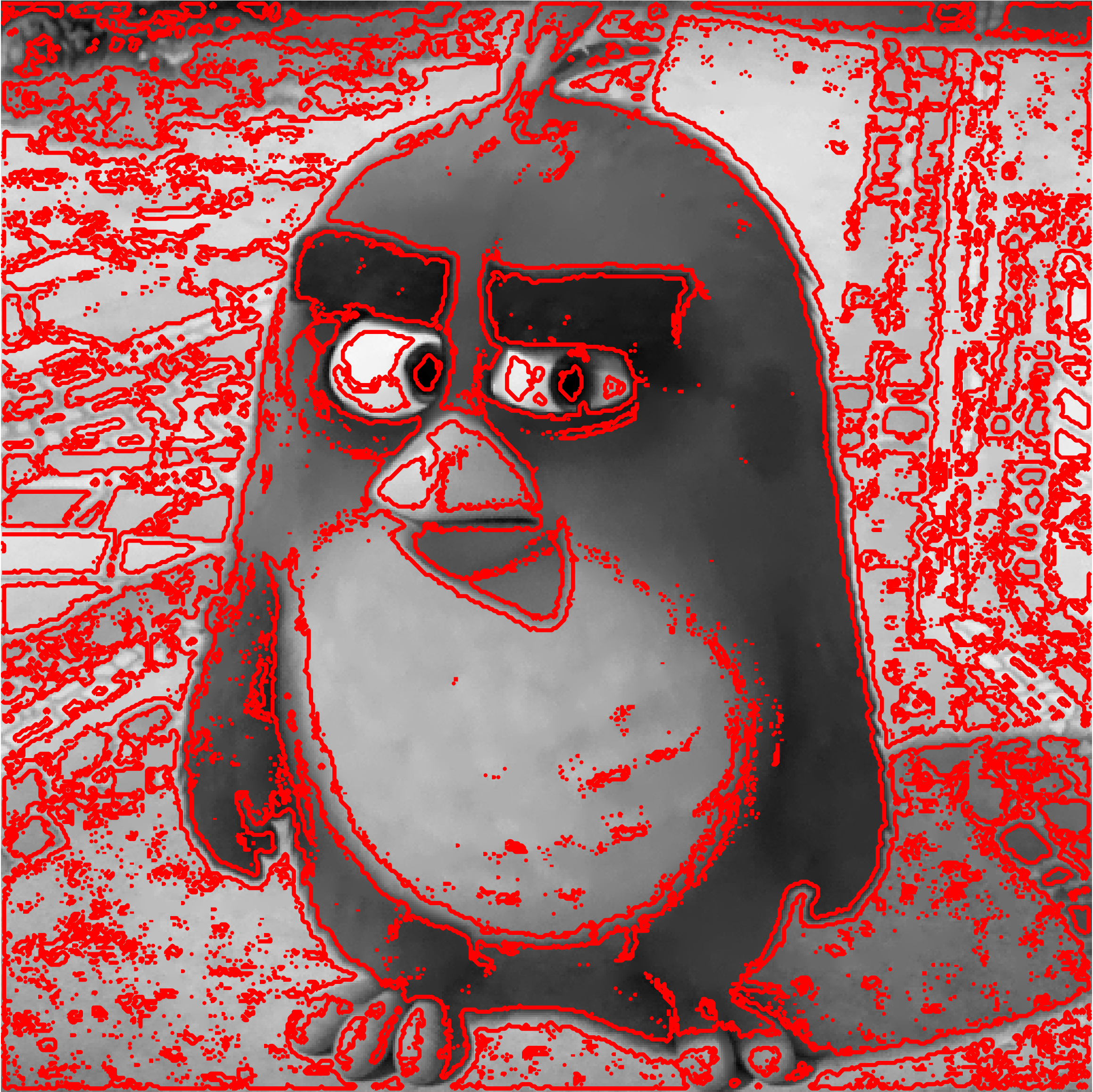}
\end{minipage}\vspace{0.25em}\\
\begin{minipage}{4.2cm}
\includegraphics[width=4.2cm]{PeppersCorrupted.pdf}
\end{minipage}
\begin{minipage}{4.2cm}
\includegraphics[width=4.2cm]{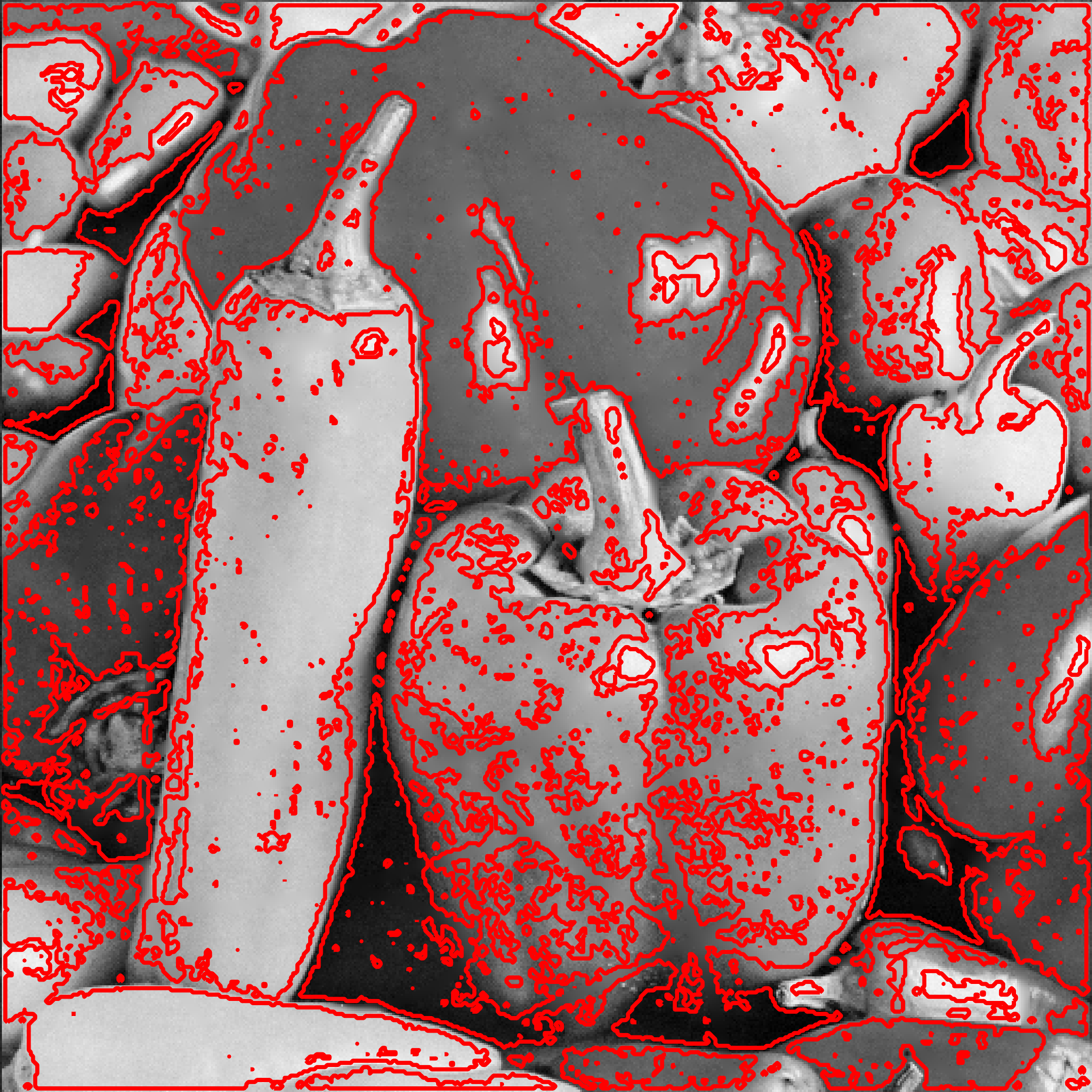}
\end{minipage}
\begin{minipage}{4.2cm}
\includegraphics[width=4.2cm]{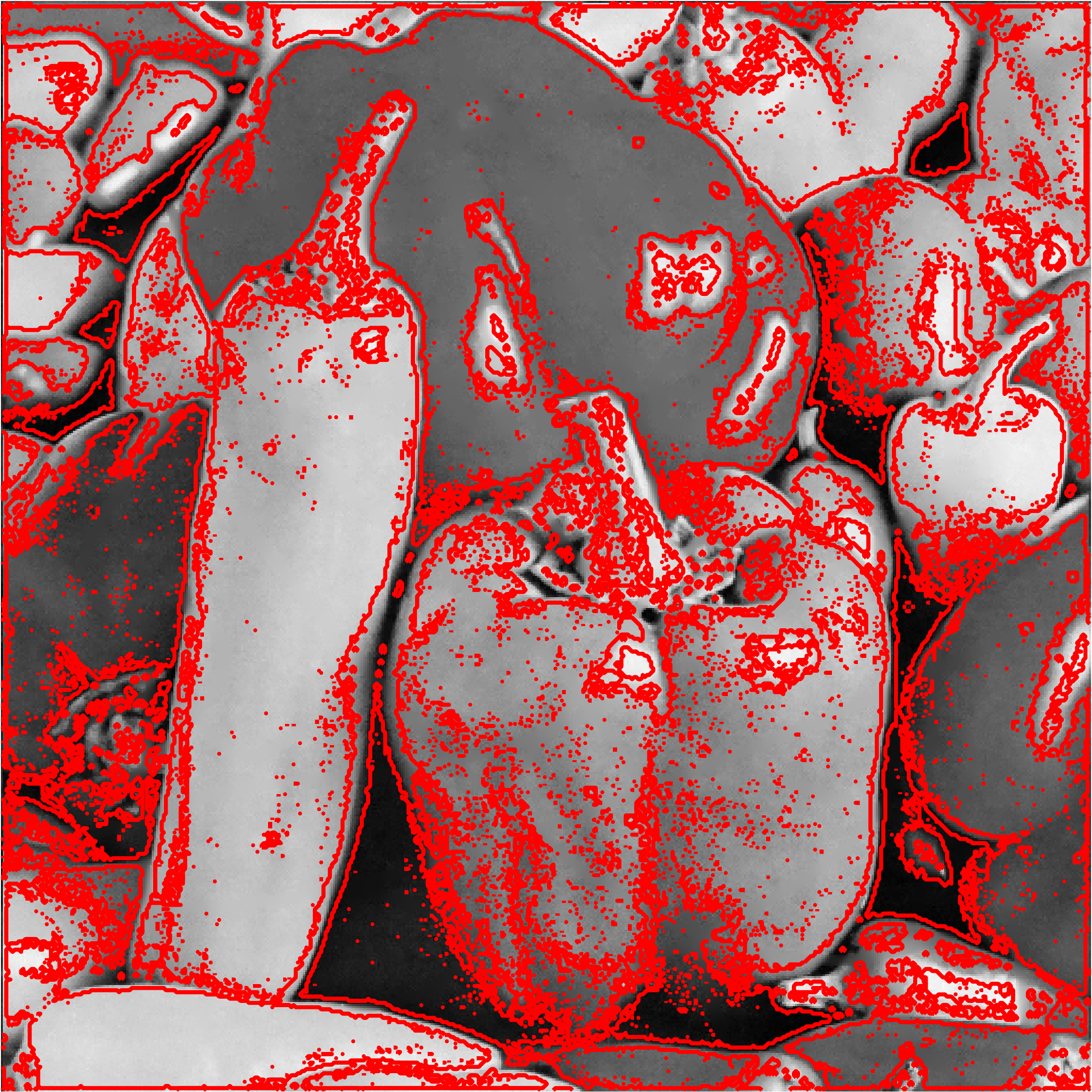}
\end{minipage}\vspace{0.25em}\\
\begin{minipage}{4.2cm}\begin{center}{\small{Observed}}\end{center}\end{minipage}\begin{minipage}{4.2cm}\begin{center}{\small{PS Model \cite{J.F.Cai2016}}}\end{center}\end{minipage}\begin{minipage}{4.2cm}\begin{center}{\small{Our Model \eqref{OurModel}}}\end{center}\end{minipage}
\caption{Comparison of estimated jump sets which are marked by red curves. We can see that the wavelet frame regularization on $\bsv$ can remove the discontinuities caused by the scratches and the texts, leading to the better inpainted results.}\label{fig:InpaintJumpEstiComp}
\end{center}
\end{figure}

The singularities estimated by the PS model \eqref{CDS} and our model \eqref{OurModel} are shown in Figure \ref{fig:InpaintJumpEstiComp}. We can easily see that the singularities estimated by our model contains less isolated singularities compared with the PS model. By relaxing the binary image $\one_{\Sig}$ into $\bsv$ taking values in $[0,1]$ and regularizing it by the wavelet frame system $\bsW''$, we can remove the isolated singularities which can be captured by solely comparing the wavelet frame coefficients. In particular, it is worth noting that the singularities estimated by our model do not include the texts and the scratches.

\begin{figure}[ht]
\begin{center}
\begin{minipage}{4.2cm}
\includegraphics[width=4.2cm]{SlopeCorruptedZoom.pdf}
\end{minipage}
\begin{minipage}{4.2cm}
\includegraphics[width=4.2cm]{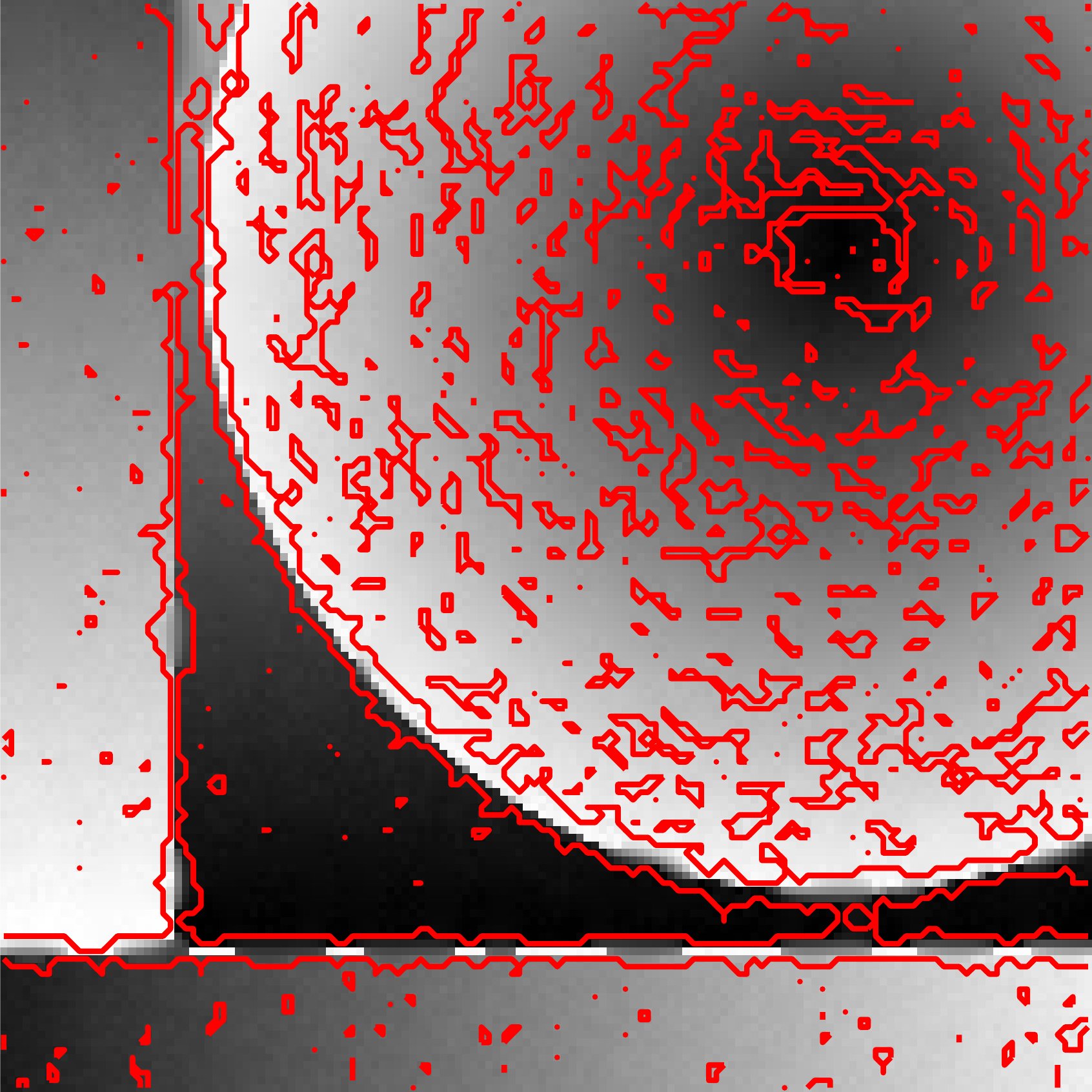}
\end{minipage}
\begin{minipage}{4.2cm}
\includegraphics[width=4.2cm]{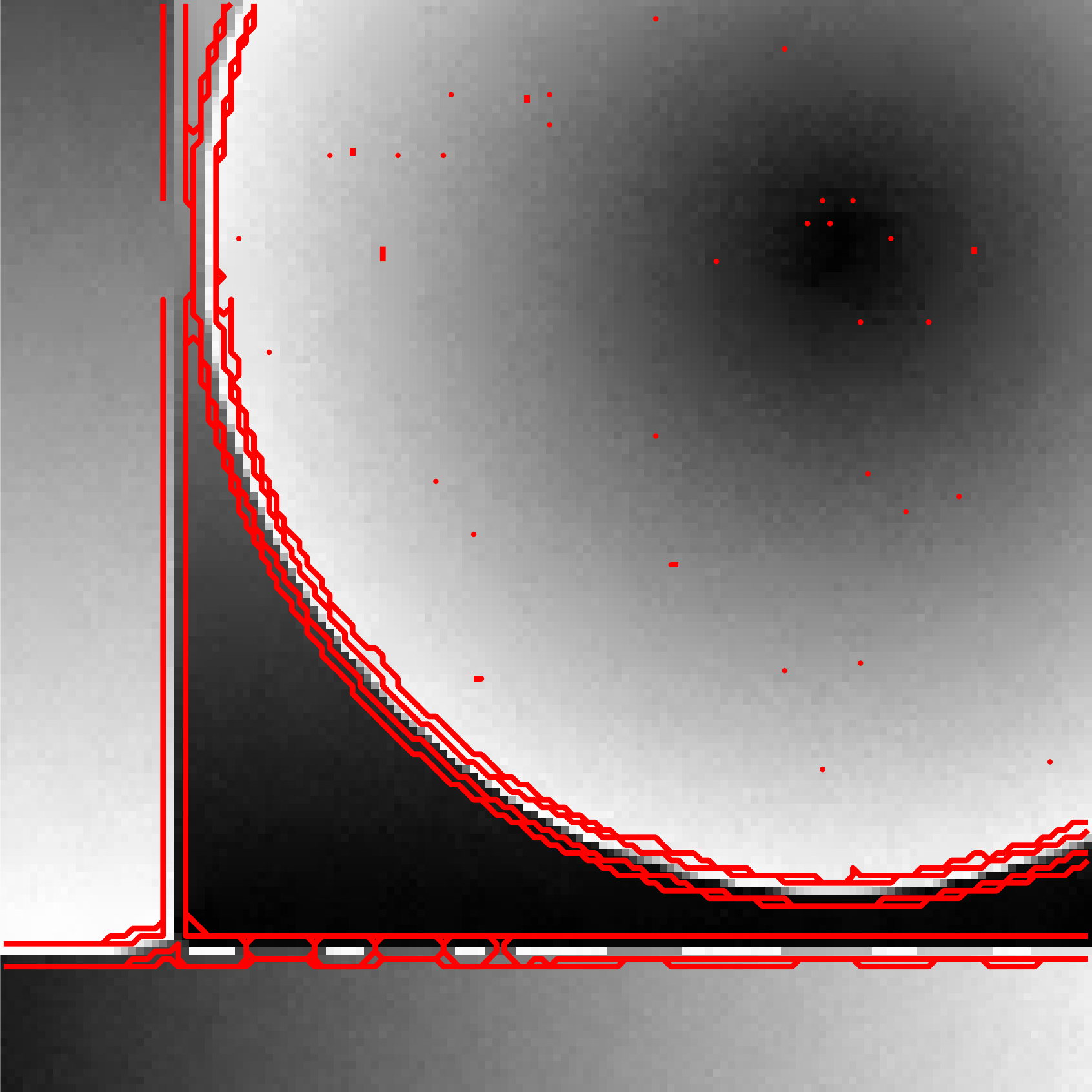}
\end{minipage}\vspace{0.25em}\\
\begin{minipage}{4.2cm}
\includegraphics[width=4.2cm]{AngryBirdsCorruptedZoom.pdf}
\end{minipage}
\begin{minipage}{4.2cm}
\includegraphics[width=4.2cm]{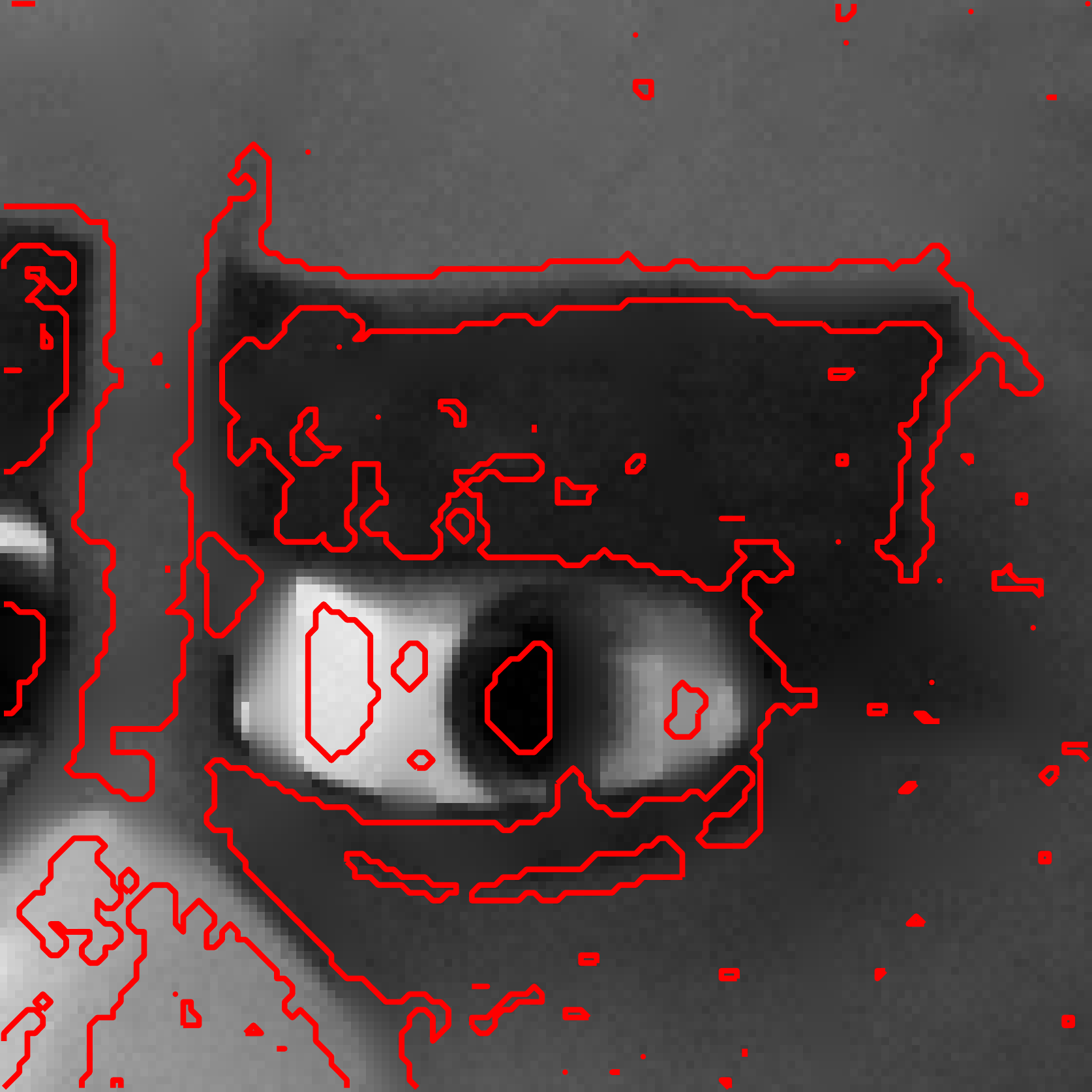}
\end{minipage}
\begin{minipage}{4.2cm}
\includegraphics[width=4.2cm]{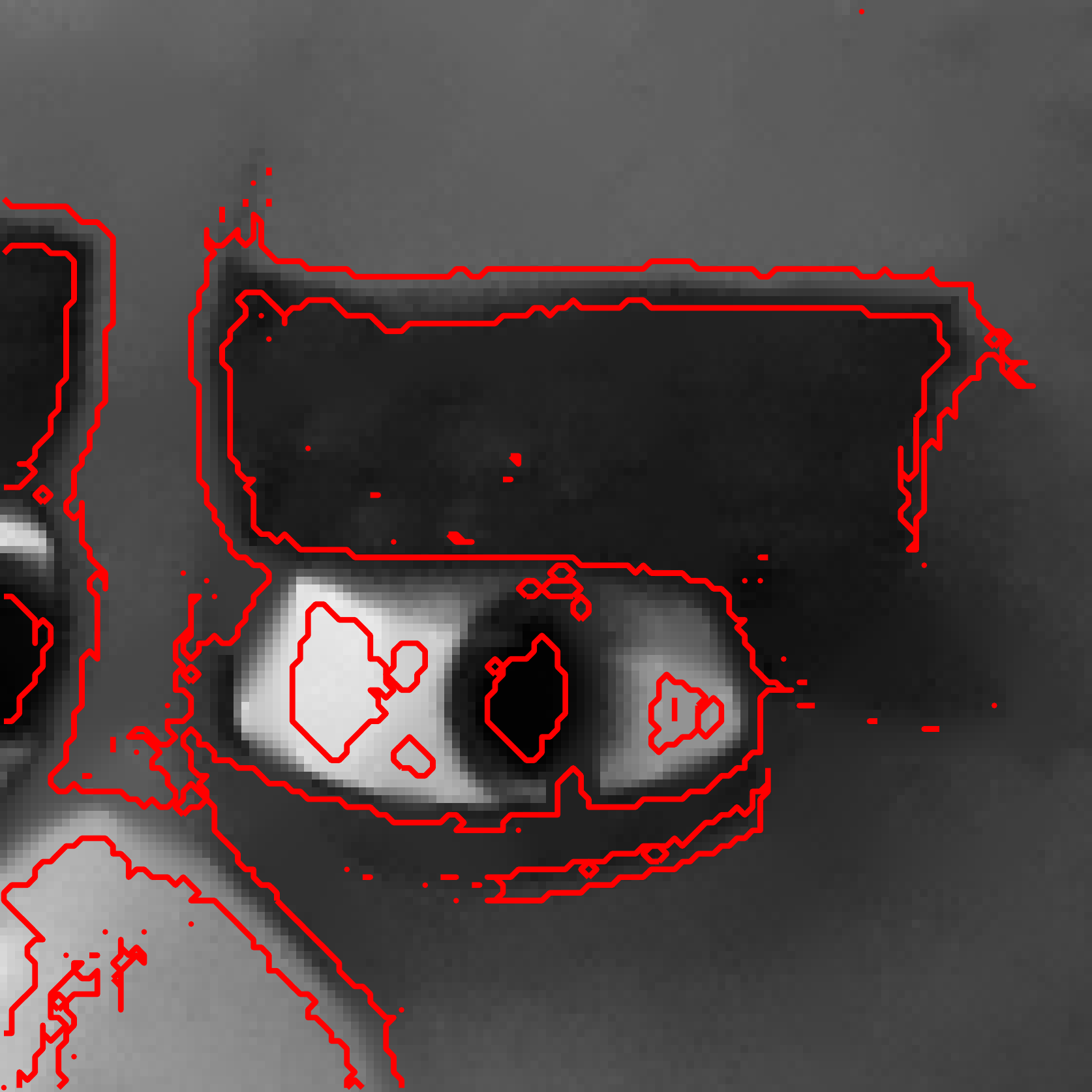}
\end{minipage}\vspace{0.25em}\\
\begin{minipage}{4.2cm}
\includegraphics[width=4.2cm]{PeppersCorruptedZoom.pdf}
\end{minipage}
\begin{minipage}{4.2cm}
\includegraphics[width=4.2cm]{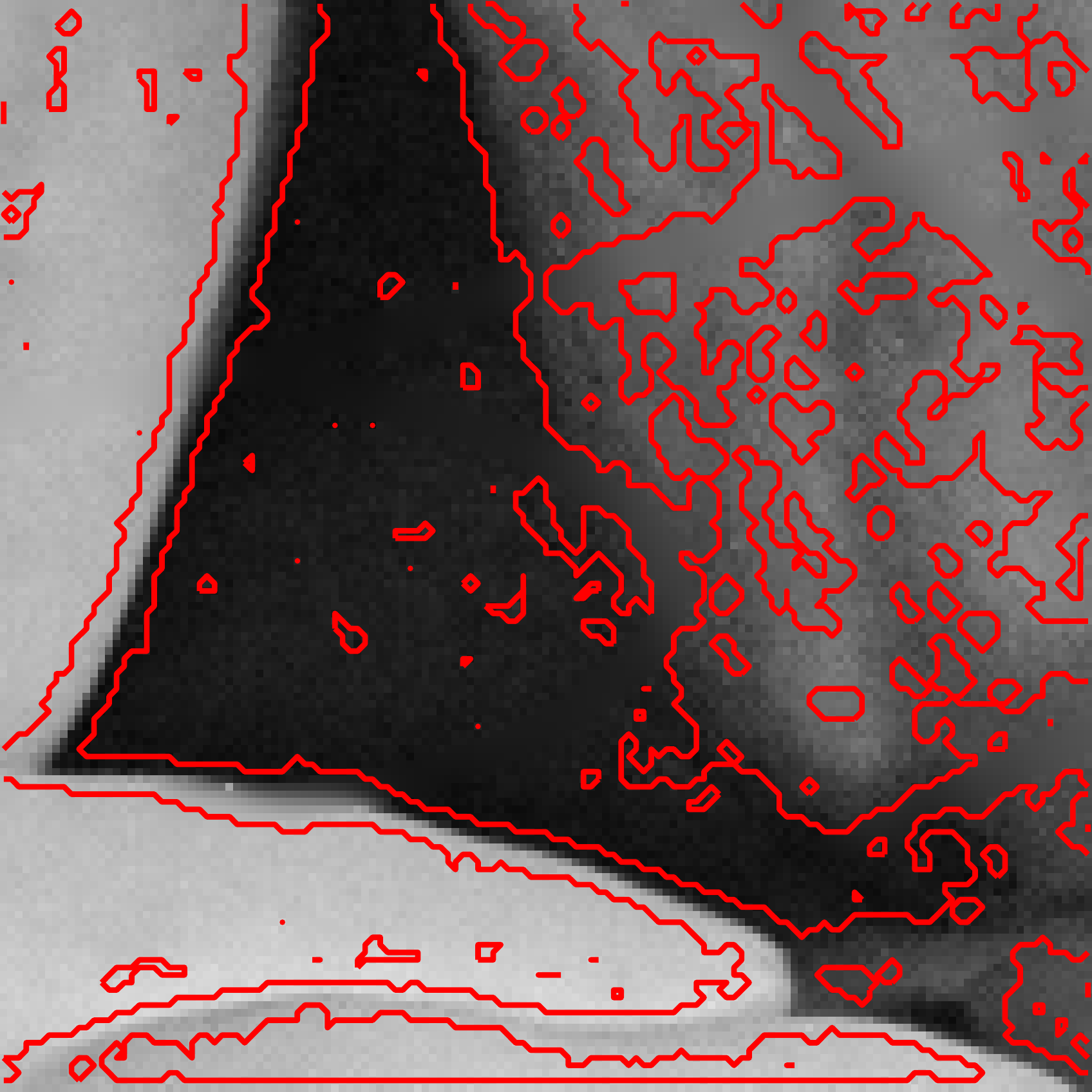}
\end{minipage}
\begin{minipage}{4.2cm}
\includegraphics[width=4.2cm]{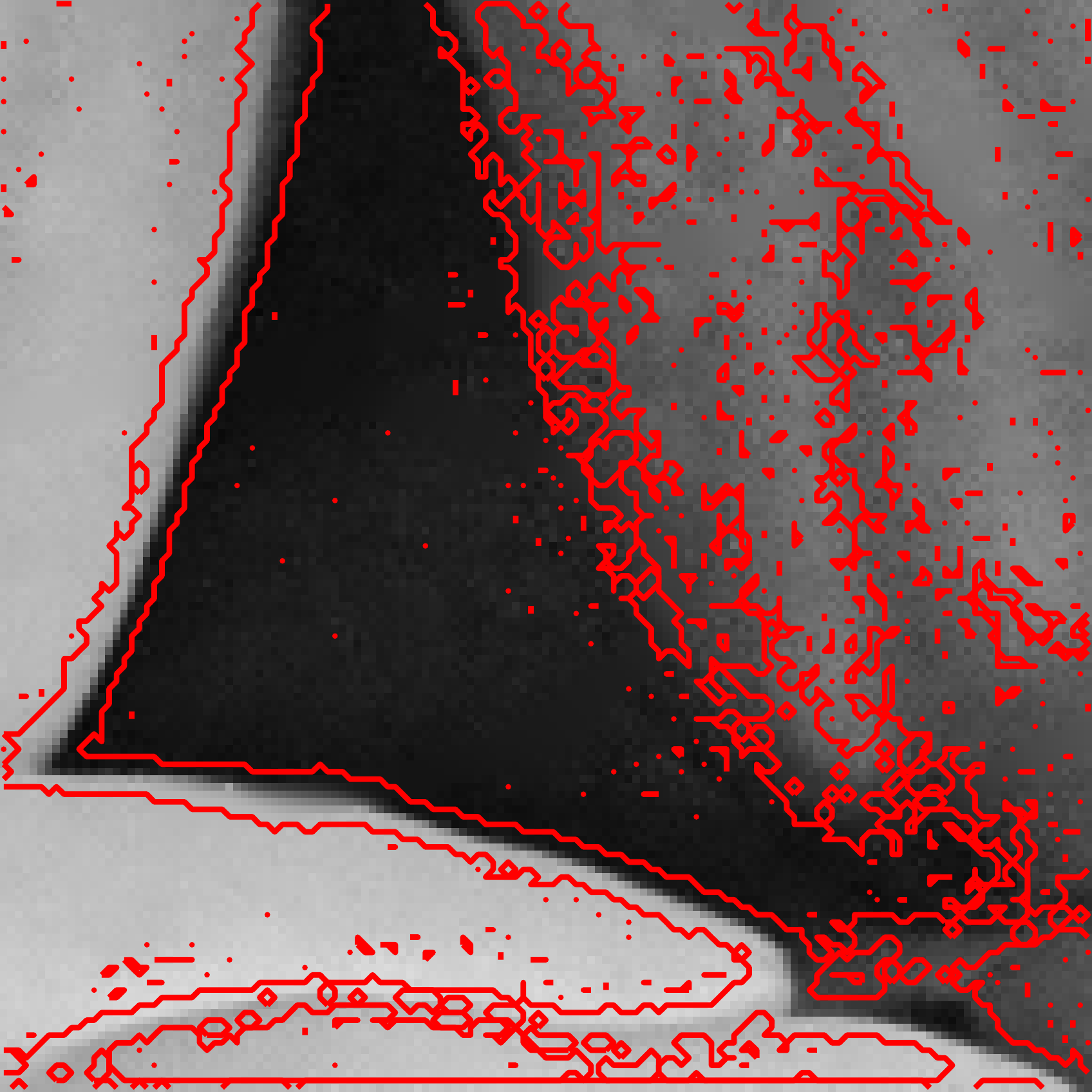}
\end{minipage}\vspace{0.25em}\\
\begin{minipage}{4.2cm}\begin{center}{\small{Observed}}\end{center}\end{minipage}\begin{minipage}{4.2cm}\begin{center}{\small{PS Model \cite{J.F.Cai2016}}}\end{center}\end{minipage}\begin{minipage}{4.2cm}\begin{center}{\small{Our Model \eqref{OurModel}}}\end{center}\end{minipage}
\caption{Zoom-in views of Figure \ref{fig:InpaintJumpEstiComp}. The estimated jump sets are marked by red curves.}\label{fig:InpaintJumpEstiCompZoom}
\end{center}
\end{figure}


\subsubsection{Image Deblurring}\label{SimulationDeblur}

\begin{figure}[ht]
\begin{center}
\begin{minipage}{3.00cm}
\includegraphics[width=3.00cm]{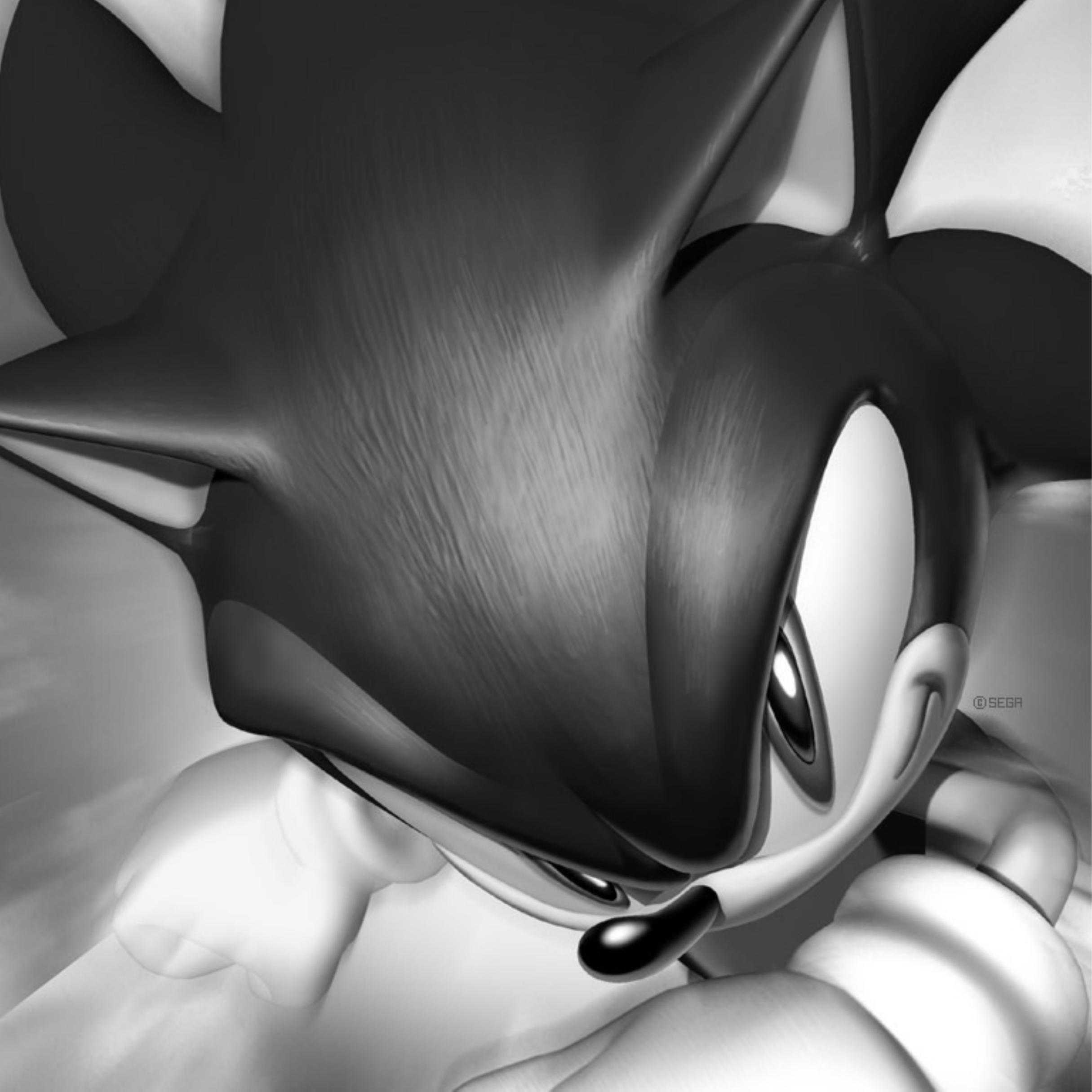}
\end{minipage}
\begin{minipage}{4.5cm}
\includegraphics[height=3.00cm]{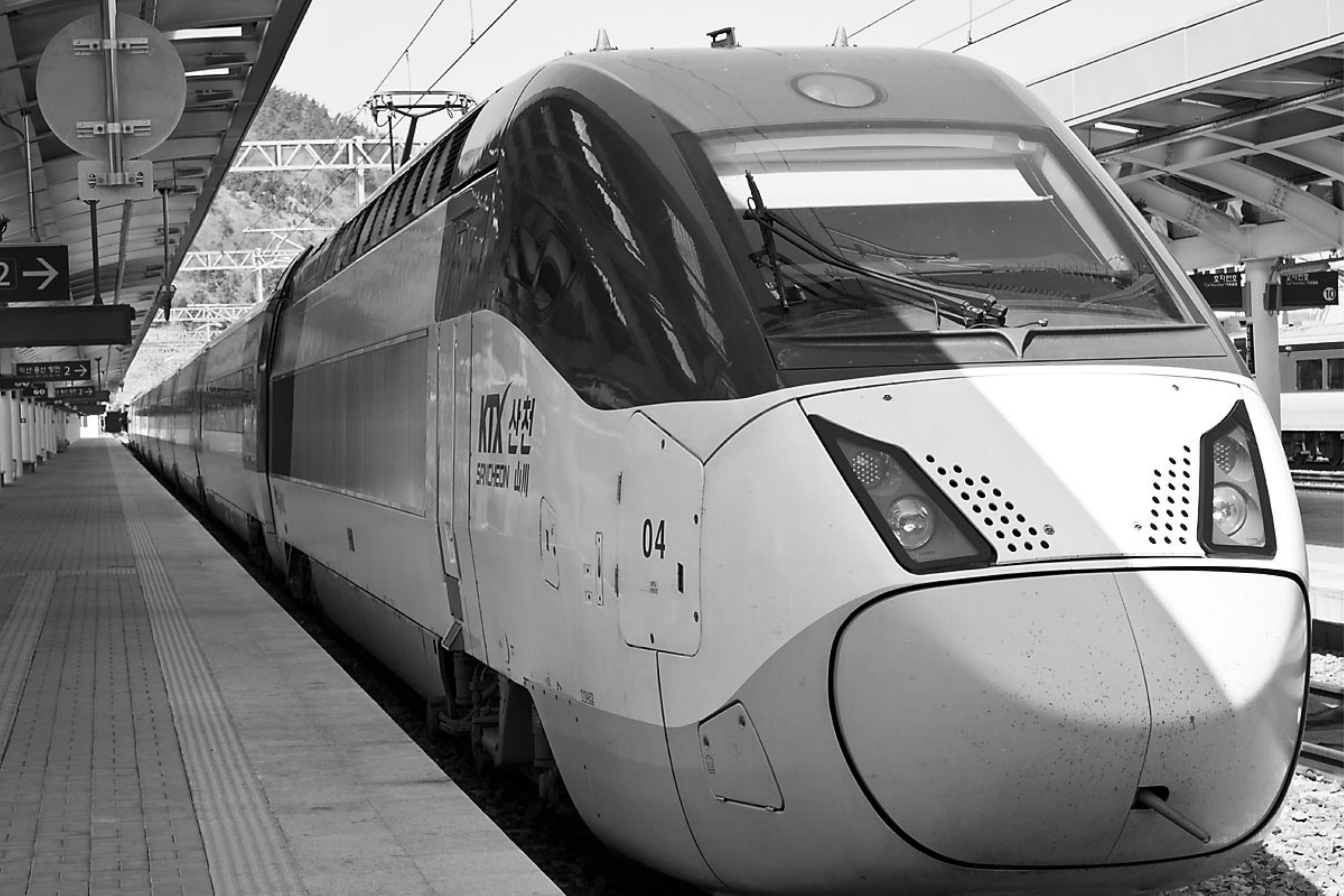}
\end{minipage}
\begin{minipage}{3.75cm}
\includegraphics[height=3.00cm]{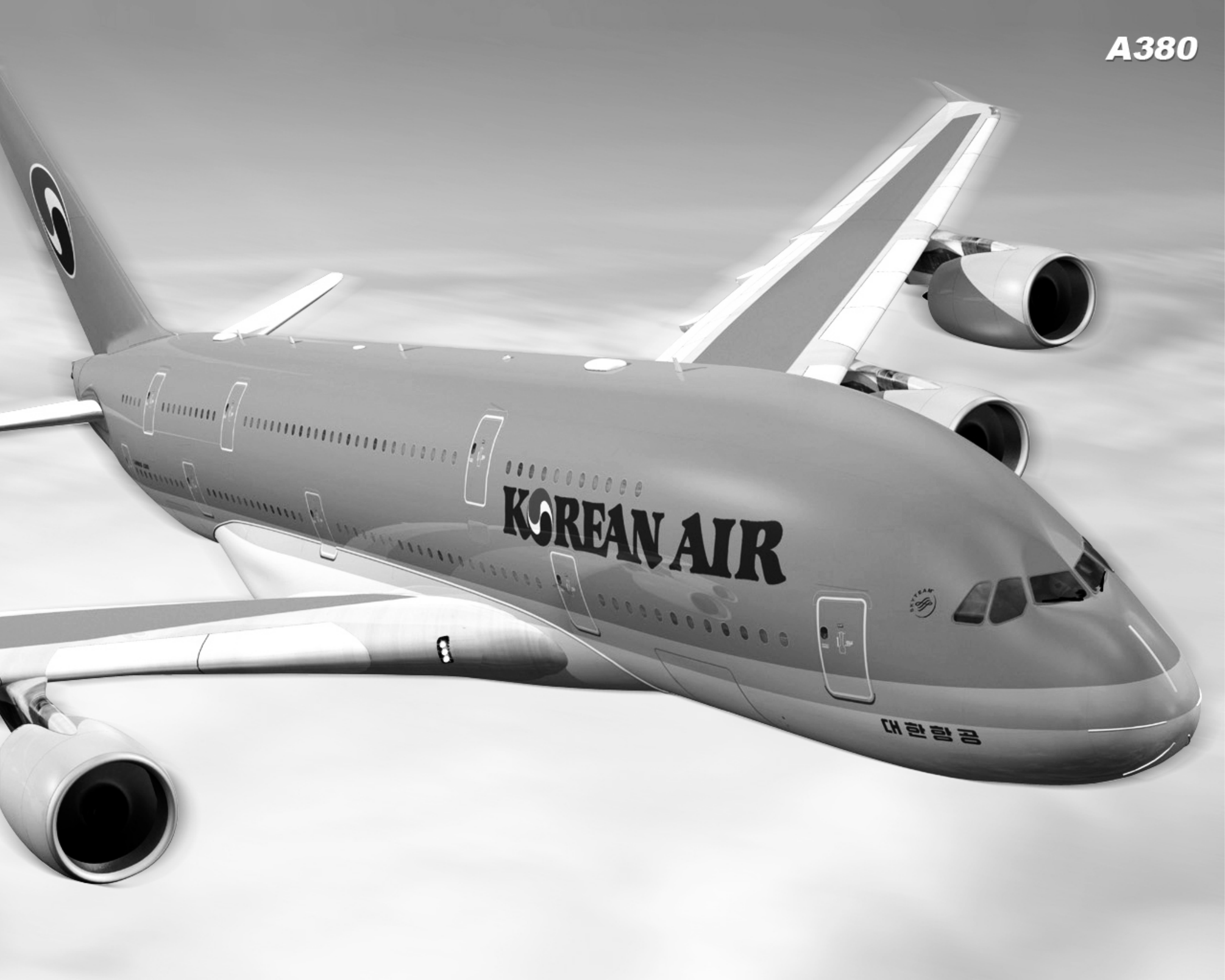}
\end{minipage}
\begin{minipage}{2.4cm}
\includegraphics[height=3.00cm]{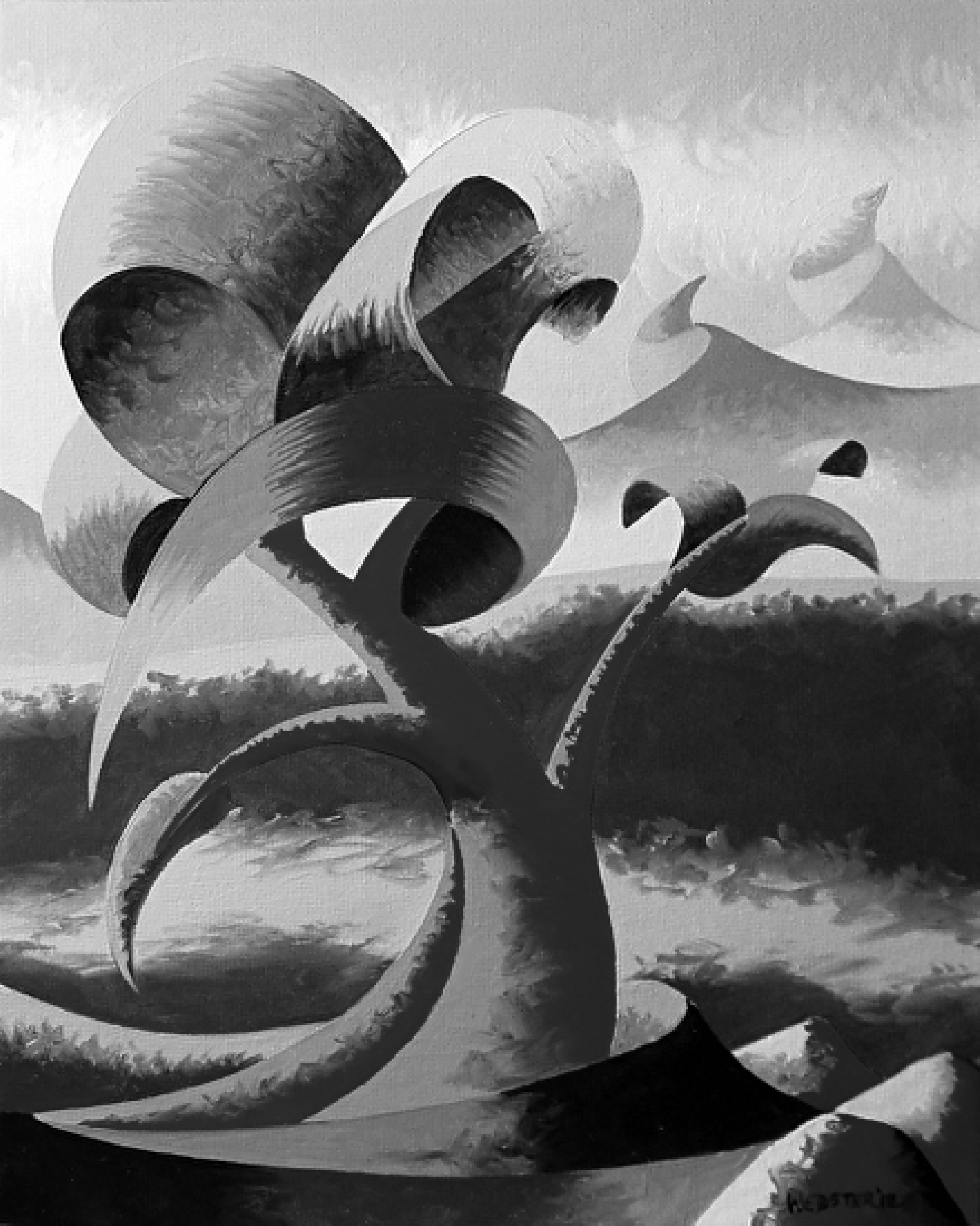}
\end{minipage}
\begin{minipage}{2.25cm}
\includegraphics[height=3.00cm]{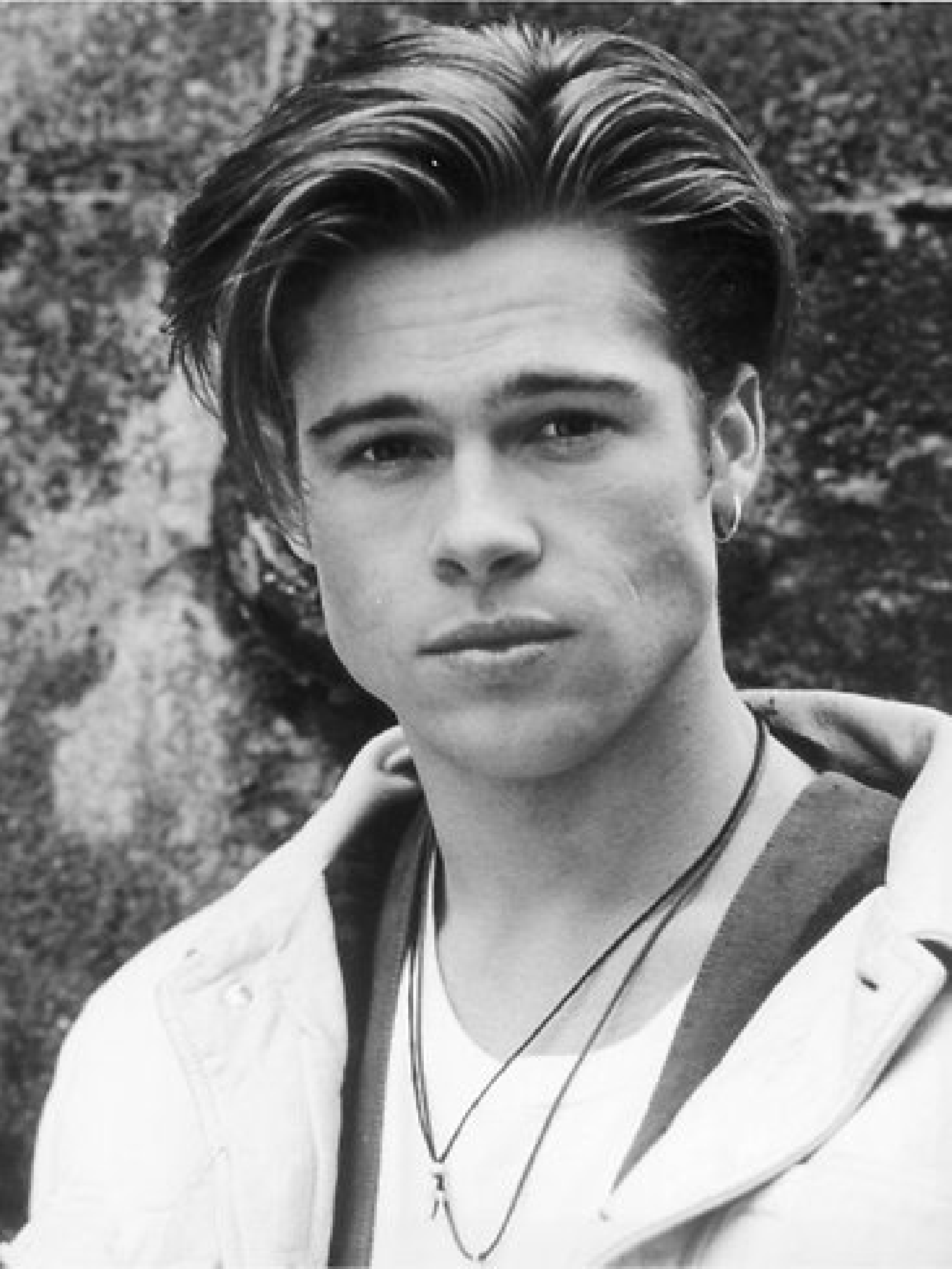}
\end{minipage}\vspace{0.25em}\\
\begin{minipage}{3.00cm}
\includegraphics[width=3.00cm]{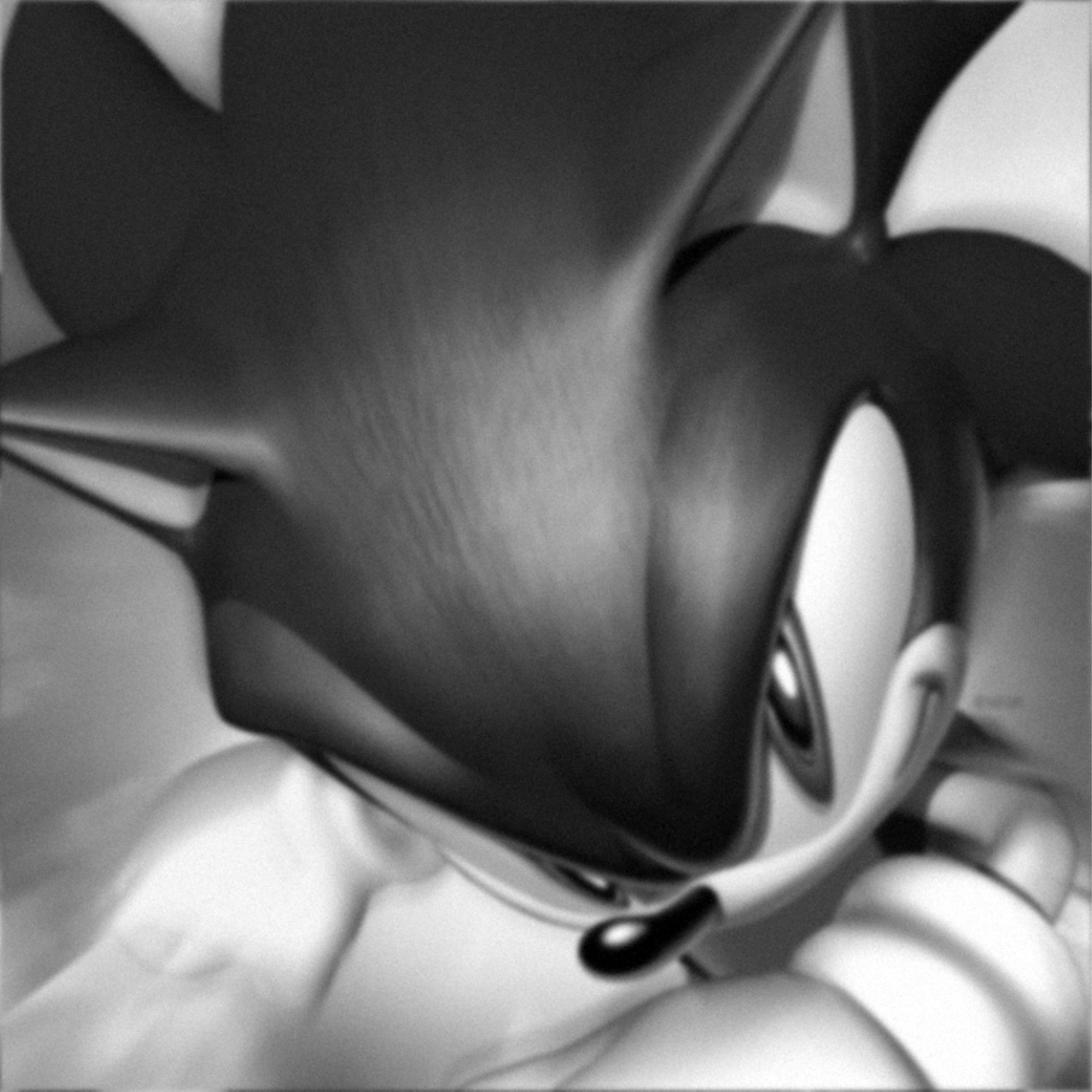}
\end{minipage}
\begin{minipage}{4.5cm}
\includegraphics[height=3.00cm]{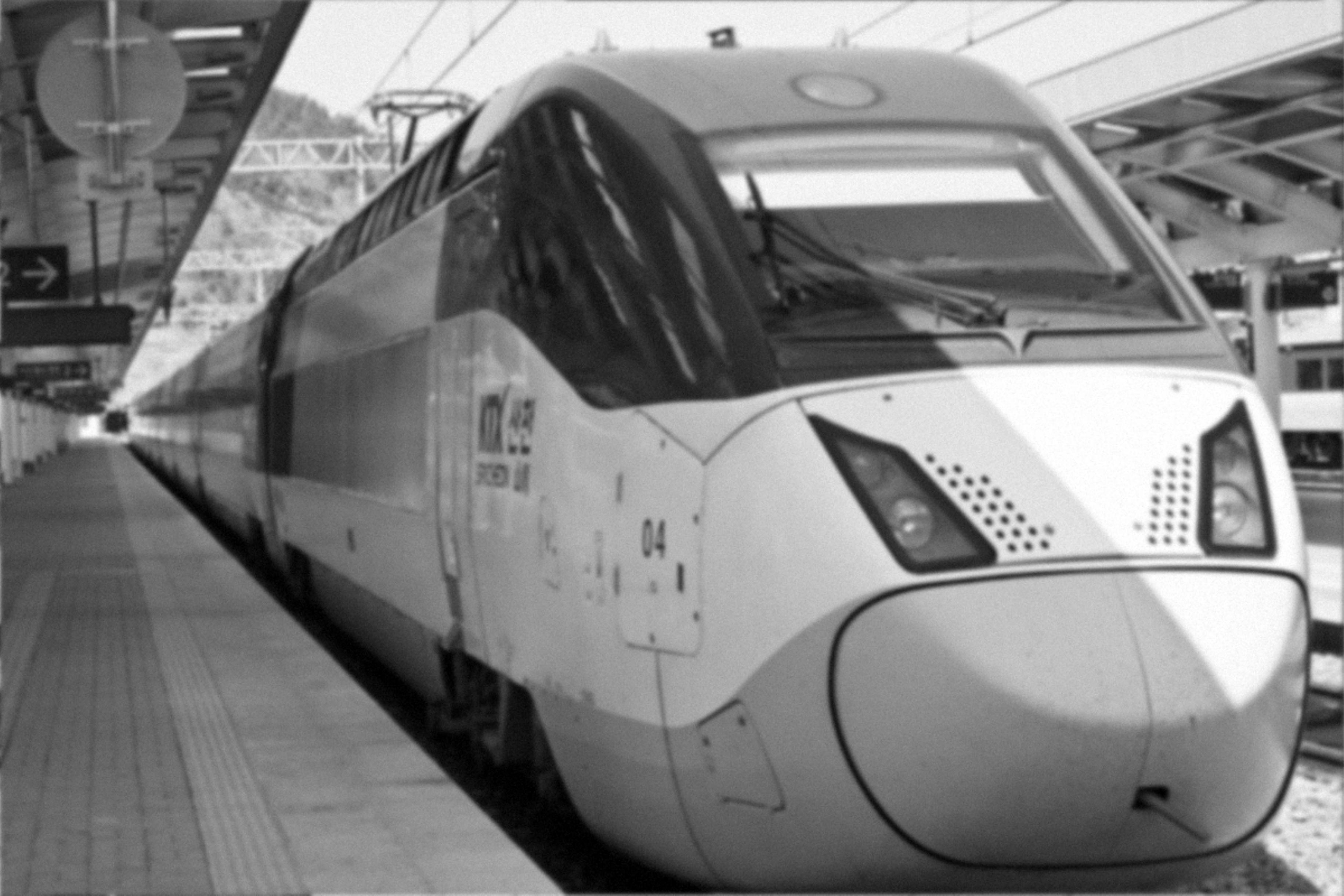}
\end{minipage}
\begin{minipage}{3.75cm}
\includegraphics[height=3.00cm]{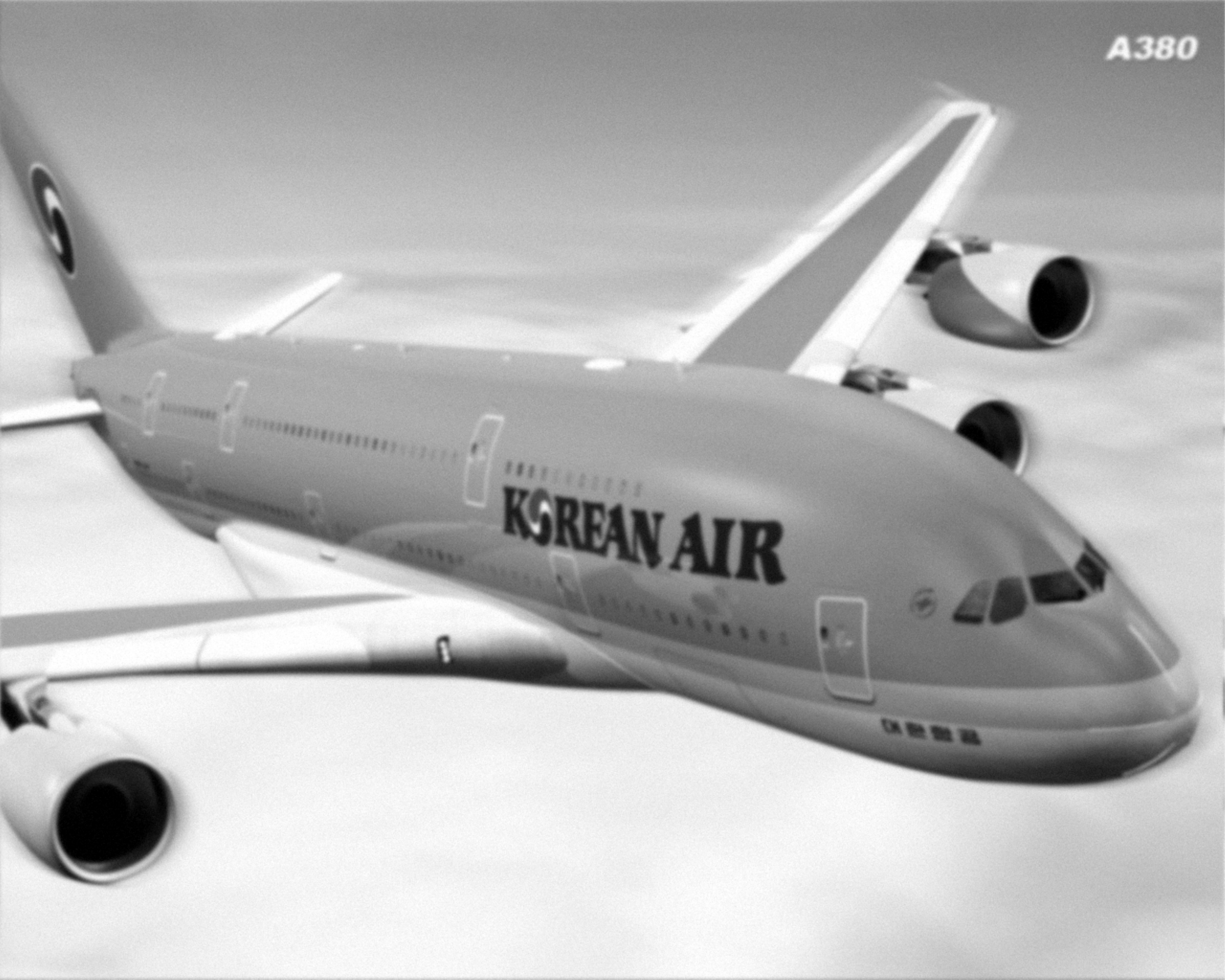}
\end{minipage}
\begin{minipage}{2.4cm}
\includegraphics[height=3.00cm]{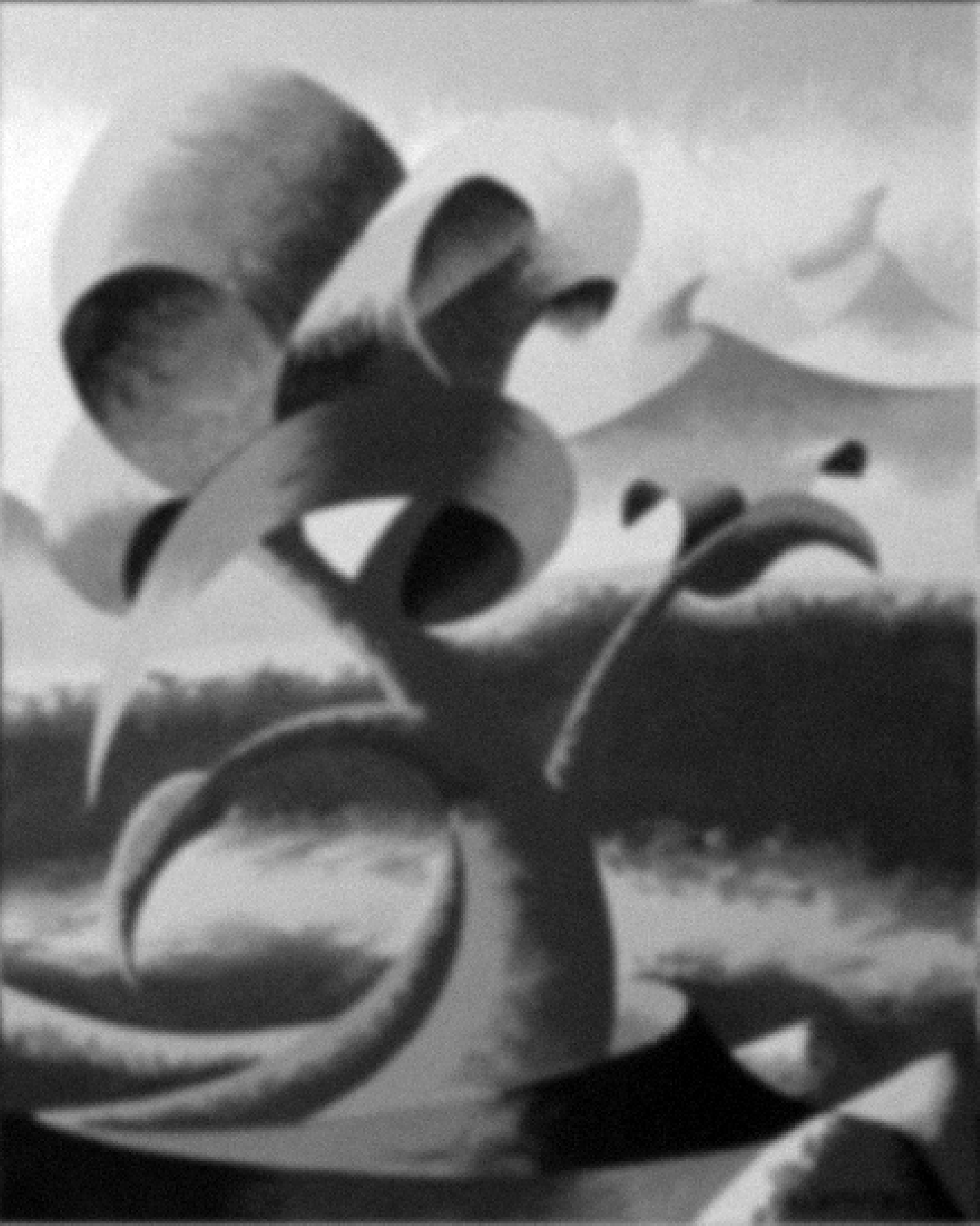}
\end{minipage}
\begin{minipage}{2.25cm}
\includegraphics[height=3.00cm]{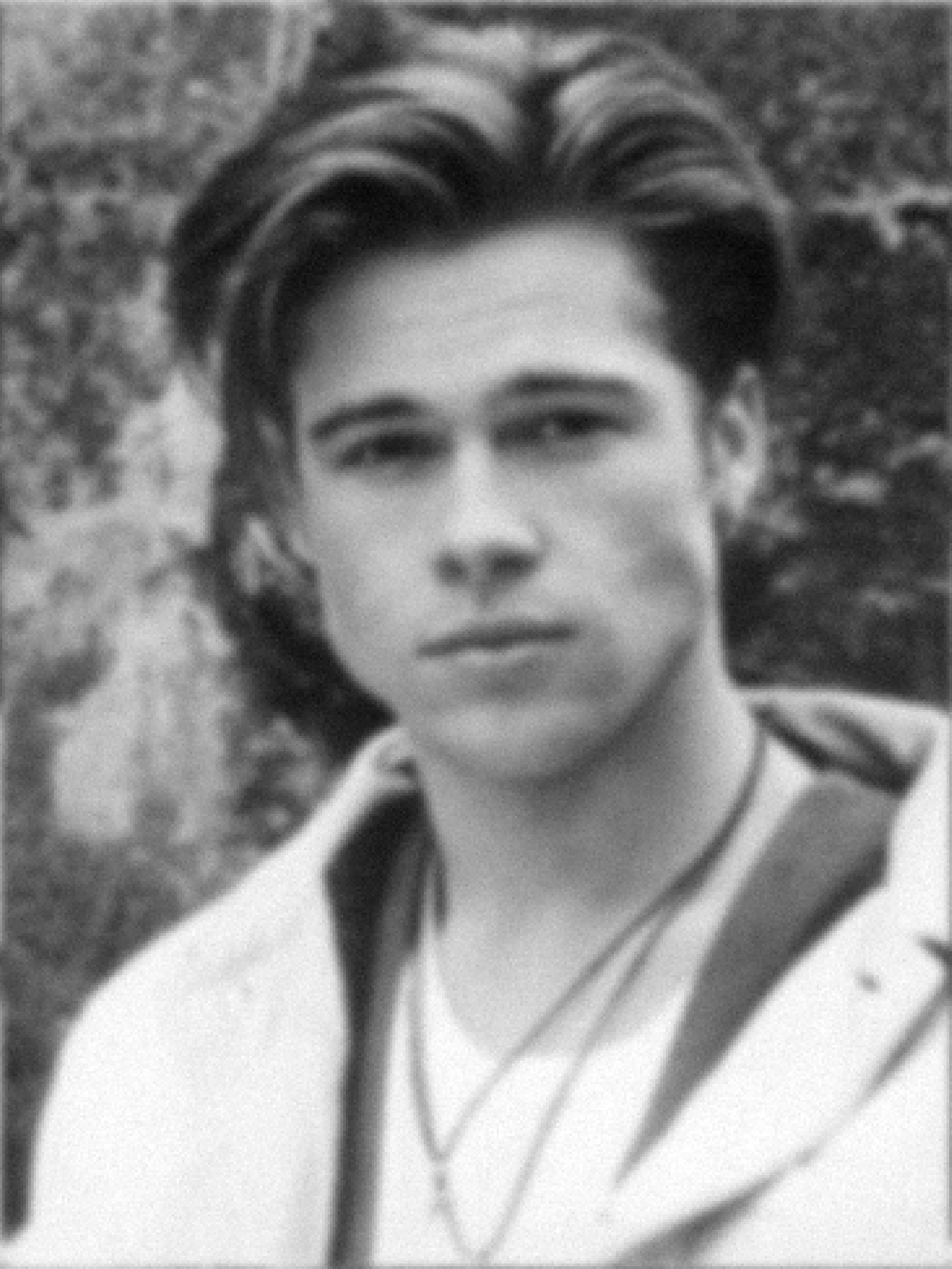}
\end{minipage}\vspace{0.25em}\\
\begin{minipage}{3.00cm}\begin{center}{\small{Sonic}}\end{center}\end{minipage}\begin{minipage}{4.5cm}\begin{center}{\small{Train}}\end{center}\end{minipage}\begin{minipage}{3.75cm}\begin{center}{\small{Airplane}}\end{center}\end{minipage}\begin{minipage}{2.4cm}\begin{center}{\small{Oil Painting}}\end{center}\end{minipage}\begin{minipage}{2.25cm}\begin{center}{\small{Pitt}}\end{center}\end{minipage}
\caption{Visualization of original images and the observed images.}\label{fig:DeblurOriginalMeasure}
\end{center}
\end{figure}

For  image deblurring, five images are tested, as shown in Figure \ref{fig:DeblurOriginalMeasure}. We refer to these images as ``Sonic'', ``Train'', ``Airplane'', ``Oil Painting'', and ``Pitt'' respectively. The algorithm is initialized by choosing $\bu^0=\0$. For $\bsv^0$, we first compute the initial guess of the singularity set from the degraded measurement $\bsf$:
\begin{align*}
\bsh_l[\bk]=\left(\sum_{\bsi\in\BB}\left|\big(\wt{\bsW}_{l,\bsi}\bsf\big)[\bk]\right|^2\right)^{\f{1}{2}}
\end{align*}
where $\wt{\bsW}$ is chosen to be the piecewise cubic B-spline wavelet frame with $2$ levels of decomposition. Then we compute $\bsv^0=(\bsv_0^0,\cdots,\bsv_{L-1}^0)$ by
\begin{align*}
\Sig_l^0=\big\{\bk\in\OO:\bsh_l[\bk]/\|\bsh_l\|_{\infty}\geq\tau_l\big\}~~~\text{so that}~~\bsv_l^0=\one_{\Sig_l^0}~~~l=0,\cdots,L-1.
\end{align*}
Throughout our numerical experiments, we set $\tau_l=\tau=0.15$ for $0\leq l\leq L-1$. (Note, however, that the reconstruction results are relatively insensitive to the choice of $\wt{\bsW}$ and $\tau$.) The level of decomposition for $\bsW$, $\bsW'$, and $\bsW''$ are all chosen to be $2$. For the PS model \eqref{CDS} and the GS model \eqref{JLS}, the piecewise linear B-spline wavelet frame with $2$ levels of decomposition are used. The parameters in \eqref{OurModel} are chosen in the same way as the image inpainting, and the parameters in all models are manually chosen for the optimal recovery results.

\begin{table}[ht]
\begin{center}
\begin{tabular}{|c||c|c|c|c|c|}\hline
Image&Observed&TGV Model \cite{K.Bredies2010}&PS Model \cite{J.F.Cai2016}&GS Model \cite{H.Ji2016}&Our Model \eqref{OurModel}\\ \hline
Sonic&$29.5403$&$35.2665$&$35.5588$&$35.2905$&$\textbf{35.9163}$\\ \hline
Train&$22.3559$&$25.5072$&$25.7154$&$25.3460$&$\textbf{25.8934}$\\ \hline
Airplane&$28.8864$&$33.2598$&$33.6136$&$33.2101$&$\textbf{33.8899}$\\ \hline
Oil Painting&$24.7937$&$28.4969$&$28.8532$&$28.0944$&$\textbf{29.0383}$\\ \hline
Pitt&$24.5077$&$27.7433$&$28.1452$&$28.0778$&$\textbf{28.3069}$\\ \hline
\end{tabular}
\caption{Comparison of the PSNR values of four models for deblurring.}\label{tab:PSNRCompareDeblur}
\end{center}
\end{table}

\begin{figure}[htp!]
\begin{center}
\begin{minipage}{3.00cm}
\includegraphics[width=3.00cm]{SonicCorrupted2.pdf}
\end{minipage}
\begin{minipage}{4.5cm}
\includegraphics[height=3.00cm]{KTXCorrupted2.pdf}
\end{minipage}
\begin{minipage}{3.75cm}
\includegraphics[height=3.00cm]{KEA388Corrupted2.pdf}
\end{minipage}
\begin{minipage}{2.4cm}
\includegraphics[height=3.00cm]{PaintCorrupted2.pdf}
\end{minipage}
\begin{minipage}{2.25cm}
\includegraphics[height=3.00cm]{Pitt2Corrupted2.pdf}
\end{minipage}\vspace{0.25em}\\
\begin{minipage}{3.00cm}
\includegraphics[width=3.00cm]{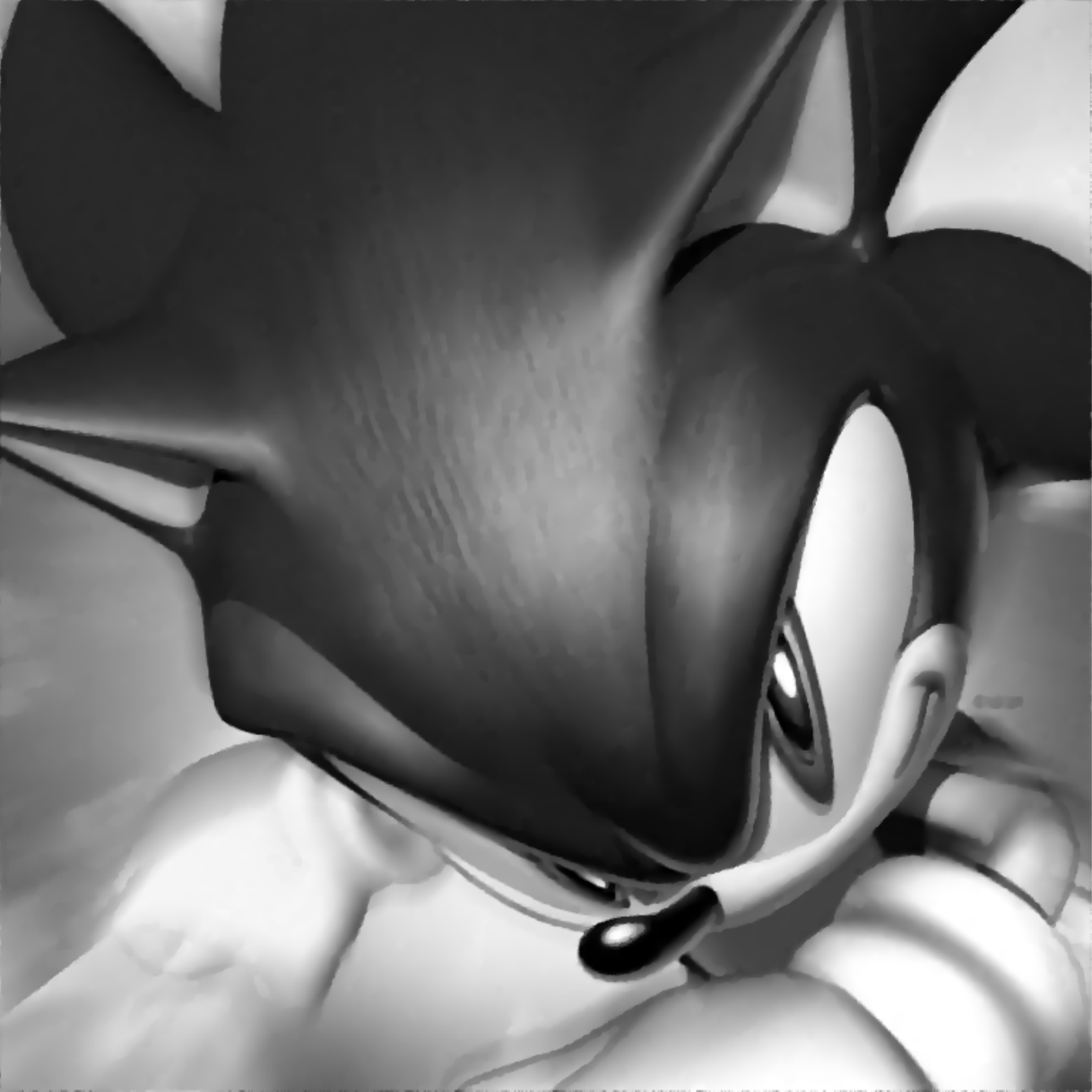}
\end{minipage}
\begin{minipage}{4.5cm}
\includegraphics[height=3.00cm]{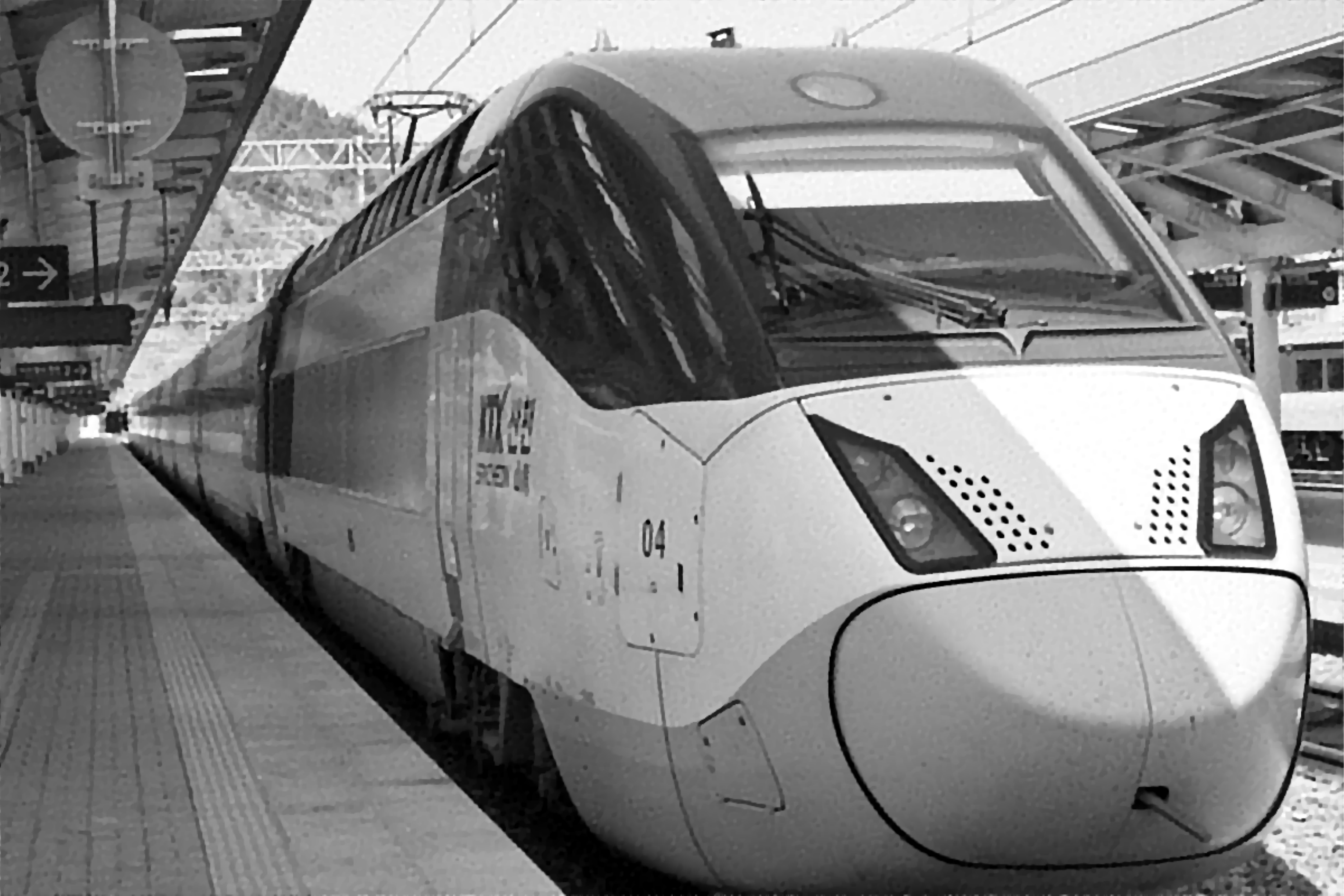}
\end{minipage}
\begin{minipage}{3.75cm}
\includegraphics[height=3.00cm]{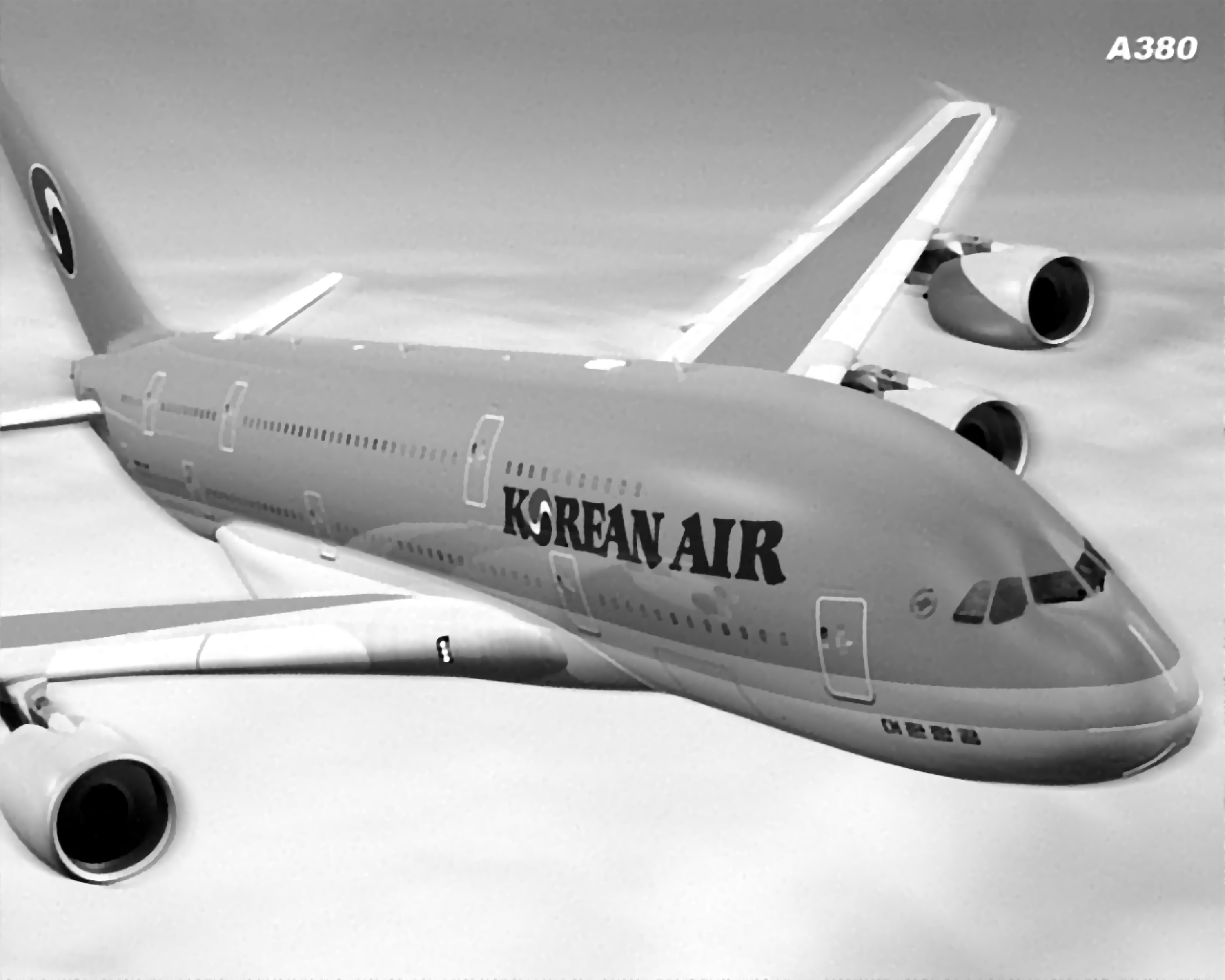}
\end{minipage}
\begin{minipage}{2.4cm}
\includegraphics[height=3.00cm]{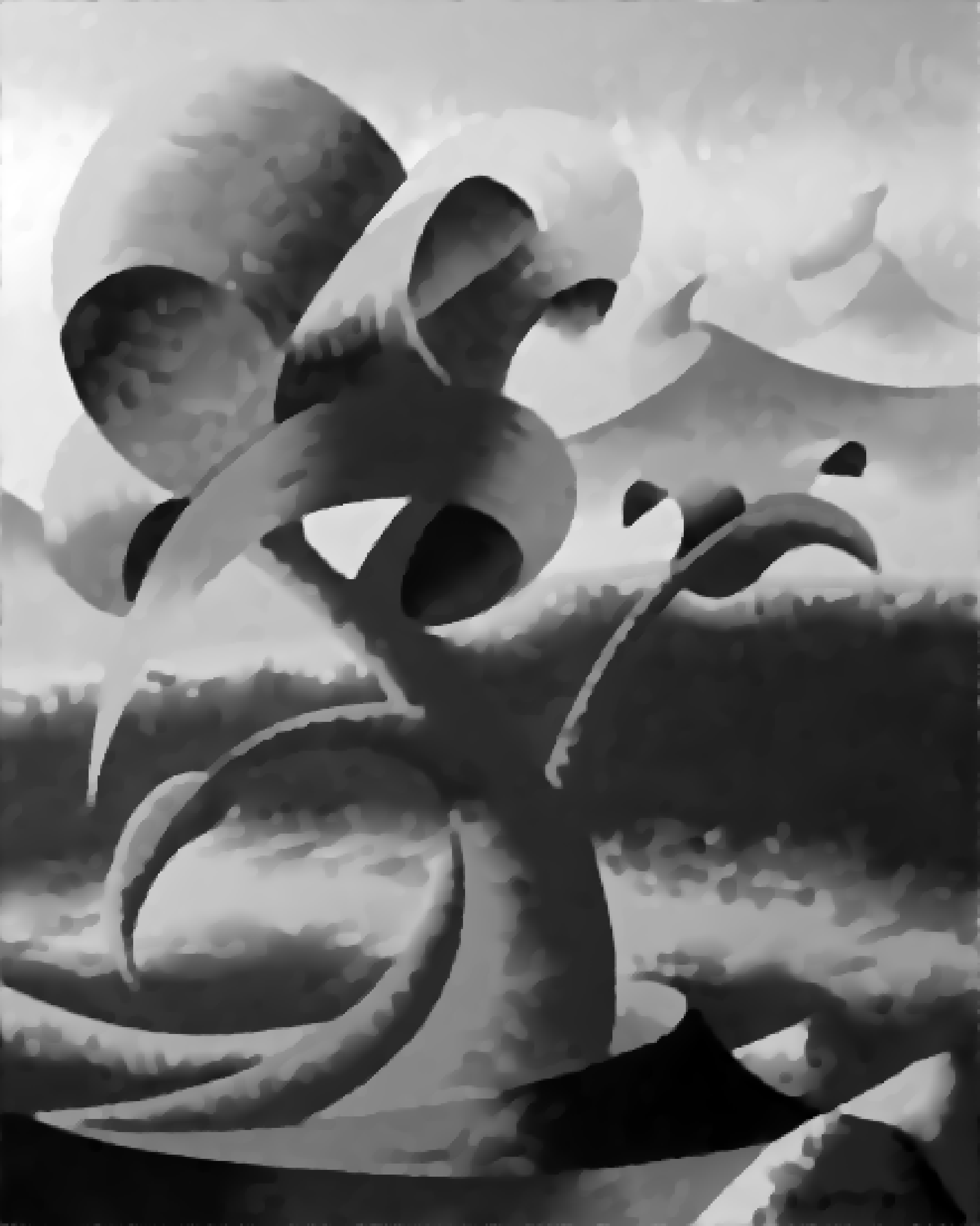}
\end{minipage}
\begin{minipage}{2.25cm}
\includegraphics[height=3.00cm]{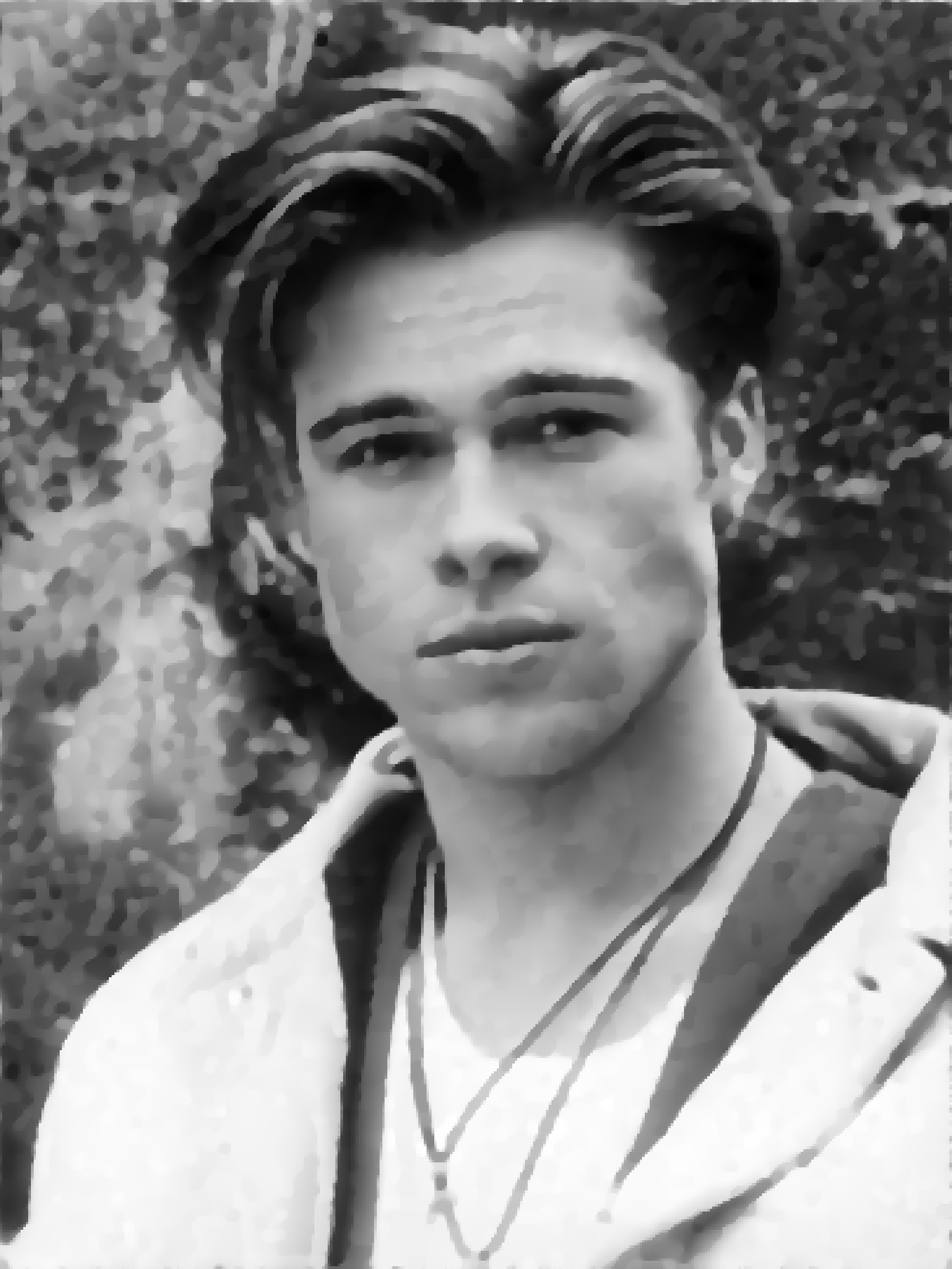}
\end{minipage}\vspace{0.25em}\\
\begin{minipage}{3.00cm}
\includegraphics[width=3.00cm]{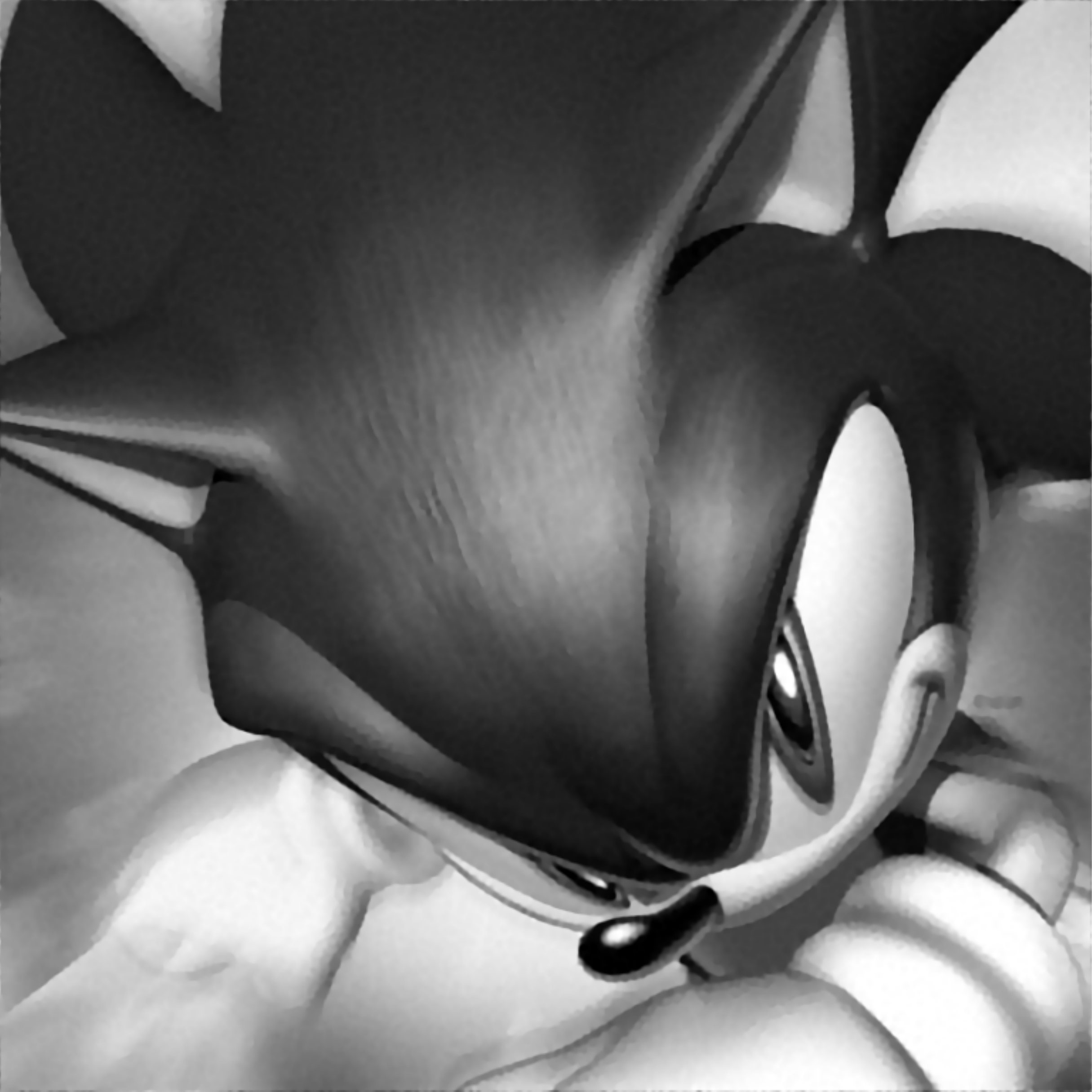}
\end{minipage}
\begin{minipage}{4.5cm}
\includegraphics[height=3.00cm]{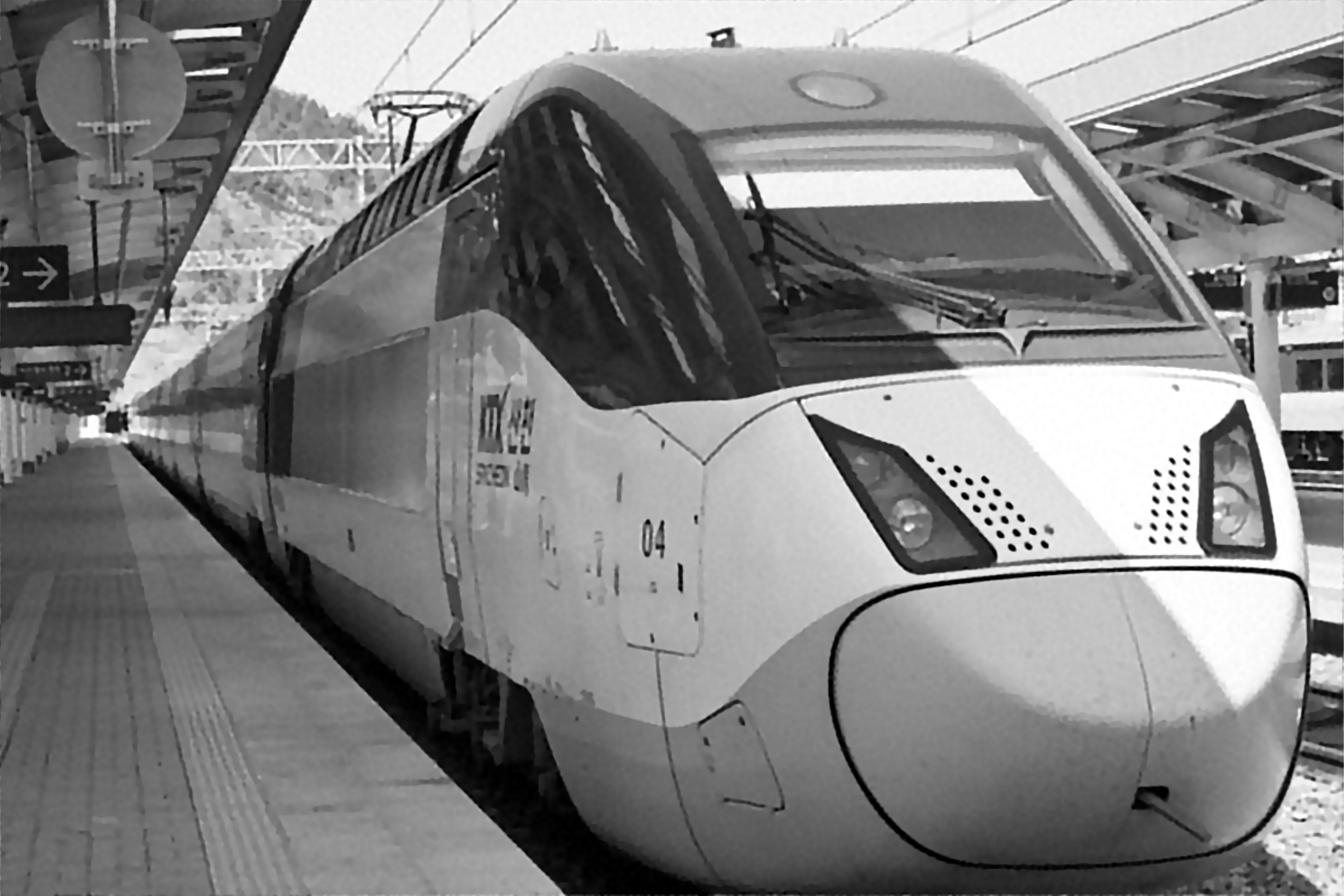}
\end{minipage}
\begin{minipage}{3.75cm}
\includegraphics[height=3.00cm]{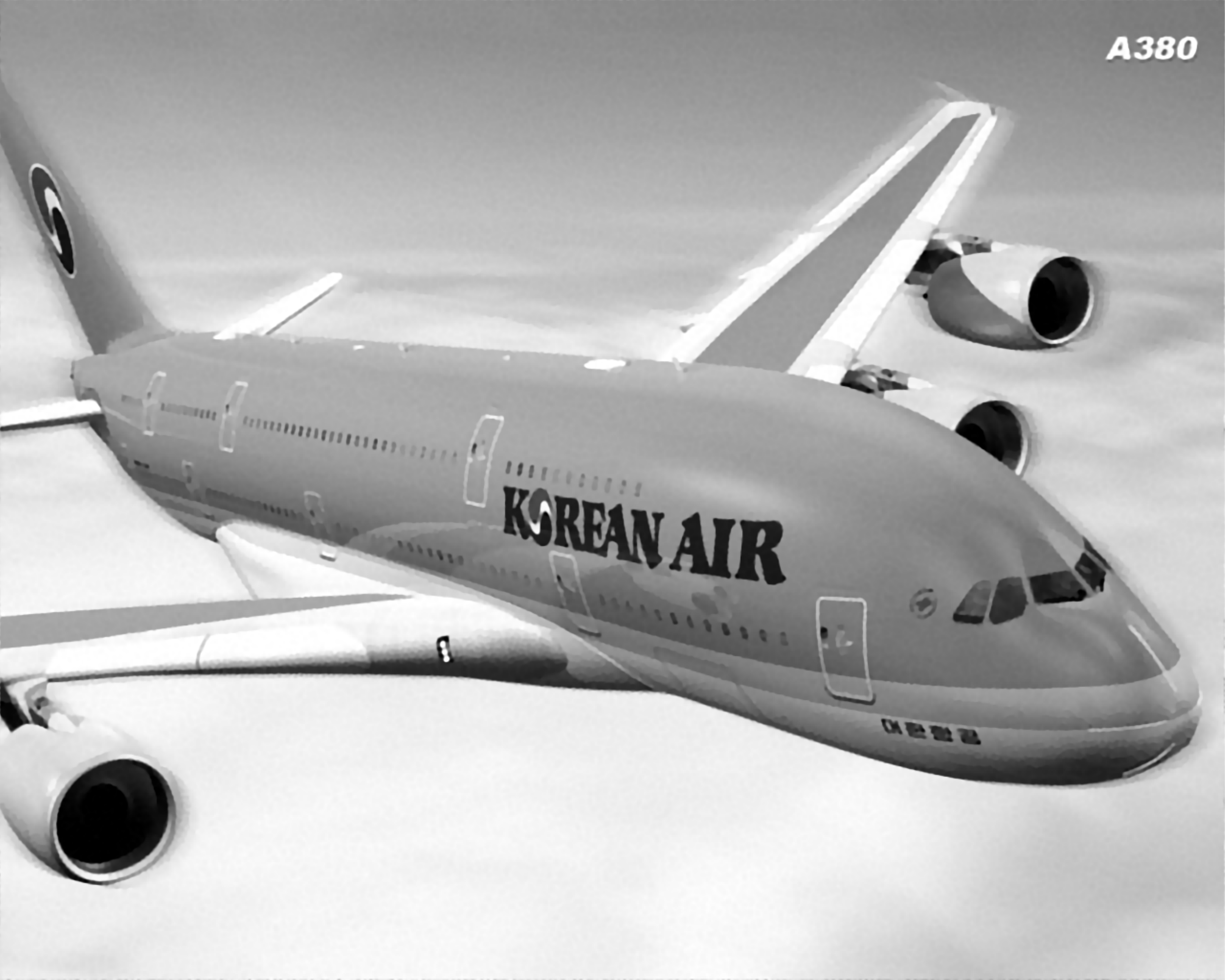}
\end{minipage}
\begin{minipage}{2.4cm}
\includegraphics[height=3.00cm]{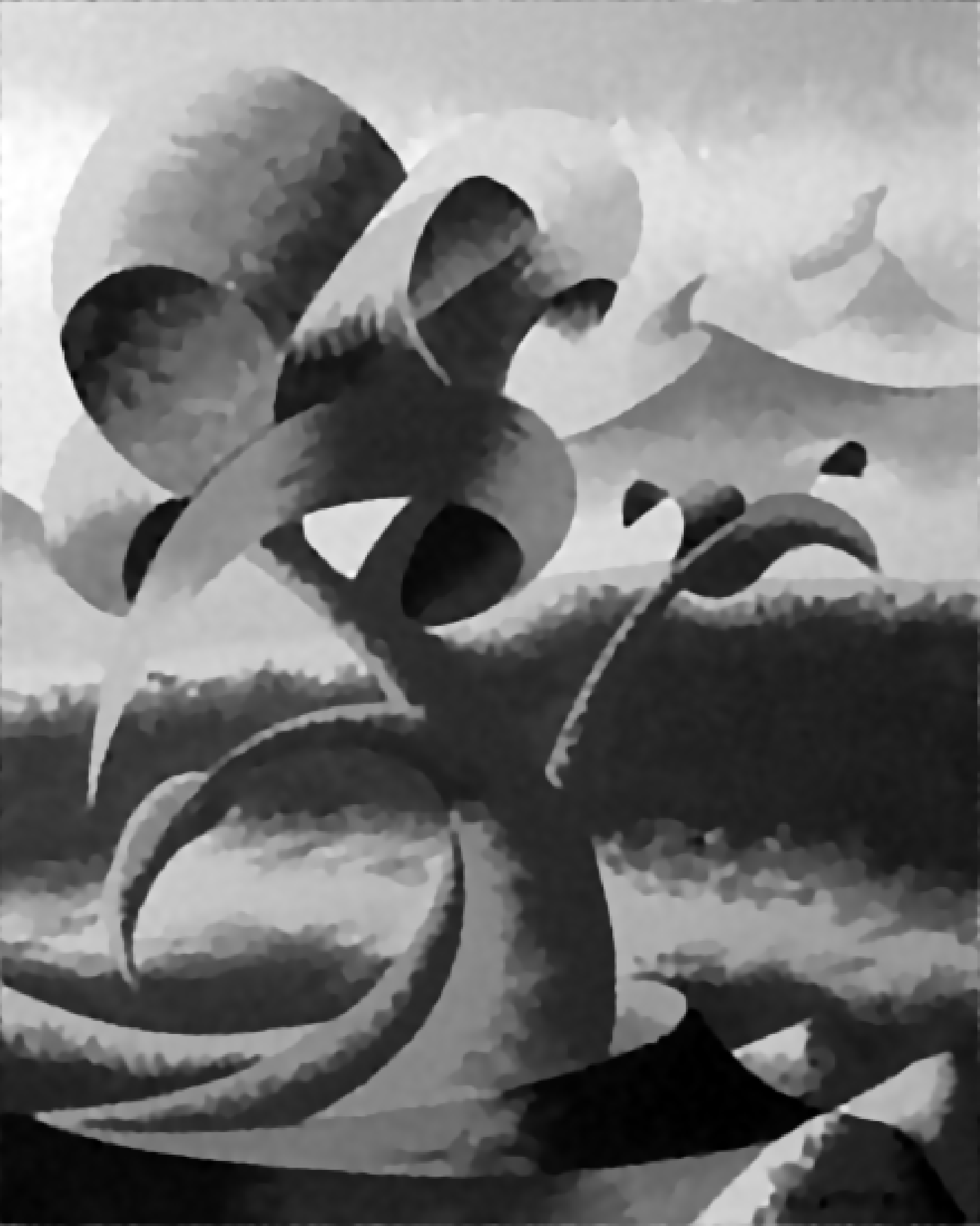}
\end{minipage}
\begin{minipage}{2.25cm}
\includegraphics[height=3.00cm]{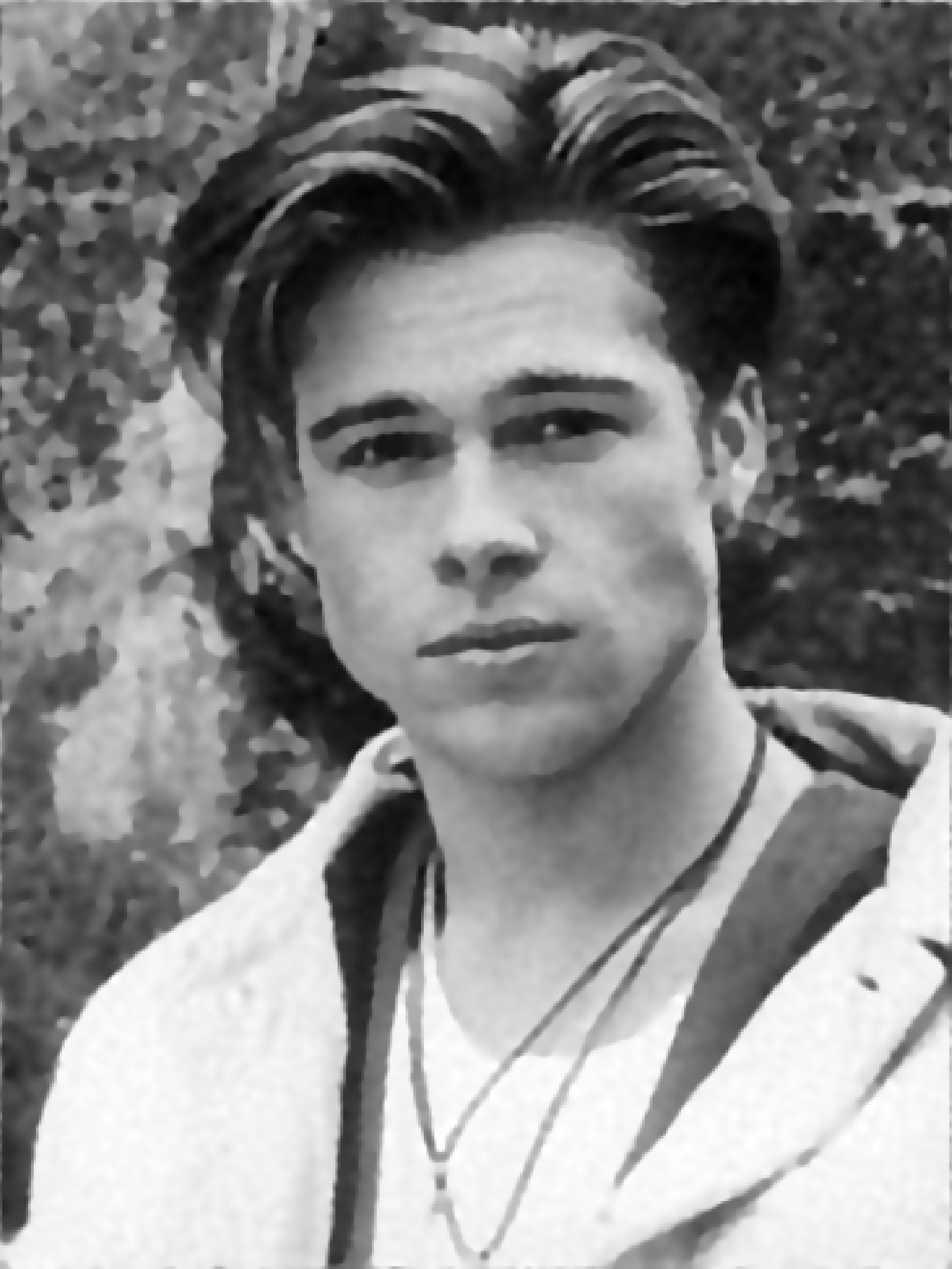}
\end{minipage}\vspace{0.25em}\\
\begin{minipage}{3.00cm}
\includegraphics[width=3.00cm]{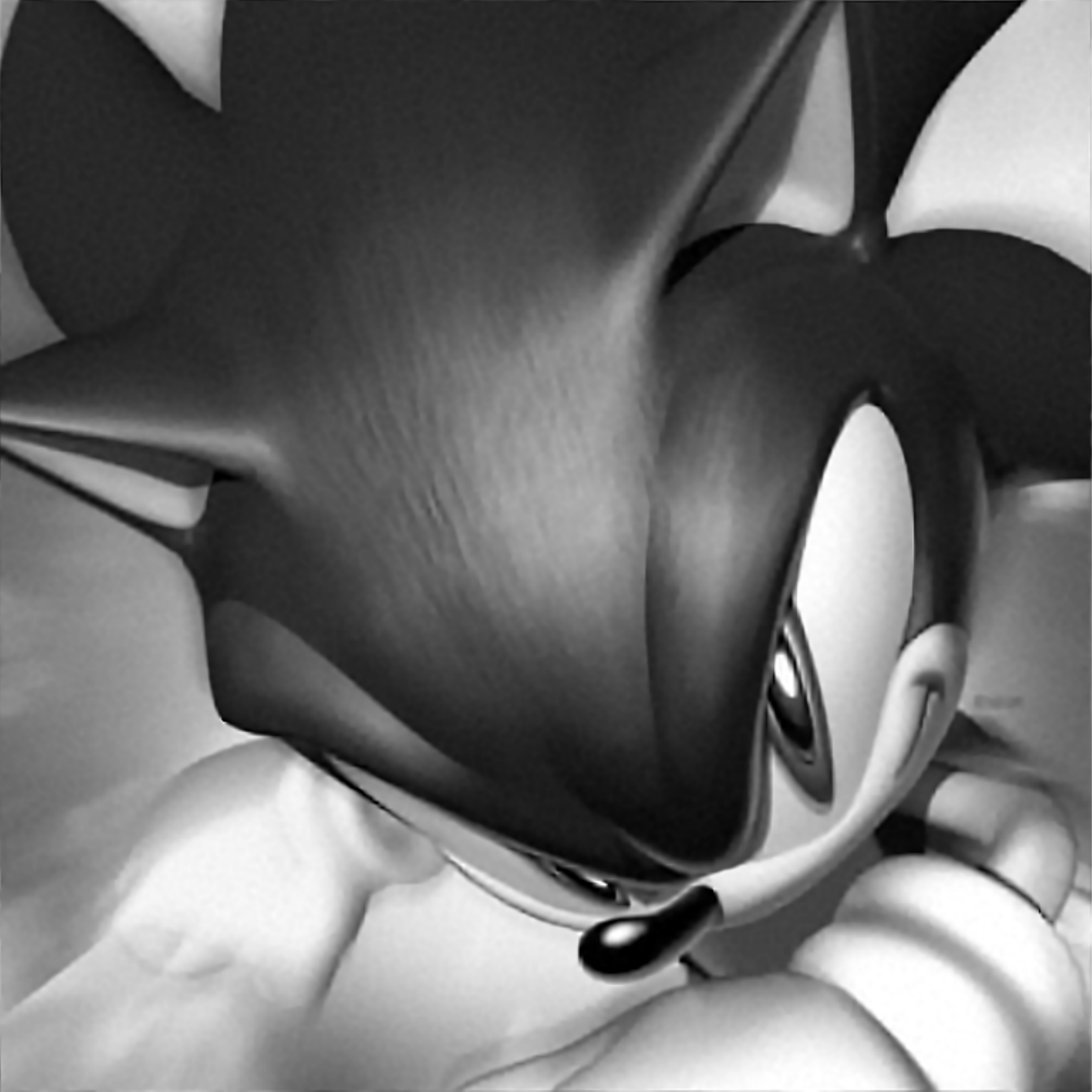}
\end{minipage}
\begin{minipage}{4.5cm}
\includegraphics[height=3.00cm]{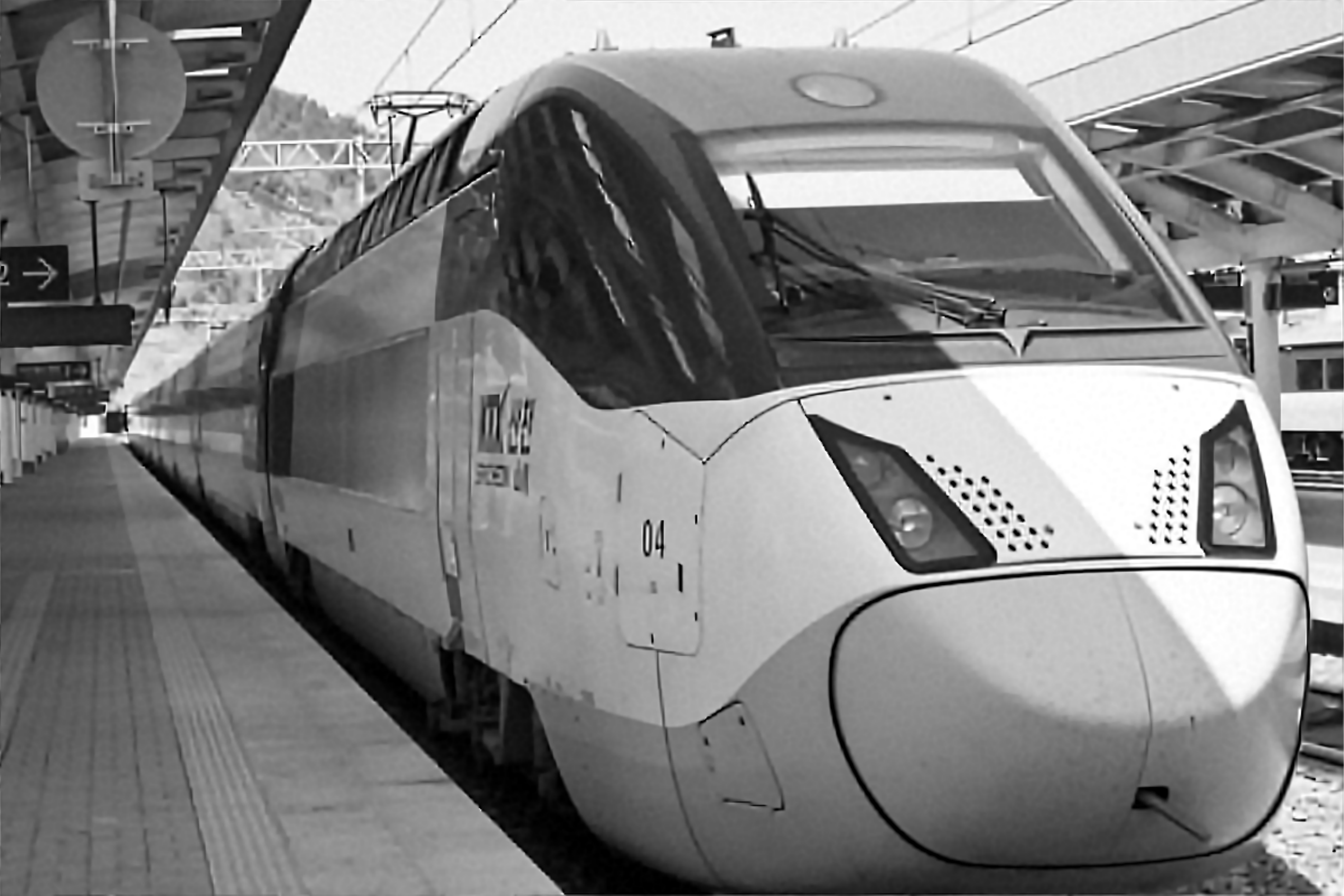}
\end{minipage}
\begin{minipage}{3.75cm}
\includegraphics[height=3.00cm]{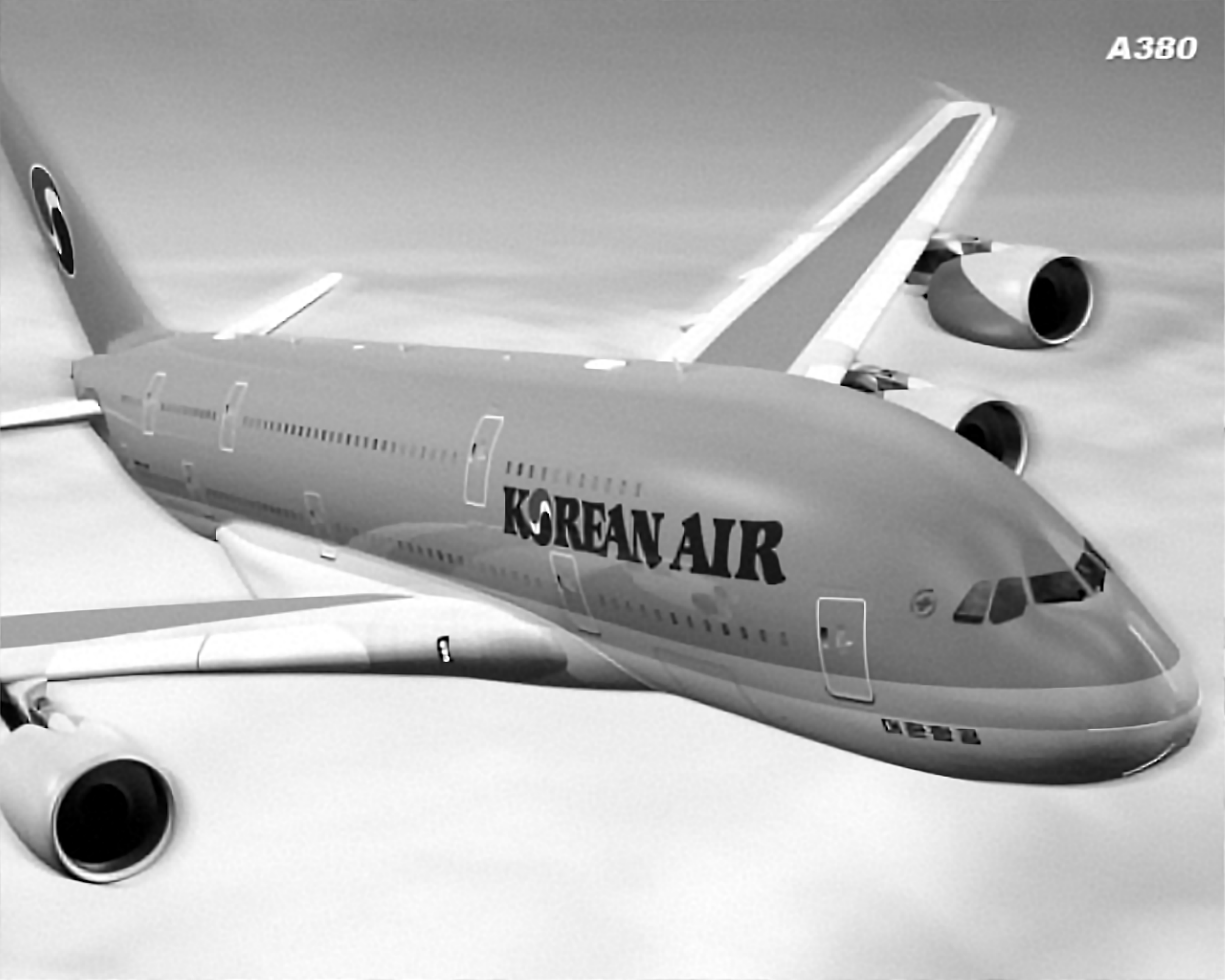}
\end{minipage}
\begin{minipage}{2.4cm}
\includegraphics[height=3.00cm]{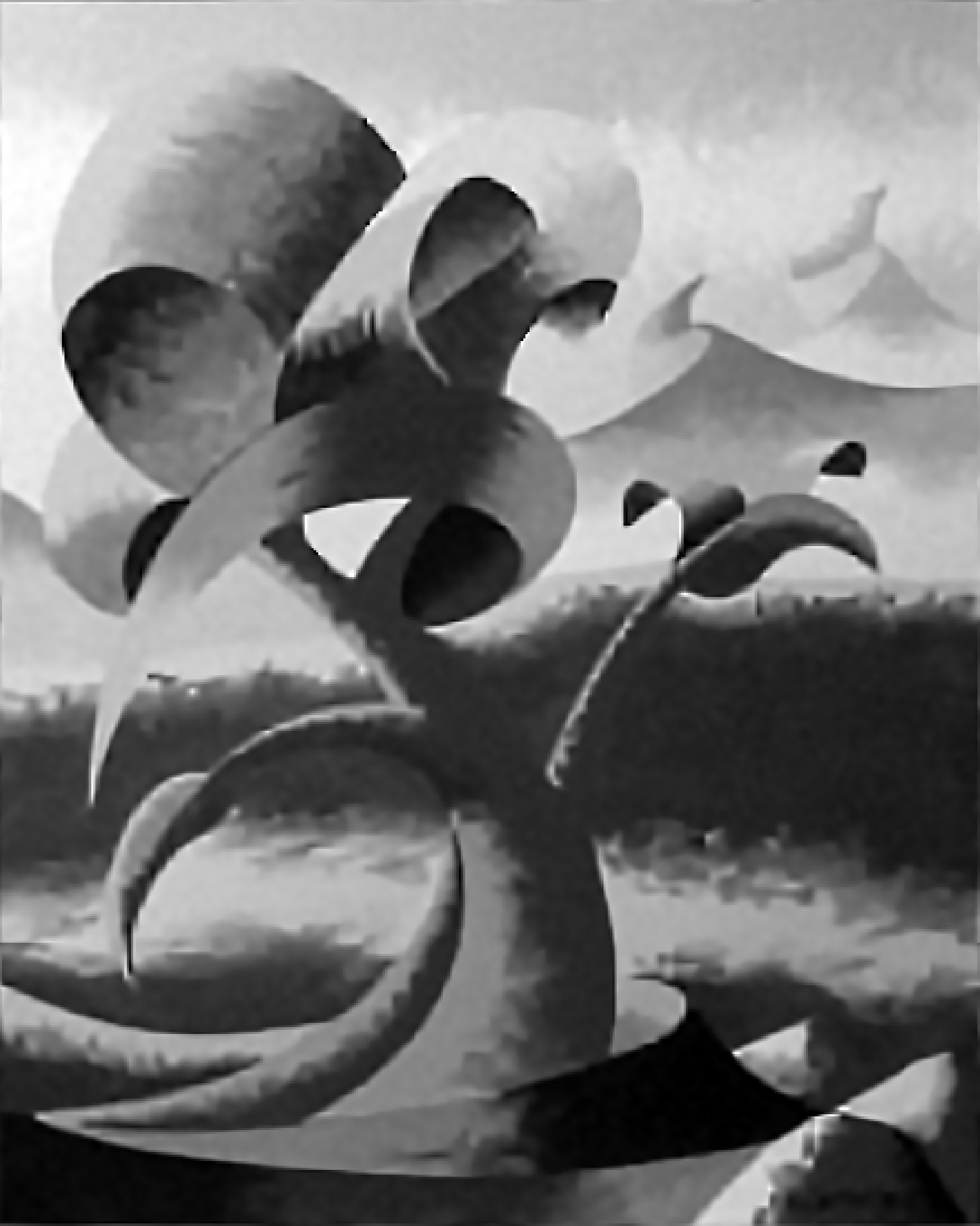}
\end{minipage}
\begin{minipage}{2.25cm}
\includegraphics[height=3.00cm]{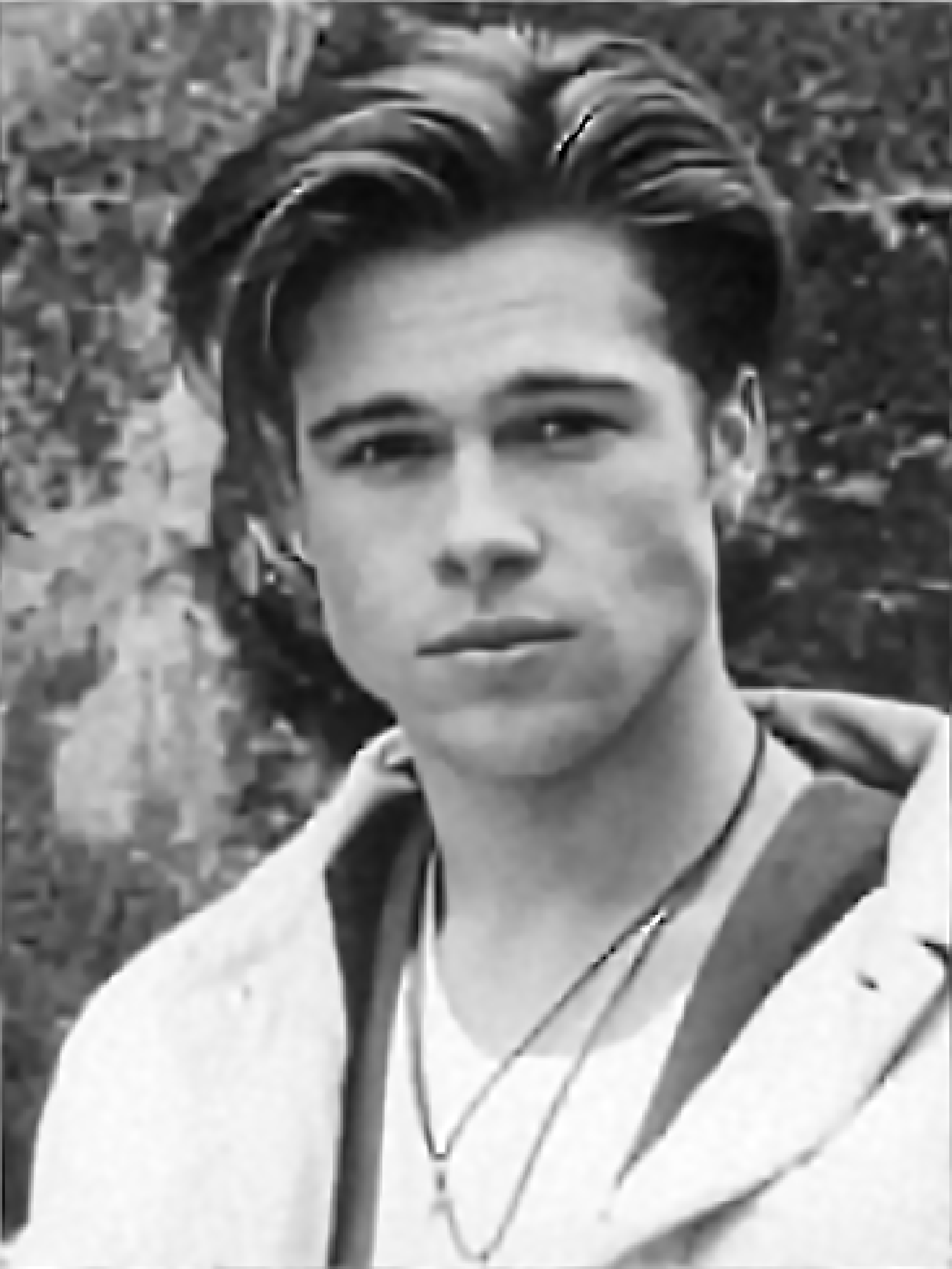}
\end{minipage}\vspace{0.25em}\\
\begin{minipage}{3.00cm}
\includegraphics[width=3.00cm]{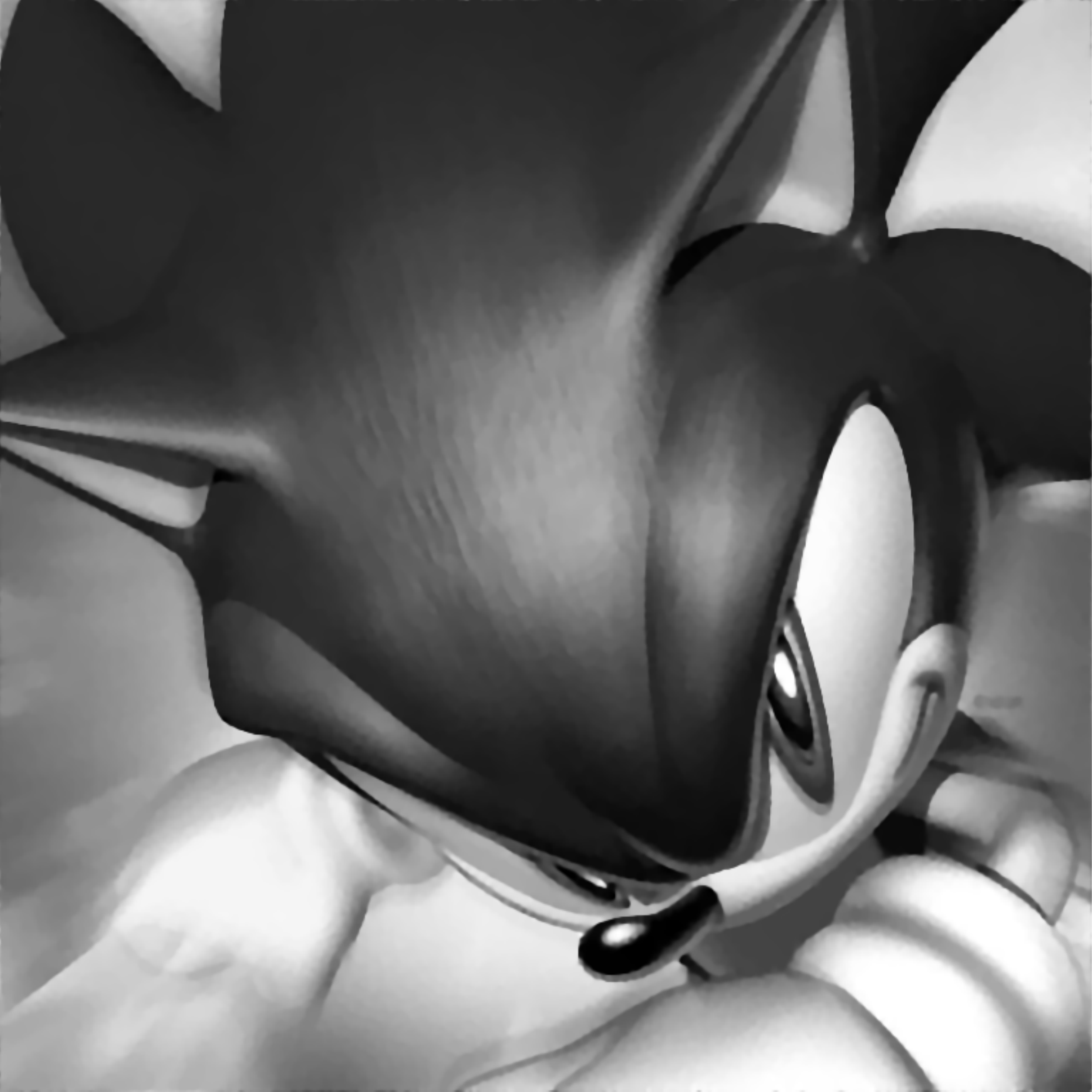}
\end{minipage}
\begin{minipage}{4.5cm}
\includegraphics[height=3.00cm]{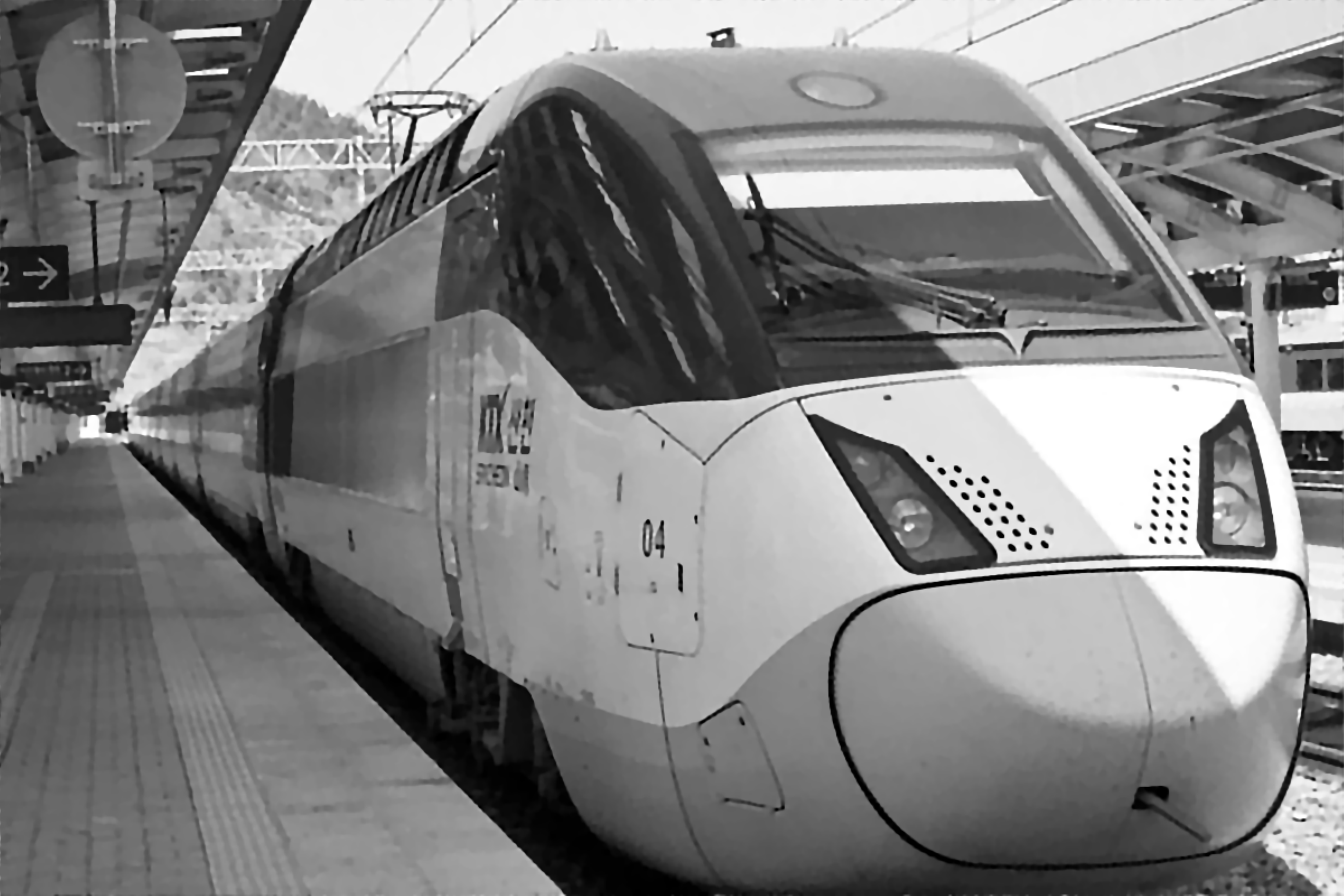}
\end{minipage}
\begin{minipage}{3.75cm}
\includegraphics[height=3.00cm]{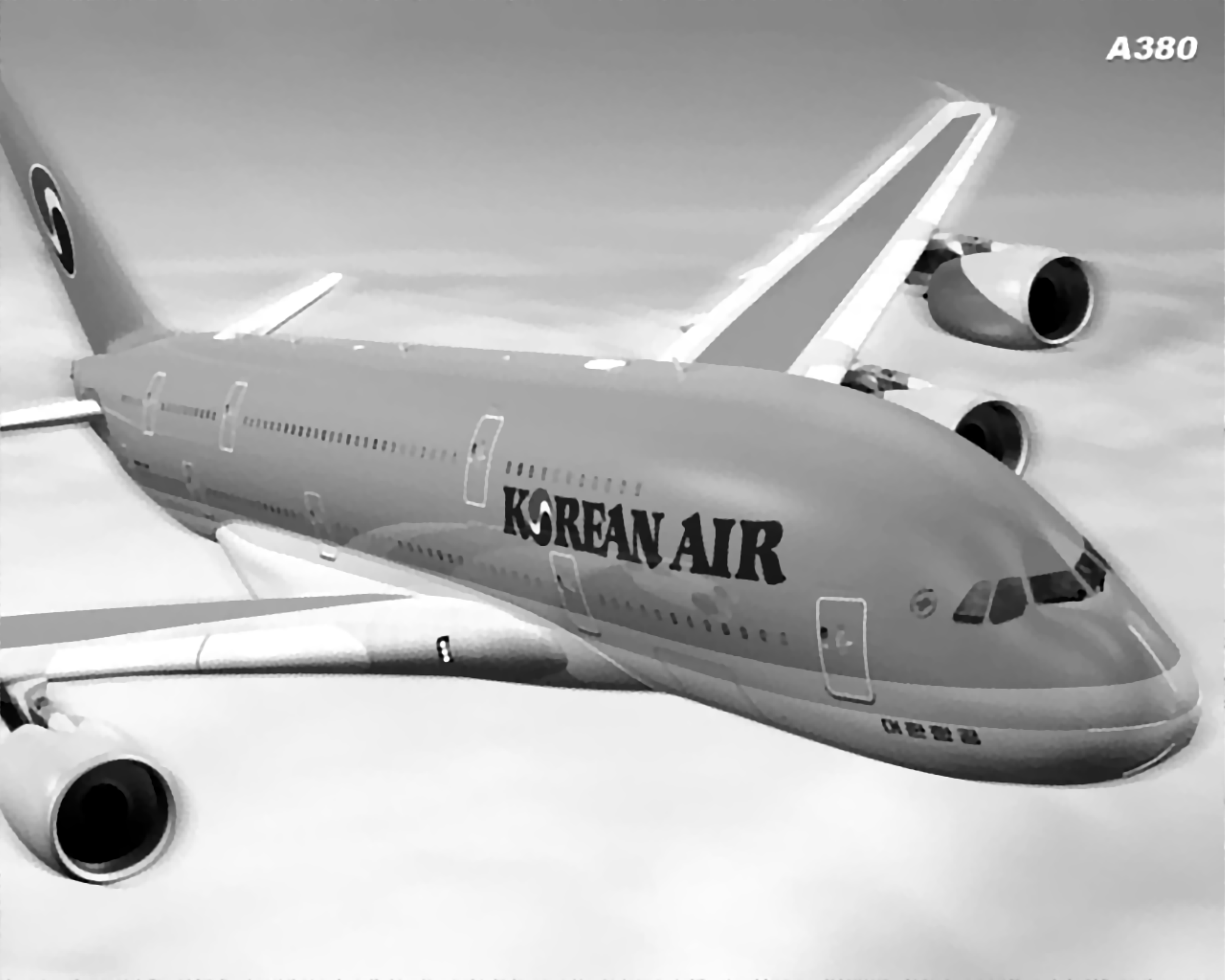}
\end{minipage}
\begin{minipage}{2.4cm}
\includegraphics[height=3.00cm]{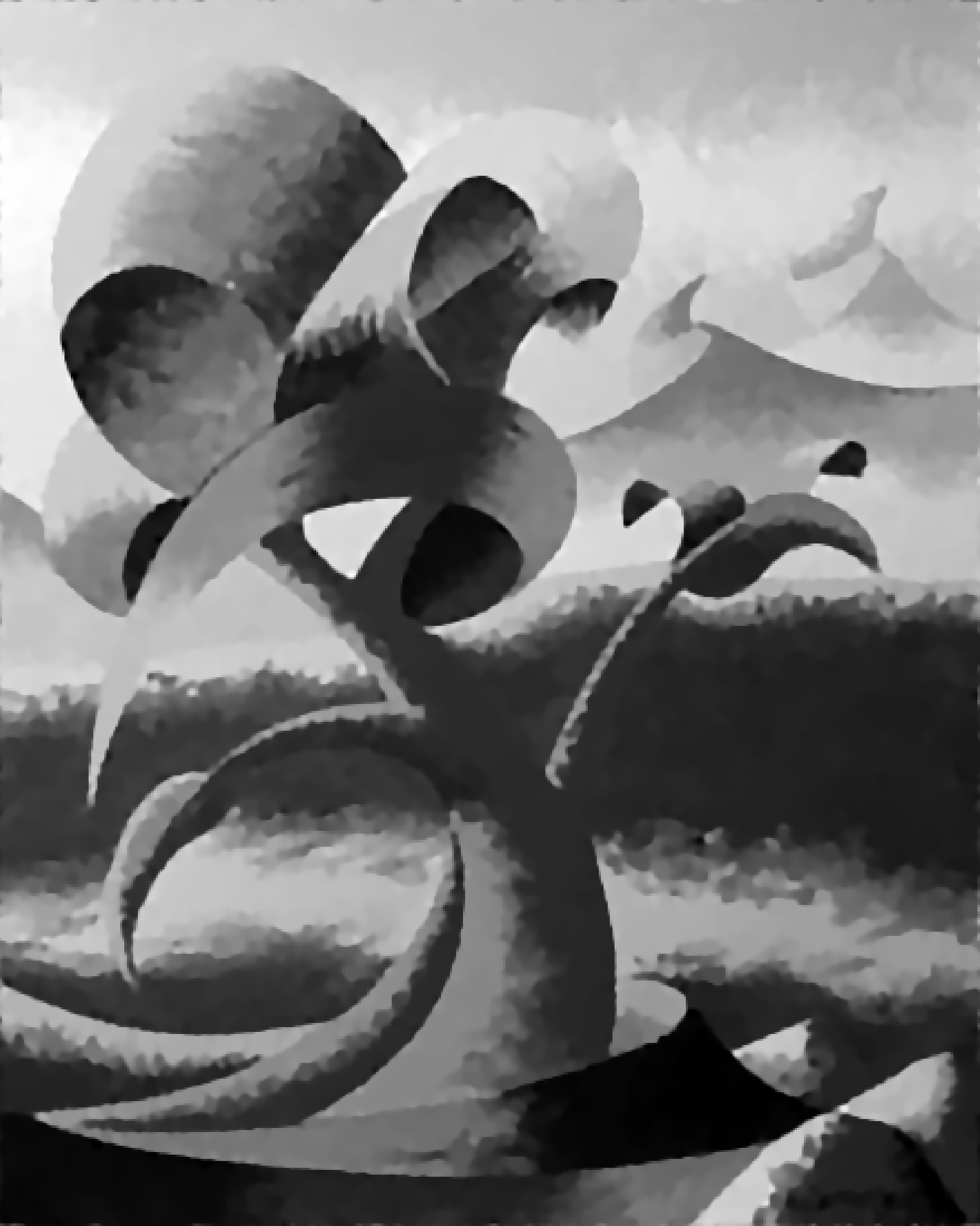}
\end{minipage}
\begin{minipage}{2.25cm}
\includegraphics[height=3.00cm]{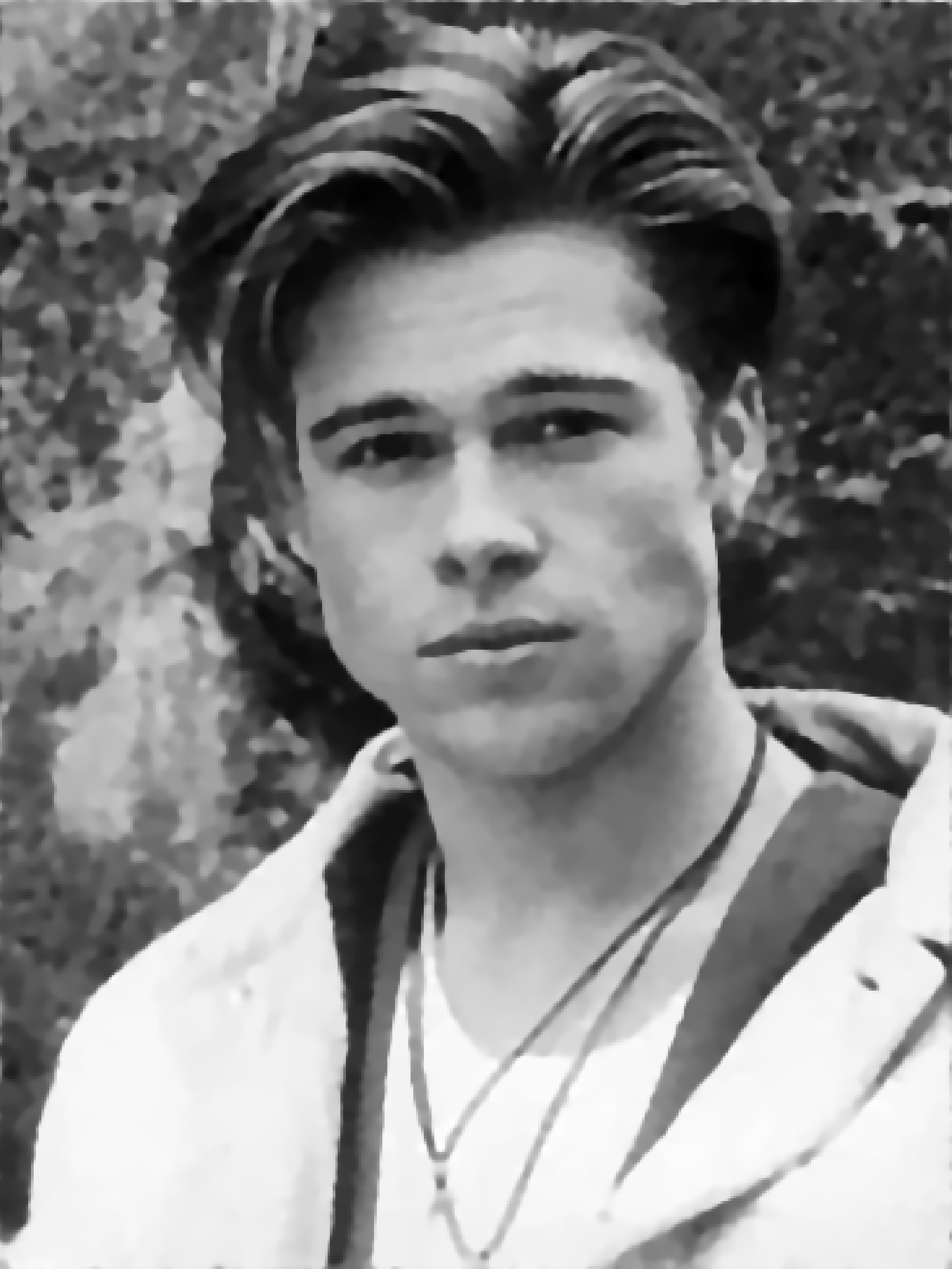}
\end{minipage}\\
\begin{minipage}{3.00cm}\begin{center}{\small{Sonic}}\end{center}\end{minipage}\begin{minipage}{4.5cm}\begin{center}{\small{Train}}\end{center}\end{minipage}\begin{minipage}{3.75cm}\begin{center}{\small{Airplane}}\end{center}\end{minipage}\begin{minipage}{2.4cm}\begin{center}{\small{Oil Painting}}\end{center}\end{minipage}\begin{minipage}{2.25cm}\begin{center}{\small{Pitt}}\end{center}\end{minipage}
\caption{Visualization of restoration results deblurred by four models. The first row describes the observed images, followed by the results of TGV model, PS model, GS model, and our model, respectively.}\label{fig:DeblurResults}
\end{center}
\end{figure}

\begin{figure}[htp!]
\begin{center}
\begin{minipage}{3.150cm}
\includegraphics[width=3.150cm]{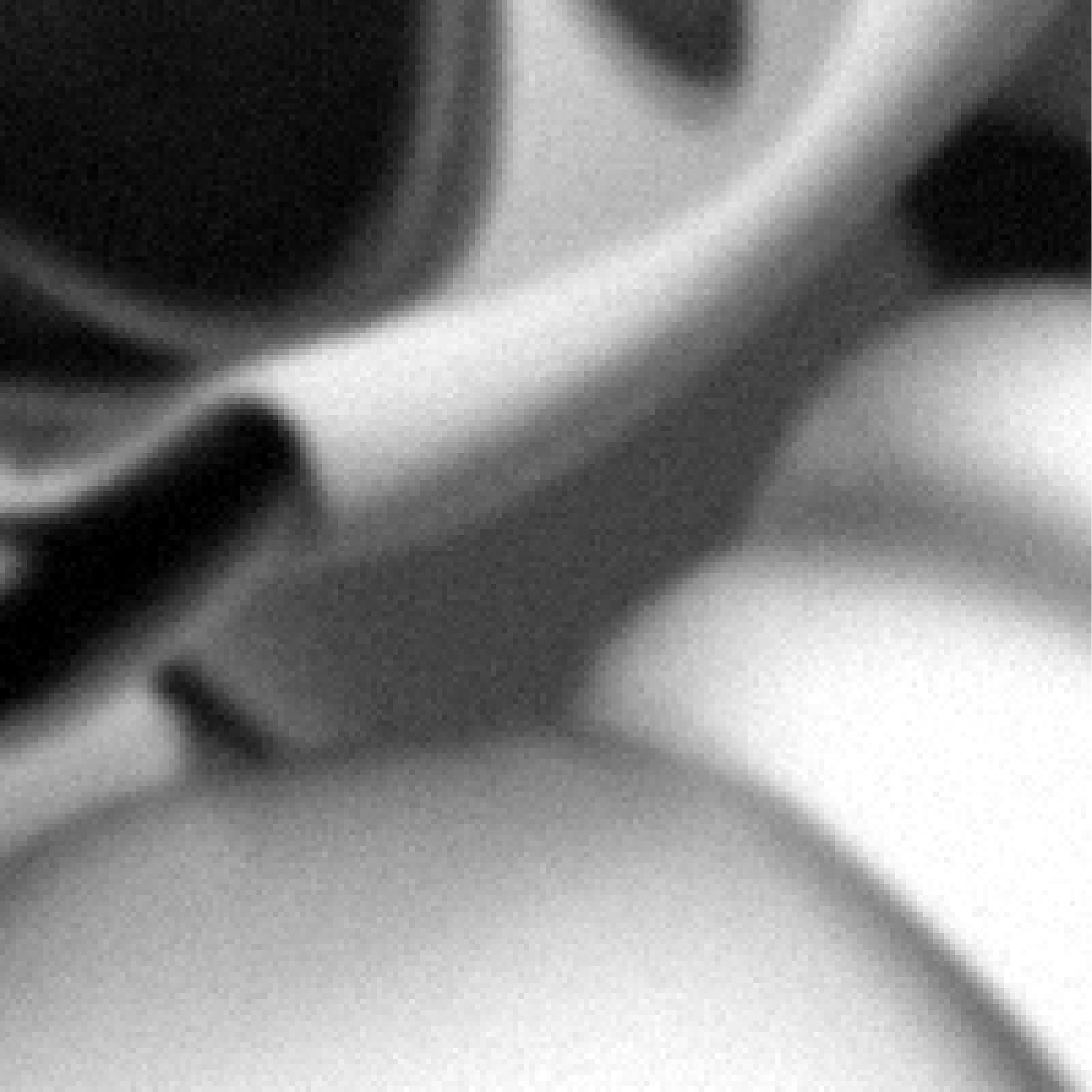}
\end{minipage}
\begin{minipage}{3.150cm}
\includegraphics[width=3.150cm]{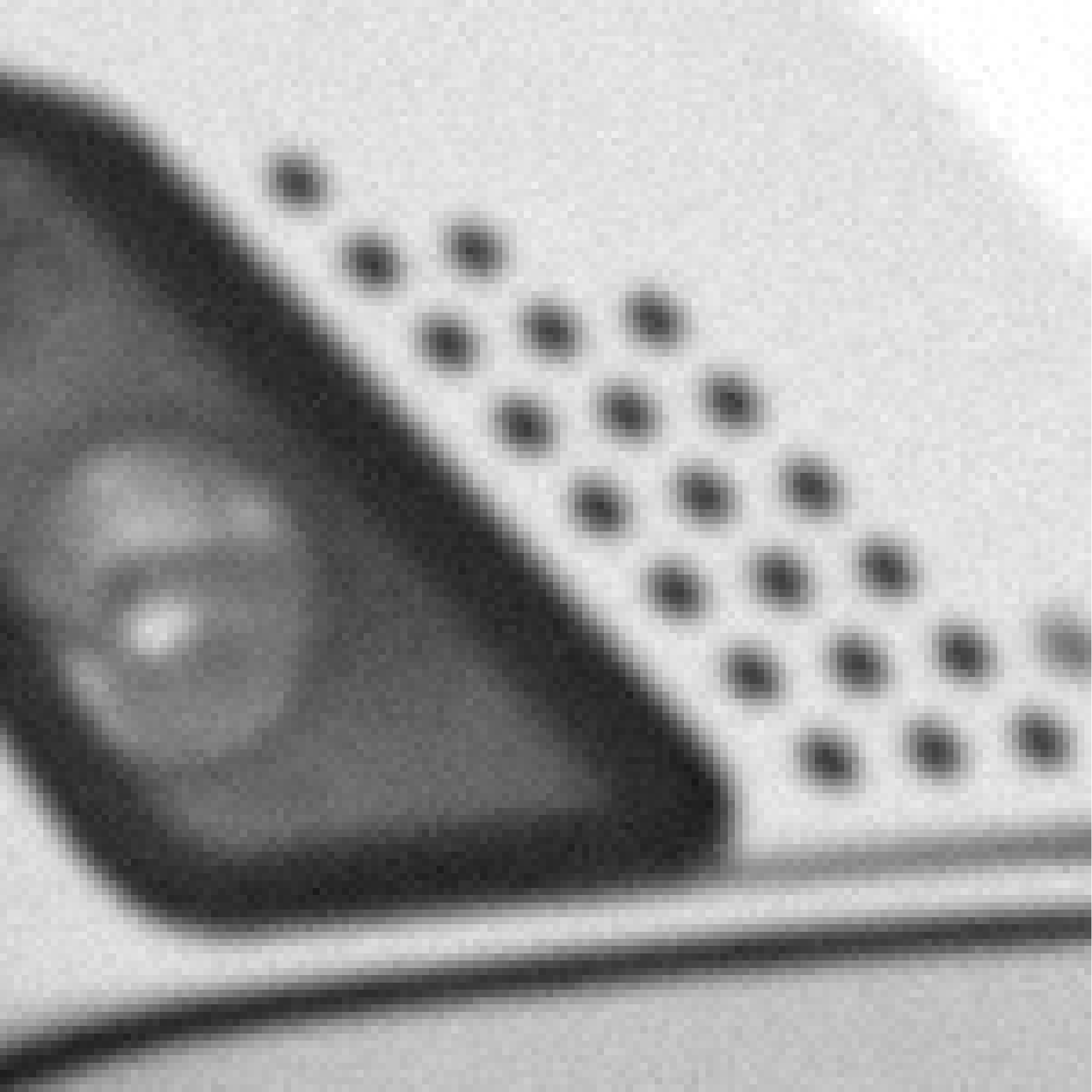}
\end{minipage}
\begin{minipage}{3.150cm}
\includegraphics[width=3.150cm]{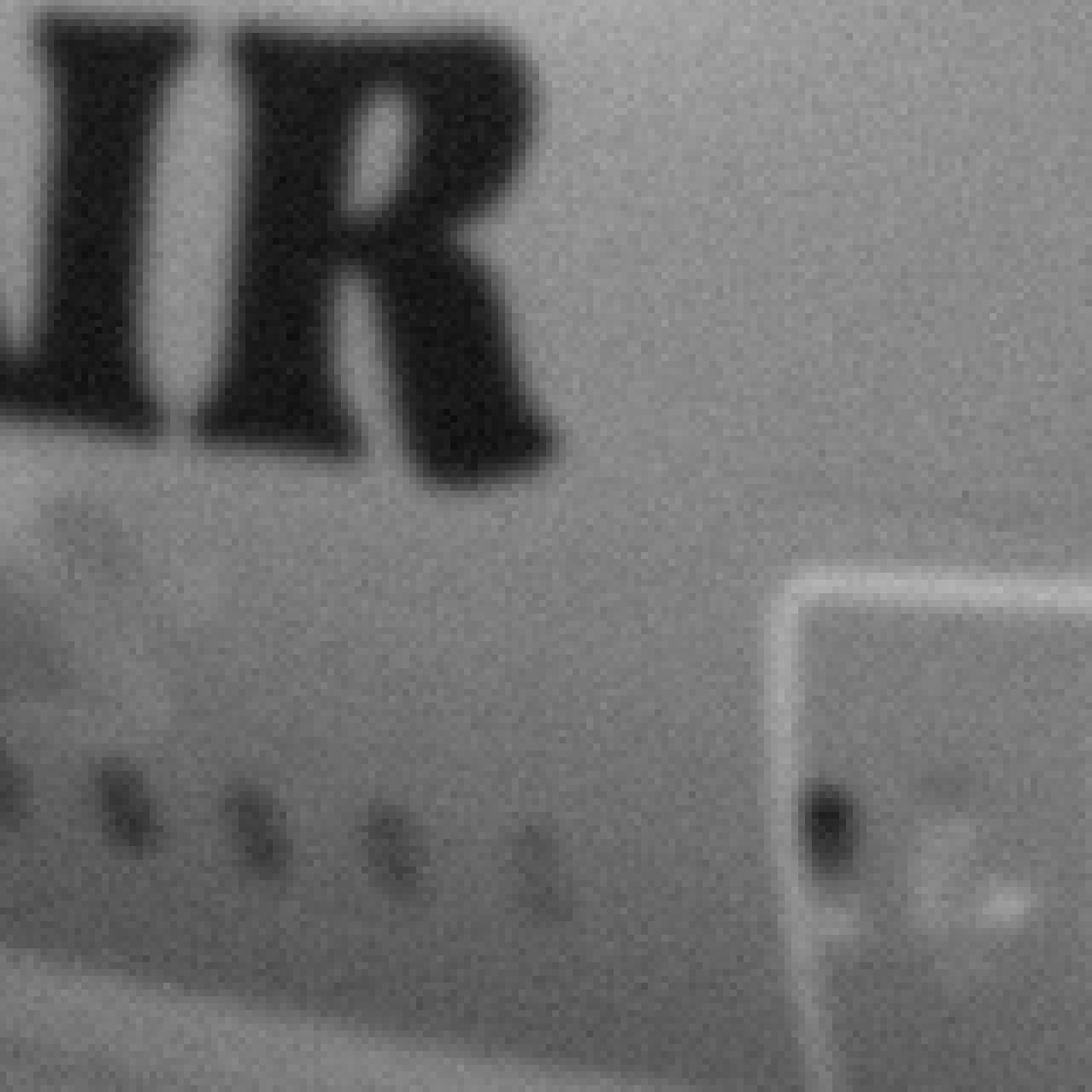}
\end{minipage}
\begin{minipage}{3.150cm}
\includegraphics[width=3.150cm]{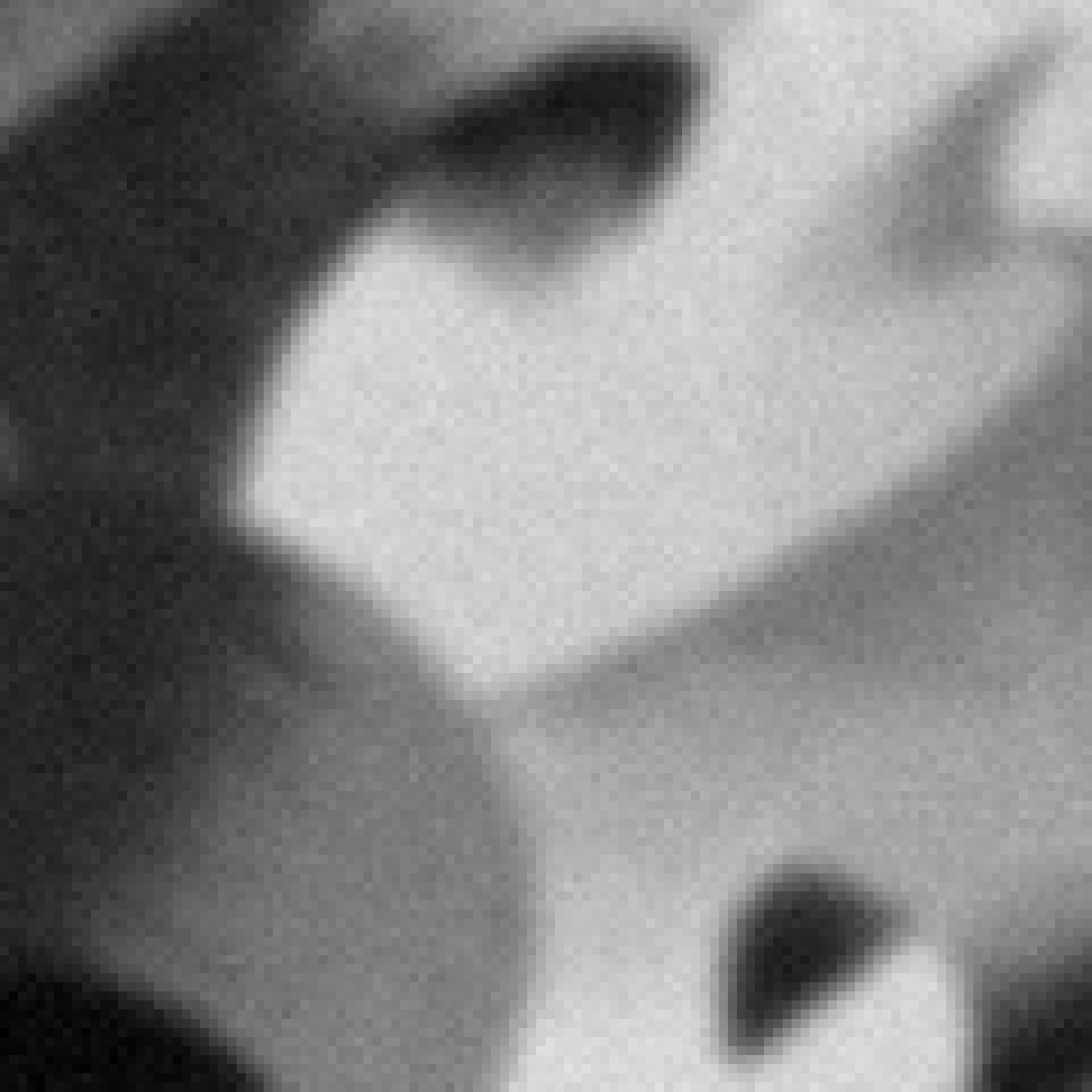}
\end{minipage}
\begin{minipage}{3.150cm}
\includegraphics[width=3.150cm]{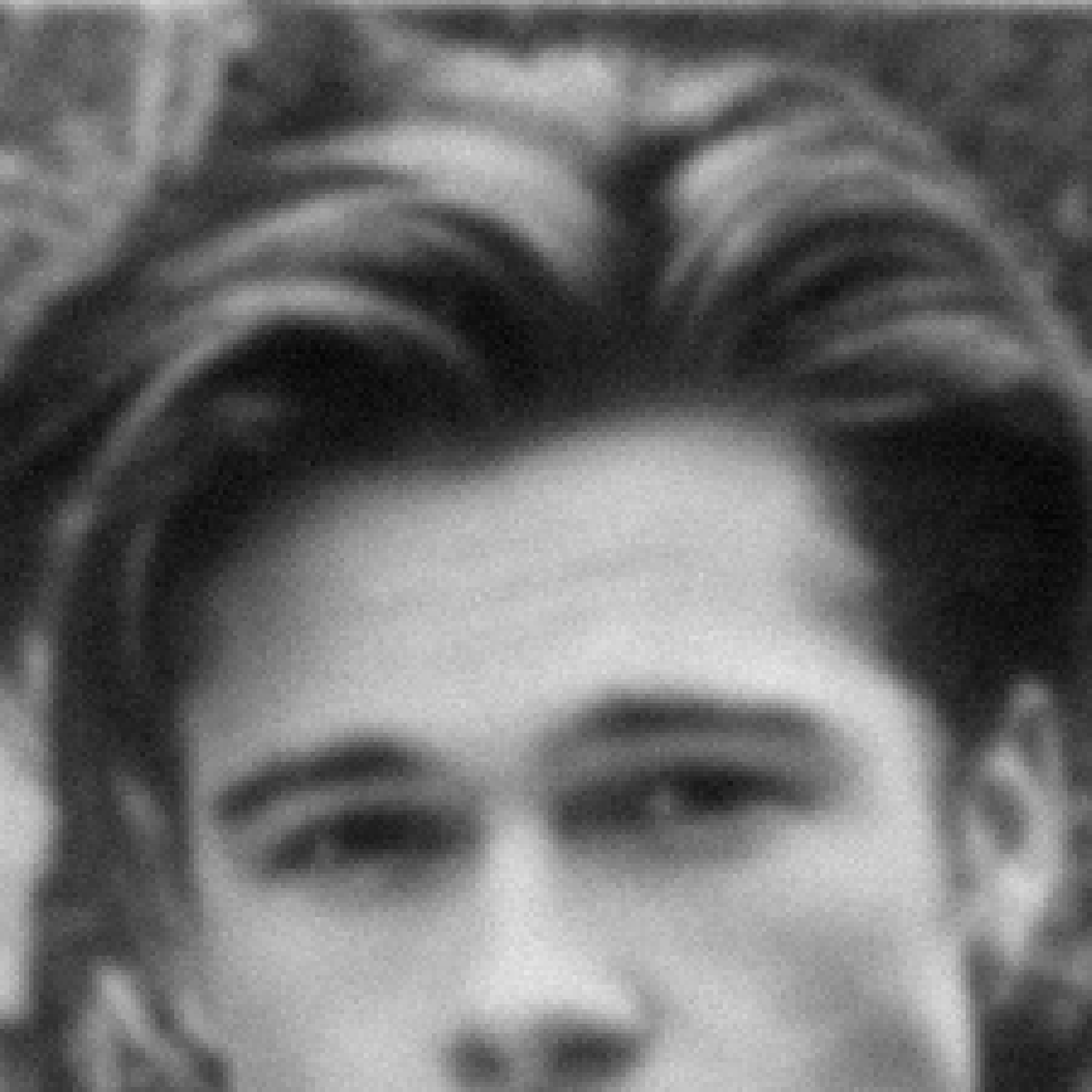}
\end{minipage}\vspace{0.25em}\\
\begin{minipage}{3.150cm}
\includegraphics[width=3.150cm]{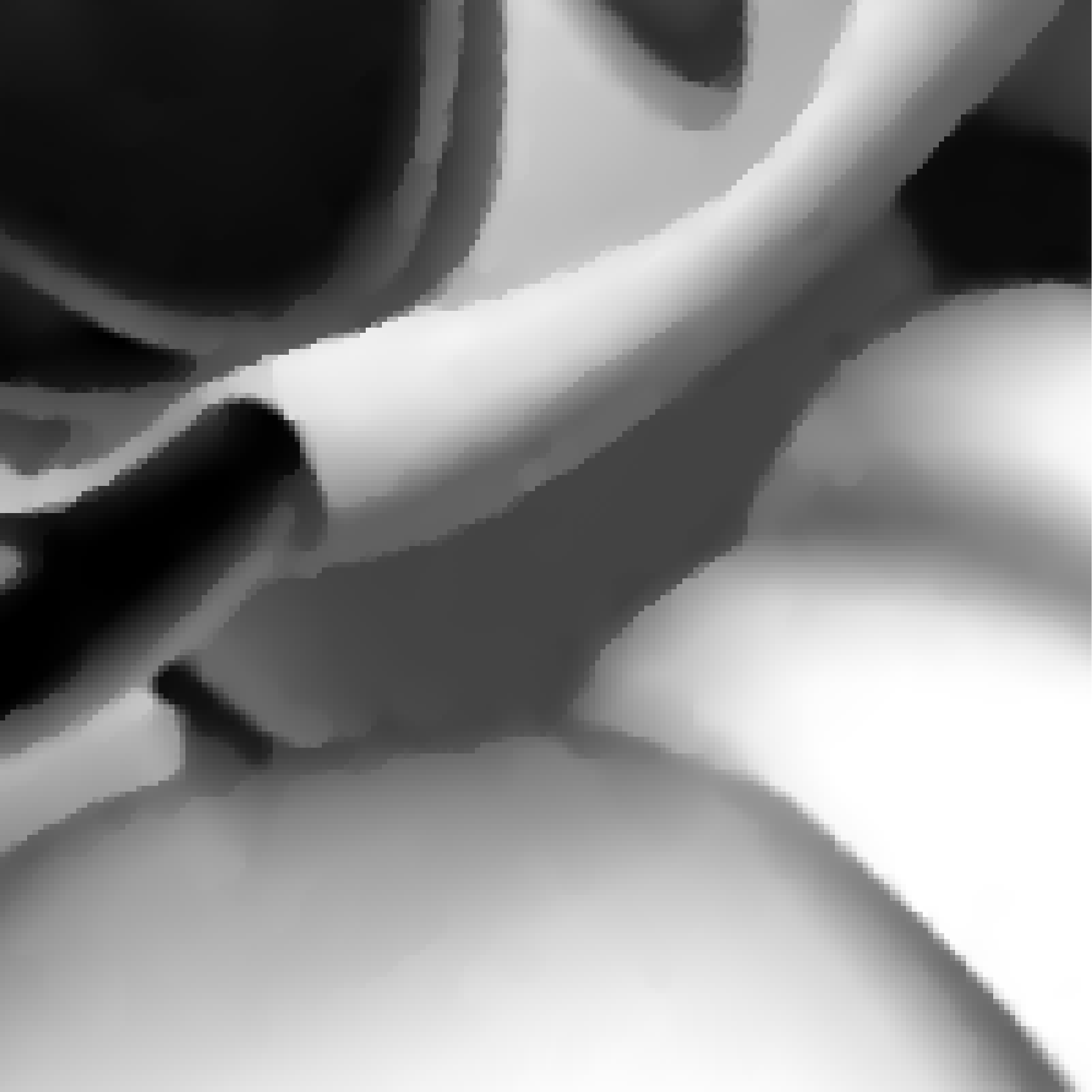}
\end{minipage}
\begin{minipage}{3.150cm}
\includegraphics[width=3.150cm]{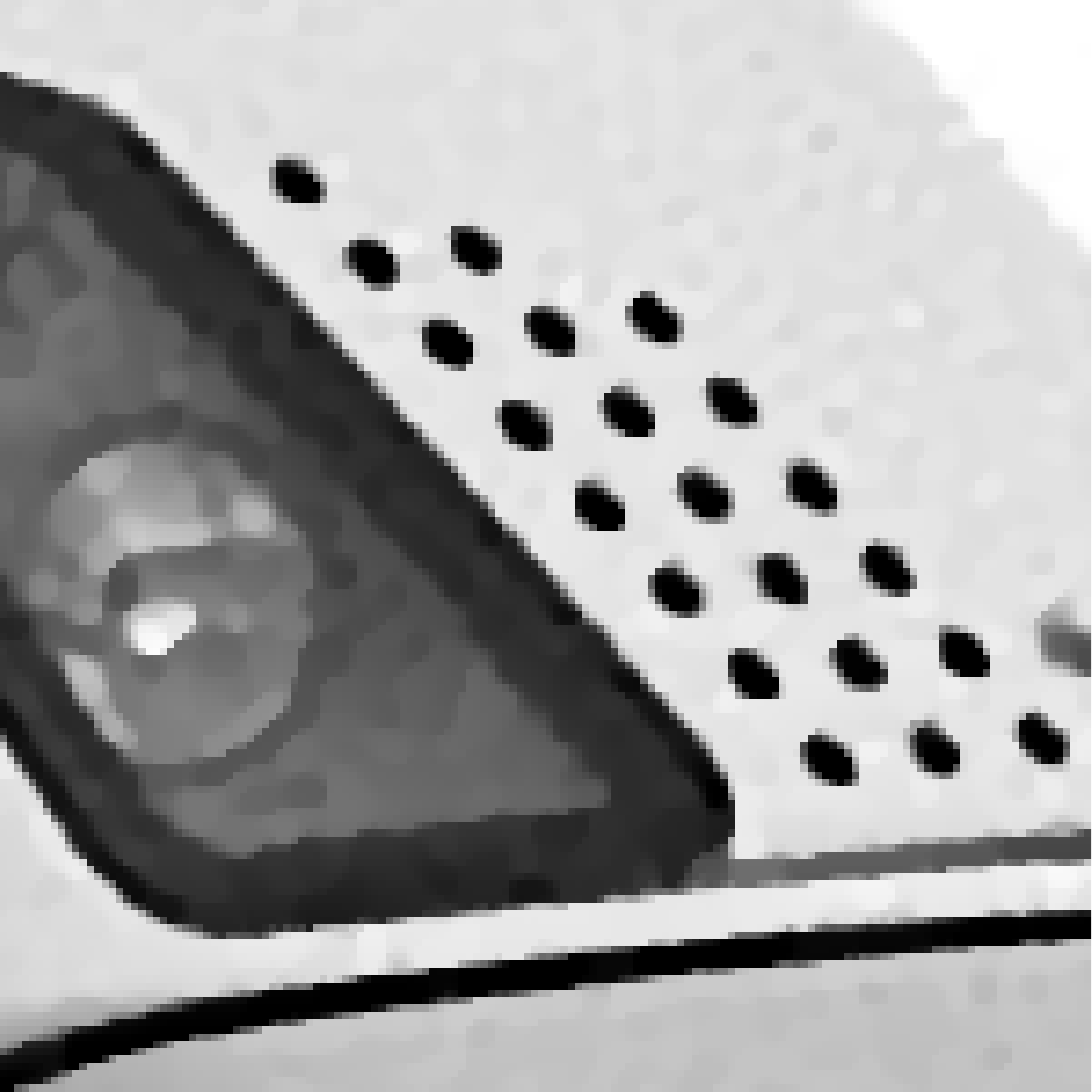}
\end{minipage}
\begin{minipage}{3.150cm}
\includegraphics[width=3.150cm]{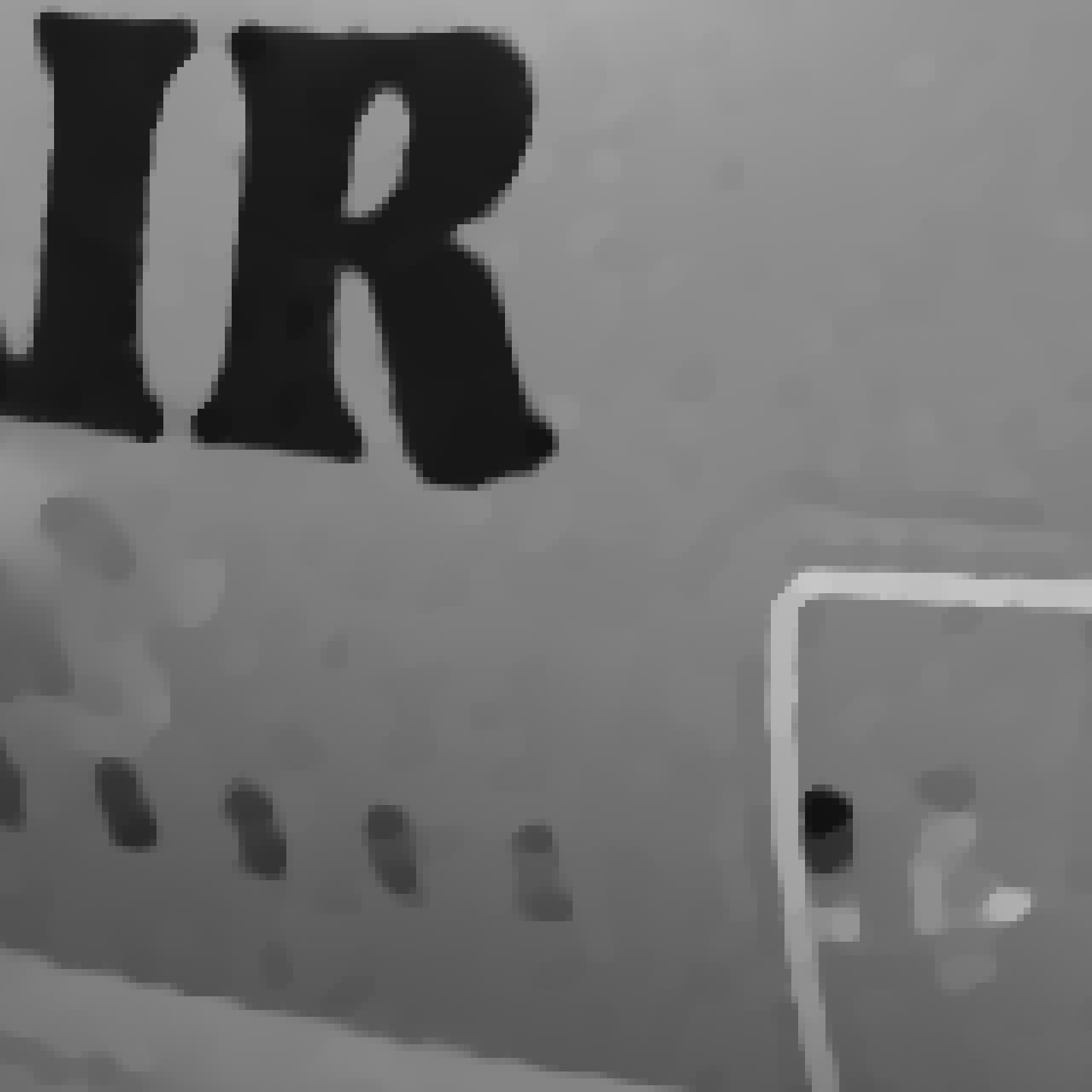}
\end{minipage}
\begin{minipage}{3.150cm}
\includegraphics[width=3.150cm]{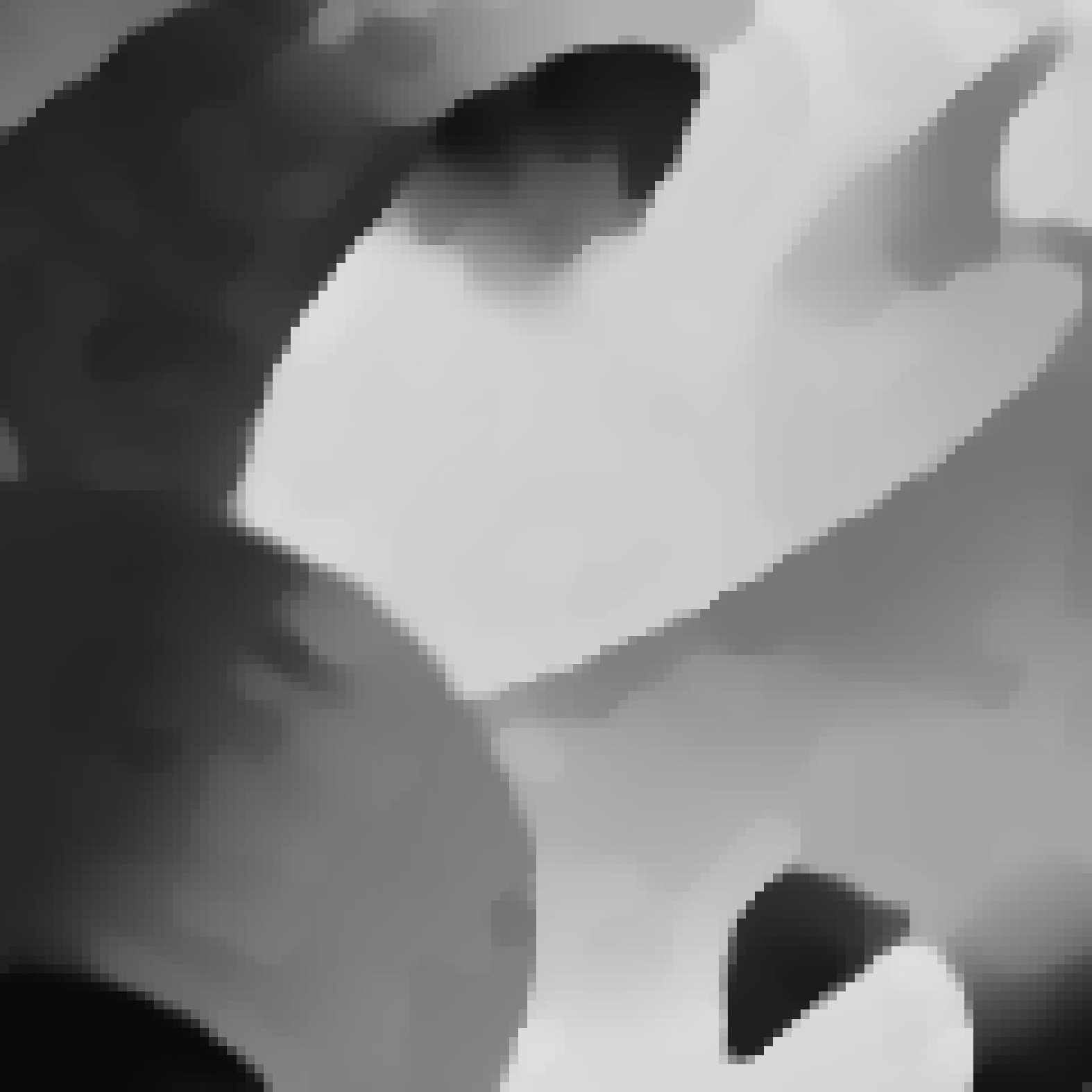}
\end{minipage}
\begin{minipage}{3.150cm}
\includegraphics[width=3.150cm]{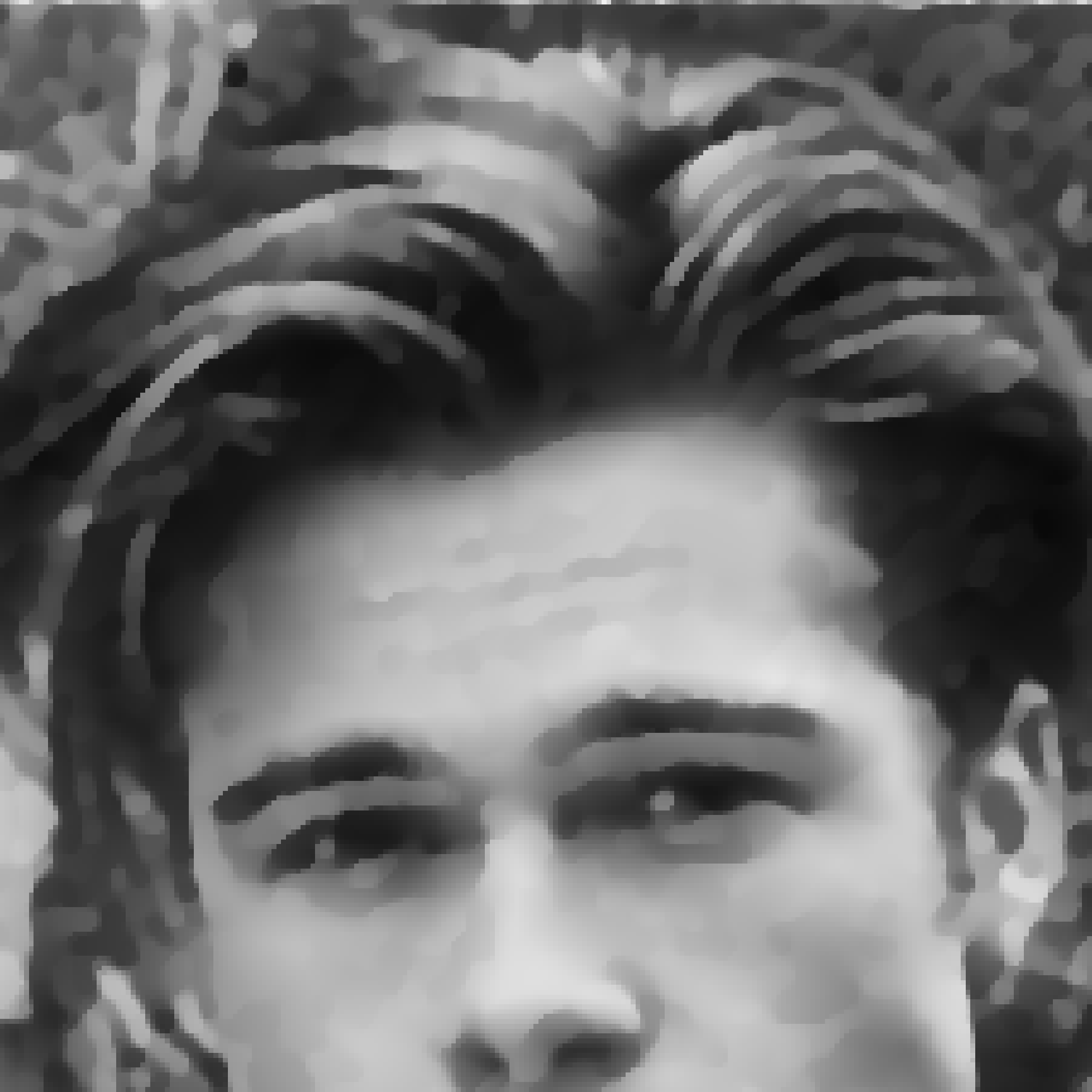}
\end{minipage}\vspace{0.25em}\\
\begin{minipage}{3.150cm}
\includegraphics[width=3.150cm]{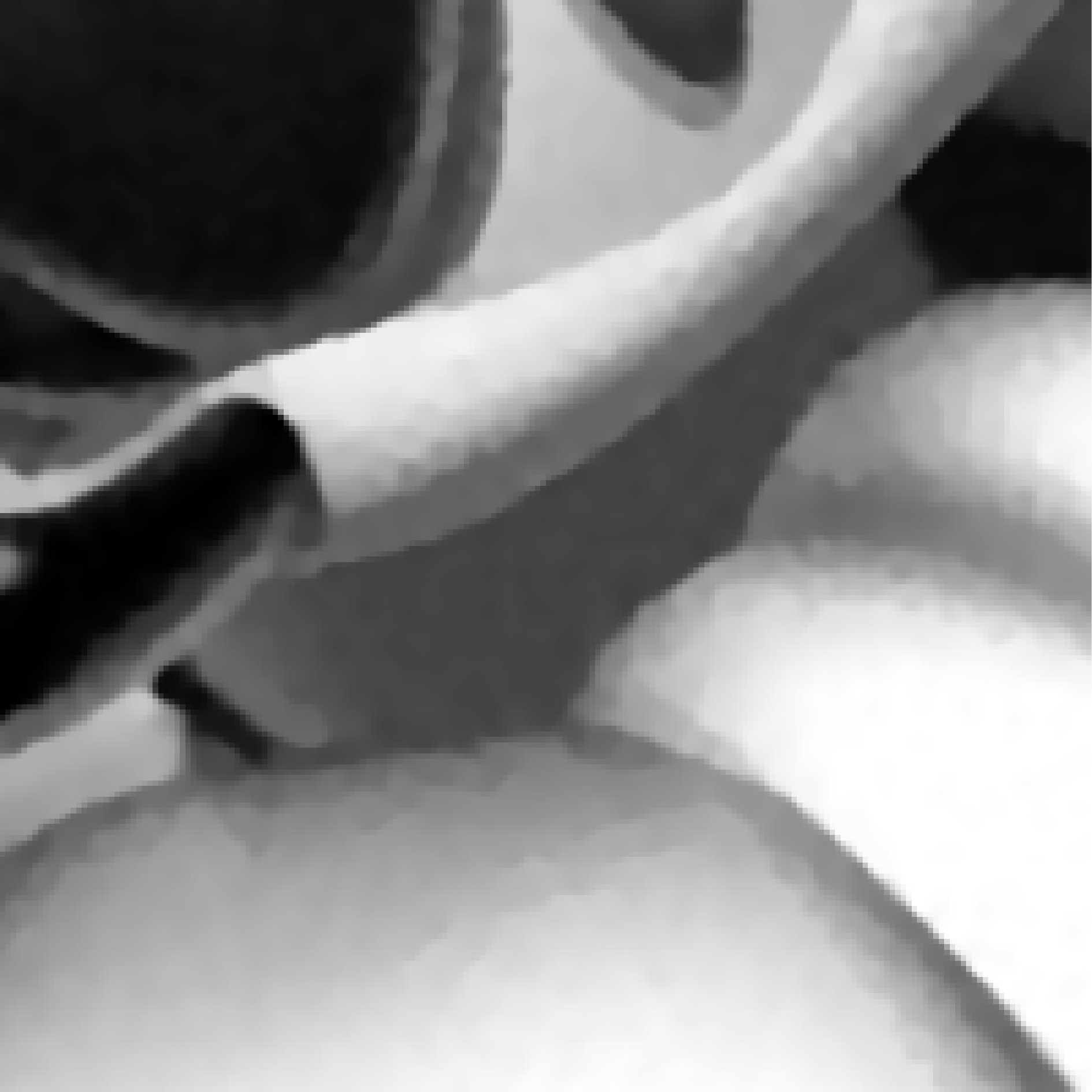}
\end{minipage}
\begin{minipage}{3.150cm}
\includegraphics[width=3.150cm]{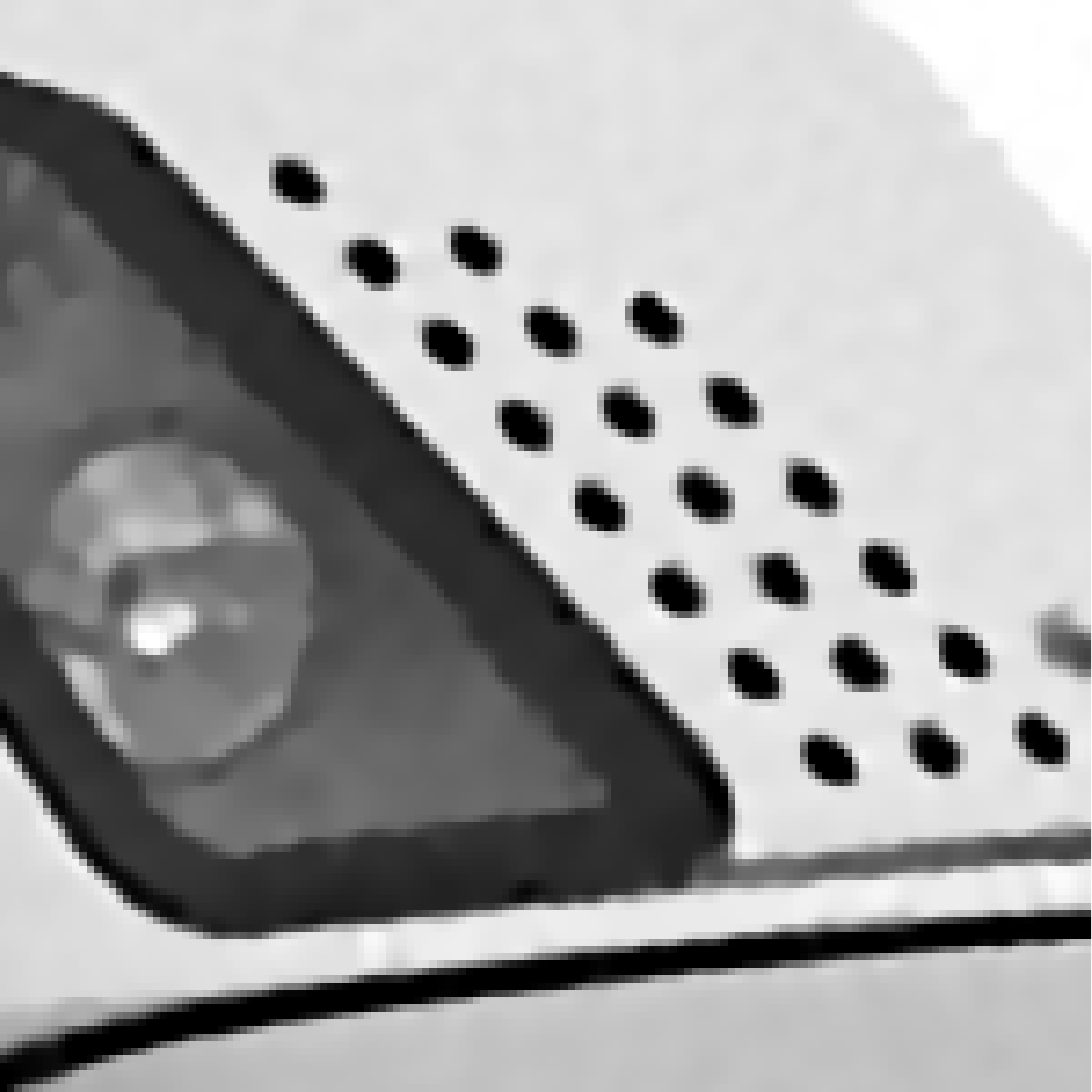}
\end{minipage}
\begin{minipage}{3.150cm}
\includegraphics[width=3.150cm]{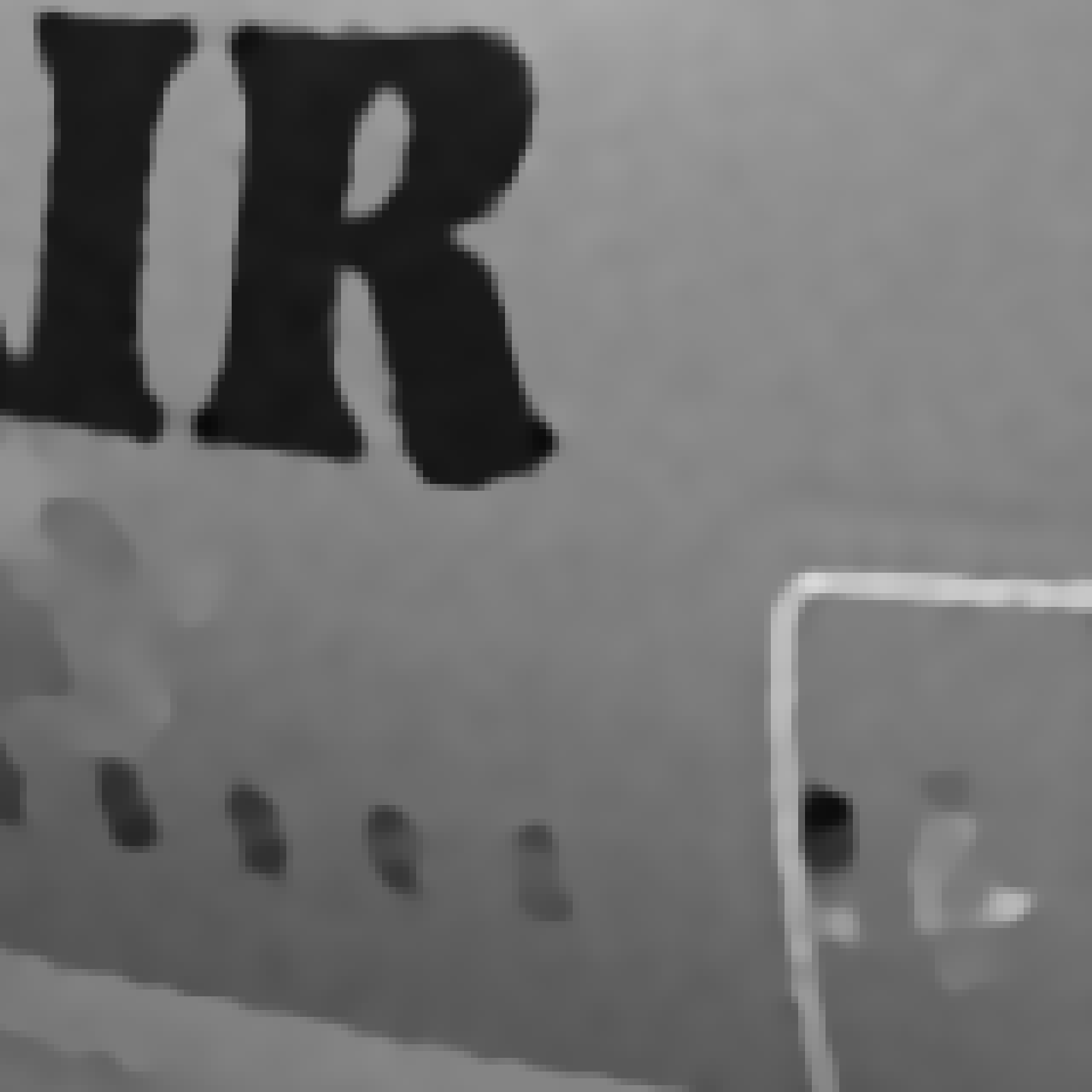}
\end{minipage}
\begin{minipage}{3.150cm}
\includegraphics[width=3.150cm]{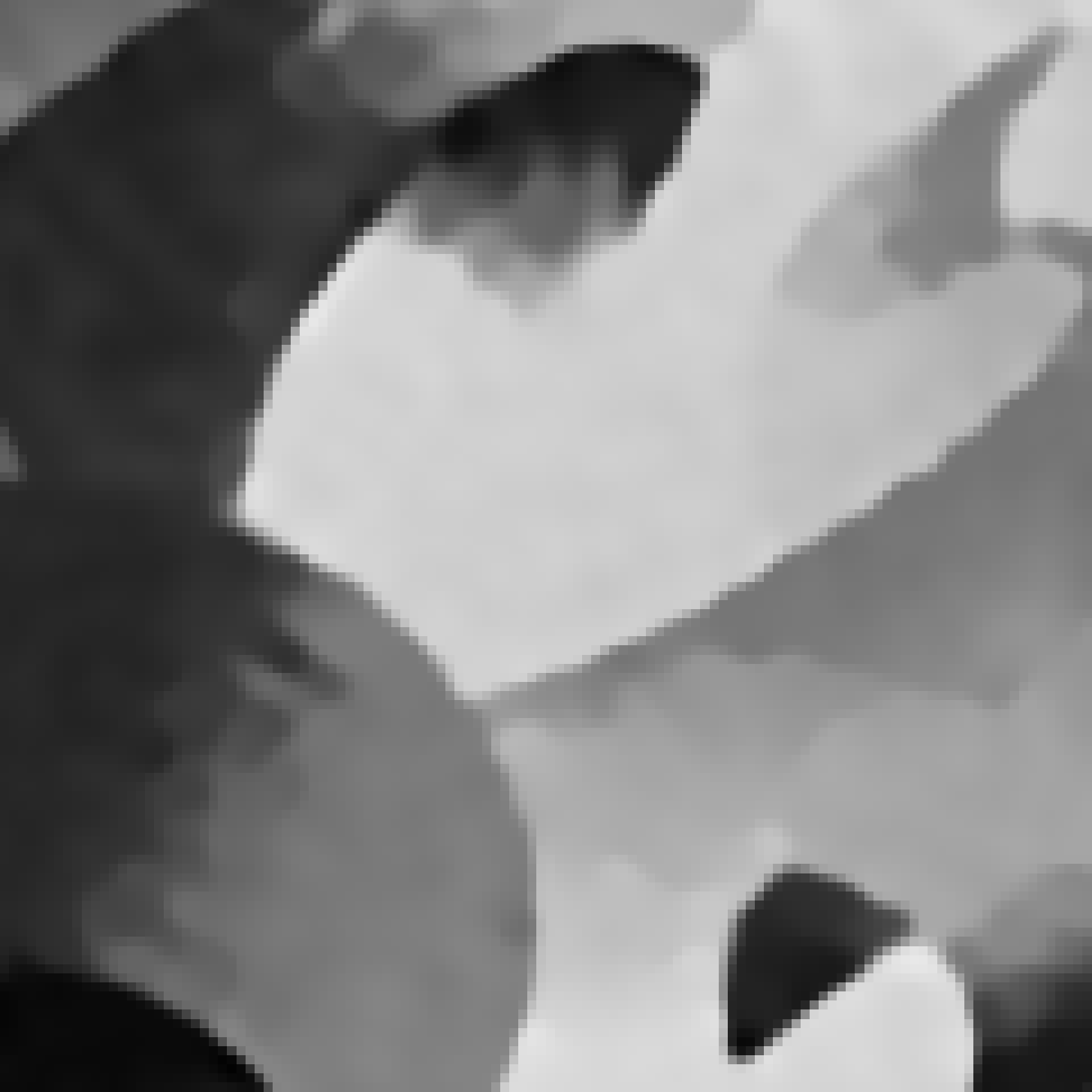}
\end{minipage}
\begin{minipage}{3.150cm}
\includegraphics[width=3.150cm]{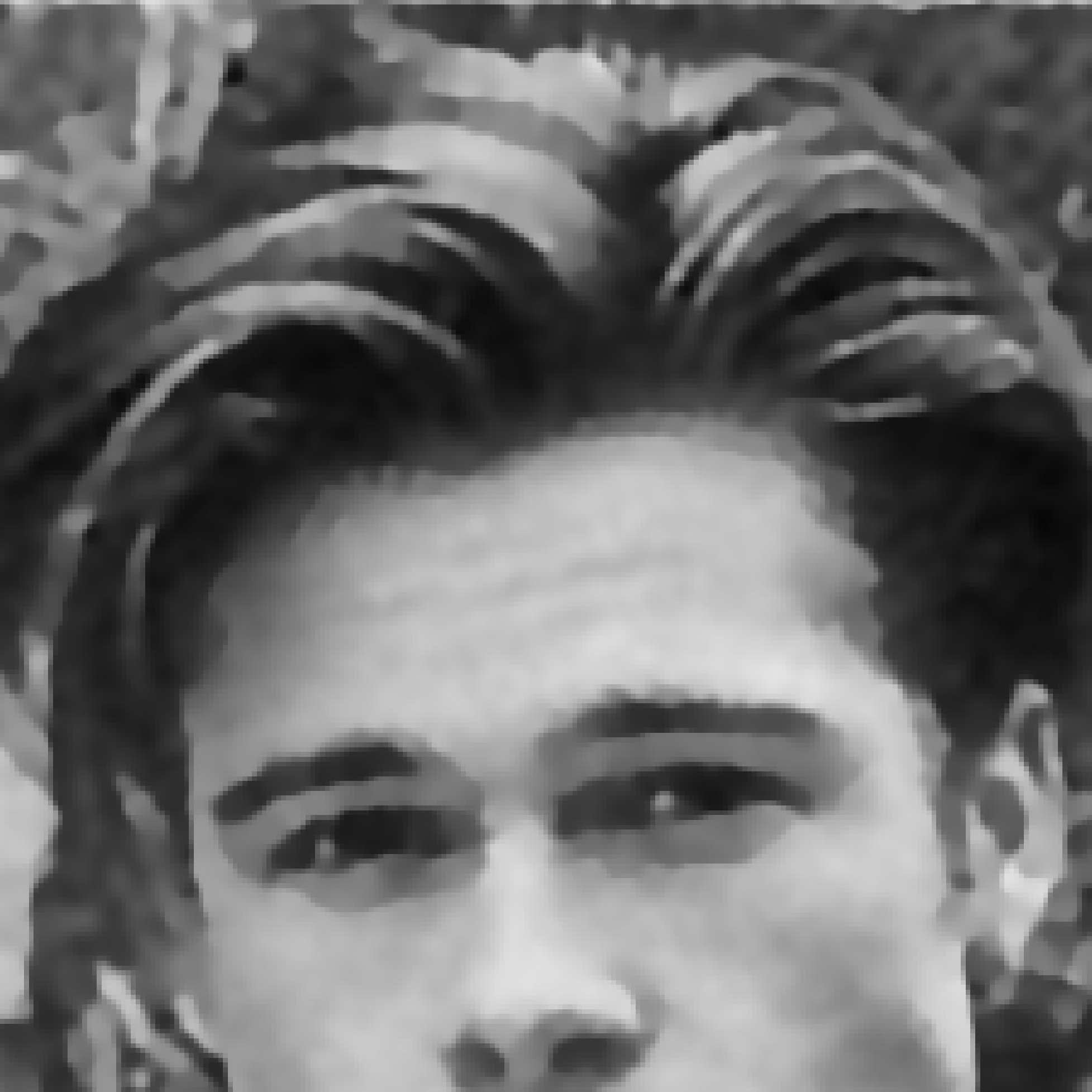}
\end{minipage}\vspace{0.25em}\\
\begin{minipage}{3.150cm}
\includegraphics[width=3.150cm]{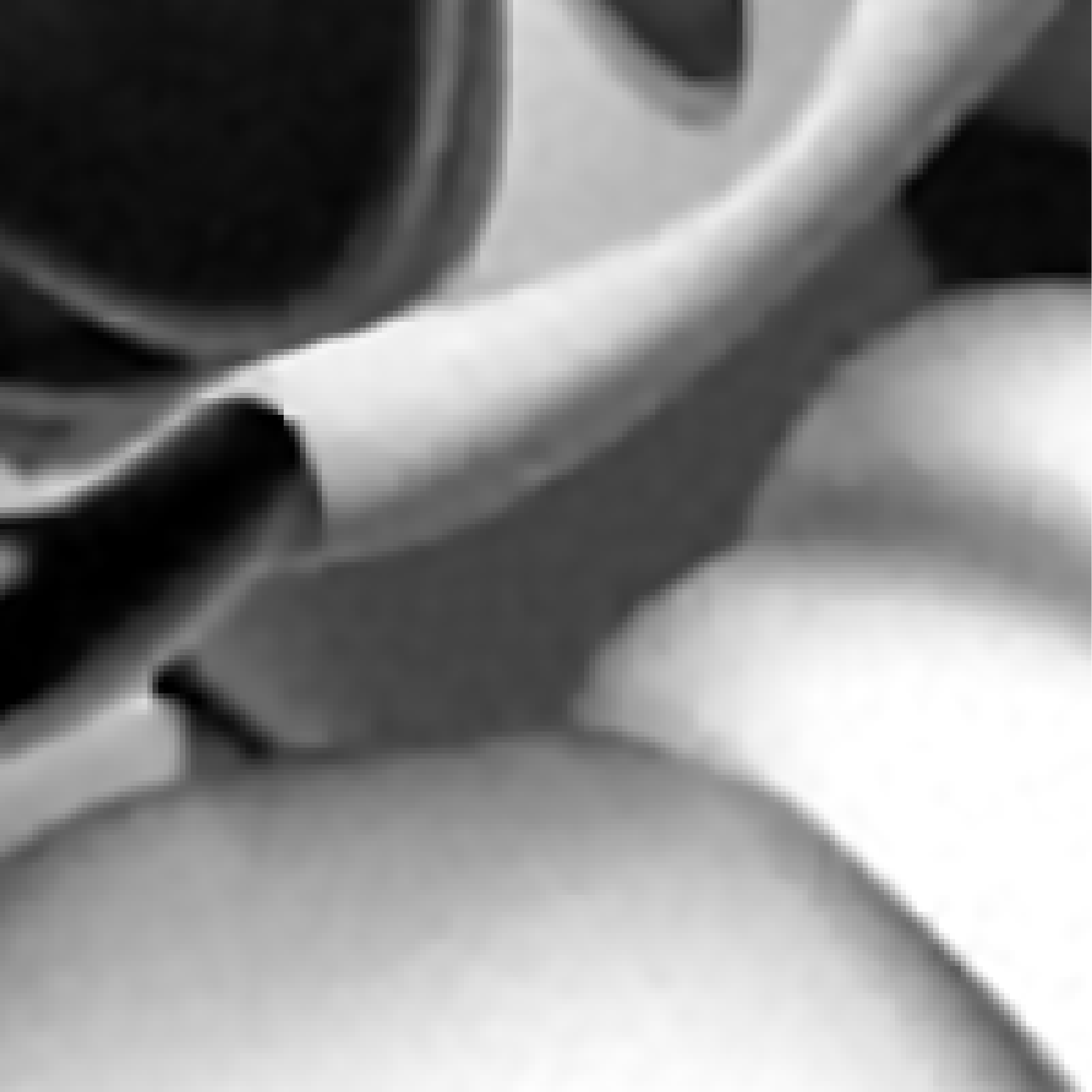}
\end{minipage}
\begin{minipage}{3.150cm}
\includegraphics[width=3.150cm]{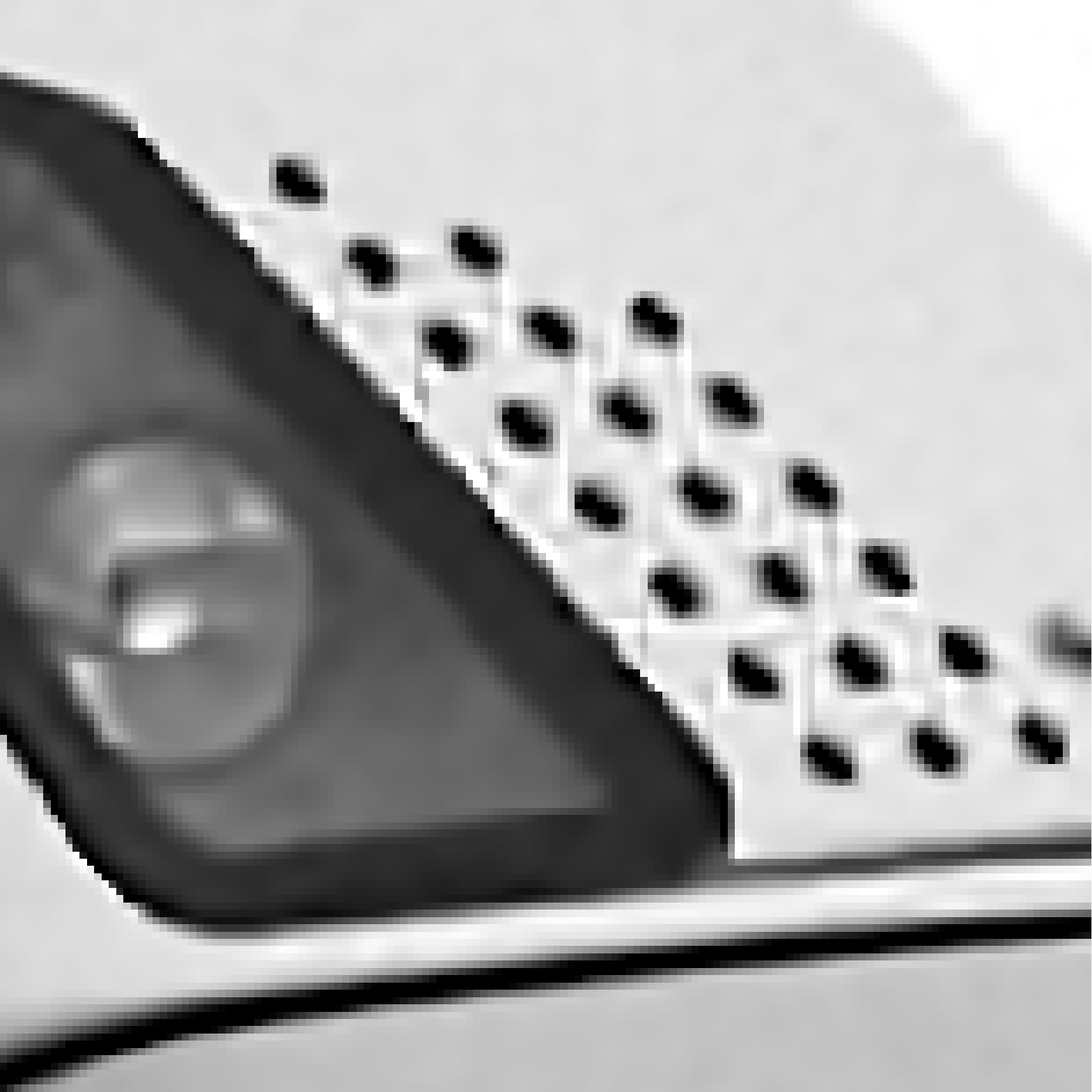}
\end{minipage}
\begin{minipage}{3.150cm}
\includegraphics[width=3.150cm]{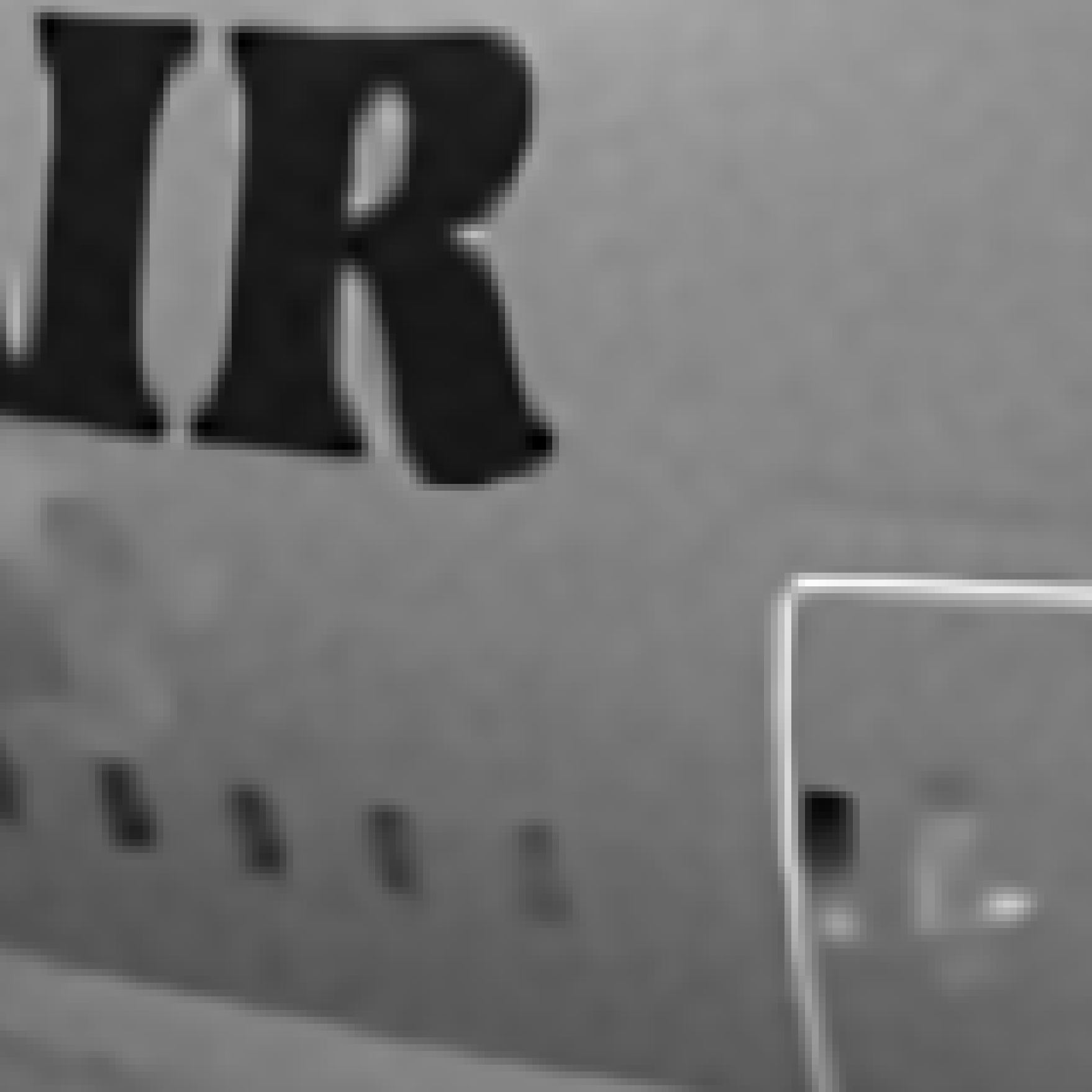}
\end{minipage}
\begin{minipage}{3.150cm}
\includegraphics[width=3.150cm]{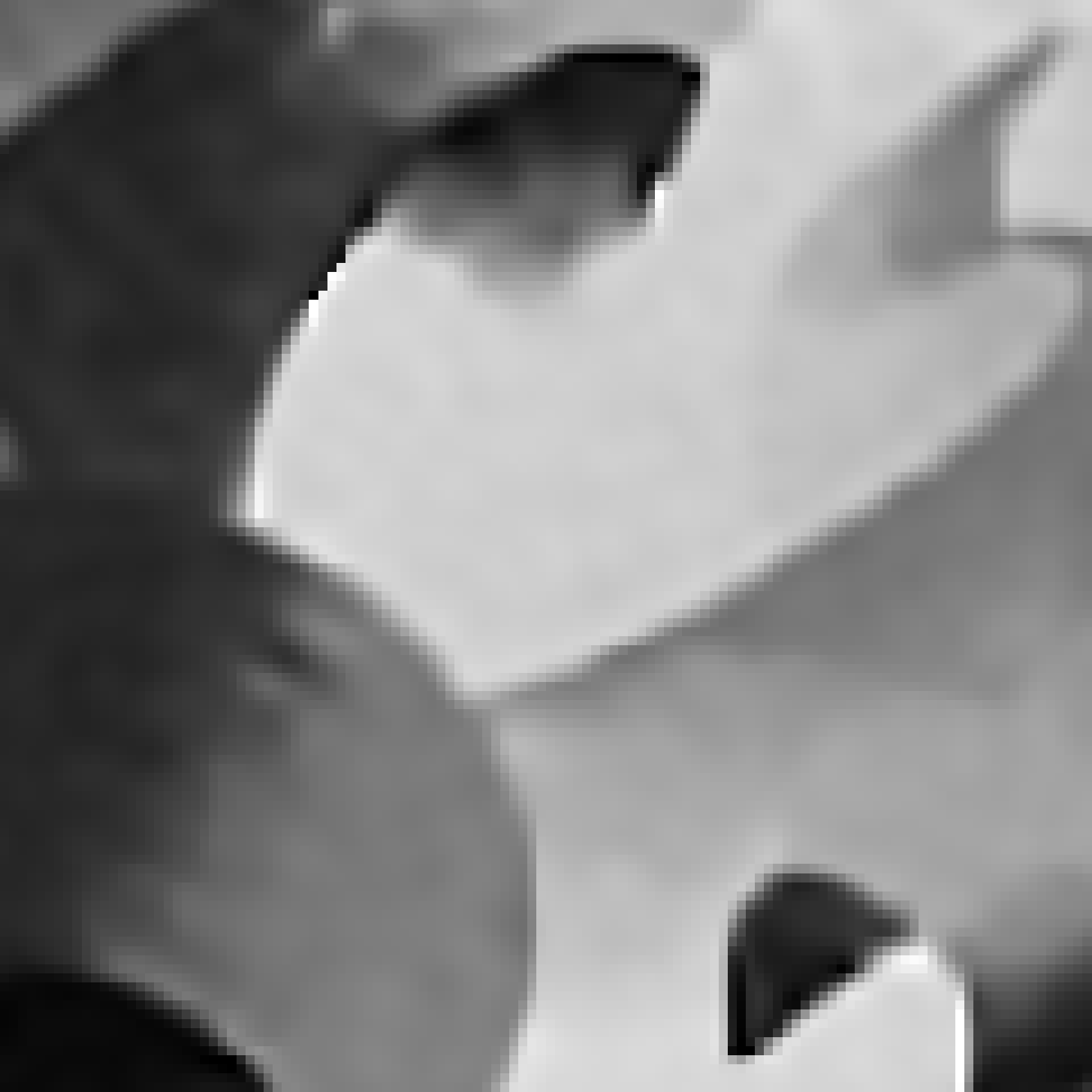}
\end{minipage}
\begin{minipage}{3.150cm}
\includegraphics[width=3.150cm]{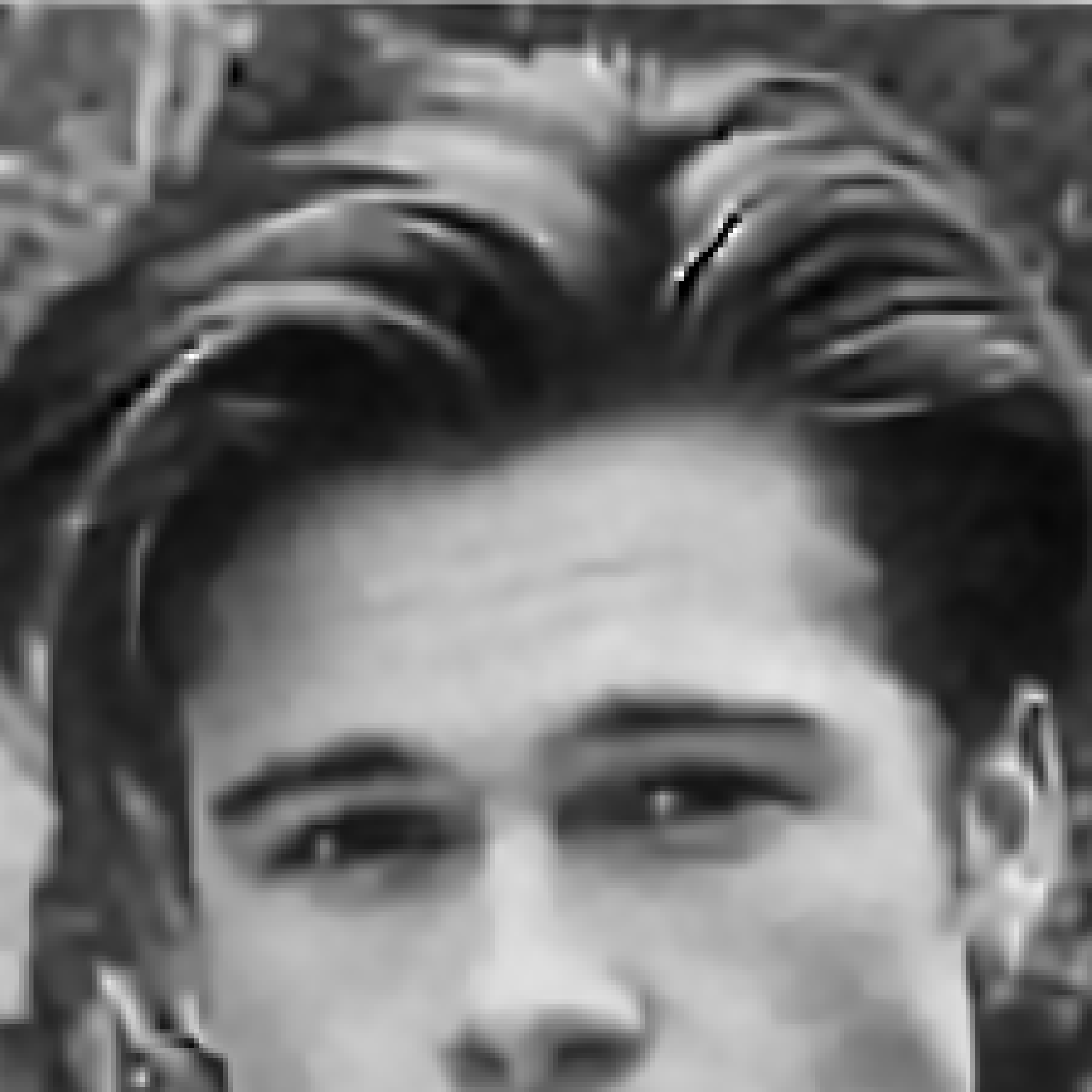}
\end{minipage}\vspace{0.25em}\\
\begin{minipage}{3.150cm}
\includegraphics[width=3.150cm]{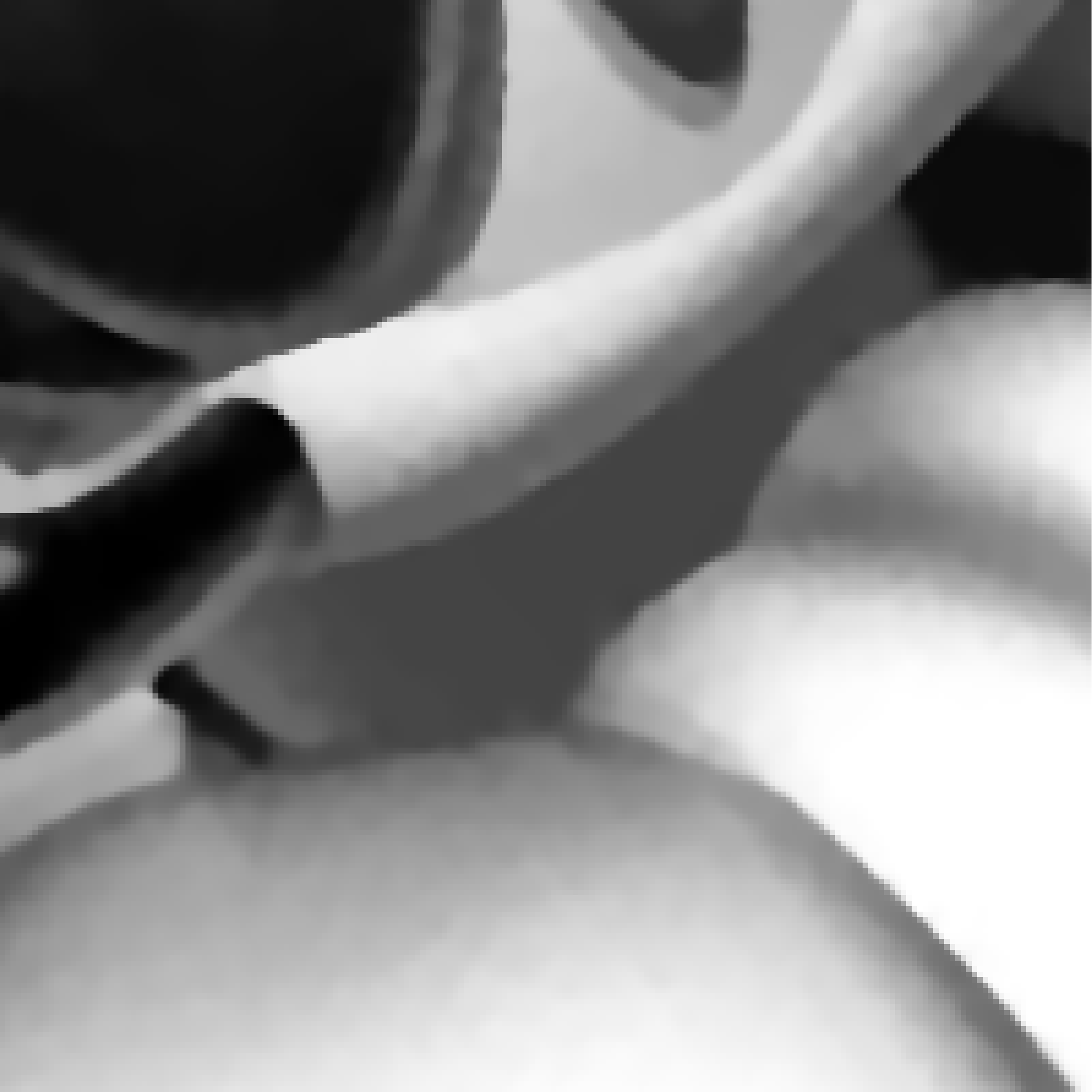}
\end{minipage}
\begin{minipage}{3.150cm}
\includegraphics[width=3.150cm]{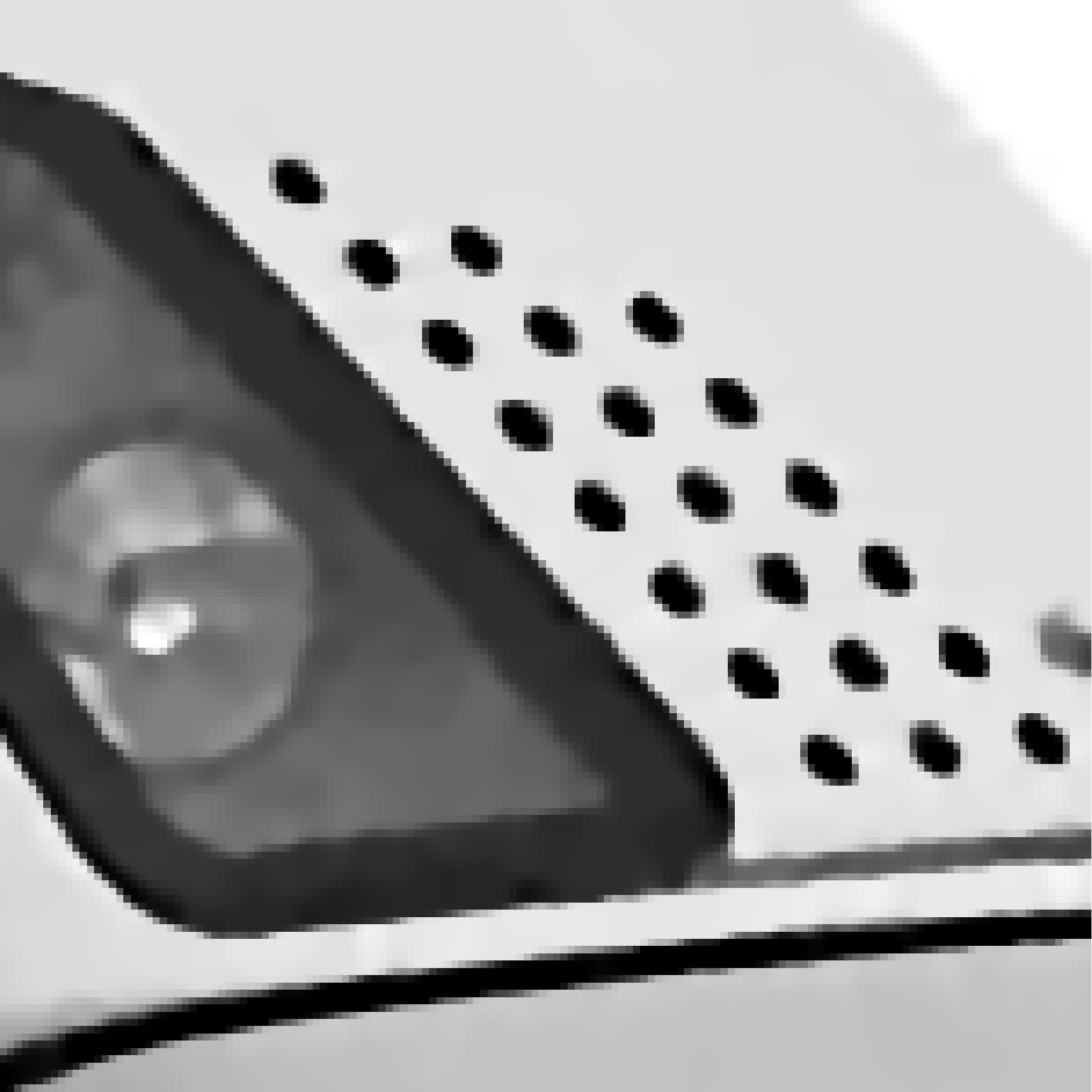}
\end{minipage}
\begin{minipage}{3.150cm}
\includegraphics[width=3.150cm]{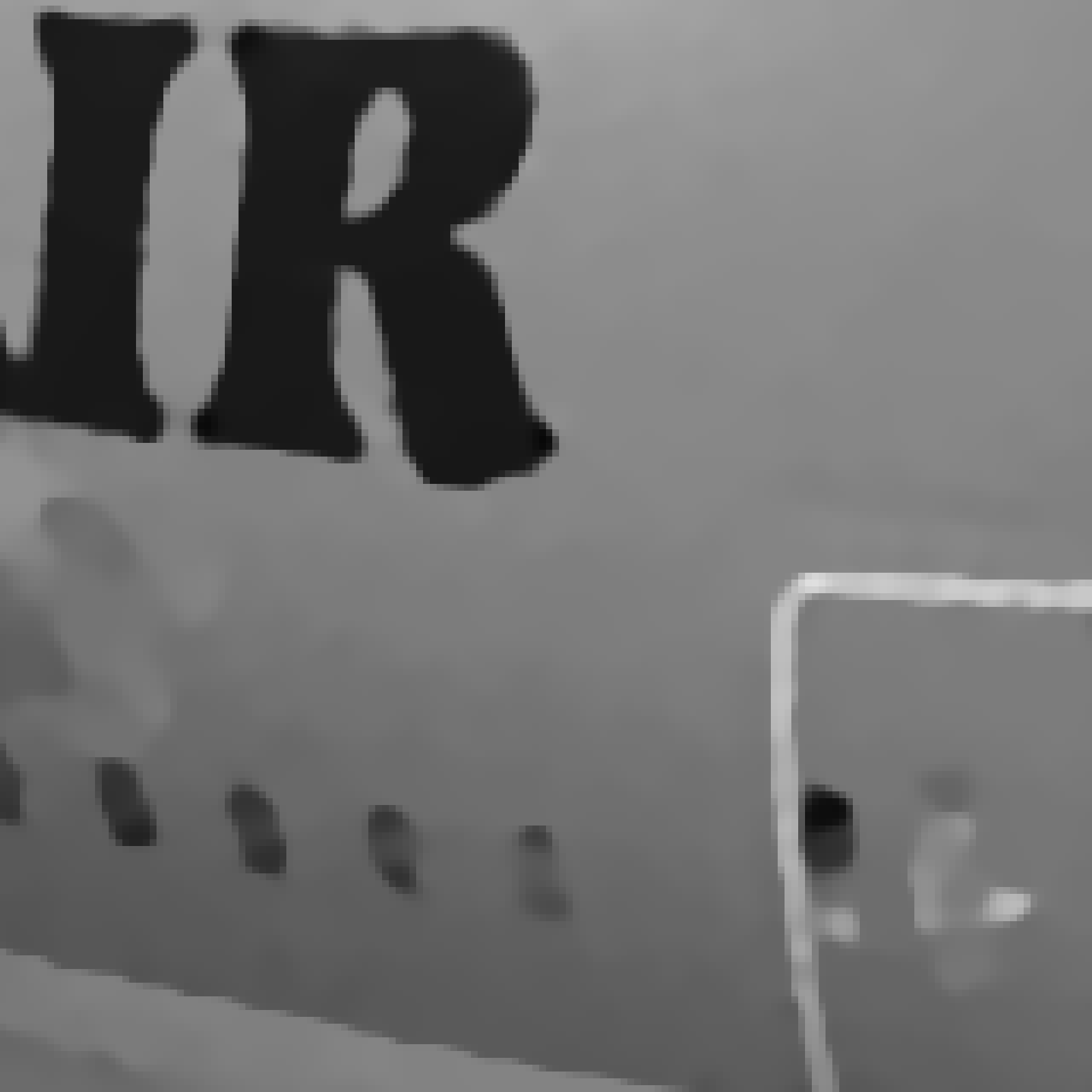}
\end{minipage}
\begin{minipage}{3.150cm}
\includegraphics[width=3.150cm]{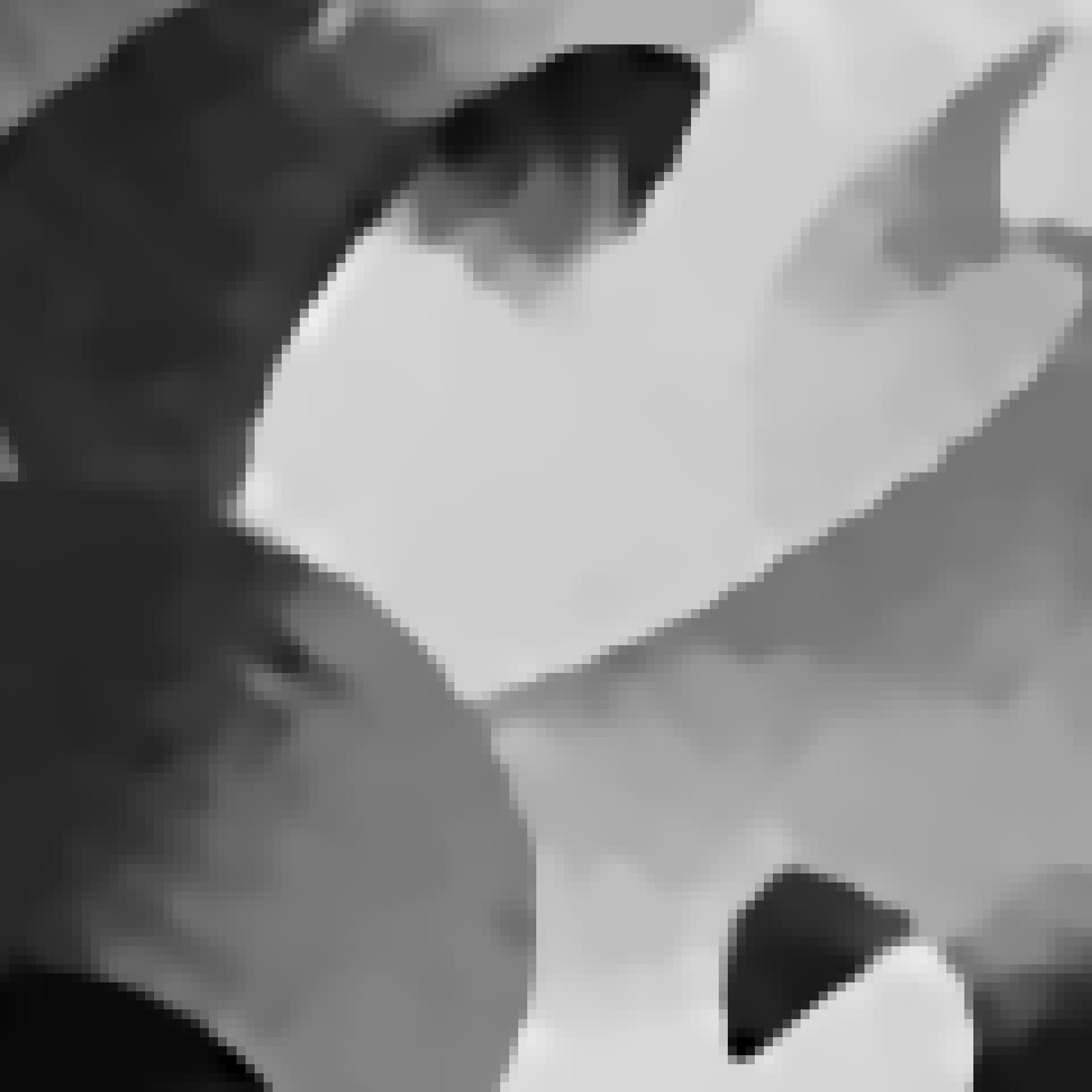}
\end{minipage}
\begin{minipage}{3.150cm}
\includegraphics[width=3.150cm]{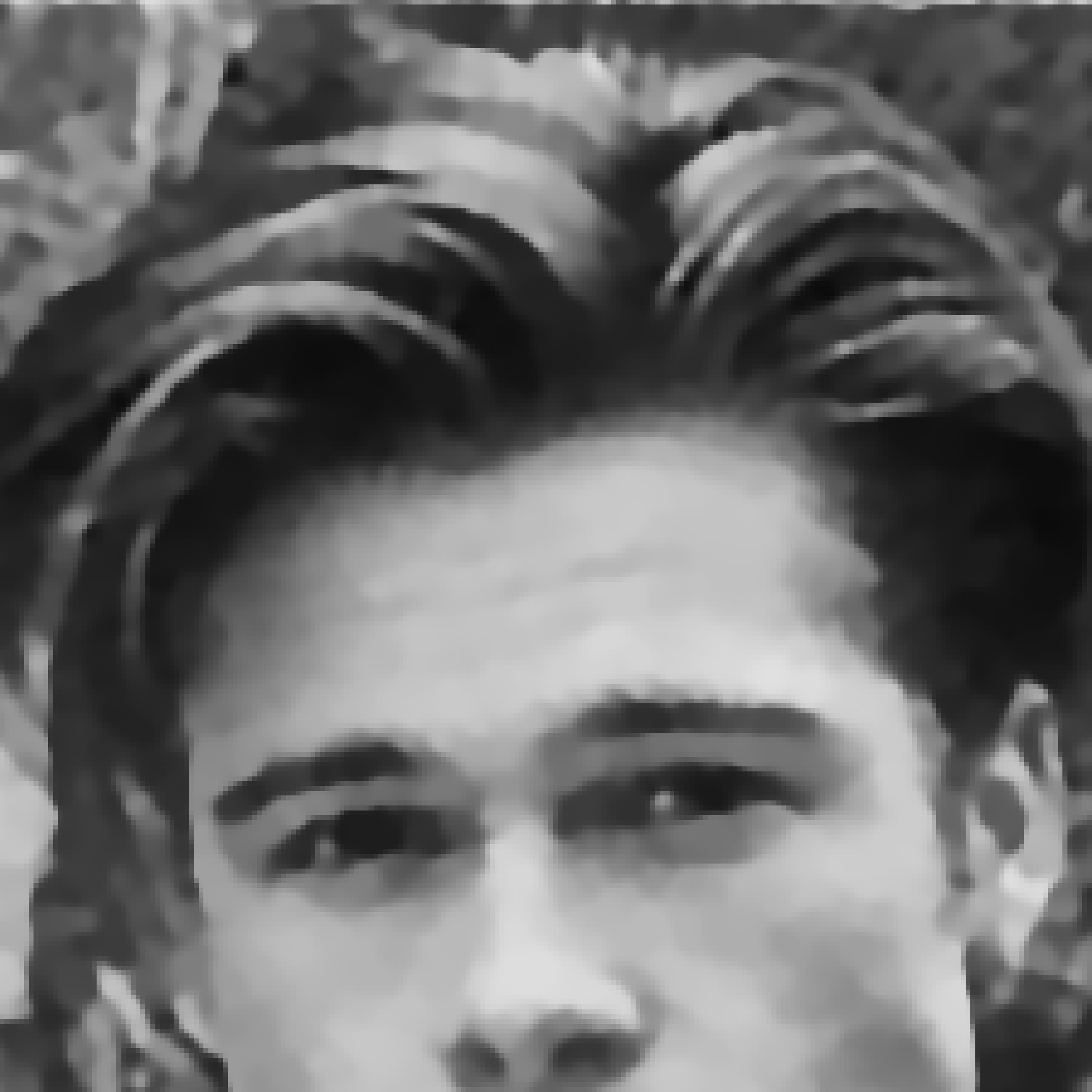}
\end{minipage}\\
\begin{minipage}{3.150cm}\begin{center}{\small{Sonic}}\end{center}\end{minipage}\begin{minipage}{3.150cm}\begin{center}{\small{Train}}\end{center}\end{minipage}\begin{minipage}{3.150cm}\begin{center}{\small{Airplane}}\end{center}\end{minipage}\begin{minipage}{3.150cm}\begin{center}{\small{Oil Painting}}\end{center}\end{minipage}\begin{minipage}{3.150cm}\begin{center}{\small{Pitt}}\end{center}\end{minipage}
\caption{Zoom-in views of Figure \ref{fig:DeblurResults}. The first row describes the degraded measurements, followed by the results of TGV model, PS model, GS model, and our model, respectively.}\label{fig:DeblurResultsZoom}
\end{center}
\end{figure}

The deblurring results of the four models are summarized in Table \ref{tab:PSNRCompareDeblur}, and presented in Figure \ref{fig:DeblurResults} and Figure \ref{fig:DeblurResultsZoom} for visual comparisons. First of all, we can observe that our model \eqref{OurModel} outperforms other three models in terms of PSNR values. The improvements of visual quality are also clearly observable in most cases. It is notable that our model is especially good for images that have gradual changes in intensities, as well as images that have relatively sparsely located singularities, such as the image ``Sonic''.

\section{Asymptotic Analysis}\label{VariationAsymAnal}

This section is devoted to provide an asymptotic analysis for the proposed edge driven model \eqref{OurModel}. We will present a new variational model \eqref{CorrespondingVariational}, and then show that \eqref{OurModel} can be regarded as a discrete approximation to the variational model through $\Gamma$-convergence \cite{Maso1993}. Relations among approximate minimizers of the discrete model and the corresponding variational model are also investigated. Some technical details are postponed to \ref{ProofTh1} and \ref{ProofProp2}.

\subsection{Variational Model and Properties}\label{VariationalModel}


As we will prove in later subsections, the variational model corresponding to our edge driven model \eqref{OurModel} takes the form:
\begin{align}\label{CorrespondingVariational}
\begin{split}
\min_{u,0\leq v\leq1}~\lambda\int_{\Om}(1-v)\left(\sum_{\aal\in\II}|\p^{\aal}u|^2\right)^{\f{1}{2}}\rd\x&+\gamma\int_{\Om}v\left(\sum_{\aal\in\II'}|\p^{\aal}u|^2\right)^{\f{1}{2}}\rd\x\\
&+\rho\int_{\Om}\left(\sum_{\aal\in\II''}|\p^{\aal}v|^2\right)^{\f{1}{2}}\rd\x+\f{1}{2}\big\|Au-f\big\|_{L_2(\Om)}^2,
\end{split}
\end{align}
where $\II$, $\II'$, and $\II''$ are three index sets. Since the first two terms are exchangeable, we impose some restriction on $\II$ and $\II'$ for clarity. Noting that the key features such as edges, ridges  can be well extracted after lower order differentiations, we choose the index sets $\II$ and $\II'$ so that there exists $\aal\in\II$ such that $\aal>\bbe$ for all $\bbe\in\II'$.

To better understand \eqref{CorrespondingVariational}, we consider a special case of it. Letting $\II''=\big\{\aal:|\aal|=1\big\}$, \eqref{CorrespondingVariational} is reduced to the following model:
\begin{align*}
\begin{split}
\min_{u,0\leq v\leq1}~\lambda\int_{\Om}(1-v)\left(\sum_{\aal\in\II}|\p^{\aal}u|^2\right)^{\f{1}{2}}\rd\x+\gamma\int_{\Om}v\left(\sum_{\aal\in\II'}|\p^{\aal}u|^2\right)^{\f{1}{2}}\rd\x+\rho\int_{\Om}|\na v|\rd\x+\f{1}{2}\big\|Au-f\big\|_{L_2(\Om)}^2,
\end{split}
\end{align*}
which can be viewed as a relaxation of
\begin{align}\label{Energy:MS}
\begin{split}
\min_{u,\Sigma\subseteq\Om}~\underbrace{\lambda\int_{\Om\setminus\Sigma}\left(\sum_{\aal\in\II}|\p^{\aal}u|^2\right)^{\f{1}{2}}\rd\x+\gamma\int_{\Sigma}\left(\sum_{\aal\in\II'}|\p^{\aal}u|^2\right)^{\f{1}{2}}\rd\x+\rho\Per(\Sigma;\Om)+\f{1}{2}\big\|Au-f\big\|_{L_2(\Om)}^2}_{:=\wt{E}(u,\Sigma)}
\end{split}
\end{align}
with $\Sigma$ being the estimated region of singularities having positive measure and an interior. Here, $\Per(\Sigma;\Om)$ is the perimeter of a Borel measurable set $\Sigma$ in $\Om$ \cite{H.Attouch2014}. Following \cite{Strang1983,Strang1983a}, we arrive at the following proposition which relates $\bsv$ from the subproblem \eqref{vsubprob} of our wavelet frame model \eqref{OurModel} to the regions with singularities.

\begin{proposition}\label{Prop3} For any given fixed $u$, we can find the global minimizer of $\wt{E}(u,\cdot)$ (given by \eqref{Energy:MS}) by solving the convex minimization problem
\begin{align*}
\min_{0\leq v\leq1}~\underbrace{\lambda\int_{\Om}(1-v)\left(\sum_{\aal\in\II}|\p^{\aal}u|^2\right)^{\f{1}{2}}\rd\x+\gamma\int_{\Om}v\left(\sum_{\aal\in\II'}|\p^{\aal}u|^2\right)^{\f{1}{2}}\rd\x+\rho\int_{\Om}|\na v|\rd\x}_{:=E_u(v)}
\end{align*}
and setting $\Sigma=\big\{\x\in\Om:v(\x)>t\big\}$ for almost every $t\in[0,1]$.
\end{proposition}

\begin{pf} The proof is similar to \cite[Theorem 2]{T.F.Chan2006}. However, for completeness, we include the proof. Since $v$ takes its values in $[0,1]$, the co-area formula \cite{H.Attouch2014} tells us that
\begin{align*}
\int_{\Om}|\na v|\rd\x=\int_0^1\Per\left(\big\{\x\in\Om:v(\x)>t\big\};\Om\right)\rd t.
\end{align*}
Let $\Sigma(t):=\big\{\x\in\Om:v(\x)>t\big\}$. For a fixed $u$, we have
\begin{align*}
\int_{\Om}v(\x)\left(\sum_{\aal\in\II'}|\p^{\aal}u(\x)|^2\right)^{\f{1}{2}}\rd\x&=\int_{\Om}\int_0^1 \chi_{[0,v(\x)]}(t)\left(\sum_{\aal\in\II'}|\p^{\aal}u(\x)|^2\right)^{\f{1}{2}}\rd t\rd\x\\
&=\int_0^1\int_{\Om}\chi_{\Sigma(t)}(\x)\left(\sum_{\aal\in\II'}|\p^{\aal}u(\x)|^2\right)^{\f{1}{2}}\rd\x \rd t=\int_0^1\int_{\Sigma(t)}\left(\sum_{\aal\in\II'}|\p^{\aal}u(\x)|^2\right)^{\f{1}{2}}\rd\x \rd t.
\end{align*}
where $\chi_{\Sigma}$ is the characteristic function of a set $\Sigma$; $\chi_{\Sigma}(\x)=1$ if $\x\in\Sigma$ and $\chi_{\Sigma}(\x)=0$ otherwise. Similarly, we have
\begin{align*}
\int_{\Om}(1-v(\x))\left(\sum_{\aal\in\II}|\p^{\aal}u(\x)|^2\right)^{\f{1}{2}}\rd\x=\int_0^1\int_{\Om\setminus\Sigma(t)}\left(\sum_{\aal\in\II}|\p^{\aal}u(\x)|^2\right)^{\f{1}{2}}\rd\x \rd t.
\end{align*}
Combining the above three equalities, we have
\begin{align*}
E_u(v)&=\int_0^1\left[\lambda\int_{\Om\setminus\Sigma(t)}\left(\sum_{\aal\in\II}|\p^{\aal}u|^2\right)^{\f{1}{2}}\rd\x+\gamma\int_{\Sigma(t)}\left(\sum_{\aal\in\II'}|\p^{\aal}u|^2\right)^{\f{1}{2}}\rd\x+\rho\Per(\Sigma(t);\Om)\right]\rd t\\
&=\int_0^1 \wt{E}(u,\Sigma(t))\rd t.
\end{align*}
It follows that if $v$ is a minimizer of $E_u$, then for a.e. $t\in[0,1]$, $\Sigma(t)$ has to be a minimizer of $\wt{E}(u,\cdot)$.\qquad$\square$
\end{pf}

Now, we consider the $u$-subproblem of \eqref{CorrespondingVariational} when $\II''=\big\{\aal:|\aal|=1\big\}$. By virtue of Proposition \ref{Prop3}, it suffices to consider the following problem:
\begin{align}\label{Varusubprob}
\min_u~\lambda\int_{\Om\setminus\Sigma}\left(\sum_{\aal\in\II}|\p^{\aal}u|^2\right)^{\f{1}{2}}\rd\x+\gamma\int_{\Sigma}\left(\sum_{\aal\in\II'}|\p^{\aal}u|^2\right)^{\f{1}{2}}\rd\x+\f{1}{2}\big\|Au-f\big\|_{L_2(\Om)}^2
\end{align}
for a fixed $\Sigma\subseteq\Om$. Then we can see how \eqref{Varusubprob} is related to several existing variational and PDE models for image restoration:
\begin{enumerate}
\item When $\II=\big\{\aal:|\aal|=2\big\}$ and $\II'=\big\{\aal:|\aal|=1\big\}$, \eqref{Varusubprob} is reduced to
\begin{align}\label{Varusubspecial}
\min_u~\lambda\int_{\Om\setminus\Sigma}|\na^2 u|\rd\x+\gamma\int_{\Sigma}|\na u|\rd\x+\f{1}{2}\big\|Au-f\big\|_{L_2(\Om)}^2,
\end{align}
which is a special type of the combined first and second order total variation (TV) model \cite{M.Bergounioux2010,J.Liang2015,K.Papafitsoros2014}. More precisely, let $\alpha(\x)=\lambda\chi_{\Om\setminus\Sigma}(\x)$ and $\beta(\x)=\gamma\chi_{\Sigma}(\x)$. Then we have the following combined first and second order TV model with spatially varying parameters \cite{K.Papafitsoros2014}
\begin{align*}
\min_u~\int_{\Om}\alpha(\x)|\na^2u(\x)|\rd\x+\int_{\Om}\beta(\x)|\na u(\x)|\rd\x+\f{1}{2}\big\|Au-f\big\|_{L_2(\Om)}^2.
\end{align*}

\item In \cite{G.Aubert2006}, the gradient descent flow of \eqref{Varusubspecial} is studied:
\begin{align}\label{PDE}
\f{\p u}{\p t}=-\lambda\div^2\left(\chi_{\Om\setminus\Sigma}\f{\na^2u}{|\na^2u|}\right)+\gamma\div\left(\chi_{\Sigma^{\circ}}\f{\na u}{|\na u|}\right)-A^T(Au-f).
\end{align}
We can easily see that there are two different nonlinear diffusions in region $\Om\setminus\Sigma$ and $\Sigma^{\circ}$, where $\Sigma^{\circ}$ stands for the interior of $\Sigma$. The second order nonlinear diffusion in $\Sigma^{\circ}$ plays a role of edge-enhancing, while the fourth order nonlinear diffusion in $\Om\setminus\Sigma$ plays a role of preventing smooth regions from being blocky \cite{B.Dong2013a,Y.L.You2000}.

\item The $u$-subproblem \eqref{Varusubprob} can be viewed (formally) as a generalized inf-convolution model \cite{A.Chambolle1997} as well; we define
\begin{align*}
u_1=u\chi_{\Om\setminus\Sigma}~~~~\text{and}~~~~u_2=u\chi_{\Sigma^{\circ}},
\end{align*}
and we set $\II$ and $\II'$ as in \eqref{Varusubspecial}. Then $u=u_1+u_2$ almost everywhere in $\Om$, and \eqref{Varusubprob}, namely \eqref{Varusubspecial} reduces to the following inf-convolution model:
\begin{align}\label{InfConv}
\min_{u_1,u_2}~\lambda\int_{\Om}|\na^2u_1|\rd\x+\gamma\int_{\Om}|\na u_2|\rd\x+\f{1}{2}\big\|A(u_1+u_2)-f\big\|_{L_2(\Om)}^2.
\end{align}
Moreover, \eqref{InfConv} can be rewritten as
\begin{align}\label{InfConv2}
\min_{u,u_1}~\lambda\int_{\Om}|\na(\na u_1)|\rd\x+\gamma\int_{\Om}|\na u-\na u_1|\rd\x+\f{1}{2}\big\|Au-f\big\|_{L_2(\Om)}^2,
\end{align}
which is a special case of the following (unsymmetrized) TGV model
\begin{align}\label{TGVSpecial}
\min_{u,w}~\lambda\int_{\Om}|\na w|\rd\x+\gamma\int_{\Om}|\na u-w|\rd\x+\f{1}{2}\big\|Au-f\big\|_{L_2(\Om)}^2.
\end{align}
\end{enumerate}

As we can see from the above discussions, the variational model \eqref{CorrespondingVariational} is an edge driven variational model which restores piecewise smooth functions by inflicting varied strength of regularization in smooth and sharp image regions and simultaneously restoring image singularities. Since the proposed discrete model \eqref{OurModel} approximates the variational model \eqref{CorrespondingVariational} as will be shown in the next subsection, we can make the same assertion on \eqref{OurModel}. Furthermore, the proposed model \eqref{OurModel} can be viewed as a more general image restoration model than the aforementioned variational models.

\subsection{Analysis}\label{Connection}

In this subsection, we find a connection between the model \eqref{OurModel} and the variational model \eqref{CorrespondingVariational}. As will be revealed in our analysis, $\lam\cdot\bsW$ can approximate various differential operators by choosing an appropriate weight for each of framelet bands. Therefore, for simplicity, we shall restrict $\bsW=\bsW'=\bsW''$ in \eqref{OurModel} and analyze the following problem
\begin{align}\label{GenericModel}
\min_{\bu,0\leq\bsv\leq1}~\left\|(\one-\bsv)\cdot\big(\lam\cdot\bsW\bu\big)\right\|_1+\left\|\bsv\cdot\big(\gga\cdot\bsW\bu\big)\right\|_1+\big\|\rrh\cdot\bsW\bsv\big\|_1+\f{1}{2}\big\|\bsA\bu-\bsf\big\|_2^2
\end{align}
with $\big\{\lam\big\}$, $\big\{\gga\big\}$, and $\big\{\rrh\big\}$ chosen differently for different framelet bands. We further assume, for simplicity, that $\bsW$ is the wavelet frame transform of piecewise B-spline wavelet frame systems. By virtue of Proposition \ref{Prop1}, it is not hard to see that our analysis can be generalized to the more general case \eqref{OurModel}.

We start with introducing some symbols and notation that will be used throughout the rest of the paper.

\begin{notation}\label{Notation1} We focus our analysis on $\R^2$, i.e. the two-dimensional cases.  All the two-dimensional refinable functions and framelets are assumed to be constructed by tensor products of univariate B-splines and the associated framelets obtained from the UEP \cite{A.Ron1997}.
\begin{enumerate}
\item All functions we consider are defined on $\Om=(0,1)^2\subseteq\R^2$, and that their discrete versions, i.e. digital images are defined on an $N\times N$ cartesian grid on $\overline{\Om}=[0,1]^2$ with $N=2^n+1$ for $n\geq 0$. We denote by $h=2^{-n}$ the meshsize of the $N\times N$ grid.

\item The bold-face letters ($\aal$, $\bbe$, $\bsi$, $\bsj$, $\bk$, etc.) are used to denote the double indices in $\Z^2$. We denote
\begin{align*}
\OO_n=\big\{\bk\in\Z^2:2^{-n}\bk\in\overline{\Om}\big\}
\end{align*}
as the set of indices of the $N\times N$ Cartesian grid.


\item Given a wavelet frame system and its corresponding refinable function $\phi$, we define
\begin{align*}
\MM_n=\big\{\bk\in\OO_n:\Lambda_{n,\bk}:=\su(\phi_{n,\bk})\subseteq\overline{\Om}\big\}.
\end{align*}
Note that since piecewise B-spline wavelet frame systems are used, we have $\su(\psi_{\aal})=\su(\phi)$ for all $\aal\in\BB=\big\{0,\cdots,r\big\}^2\setminus\big\{\0\big\}$, so that
\begin{align*}
\su(\psi_{\aal,n,\bk})=\su(\varphi_{\aal,n,\bk})=\Lambda_{n,\bk}
\end{align*}
for all $n\in\N$ and $\bk\in\Z^2$.

\item The spaces to which $\bu$ and the components of $\bsv$ belong are respectively denoted as $\R^{\MM_n}$ and $[0,1]^{\MM_n}$. Here, for given sets $A$ and $B$, $B^A=\big\{f:A\rightarrow B\big\}$ denotes the space of all functions mapping from $A$ to $B$. Note that since $\MM_n$ is a finite set, we have $\R^{\MM_n}\simeq\R^{|\MM_n|}$ and $[0,1]^{\MM_n}\simeq[0,1]^{|\MM_n|}$.

\item For the simplicity, we assume that the level of decomposition is $1$, i.e. $L=L''=1$, while it is not hard to extend our analysis to $L$, $L''>1$ as mentioned in \cite{B.Dong2016}. Note that if $L=L''=1$, then $\bsv\in[0,1]^{\MM^2}$.

\item We define the index set $\KK_n\subseteq\MM_n$ by
\begin{align*}
\KK_n:=\left\{\bk\in\MM_n:\bk+S_{\aal}\subseteq\MM_n~~\text{for all}~~\aal\in\BB\cup\{\0\}\right\}
\end{align*}
where $S_{\aal}$ is the support of $\bq_{\aal}$. In other words, $\KK_n$ consists of double indices such that the boundary condition of $\bq_{\aal}[-\cdot]\circledast\bu$ is inactive for all $\aal\in\BB\cup\{\0\}$, so that $\bq_{\aal}\ast\bu$ is well defined, and $\bsW_{\aal}:\R^{\MM_n}\rightarrow\R^{\KK_n}$ for all $\aal\in\BB\cup\{\0\}$. In addition, note that $\OO_n$, $\MM_n$, and $\KK_n$ all depend on the resolution $n$.

\item In order to link the continuous and the discrete settings, we need to take resolution into account. Hence, for any $\bu\in\R^{\MM_n}$, the discrete $\ell_p$ norm we are using is defined as
\begin{align*}
\|\bu\|_p^p:=h^2\sum_{\bk\in\MM_n}\big|\bu[\bk]\big|^p.
\end{align*}
\end{enumerate}
\end{notation}

Using the above notation, we can take image resolution into account in model \eqref{GenericModel}. Namely, the first three terms in \eqref{GenericModel} are respectively defined as
\begin{align*}
\left\|(\one-\bsv)\cdot\big(\lam\cdot\bsW\bu\big)\right\|_1&=h^2\sum_{\bk\in\KK_n}(\one-\bsv[\bk])\left(\sum_{\aal\in\BB}\lambda_{\aal}[\bk]\bigg|\big(\bsW_{\aal}\bu\big)[\bk]\bigg|^2\right)^{\f{1}{2}},\\
\left\|\bsv\cdot\big(\gga\cdot\bsW\bu\big)\right\|_1&=h^2\sum_{\bk\in\KK_n}\bsv[\bk]\left(\sum_{\aal\in\BB}\gamma_{\aal}[\bk]\bigg|\big(\bsW_{\aal}\bu\big)[\bk]\bigg|^2\right)^{\f{1}{2}},\\
\big\|\rrh\cdot\bsW\bsv\big\|_1&=h^2\sum_{\bk\in\KK_n}\left(\sum_{\aal\in\BB}\rho_{\aal}[\bk]\bigg|\big(\bsW_{\aal}\bsv\big)[\bk]\bigg|^2\right)^{\f{1}{2}}.
\end{align*}

To analyze the relation between \eqref{GenericModel} and \eqref{CorrespondingVariational}, we first reformulate the objective function \eqref{GenericModel} to a functional defined on the same function spaces as that of \eqref{CorrespondingVariational}. Denote the energy functional of the variational model \eqref{CorrespondingVariational} as
\begin{align}\label{CDZ11Variational}
\begin{split}
E(u,v)=\lambda\int_{\Om}(1-v)\left(\sum_{\aal\in\II}|\p^{\aal}u|^2\right)^{\f{1}{2}}\rd\x&+\gamma\int_{\Om}v\left(\sum_{\aal\in\II'}|\p^{\aal}u|^2\right)^{\f{1}{2}}\rd\x\\
 &+\rho\int_{\Om}\left(\sum_{\aal\in\II''}|\p^{\aal}v|^2\right)^{\f{1}{2}}\rd\x+\f{1}{2}\big\|Au-f\big\|_{L_2(\Om)}^2,
\end{split}
\end{align}
where $\II$ and $\II'$ are chosen such that there exists $\aal\in\II$ such that $\aal>\bbe$ for all $\bbe\in\II'$, and $u\in W_1^s(\Om)$ and $v\in W_1^r(\Om,[0,1])$. Here, $W_1^s(\Om)$ is the Sobolev space defined as \eqref{W_1^sSobolev}  and $W_1^r(\Om,[0,1])$ is defined as
\begin{align*}
W_1^r(\Om,[0,1])=\big\{v\in W_1^r(\Om):0\leq v\leq 1~~\text{a.e. in}~\Om\big\}
\end{align*}
with $s=\max_{\aal\in\II\cup\II'}|\aal|$ and $r=\max_{\aal\in\II'}|\aal|$. Then by Sobolev imbedding theorem \cite{Adams1975,H.Attouch2014}, $W_1^r(\Om,[0,1])\subseteq W_1^r(\Om)\subseteq L_2(\Om)$.

Let $\phi$ be the refinable function corresponding to $\bsW$. Define a linear operator $\bsT_n$ on $L_2(\Om)$ by
\begin{align*}
\bsT_n u=\big\{2^n\la u,\phi_{n,\bk}\ra:\bk\in\MM_n\big\}\in\R^{\MM_n}.
\end{align*}
Then we define
\begin{align}\label{CDZ11_SemiDiscrete}
\begin{split}
E_n(u,v)=\left\|(\one-\bsT_n v)\cdot\big(\lam_n\cdot\bsW_n\bsT_nu\big)\right\|_1+\left\|\bsT_n v\cdot\big(\gga_n\cdot\bsW_n\bsT_nu\big)\right\|_1+\big\|\rrh_n\cdot\bsW_n\bsT_n v\big\|_1+\f{1}{2}\big\|\bsA_n\bsT_nu-\bsT_nf\big\|_2^2.
\end{split}
\end{align}
For notational simplicity, we will denote the energy functional in \eqref{GenericModel} by $F_n$:
\begin{align}\label{CDZ11_Discrete}
F_n(\bu_n,\bsv_n)=\left\|(\one-\bsv_n)\cdot\big(\lam_n\cdot\bsW_n\bu_n\big)\right\|_1+\left\|\bsv_n\cdot\big(\gga_n\cdot\bsW_n\bu_n\big)\right\|_1+\big\|\rrh_n\cdot\bsW_n\bsv_n\big\|_1+\f{1}{2}\big\|\bsA_n\bu_n-\bsf_n\big\|_2^2
\end{align}
where the subscript $n$ is used to emphasize the dependence of $\bsW$ and $\bsA$ on the image resolution $n$. We first consider
\begin{align*}
P_F&=\inf\big\{F_n(\bu_n,\bsv_n):\bu_n\in\R^{\MM_n},~\bsv_n\in[0,1]^{\MM_n}\big\}\\
P_E&=\inf\big\{E_n(u,v):u\in W_1^s(\Om),~v\in W_1^r(\Om,[0,1])\big\}.
\end{align*}
Then it is obvious that $P_F\leq P_E$ because for every $(u,v)\in W_1^s(\Om)\times W_1^r(\Om,[0,1])$,
\begin{align*}
F_n(\bsT_n u,\bsT_n v)=E_n(u,v).
\end{align*}
Note that in general, we do not have $P_F=P_E$ because $\bsv_n\in[0,1]^{\MM_n}$ may not necessarily lie in $\bsT_n(W_1^r(\Om,[0,1])$. Indeed, $\bsT_n(L_2(\Om,[0,1])$ where $L_2(\Om,[0,1])=\big\{u\in L_2(\Om):0\leq u\leq 1~\text{a.e. in}~\Om\big\}$ is a proper subset of $[0,1]^{\MM_n}$.


\begin{rmk}\label{rmk} We further mention that in fact it is not necessary to impose the restriction on $\bsW''$. Using the refinable function $\phi''$ corresponding to the piecewise B-spline wavelet frame system $\bsW''$ and defining corresponding index sets appropriately, we can establish the relation between (the reformulation of) the following model
\begin{align*}
\min_{\bu,0\leq\bsv\leq1}~\left\|(\one-\bsv)\cdot\big(\lam\cdot\bsW\bu\big)\right\|_1+\left\|\bsv\cdot\big(\gga\cdot\bsW\bu\big)\right\|_1+\big\|\rrh\cdot\bsW''\bsv\big\|_1+\f{1}{2}\big\|\bsA\bu-\bsf\big\|_2^2
\end{align*}
and the variational model \eqref{CorrespondingVariational}. Nevertheless, for simplicity, we focus on analyzing the relation between \eqref{CDZ11_SemiDiscrete} and \eqref{CorrespondingVariational}.
\end{rmk}

For convenience, we write $E_n$ and $E$ respectively as
\begin{align*}
E_n(u,v)&=\left\|(\one-\bsT_n v)\cdot\big(\lam_n\cdot\bsW_n\bsT_nu\big)\right\|_1+\left\|\bsT_n v\cdot\big(\gga_n\cdot\bsW_n\bsT_nu\big)\right\|_1+\big\|\rrh_n\cdot\bsW_n\bsT_n v\big\|_1+\f{1}{2}\big\|\bsA_n\bsT_nu-\bsT_nf\big\|_2^2\\
&=E_n^{(1)}(u,v)+E_n^{(2)}(u,v)+E_n^{(3)}(v)+E_n^{(4)}(u),
\end{align*}
and
\begin{align*}
E(u,v)&=\lambda\int_{\Om}(1-v)\left(\sum_{\aal\in\II}|\p^{\aal}u|^2\right)^{\f{1}{2}}\rd\x+\gamma\int_{\Om}v\left(\sum_{\aal\in\II'}|\p^{\aal}u|^2\right)^{\f{1}{2}}\rd\x+\rho\int_{\Om}\left(\sum_{\aal\in\II''}|\p^{\aal}v|^2\right)^{\f{1}{2}}\rd\x+\f{1}{2}\big\|Au-f\big\|_{L_2(\Om)}^2\\
&=E^{(1)}(u,v)+E^{(2)}(u,v)+E^{(3)}(v)+E^{(4)}(u).
\end{align*}
Here, without loss of generality, we assume that $\lambda=\gamma=\rho=1$ for $E(u,v)$. To draw an asymptotic relation between $E_n$ and $E$, we need the assumptions on the operator $A$ and its discretization $\bsA_n$, and the parameters $\big\{\lam_n\big\}$, $\big\{\gga_n\big\}$, and $\big\{\rrh_n\big\}$:

\begin{enumerate}
\item[A1.] $A$ is a continuous linear operator mapping $L_2(\Om)$ into itself, and its discretization $\bsA_n$ satisfies
\begin{align}\label{Assumption_A}
\lim_{n\rightarrow\infty}\|\bsT_nAu-\bsA_n\bsT_nu\|_2=0~~~~~\text{for all}~~~u\in L_2(\Om).
\end{align}
Note that $A$ which corresponds to denoising, deblurring, and inpainting satisfies the above assumption \cite{J.F.Cai2012,J.F.Cai2016,B.Dong2016}.
\item[A2.] We split the framelet band $\BB$ into $\BB=\II\cup\JJ$ where $\II$ is the index set in $E^{(1)}(u,v)$. For $\aal\in\II$, we set $\lambda_{\aal}=\big(c_{\aal}^{-1}2^{|\aal|(n-1)}\big)^2$, where $c_{\aal}$ is given in Proposition \ref{Prop1}. For $\aal\in\JJ$, we set $0\leq\lambda_{\aal}\leq O(2^{2|\bbe|(n-1)})$ for some $\bbe\in\BB\cup\big\{\0\big\}$ such that $\0\leq\bbe<\aal$ and $|\bbe|\leq s$. The remaining parameters $\big\{\gga_n\big\}$ and $\big\{\rrh_n\big\}$ are defined as in the similar way except for changing $\II$ with $\II'$ in $E^{(2)}(u,v)$ and $\II''$ in $E^{(3)}(v)$ respectively. In particular, we replace $s$ with $r$ when we set $\big\{\rrh_n\big\}$.
\end{enumerate}

It remains to impose an appropriate topology on $W_1^r(\Om,[0,1])$ which makes it complete. To do this, we define $\msX=W_1^r(\Om)\cap L_{\infty}(\Om)$ equipped with the norm defined by
\begin{align*}
\|v\|_{\msX}=\|v\|_{W_1^r(\Om)}+\|v\|_{L_{\infty}(\Om)}.
\end{align*}
Note that $\msX$ equipped with the norm defined above is a Banach space, and $W_1^r(\Om,[0,1])$ is closed in $\msX$. Hence, in what follows, by a topology on $W_1^r(\Om,[0,1])$, we mean the subspace topology inherited from $\msX$.

The first relation between $E_n$ and $E$ that we want to present is the pointwise convergence of $E_n(u,v)$ to $E(u,v)$ for each $(u,v)$. Since the proof is long and technical, it is postponed to \ref{ProofTh1}.

\begin{thm}[Pointwise Convergence]\label{Th1} Assume that A1 and A2 are satisfied. Then for any $(u,v)\in W_1^s(\Om)\times W_1^r(\Om,[0,1])$, we have
\begin{align}\label{PointwiseConvergence}
\lim_{n\rightarrow\infty}E_n(u,v)=E(u,v).
\end{align}
\end{thm}

With Theorem \ref{Th1}, we can show that the sequence $\big\{E_n:n\in\N\big\}$ is equicontinuous.

\begin{proposition}\label{Prop2} Assume that A1 and A2 are satisfied. Let $(u,v)\in W_1^s(\Om)\times W_1^r(\Om,[0,1])$ be given. Then for every $\eps>0$, there exist $\delta>0$ and $\mN\in\N$ both of which are independent of $n$ such that for any $(u',v')\in W_1^s(\Om)\times W_1^r(\Om,[0,1])$ with $\|u'-u\|_{W_1^s(\Om)}+\|v'-v\|_{\msX}<\delta$ and $n>\mN$, we have $|E_n(u',v')-E_n(u,v)|<\eps$.
\end{proposition}

\begin{pf} See \ref{ProofProp2}.\qquad$\square$
\end{pf}

With the aid of Theorem \ref{Th1} and Proposition \ref{Prop2}, we have the following theorem showing that the convergence of $E_n$ to $E$ is stronger than pointwise convergence. A direct consequence of such convergence is the $\Gamma$-convergence of $E_n$ to $E$ in $W_1^s(\Om)\times W_1^r(\Om,[0,1])$ with the subspace topology inherited from $W_1^s(\Om)\times\msX$. The proof is almost the same as \cite[Theorem 3.1]{B.Dong2016} provided that Theorem \ref{Th1} and Proposition \ref{Prop2} are established. Therefore, we shall omit the proof of Theorem \ref{Th2}.

\begin{definition} Let $\msY$ be a topological space. Given $E_n$, $E:\msY\rightarrow\overline{\R}$, we say that $E_n$ $\Gamma$-converges to $E$ in $\msY$ if
\begin{enumerate}
\item for every sequence $u_n\rightarrow u$ in $\msY$, $E(u)\leq\liminf_{n\rightarrow\infty}E_n(u_n)$,
\item for every $u\in\msY$, there is a sequence $u_n\rightarrow u$ in $\msY$ such that $E(u)\geq\limsup_{n\rightarrow\infty}E_n(u_n)$.
\end{enumerate}
\end{definition}

\begin{thm}\label{Th2} Suppose that the assumptions A1 and A2 are satisfied. For every $(u_n,v_n)$, $(u,v)\in W_1^s(\Om)\times W_1^r(\Om,[0,1])$ with
\begin{align*}
\lim_{n\rightarrow\infty}\left(\|u_n-u\|_{W_1^s(\Om)}+\|v_n-v\|_{\msX}\right)=0,
\end{align*}
we have
\begin{align*}
\lim_{n\rightarrow\infty} E_n(u_n,v_n)=E(u,v).
\end{align*}
Consequently, $E_n$ $\Gamma$-converges to $E$ in $W_1^s(\Om)\times W_1^r(\Om,[0,1])$ with the subspace topology inherited from $W_1^s(\Om)\times\msX$.
\end{thm}

From a practical point of view, it is more important to relate the (approximate) solutions of the optimizations problems. Recall that $(u^*,v^*)$ is the $\eps$-minimizer of $E(u,v)$ if
\begin{align*}
E(u^*,v^*)\leq\inf_{u,v}~E(u,v)+\eps~~~~\text{for some}~~~\eps>0.
\end{align*}
In particular, $(u^*,v^*)$ is the minimizer of $E$ if $E(u^*,v^*)=\inf_{u,v}~E(u,v)$. Theorem \ref{Th2} implies the following relation between the ($\eps$-)minimizers of the original discrete model $F_n$ in \eqref{CDZ11_Discrete} and the variational model $E$ in \eqref{CDZ11Variational}.

\begin{corollary} Let $(\bu_n^*,\bsv_n^*)$ be an $\eps$-minimizer of $F_n$ for a given $\eps>0$ and for all $n$. Then we have
\begin{align*}
\limsup_{n\rightarrow\infty} F_n(\bu_n^*,\bsv_n^*)\leq \inf_{u,v}~E(u,v)+\eps.
\end{align*}
In particular, when $(\bu_n^*,\bsv_n^*)$ is a minimizer of $F_n$, then
\begin{align*}
\limsup_{n\rightarrow\infty} F_n(\bu_n^*,\bsv_n^*)\leq \inf_{u,v}~E(u,v).
\end{align*}
\end{corollary}

\begin{pf} For a given $(u,v)\in W_1^s(\Om)\times W_1^r(\Om,[0,1])$, let $(u_n,v_n)$ be the sequence as given in item 2 of the definition of $\Gamma$-convergence. Together with $\inf_{\bu_n,\bsv_n}~F_n(\bu_n,\bsv_n)\leq\inf_{u,v}~E_n(u,v)$, we have
\begin{align*}
E(u,v)&\geq\limsup_{n\rightarrow\infty}E_n(u_n,v_n)\geq\limsup_{n\rightarrow\infty}\left(\inf_{u,v}~E_n(u,v)\right)\\
&\geq\limsup_{n\rightarrow\infty}\left(\inf_{\bu_n,\bsv_n}~F_n(\bu_n,\bsv_n)\right)\geq\limsup_{n\rightarrow\infty} F_n(\bu_n^*,\bsv_n^*)-\eps,
\end{align*}
which completes the proof.\qquad$\square$
\end{pf}

\section{Conclusion}

In this paper, we proposed a new edge driven wavelet frame based image restoration model by approximating images as piecewise smooth functions. The proposed model inflicts different strength of regularization in smooth image regions and near image singularities such as edges, and actively regularize image singularities at the same time. The performance gain of the proposed model over the existing piecewise smooth image restoration models is mainly due to its robustness to the estimation of image singularities and better regularization on the singularity set.  Finally, the formulation of using an implicit representation of the singularities set also enables an asymptotic analysis of the proposed edge driven model and a rigorous connection between the discrete  model and a  general variational model in the continuum setting.

{\appendix
\section{Proof of Proposition \ref{Prop1}}\label{ProofProp1}

Since $\psi_{\aal}\in L_2(\R^2)$ is constructed by the tensor product of the univariate framelets, we first consider one-dimensional case. Let $\psi_{\alpha}\in L_2(\R)$ have vanishing moments of order $\alpha$, and let $K=\su(\psi_{\alpha})$. From the assumption, $K$ is a closed interval. We also denote by $H_K$ the supporting function on $K$:
\begin{align*}
H_K(\xi)=\sup_{x\in K}x\xi.
\end{align*}
Since $\psi_{\alpha}$ has vanishing moments of order $\alpha$, it follows that
\begin{align}\label{Vanishing}
\int_{-\infty}^{\infty}x^{\beta}\psi_{\alpha}(x)\rd x=i^{\beta}\wh{\psi}_{\alpha}^{(\beta)}(0)=0
\end{align}
for all $\beta<\alpha$, but $\int_{\infty}^{\infty}x^{\alpha}\psi(x)\rd x=i^{\alpha}\wh{\psi}_{\alpha}^{(\alpha)}(0)\neq0$. Since $\psi_{\alpha}$ is compactly supported, its Fourier transform
\begin{align*}
\wh{\psi}_{\alpha}(\xi)=\int_{-\infty}^{\infty}\psi_{\alpha}(x)e^{-i\xi x}\rd x
\end{align*}
can be extended to an entire function of $\zeta\in\C$, called Fourier-Laplace transform, which satisfies \eqref{Vanishing}. Then the Taylor series expansion of $\wh{\psi}_{\alpha}$ at $0$ satisfies
\begin{align*}
\wh{\psi}_{\alpha}(\zeta)=\sum_{\beta=0}^{\infty}\f{\wh{\psi}_{\alpha}^{(\beta)}(0)}{\beta !}\zeta^{\beta}=\sum_{\beta=\alpha}^{\infty}\f{\wh{\psi}_{\alpha}^{(\beta)}(0)}{\beta !}\zeta^{\beta}=\zeta^{\alpha}\sum_{\beta=0}^{\infty}\f{\wh{\psi}_{\alpha}^{(\alpha+\beta)}(0)}{(\alpha+\beta) !}\zeta^{\beta}.
\end{align*}
In other words, there exists an entire function $g_{\alpha}$ such that
\begin{align}\label{Entire_Eq}
g_{\alpha}(0)\neq0~~~~\text{and}~~~~\wh{\psi}_{\alpha}(\zeta)=\zeta^{\alpha}g_{\alpha}(\zeta)~~~\text{for}~~\zeta\in\C.
\end{align}
For a given $\zeta\in\C$, we define
\begin{align*}
p_{\zeta}(w)=(1+\overline{\zeta}w)^{\alpha}=w^{\alpha}(w^{-1}+\overline{\zeta})^{\alpha},~~~~~~w\in\C.
\end{align*}
Note that $p_{\zeta}(0)=1$ and $|p_{\zeta}(w)|=|(\zeta+w)^{\alpha}|$ for $|w|=1$. Then by maximum modulus principle (e.g. \cite{E.M.Stein2003}), we have
\begin{align}\label{MaximumModulus}
|g_{\alpha}(\zeta)|=|p_{\zeta}(0)g_{\alpha}(\zeta)|\leq\sup_{|w|=1}|p_{\zeta}(w)g_{\alpha}(\zeta+w)|=\sup_{|w|=1}|(\zeta+w)^{\alpha}g_{\alpha}(\zeta+w)|=\sup_{|w|=1}|\psi_{\alpha}(\zeta+w)|,
\end{align}
and by Paley-Wiener-Schwartz theorem \cite{Hormander1983}, there exist $N\in\N$ and $C>0$ such that
\begin{align*}
|\wh{\psi}_{\alpha}(\zeta+w)|\leq C(1+|\zeta+w|)^Ne^{H_K(\im(\zeta+w))}.
\end{align*}
Since $K$ is a closed interval, we can find $R>0$ such that $K\subseteq[-R,R]$. Then for $|w|=1$, $\zeta\in\C$ and $x\in K$, we have
\begin{align*}
x\im(\zeta+w)=x\im w+x\im\zeta\leq|x||w|+x\im\zeta\leq R+x\im\zeta\leq R+H_K(\im\zeta),
\end{align*}
which means that
\begin{align}\label{InEq1}
H_K(\im(\zeta+w))\leq R+H_K(\im\zeta).
\end{align}
In addition, we note that for $\zeta\in\C$ and $|w|=1$,
\begin{align}\label{InEq2}
1+|\zeta+w|\leq 1+|\zeta|+|w|=2+|\zeta|\leq 2(1+|\zeta|).
\end{align}
Combining the above two inequalities \eqref{InEq1} and \eqref{InEq2}, we have
\begin{align*}
(1+|\zeta+w|)^Ne^{H_K(\im(\zeta+w))}\leq 2^Ne^{R}(1+|\zeta|)^Ne^{H_K(\im\zeta)}.
\end{align*}
Consequently, \eqref{MaximumModulus} leads to
\begin{align*}
|g_{\alpha}(\zeta)|\leq 2^NCe^{R}(1+|\zeta|)^Ne^{H_K(\im\zeta)}.
\end{align*}
Again by Paley-Wiener-Schwartz theorem, there exists a distribution $\varphi_{\alpha}$ on $\R$ such that $\su(\varphi_{\alpha})\subseteq K$ and $\wh{\varphi}_{\alpha}(\zeta)=i^{-\alpha}g_{\alpha}(\zeta)$ for $\zeta\in\C$. The uniqueness of $\varphi_{\alpha}$ is obvious since there exists at most one entire function $g_{\alpha}$ satisfying \eqref{Entire_Eq}. In addition, $c_{\alpha}=\wh{\varphi}_{\alpha}(0)\neq0$ because $g_{\alpha}(0)\neq0$. From the construction of $\varphi_{\alpha}$, we have
\begin{align*}
(-1)^{\alpha}\big\la\varphi_{\alpha},\wh{f}^{(\alpha)}\big\ra&=\big\la\varphi_{\alpha},\msF((i\xi)^{\alpha}f)\big\ra=\int_{-\infty}^{\infty}\wh{\varphi}_{\alpha}(\xi)(i\xi)^{\alpha}f(\xi)\rd\xi=\int_{-\infty}^{\infty}i^{-\alpha}g_{\alpha}(\xi)(i\xi)^{\alpha}f(\xi)\rd\xi\\
&=\int_{-\infty}^{\infty}\xi^{\alpha}g_{\alpha}(\xi)f(\xi)\rd\xi=\int_{-\infty}^{\infty}\wh{\psi}_{\alpha}(\xi)f(\xi)\rd\xi~~~~~~f\in\msS,
\end{align*}
where $\msS$ denotes the space of rapidly decaying smooth functions. In other words, $\msF(\varphi_{\alpha}^{(\alpha)})=\wh{\psi}_{\alpha}$ in $\msS'$, where $\msS'$ stands for the space of tempered distributions (i.e. continuous linear functionals on $\msS$) and $\varphi_{\alpha}^{(\alpha)}$, the $\alpha$th derivative of $\varphi_{\alpha}$, is the distribution derivative. Based on the fact that the Fourier transform is a linear isomorphism on $\msS'$ (e.g. \cite{Hormander1983}), it follows that
\begin{align}\label{Diff}
\varphi_{\alpha}^{(\alpha)}=\psi_{\alpha}
\end{align}
in the sense of distribution. Then $K\subseteq\su(\varphi_{\alpha})$, whence $\su(\varphi_{\alpha})=K$. In addition, since $(\varphi_{\alpha}^{(\beta)})^{(\alpha-\beta)}=\psi_{\alpha}$ and $\su(\varphi_{\alpha}^{(\beta)})\subseteq\su(\varphi_{\alpha})=K$ for $0<\beta<\alpha$, we have $\su(\varphi_{\alpha}^{(\beta)})=K$ for all $0\leq\beta\leq\alpha$.

For the regularity of $\varphi_{\alpha}$, first note that since $\psi_{\alpha}\in L_2(\R)$ and $\varphi_{\alpha}$ satisfy \eqref{Diff}, we have $\varphi_{\alpha}\in H^{\alpha}(\R)$ by the elliptic regularity theorem \cite{Folland1999}, where $H^{\alpha}(\R)$ is the Sobolev space defined as
\begin{align*}
H^{\alpha}(\R)&=\big\{u\in\msS':(1+|\xi|^2)^{\alpha/2}\wh{u}\in L_2(\R)\big\}=\big\{u\in L_2(\R):u^{(\beta)}\in L_2(\R)~~\mbox{for}~~\beta\leq\alpha\big\}.
\end{align*}
This means that $\varphi_{\alpha}$ has weak derivatives up to order $\alpha$, and its $\alpha$th weak derivative equals $\psi_{\alpha}$. Moreover, by Sobolev Lemma \cite{Folland1999,Rudin1991}, $\varphi_{\alpha}\in H^{\alpha}(\R)\subseteq C_0^{\alpha-1}(\R)$ where
\begin{align*}
C_0^k(\R)=\big\{f\in C^k(\R):f^{(\alpha)}\in C_0(\R)~~\text{for}~~\alpha\leq k\big\},~~~~k\in\N.
\end{align*}
Hence, $\varphi_{\alpha}$ is $\alpha-1$ differentiable in the classic sense. For the $\alpha$th derivative of $\varphi_{\alpha}$, since $\su(\varphi_{\alpha}^{(\beta)})=K$ for $0\leq\beta\leq\alpha$, it suffices to consider the restrictions of $\varphi_{\alpha}$, $\varphi_{\alpha}',\cdots,\varphi_{\alpha}^{(\alpha-1)},\varphi_{\alpha}^{(\alpha)}=\psi_{\alpha}$ on $K$, with an abuse of notation. Note that $\varphi_{\alpha}\in H^{\alpha}(K^{\circ})$ implies $\varphi_{\alpha}^{(\alpha-1)}\in H^1(K^{\circ})$. Then since $\varphi_{\alpha}^{(\alpha-1)}$ is differentiable a.e. in $K^{\circ}$ and its derivative agrees with its weak derivative a.e. in $K^{\circ}$ \cite{Evans2010}, it follows that \eqref{Diff} holds in the classic sense a.e. in $K^{\circ}$, and thus, a.e. in $\R$.

To complete the proof, we write $\psi_{\aal}$ as
\begin{align*}
\psi_{\aal}(\x)=\psi_{\alpha_1}(x_1)\psi_{\alpha_2}(x_2)
\end{align*}
with $\psi_{\alpha_j}$ being the univariate framelet function having vanishing moments of order $\alpha_j$, and being supported in $[a_j,b_j]$. Then for each $j=1,2$, there exists the unique $\varphi_{\alpha_j}\in H^{\alpha_j}(\R)\subseteq C_0^{\alpha_j-1}(\R)$ such that $\su(\varphi_{\alpha_j})=[a_j,b_j]$,
\begin{align*}
c_{\alpha_j}=\int_{-\infty}^{\infty}\varphi_{\alpha_j}(x)\rd x\neq0~~~~\text{and}~~~~\psi_{\alpha_j}=\varphi_{\alpha_j}^{(\alpha_j)}~~~~\text{a.e. in}~~\R.
\end{align*}
Define $\varphi_{\aal}(\x)=\varphi_{\alpha_1}(x_1)\varphi_{\alpha_2}(x_2)$. Then from the construction of $\varphi_{\aal}$, $\varphi_{\aal}\in L_2(\R^2)$, and
\begin{align*}
c_{\aal}=\int_{\R^2}\varphi_{\aal}(\x)\rd\x=\prod_{j=1}^2\left(\int_{-\infty}^{\infty}\varphi_{\alpha_j}(x_j)\rd x_j\right)=c_{\alpha_1}c_{\alpha_2}\neq0.
\end{align*}
Since each $\varphi_{\alpha_j}$ is uniquely determined from $\psi_{\alpha_j}$, $\varphi_{\aal}$ is uniquely determined from $\psi_{\aal}$ as well. Finally, since each $\varphi_{\alpha_j}$ is $\alpha_j$ differentiable a.e. in $\R$, $\varphi_{\aal}$ is differentiable up to order $\aal$ a.e. in $\R^2$, and we have
\begin{align*}
\p^{\aal}\varphi_{\aal}(\x)=\p^{\aal}\left(\prod_{j=1}^2\varphi_{\alpha_j}(x_j)\right)=\prod_{j=1}^2\p_j^{\alpha_j}\varphi_{\alpha_j}(x_j)=\prod_{j=1}^2\psi_{\alpha_j}(x_j)=\psi_{\aal}(\x).
\end{align*}
This completes the proof of Proposition \ref{Prop1}.

\section{Proof of Theorem \ref{Th1}}\label{ProofTh1}

Note that $E_n^{(3)}(v)\rightarrow E^{(3)}(v)$ and $E_n^{(4)}(u)\rightarrow E^{(4)}(u)$ are already proven in \cite{J.F.Cai2012}, since $W_1^s(\Om)\subseteq L_2(\Om)$ by Sobolev imbedding theorem \cite{Adams1975,H.Attouch2014} and $W_1^r(\Om,[0,1])\subseteq W_1^r(\Om)$. Therefore, we focus on $E_n^{(i)}(u,v)\rightarrow E^{(i)}(u,v)$ for $i=1$, $2$. To prove this, we note that if $v\in W_1^r(\Om,[0,1])$, then so is $1-v$, and $\one-\bsT_n v=\bsT_n(1-v)$. In other words, it is sufficient to prove that for every $(u,v)\in W_1^s(\Om)\times W_1^r(\Om,[0,1])$, we have
\begin{align*}
\left\|\bsT_n v\cdot\big(\lam_n\cdot\bsW_n\bsT_nu\big)\right\|_1:=h^2\sum_{\bk\in\KK_n}\bsT_n v[\bk]\left(\sum_{\aal\in\BB}\lambda_{\aal}\left|\big(\bq_{\aal}[-\cdot]\ast\bsT_nu\big)[\bk]\right|^2\right)^{\f{1}{2}}\longrightarrow\int_{\Om}v\left(\sum_{\aal\in\II}|\p^{\aal}u|^2\right)^{\f{1}{2}}\rd\x
\end{align*}
under a properly chosen $\big\{\lam_n\big\}$.

We split the framelet band $\BB$ into the following two parts:
\begin{align*}
\BB=\II\cup\JJ.
\end{align*}
For $\aal\in\II$, we set $\lambda_{\aal}=\left(c_{\aal}^{-1}2^{|\aal|(n-1)}\right)^2$ where $c_{\aal}$ is given in Proposition \ref{Prop1}. For $\aal\in\JJ$, we set $0\leq\lambda_{\aal}\leq O(2^{2|\bbe|(n-1)})$ for some $\bbe\in\BB\cup\big\{\0\big\}$ such that $\0\leq\bbe<\aal$ and $|\bbe|\leq s$. First we consider $\JJ=\emptyset$. By \eqref{WCMeaning} in Proposition \ref{Prop4}, we have
\begin{align*}
(\bq_{\aal}[-\cdot]\ast\bsT_nu)[\bk]&=\sum_{\bsj\in S_{\aal}+\bk}\bq_{\aal}[\bsj-\bk]\big(\bsT_nu\big)[\bsj]=2^n\left\la u,\sum_{\bsj\in S_{\aal}+\bk}\bq_{\aal}[\bsj-\bk]\phi_{n,\bsj}\right\ra\\
&=2^n\big\la u,\psi_{\aal,n-1,\bk}\big\ra=(-1)^{|\aal|}2^{|\aal|(1-n)+n}\big\la\p^{\aal}u,\varphi_{\aal,n-1,\bk}\big\ra
\end{align*}
for $\bk\in\KK_n$. Hence, it follows that
\begin{align*}
\left\|\bsT_n v\cdot\big(\lam_n\cdot\bsW_n\bsT_nu\big)\right\|_1=2^{-n}\sum_{\bk\in\KK_n}\bsT_n v[\bk]\left(\sum_{\aal\in\II}\left|\big\la\p^{\aal}u,c_{\aal}^{-1}\varphi_{\aal,n-1,\bk}\big\ra\right|^2\right)^{\f{1}{2}}.
\end{align*}
Let $I_{n,\bk}=[\f{k_1}{2^n},\f{k_1+1}{2^n}]\times[\f{k_2}{2^n},\f{k_2+1}{2^n}]$ for $\bk=(k_1,k_2)$. Then we have
\begin{align*}
&\left|\int_{\Om}v\left(\sum_{\aal\in\II}|\p^{\aal}u|^2\right)^{\f{1}{2}}\rd\x-\left\|\bsT_n v\cdot\big(\lam_n\cdot\bsW_n\bsT_nu\big)\right\|_1\right|\\
&=\left|\sum_{\bk\in\OO_n}\int_{I_{n,\bk}}v\left(\sum_{\aal\in\II}|\p^{\aal}u|^2\right)^{\f{1}{2}}\rd\x-2^{-n}\sum_{\bk\in\KK_n}\bsT_n v[\bk]\left(\sum_{\aal\in\II}\left|\big\la\p^{\aal}u,c_{\aal}^{-1}\varphi_{\aal,n-1,\bk}\big\ra\right|^2\right)^{\f{1}{2}}\right|\\
&\leq\left|\sum_{\bk\in\KK_n}\int_{I_{n,\bk}}v\left(\sum_{\aal\in\II}|\p^{\aal}u|^2\right)^{\f{1}{2}}\rd\x-2^{-n}\sum_{\bk\in\KK_n}\bsT_n v[\bk]\left(\sum_{\aal\in\II}\left|\big\la\p^{\aal}u,c_{\aal}^{-1}\varphi_{\aal,n-1,\bk}\big\ra\right|^2\right)^{\f{1}{2}}\right|+\int_{\mathfrak{S}_n}v\left(\sum_{\aal\in\II}|\p^{\aal}u|^2\right)^{\f{1}{2}}\rd\x
\end{align*}
where $\mathfrak{S}_n=\cup_{\bk\in\OO_n\setminus\KK_n}I_{n,\bk}$. Note that the Lebesgue measure $\mathfrak{L}$ of $\mathfrak{S}_n$ satisfies
\begin{align*}
\mathfrak{L}\left(\mathfrak{S}_n\right)\leq 4\left(\f{\diam(\Lambda_{n,\0})}{2^{-n}}+1\right)\left(2^{-n}\right)^2=4c\cdot 2^{-n}
\end{align*}
where $\diam(\Lambda_{n,\0})$ denotes the diameter of $\Lambda_{n,\0}$. Hence, the Lebesgue dominated convergence theorem (e.g. \cite{Folland1999}) leads to
\begin{align*}
\lim_{n\rightarrow\infty}\int_{\mathfrak{S}_n}v\left(\sum_{\aal\in\II}|\p^{\aal}u|^2\right)^{\f{1}{2}}\rd\x=0,
\end{align*}
since $v\in L_{\infty}(\Om)$ and $u\in W_1^s(\Om)$, i.e., $\p^{\aal}u\in L_1(\Om)$ for all $\aal\in\II$, so that the integrand is in $L_1(\Om)$ by the H\"{o}lder's inequality (e.g. \cite{Folland1999}).

For the remaining term, since $0\leq\bsT_nv[\bk]\leq1$ for all $\bk\in\MM_n$, we have
\begin{align*}
&\left|\sum_{\bk\in\KK_n}\int_{I_{n,\bk}}v\left(\sum_{\aal\in\II}|\p^{\aal}u|^2\right)^{\f{1}{2}}\rd\x-2^{-n}\sum_{\bk\in\KK_n}\bsT_n v[\bk]\left(\sum_{\aal\in\II}\left|\big\la\p^{\aal}u,c_{\aal}^{-1}\varphi_{\aal,n-1,\bk}\big\ra\right|^2\right)^{\f{1}{2}}\right|\\
&\leq\sum_{\bk\in\KK_n}\int_{I_{n,\bk}}|v-\bsT_n v[\bk]|\left(\sum_{\aal\in\II}|\p^{\aal}u|^2\right)^{\f{1}{2}}\rd\x\\
&\hspace{11em}+\sum_{\bk\in\KK_n}\int_{I_{n,\bk}}\bsT_n v[\bk]\left|\left(\sum_{\aal\in\II}|\p^{\aal}u|^2\right)^{\f{1}{2}}-\left(\sum_{\aal\in\II}\left|2^n\big\la\p^{\aal}u,c_{\aal}^{-1}\varphi_{\aal,n-1,\bk}\big\ra\right|^2\right)^{\f{1}{2}}\right|\rd\x\\
&\leq\sum_{\bk\in\KK_n}\int_{I_{n,\bk}}|v-\bsT_n v[\bk]|\left(\sum_{\aal\in\II}|\p^{\aal}u|^2\right)^{\f{1}{2}}\rd\x+\sum_{\bk\in\KK_n}\int_{I_{n,\bk}}\left(\sum_{\aal\in\II}\left|\p^{\aal}u-2^n\big\la\p^{\aal}u,c_{\aal}^{-1}\varphi_{\aal,n-1,\bk}\big\ra\right|^2\right)^{\f{1}{2}}\rd\x\\
&\leq\sum_{\bk\in\OO_n}\int_{I_{n,\bk}}\left|v-2^n\big\la v,\phi_{n,\bk}\big\ra\right|\left(\sum_{\aal\in\II}|\p^{\aal}u|^2\right)^{\f{1}{2}}\rd\x+\sum_{\bk\in\OO_n}\int_{I_{n,\bk}}\sum_{\aal\in\II}\left|\p^{\aal}u-2^n\big\la\p^{\aal}u,c_{\aal}^{-1}\varphi_{\aal,n-1,\bk}\big\ra\right|\rd\x\\
&=\int_{\Om}\left|v-\sum_{\bk\in\OO_n}2^n\big\la v,\phi_{n,\bk}\big\ra\chi_{I_{n,\bk}}\right|\left(\sum_{\aal\in\II}|\p^{\aal}u|^2\right)^{\f{1}{2}}\rd\x+\sum_{\aal\in\II}\int_{\Om}\left|\p^{\aal}u-\sum_{\bk\in\OO_n}2^n\big\la\p^{\aal}u,c_{\aal}^{-1}\varphi_{\aal,n-1,\bk}\big\ra\chi_{I_{n,\bk}}\right|\rd\x\\
&\leq\left\|v-\sum_{\bk\in\OO_n}2^n\big\la v,\phi_{n,\bk}\big\ra\chi_{I_{n,\bk}}\right\|_{L_{\infty}(\Om)}\int_{\Om}\left(\sum_{\aal\in\II}|\p^{\aal}u|^2\right)^{\f{1}{2}}\rd\x+\sum_{\aal\in\II}\left\|\p^{\aal}u-\sum_{\bk\in\OO_n}2^n\big\la\p^{\aal}u,c_{\aal}^{-1}\varphi_{\aal,n-1,\bk}\big\ra\chi_{I_{n,\bk}}\right\|_{L_1(\Om)},
\end{align*}
where the last inequality comes from applying the H\"{o}lder's inequality to the first term, and the second to the last equality follows from the fact that
\begin{align*}
\cup_{\bk\in\OO_n}I_{n,\bk}=\overline{\Om}~~~~\text{and}~~~\mathfrak{L}(I_{n,\bsj}\cap I_{n,\bk})=0~~~\text{for}~~\bsj\neq\bk.
\end{align*}
Note that $2^n\chi_{I_{n,\bk}}=\phi_{n,\bk}^H$ where $\phi^H=\chi_{\Om}$, i.e. the refinable function corresponding to Haar framelet which satisfies the partition of unity. Since the piecewise B-spline wavelet frame systems are used, it is obvious that $\int_{\R^2}\phi\rd\x=1$. Moreover, by \eqref{Frameletn} and \eqref{Frameletn-1}, we have
\begin{align*}
c_{\aal}^{-1}\varphi_{\aal,n-1,\bk}=\f{2^n}{4 c_{\aal}}\varphi_{\aal}(2^{n-1}\cdot-2^{-1}\bk)=\left(\f{\varphi_{\aal}(2^{-1}\cdot)}{4 c_{\aal}}\right)_{n,\bk}~~~\text{and}~~~\int_{\R^2}\f{\varphi_{\aal}(2^{-1}\x)}{4 c_{\aal}}\rd\x=1.
\end{align*}
We also note that both $\su(\phi)$ and $\su(\varphi_{\aal}(2^{-1}\cdot))$ contain $\su(\phi^H)$. Therefore, together with $v\in L_{\infty}(\Om)$ and $u\in W_1^s(\Om)$, i.e., $\p^{\aal}u\in L_1(\Om)$ for all $\aal\in\II$, we establish
\begin{align*}
&\left\|v-\sum_{\bk\in\OO_n}2^n\big\la v,\phi_{n,\bk}\big\ra\chi_{I_{n,\bk}}\right\|_{L_{\infty}(\Om)}\int_{\Om}\left(\sum_{\aal\in\II}|\p^{\aal}u|^2\right)^{\f{1}{2}}\rd\x+\sum_{\aal\in\II}\left\|\p^{\aal}u-\sum_{\bk\in\OO_n}2^n\big\la\p^{\aal}u,c_{\aal}^{-1}\varphi_{\aal,n-1,\bk}\big\ra\chi_{I_{n,\bk}}\right\|_{L_1(\Om)}\longrightarrow0
\end{align*}
by the approximation lemma \cite[Lemma 4.1]{J.F.Cai2012}.

For $\JJ\neq\emptyset$, if we show that
\begin{align}\label{Goal2}
\lim_{n\rightarrow\infty}\left(\lambda_{\aal}\right)^{\f{1}{2}}\big\|\bq_{\aal}[-\cdot]\ast\bsT_nu\big\|_1=0~~~\text{for all}~~~\aal\in\JJ,
\end{align}
then we complete the proof. Indeed, we define
\begin{align*}
E_{\II}:=h^2\sum_{\bk\in\KK_n}\bsT_n v[\bk]\left(\sum_{\aal\in\II}\lambda_{\aal}\left|\big(\bq_{\aal}[-\cdot]\ast\bsT_nu\big)[\bk]\right|^2\right)^{\f{1}{2}}.
\end{align*}
Then since $0\leq\bsT_nv[\bk]\leq 1$ for all $\bk\in\MM_n$, we have
\begin{align*}
E_{\II}\leq\left\|\bsT_n v\cdot\big(\lam_n\cdot\bsW_n\bsT_nu\big)\right\|_1\leq E_{\II}+\sum_{\aal\in\JJ}\left(\lambda_{\aal}\right)^{\f{1}{2}}\big\|\bq_{\aal}[-\cdot]\ast\bsT_nu\big\|_1.
\end{align*}
Once we have \eqref{Goal2}, then taking the limit of the above inequality leads to
\begin{align*}
\lim_{n\rightarrow\infty}E_{\II}=\lim_{n\rightarrow\infty}\left\|\bsT_n v\cdot\big(\lam_n\cdot\bsW_n\bsT_nu\big)\right\|_1.
\end{align*}
By Proposition \ref{Prop1}, there exist $\varphi_{\aal}$ and $\varphi_{\bbe}$ such that $\p^{\aal}\varphi_{\aal}=\psi_{\aal}$ and $\p^{\bbe}\varphi_{\bbe}=\psi_{\bbe}$ a.e. We set $\bbe\in\BB\cup\big\{\0\big\}$ such that $\0\leq\bbe<\aal$ and $|\bbe|\leq s$, as mentioned in the beginning of the proof. Indeed, such $\bbe$ always exists, since, for example, one may pick $\bbe=\0$. Let $\overline{\psi}_{\aal}=\p^{\aal-\bbe}\varphi_{\aal}$. Then it is obvious that $\p^{\bbe}\overline{\psi}_{\aal}=\psi_{\aal}$ due to the tensor product structure of $\varphi_{\aal}$. For $t\geq0$, we define
\begin{align*}
\wt{\varphi}_t=c_{\bbe}^{-1}\varphi_{\bbe}+t\overline{\psi}_{\aal}.
\end{align*}
Then $\wt{\varphi}_t$ is compactly supported, (i.e. $\su(\wt{\varphi}_t)\subseteq\su(\phi)$), differentiable a.e. up to order $\bbe$, and $\int_{\R^2}\wt{\varphi}_t\rd\x=1$. Together with $\p^{\bbe}\wt{\varphi}_t=c_{\bbe}^{-1}\psi_{\bbe}+t\psi_{\aal}$, we have
\begin{align*}
\big\la\p^{\bbe}u,\wt{\varphi}_{t,n-1,\bk}\big\ra=(-1)^{|\bbe|}2^{|\bbe|(n-1)}\big\la u,c_{\bbe}^{-1}\psi_{\bbe,n-1,\bk}+t\psi_{\aal,n-1,\bk}\big\ra
\end{align*}
for $u\in W_1^s(\Om)$. Therefore,
\begin{align*}
2^{|\bbe|(n-1)}\left\|\big(c_{\bbe}^{-1}\bq_{\bbe}[-\cdot]+t\bq_{\aal}[-\cdot]\big)\ast\bsT_nu\right\|_1=2^{-n}\sum_{\bk\in\KK_n}\left|\big\la\p^{\bbe}u,\wt{\varphi}_{t,n-1,\bk}\big\ra\right|,
\end{align*}
and following the similar steps as $\JJ=\emptyset$ by setting $v\equiv1$ and replacing isotropic $\ell_1$ norm by anisotropic $\ell_1$ norm, $\II$ by $\JJ$, $\varphi$ by $\wt{\varphi}_t$, $\p^{\aal}$ by $\p^{\bbe}$, and $c_{\aal}$ by $c_{\bbe}$, we have
\begin{align*}
\lim_{n\rightarrow\infty}2^{|\bbe|(n-1)}\left\|\big(c_{\bbe}^{-1}\bq_{\bbe}[-\cdot]+t\bq_{\aal}[-\cdot]\big)\ast\bsT_nu\right\|_1=\int_{\Om}\left|\p^{\bbe}u\right|\rd\x.
\end{align*}
In particular, when $t=0$, we have
\begin{align*}
\lim_{n\rightarrow\infty}2^{|\bbe|(n-1)}\left\|c_{\bbe}^{-1}\bq_{\bbe}[-\cdot]\ast\bsT_nu\right\|_1=\int_{\Om}\left|\p^{\bbe}u\right|\rd\x.
\end{align*}
These two equalities imply that
\begin{align*}
t\limsup_{n\rightarrow\infty}2^{|\bbe|(n-1)}\big\|\bq_{\aal}[-\cdot]\ast\bsT_nu\big\|_1\leq2\int_{\Om}\left|\p^{\bbe}u\right|\rd\x.
\end{align*}
Since $t\geq0$ is arbitrary, it must be
\begin{align*}
\lim_{n\rightarrow\infty}2^{|\bbe|(n-1)}\big\|\bq_{\aal}[-\cdot]\ast\bsT_nu\big\|_1=0.
\end{align*}
In view of $0\leq\lambda_{\aal}\leq O(2^{2|\bbe|(n-1)})$ for $\aal\in\JJ$, we obtain \eqref{Goal2}. This completes the proof of Theorem \ref{Th1}.

\section{Proof of Proposition \ref{Prop2}}\label{ProofProp2}

Since $W_1^s(\Om)\times W_1^r(\Om,[0,1])$ is closed in $W_1^s(\Om)\times\msX$, it suffices to prove that $E_n$ is equicontinuous as a sequence of functionals on $W_1^s(\Om)\times\msX$. First we note that the equicontinuity of $E_n^{(4)}$ is already proved in \cite[Proposition 3.2]{J.F.Cai2012}. Moreover, the proof of $E_n^{(3)}$ follows the same step as \cite[Proposition 3.2]{J.F.Cai2012} by replacing $\|\cdot\|_{W_1^s(\Om)}$ with $\|\cdot\|_{\msX}$. Hence, we shall focus on the equicontinuity of $E_n^{(i)}$ in $W_1^s(\Om)\times\msX$ for $i=1$, $2$. To do this, we note that if we extend $E_n$ and $E$ to $W_1^s(\Om)\times\msX$, then the first two terms become
\begin{align*}
E_n^{(1)}(u,v)&=h^2\sum_{\bk\in\KK_n}|\one-\bsT_n v[\bk]|\left(\sum_{\aal\in\BB}\lambda_{\aal}[\bk]|(\bq_{\aal}[-\cdot]\ast\bsT_nu)[\bk]|^2\right)^{\f{1}{2}}\\
E_n^{(2)}(u,v)&=h^2\sum_{\bk\in\KK_n}|\bsT_n v[\bk]|\left(\sum_{\aal\in\BB}\gamma_{\aal}[\bk]|(\bq_{\aal}[-\cdot]\ast\bsT_nu)[\bk]|^2\right)^{\f{1}{2}}\\
E^{(1)}(u,v)&=\lambda\int_{\Om}|1-v|\left(\sum_{\aal\in\II}|\p^{\aal}u|^2\right)^{\f{1}{2}}\rd\x\\
E^{(2)}(u,v)&=\gamma\int_{\Om}|v|\left(\sum_{\aal\in\II'}|\p^{\aal}u|^2\right)^{\f{1}{2}}\rd\x,
\end{align*}
and the pointwise convergence of $E_n^{(i)}(u,v)$ to $E^{(i)}(u,v)$ for $(u,v)\in W_1^s(\Om)\times\msX$ can be proven in the similar way. In addition, if $v\in\msX$, then so is $1-v$, and $\one-\bsT_nv=\bsT_n(1-v)$. Therefore, as in Theorem \ref{Th1}, it is sufficient to prove the equicontinuity of
\begin{align*}
\left\|\bsT_n v\cdot\big(\lam_n\cdot\bsW_n\bsT_nu\big)\right\|_1=h^2\sum_{\bk\in\KK_n}\left|\bsT_n v[\bk]\right|\left(\sum_{\aal\in\BB}\lambda_{\aal}|(\bq_{\aal}[-\cdot]\ast\bsT_nu)[\bk]|^2\right)^{\f{1}{2}}
\end{align*}
on $W_1^s(\Om)\times\msX$ under the parameter $\big\{\lam_n\big\}$ chosen as in \ref{ProofTh1}.

We define the space $\ell_{1,2}^{\star}(\Z^2):=\big\{\bb:\|\bb\|_{1,2}^{\star}<\infty\big\}$ with
\begin{align*}
\|\bb\|_{1,2}^{\star}=\sum_{\bk\in\Z^2}\left(\sum_{\aal\in\BB}|\bb_{\aal}[\bk]|^2\right)^{\f{1}{2}}.
\end{align*}
We fix $v\in\msX$. For any given $n$ and $u\in W_1^s(\Om)$, we have $\bsT_n v\cdot\big(\lam_n\cdot\bsW_n\bsT_nu\big)\in\ell_{1,2}^{\star}(\Z^2)$:
\begin{align*}
\left\|2^{-2n}\bsT_n v\cdot\big(\lam_n\cdot\bsW_n\bsT_nu\big)\right\|_{1,2}^{\star}=\left\|\bsT_n v\cdot\big(\lam_n\cdot\bsW_n\bsT_nu\big)\right\|_1.
\end{align*}
Since $\bsT_n$ is a bounded linear operator on $L_2(\Om)$ to a finite dimensional space $\R^{\MM_n}\simeq\R^{|\MM_n|}$ and $\bsW_n$ can be understood as a $(r+1)^2|\KK_n|\times|\MM_n|$ matrix, we have
\begin{align*}
\left\|2^{-2n}\bsT_n v\cdot\big(\lam_n\cdot\bsW_n\bsT_nu\big)\right\|_{1,2}^{\star}\leq A_n(v)\|u\|_{L_2(\Om)}\leq\wt{A}_n(v)\|u\|_{W_1^s(\Om)}
\end{align*}
where the last inequality follows from the Sobolev imbedding theorem \cite{Adams1975,H.Attouch2014}, and the constant is depend on $n$ and $v\in\msX$. This means that for each $v\in\msX$,
\begin{align*}
2^{-2n}\bsT_n v\cdot\big(\lam_n\cdot\bsW_n\bsT_n(\cdot)\big)\in\mB(W_1^s(\Om),\ell_{1,2}^{\star}(\Z^2)).
\end{align*}
In addition, since for any given $u\in W_1^s(\Om)$,
\begin{align}\label{Limit}
\lim_{n\rightarrow\infty}\left\|\bsT_n v\cdot\big(\lam_n\cdot\bsW_n\bsT_nu\big)\right\|_1=\int_{\Om}|v|\left(\sum_{\aal\in\II}|\p^{\aal}u|^2\right)^{\f{1}{2}}\rd\x,
\end{align}
we have
\begin{align*}
\sup_n\left\|2^{-2n}\bsT_n v\cdot\big(\lam_n\cdot\bsW_n\bsT_nu\big)\right\|_{1,2}^{\star}=\sup_n\left\|\bsT_n v\cdot\big(\lam_n\cdot\bsW_n\bsT_nu\big)\right\|_1<\infty.
\end{align*}
Recall from the uniform boundedness principle (e.g. \cite{Conway1990}) that for a sequence of bounded linear operators on a Banach space, pointwise boundedness is equivalent to uniform boundedness in operator norm. Therefore, we have
\begin{align*}
\sup_n\left\|2^{-2n}\bsT_n v\cdot\big(\lam_n\cdot\bsW_n\bsT_n(\cdot)\big)\right\|\leq A(v)<\infty
\end{align*}
for some constant $A(v)>0$ depending only on $v\in\msX$. Here, $\|\cdot\|$ stands for the operator norm.

We again define the space $\ell_1^{\star}(\Z^2)=\big\{\bsv:\|\bsv\|_1^\star<\infty\big\}$ with
\begin{align*}
\|\bsv\|_1^\star=\sum_{\bk\in\Z^2}|\bsv[\bk]|.
\end{align*}
Here, we fix $u\in W_1^s(\Om)$. For any given $n$ and $v\in\msX$, we have $\bsT_n v\cdot\big(\lam_n\cdot\bsW_n\bsT_nu\big)\in\ell_1^\star(\Z^2)$:
\begin{align*}
\left\|2^{-2n}\bsT_n v\cdot\big(\lam_n\cdot\bsW_n\bsT_nu\big)\right\|_1^\star=\left\|\bsT_n v\cdot\big(\lam_n\cdot\bsW_n\bsT_nu\big)\right\|_1.
\end{align*}
Since $\bsT_n$ is a bounded linear operator from $L_2(\Om)$ to $\R^{\MM_n}\simeq\R^{|\MM_n|}$ and the mapping $\bsv_n\mapsto\bsv_n\cdot\big(\lam_n\cdot\bsW_n\bsT_nu\big)$ can be understood as the multiplication of a diagonal matrix and a vector, we have
\begin{align*}
\left\|2^{-2n}\bsT_n v\cdot\big(\lam_n\cdot\bsW_n\bsT_nu\big)\right\|_1^\star\leq B_n(u)\|v\|_{L_2(\Om)}\leq\wt{B}_n(u)\|v\|_{\msX}
\end{align*}
where the last inequality follows from the fact that $\msX\subseteq L_2(\Om)$ due to Sobolev imbedding theorem and the boundedness of $\Om$. Again, the constant is dependent on $u\in W_1^s(\Om)$. This means that for each $u\in W_1^s(\Om)$,
\begin{align*}
2^{-2n}\bsT_n(\cdot)\cdot\big(\lam_n\cdot\bsW_n\bsT_nu\big)\in\mB(\msX,\ell_1^{\star}(\Z^2)).
\end{align*}
Since \eqref{Limit} holds for every $v\in\msX$ with a fixed $u\in W_1^s(\Om)$ as well, we have
\begin{align*}
\sup_n\left\|2^{-2n}\bsT_n v\cdot\big(\lam_n\cdot\bsW_n\bsT_nu\big)\right\|_1^\star=\sup_n\left\|\bsT_n v\cdot\big(\lam_n\cdot\bsW_n\bsT_nu\big)\right\|_1<\infty
\end{align*}
for every $v\in\msX$. Again, by the uniform boundedness principle, we have
\begin{align*}
\sup_n\left\|2^{-2n}\bsT_n(\cdot)\cdot\big(\lam_n\cdot\bsW_n\bsT_nu\big)\right\|\leq B(u)
\end{align*}
for some constant $B(u)>0$ depending only on $u\in W_1^s(\Om)$.

Let $(u,v)\in W_1^s(\Om)\times\msX$. For $(u',v')\in W_1^s(\Om)\times\msX$, we have
\begin{align*}
\bigg|&\left\|\bsT_n v'\cdot\big(\lam_n\cdot\bsW_n\bsT_nu'\big)\right\|_1-\left\|\bsT_n v\cdot\big(\lam_n\cdot\bsW_n\bsT_nu\big)\right\|_1\bigg|\\
&\leq\bigg|\left\|\bsT_n v'\cdot\big(\lam_n\cdot\bsW_n\bsT_nu'\big)\right\|_1-\left\|\bsT_n v'\cdot\big(\lam_n\cdot\bsW_n\bsT_nu\big)\right\|_1\bigg|\\
&\hspace{13.1em}+\bigg|\left\|\bsT_n v'\cdot\big(\lam_n\cdot\bsW_n\bsT_nu\big)\right\|_1-\left\|\bsT_n v\cdot\big(\lam_n\cdot\bsW_n\bsT_nu\big)\right\|_1\bigg|\\
&\leq\left\|\bsT_n v'\cdot\big[\lam_n\cdot\bsW_n\bsT_n(u'-u)\big]\right\|_1+\left\|\bsT_n(v'-v)\cdot\big(\lam_n\cdot\bsW_n\bsT_nu\big)\right\|_1\\
&\leq\|v'\|_{L_{\infty}(\Om)}\big\|\lam_n\cdot\bsW_n\bsT_n(u'-u)\big\|_1+\left\|\bsT_n(v'-v)\cdot\big(\lam_n\cdot\bsW_n\bsT_nu\big)\right\|_1\\
&=\|v'\|_{L_{\infty}(\Om)}\big\|2^{-2n}\lam_n\cdot\bsW_n\bsT_n(u'-u)\big\|_{1,2}^\star+\left\|2^{-2n}\bsT_n(v'-v)\cdot\big(\lam_n\cdot\bsW_n\bsT_nu\big)\right\|_1^\star\\
&\leq\|v'\|_{\msX}A(1)\|u'-u\|_{W_1^s(\Om)}+B(u)\|v'-v\|_{\msX}\\
&\leq A(1)(\|v\|_{\msX}+\|v'-v\|_{\msX})\|u'-u\|_{W_1^s(\Om)}+B(u)\|v'-v\|_{\msX}\\
&\leq A(1)\|v\|_{\msX}\|u'-u\|_{W_1^s(\Om)}+\f{A(1)}{2}\big(\|u'-u\|_{W_1^s(\Om)}+\|v'-v\|_{\msX}\big)^2+B(u)\|v'-v\|_{\msX}\\
&\leq C\big(\|u'-u\|_{W_1^s(\Om)}+\|v'-v\|_{\msX}\big)\left[\big(\|u'-u\|_{W_1^s(\Om)}+\|v'-v\|_{\msX}\big)+1\right]
\end{align*}
where $C=\max\left\{A(1)\|v\|_{\msX},B(u),A(1)/2\right\}$ is independent of $n$, and the third inequality follows from the stability of $\bsT_n$. For a given $\eps>0$, we choose $\mN=1$ and
\begin{align*}
\delta=\f{-C+\sqrt{C^2+4C\eps}}{2C}>0
\end{align*}
both of which are again independent of $n$. Therefore, whenever $n>\mN$ and $\|u'-u\|_{W_1^s(\Om)}+\|v'-v\|_{\msX}<\delta$, we have
\begin{align*}
\bigg|\left\|\bsT_n v'\cdot\big(\lam_n\cdot\bsW_n\bsT_nu'\big)\right\|_1-\left\|\bsT_n v\cdot\big(\lam_n\cdot\bsW_n\bsT_nu\big)\right\|_1\bigg|<\eps,
\end{align*}
which completes the proof of Proposition \ref{Prop2}.
}

\section*{References}

\bibliographystyle{siam}

\end{document}